\documentstyle[12pt]{article}
\newtheorem{df}{Definition}[section]
\newtheorem{lm}{Lemma}[section]
\newtheorem{pr}{Proposition}[section]
\newtheorem{co}{Corollary}[section]
\newtheorem{th}{Theorem}[section]
\newtheorem{rk}{Remark}[section]
\newtheorem{ex}{Example}[section]
\newtheorem{prob}{Problem}[section]

\newtheorem{conj}{Conjecture}[section]

\newcommand{\bsquare}{\hbox{\rule{6pt}{6pt}}}

\newcommand{\mapright}[1]{\smash{\mathop{
\hbox to 1cm{\rightarrowfill}}\limits^{#1}}}

\newcommand{\mapleft}[1]{\smash{\mathop{
\hbox to 1cm{\leftarrowfill}}\limits_{#1}}}

\newcommand{\mapdown}[1]{\Big\downarrow
\llap{$\vcenter{\hbox{$\scriptstyle#1\,$}}$ }}

\newcommand{\mapdownr}[1]{\Big\downarrow
\rlap{$\vcenter{\hbox{$\scriptstyle#1\,$}}$ }}

\begin{document}
\title{The Euler Characteristic Formula for Logarithmic Minimal 
Degenerations \\ of Surfaces with Kodaira Dimension Zero\\and its application 
to Calabi-Yau Threefolds with a pencil
\thanks{This is an extended version of the preprint titled \lq\lq 
The Euler Characteristic Formula for Logarithmic Minimal Degenerations of 
Surfaces with Kodaira Dimension Zero''circulated earlier.} }
\author{Koji Ohno
\\ Department of Mathematics Graduate School of Science \\Osaka University}

\maketitle 

\begin{abstract} In this paper, the Euler characteristic formula for projective logarithmic minimal
degenerations of surfaces with Kodaira dimension zero over a 1-dimensional complex disk is proved under a reasonable assumption and as its application, the singularity of logarithmic
minimal degenerations are determined in the abelian or hyperelliptic case. By globalizing this local
analysis of singular fibres via generalized canonical bundle formulae due to Fujino-Mori, we bound the number of singular fibres of abelian fibred Calabi-Yau threefolds from above,
which was previously done by Oguiso in the potentially good reduction case.

\end{abstract}

\par
\tableofcontents

\section{Introduction}
Based on the 2-dimensional minimal model theory, Kodaira classified the singular fibres of degenerations of elliptic curves ( \cite{kodaira}, Theorem 6.2 ). It is quite natural that
many people have been interested in the degenerations of surfaces with Kodaira dimension zero as a next problem. The first effort began by his student Iitaka and Ueno who studied the
first kind degeneration ( i.e., degeneration with the finite monodromy ) of abelian surfaces with a principal polarization ( \cite{ueno1} and \cite{ueno2} ) while in that time,
3-dimensional minimal model theory had not been known. After Kulikov succeeded to construct the minimal models of degenerations of algebraic K3 surfaces in the analytic category from
semistable degenerations and to classify  their singular fibres (\cite{kulikov}), extension to the case of the other surfaces with Kodaira dimension zero has been 
done ( see for example, \cite{persson-pinkham}, \cite{morrison}, \cite{shokurov2} ). As for the non-semistable case, there are works due to Crauder and
Morrison who classified triple point free degeneration (\cite{crauder}, \cite{crauder2}). 3-dimensional minimal model theory in the projective category was established by
Mori ( \cite{mori} ) but we can not start studying the degenerations from minimal models because of their complexity while it has been known that log minimal models of degenerations of 
elliptic curves behaves nicely (see \cite{reid}, (8.9) Added in Proof.).  After the establishment of 3-dimensional log minimal model theory, we introduced the notion of a {\it
logarithmic minimal degeneration} in \cite{ohno} as a good intermediate model to a minimal model which acts like a \lq\lq quotient" of minimal semistable degeneration by the
transformation group induced from a semistable reduction. Of course, because of the non-uniqueness of minimal models, the transformation group does not act holomorphically on the total
space in general.

\begin{df}\label{df:lmf}\begin{em} 
Let $f:X\rightarrow B$ be a proper connected morphism from a normal {\bf Q}-factorial variety defined over the complex number
field { \mbox{\boldmath $C$} } (resp. a normal {\bf Q}-factorial complex analytic space)
$X$ onto a smooth projective curve (resp. a unit disk ${\cal D}:=\{z\in \mbox{\boldmath $C$};|z|<1\}$ ) $B$ such that a general fibre $f^{\ast}(p)$
(resp. any fibre $f^{\ast}(p)$ where $p$ is not the origin) is a normal algebraic variety with only terminal singularity. Let $\Sigma$ be a set of
points in $B$ (resp. the origin $0$) such that the fibre $f^{\ast}(p)$ is not a normal algebraic variety with only terminal singularity. Put $\Theta_p:=f^{\ast}(p)_{\rm
red}$ and $\Theta:=\sum_{p\in\Sigma}\Theta_p$.  

\begin{description}
\item[ ($1$) ]$f:X\rightarrow B$ is called a {\it minimal fibration $($resp. degeneration$)$} if 
$X$ has only terminal singularity and
$K_X$ is
$f$-nef (i.e., The intersection number of $K_X$ and any complete curve contained in a fibre of $f$ is non-negative). 
\item[ ($2$) ]$f:X\rightarrow B$ is called a {\it logarithmic minimal $($or abbreviated, log minimal$)$ fibration $($resp. degeneration$)$} if 
$(X,\Theta)$ is divisorially log terminal and
$K_X+\Theta$ is
$f$-nef. 
\item[ ($3$) ]$f:X\rightarrow B$ is called a {\it strictly logarithmic minimal $($or abbreviated, strictly log minimal$)$ fibration $($resp.
degeneration$)$} if $(X,\Theta)$ is log canonical with $K_X$, $\Theta$ being both $f$-nef.
\end{description}
\end{em}\end{df}

\begin{rk}\begin{em}\label{rk:construction of slmm} We note that any fibrations (resp. degenerations) of algebraic surfaces with Kodaira dimension zero over a smooth projective curve
(resp. 1-dimensional unit disk)
$B$ are birationally (resp. bimeromorphically) equivalent over $B$ to a projective log minimal fibration (resp. degeneration) and also to a projective strictly log minimal
fibration (resp. degeneration). In fact, firstly we can take a birational (resp. bimeromorphic) model
$g:Y\rightarrow B$, where $g$ is a relatively projective connected morphism from a smooth variety (resp. complex analytic space) $Y$ such that for any singular fibre of $g$, its
support has only simple normal crossing singularity with each component smooth by the Hironaka's theorem (\cite{hironaka}). Let $\Sigma\subset B$ be a set of all the
points such that $g$ is not smooth over $p\in B$ (resp. the origin $0$). By the existence theorem of log minimal models established in \cite{shokurov}, \cite{kawamata term.}, Theorem
2, \cite{kollar/utah}, Theorem 1.4, we can run the log minimal program with respect to $K_Y+\sum _{p\in
\Sigma}g^{\ast}(p)_{\rm red}$ starting from
$Y$ to get a log minimal model $f:X\rightarrow B$ which is a log minimal fibration (resp. degenerations) in the sense of Definition~\ref{df:lmf}. Here we note that by the Base Point Free
Theorem in
\cite{nakayama}, we infer that $K_X+\Theta\sim_{ \mbox{\boldmath $Q$} }f^{\ast}D$ for some {\bf Q}-divisor $D$ on $B$. By applying the log minimal program
with respect to
$K_X$ starting from $X$, we obtain a model $f^s:X^s\rightarrow B$ which obviously turns out to be a strictly log minimal fibration (resp. degenerations).
\end{em}
\end{rk}

\begin{df}\begin{em}
Let $G$ be a finite group and $\rho:G\rightarrow
\mbox{GL}(3,\mbox{\boldmath $C$})$ be a faithful representation. Let
$\mbox{\boldmath
$C$}^3/(G,\rho)$ denote the quotient of $\mbox{\boldmath $C$}^3$ by the
action of $G$ defined by $\rho$. We assume that the quotient map
$\mbox{\boldmath $C$}^3\rightarrow \mbox{\boldmath
$C$}^3/(G,\rho)$ is \'etale in codimension one. A pair
$(X, D)$ which consists of a normal complex analytic space $X$ and a
reduced divisor $D$ on $X$ is said to {\it have singularity of type
$V_1(G,\rho)$
$($resp.
$V_2(G,\rho)$
$)$ at
$p\in X$} if there exists an analytic isomorphism
$\varphi:(X,p)\rightarrow (\mbox{\boldmath
$C$}^3/(G,\rho),0)$ between germs and a hypersurface $H$ in
$\mbox{\boldmath $C$}^3$ defined by the equation $z=0$ (resp. $xy=0$),
where
$x$, $y$ and $z$ are semi-invariant coordinates of $\mbox{\boldmath
$C$}^3$ at
$0$ such that
$D=\varphi^{\ast}(H/(G,\rho))$. In particular, if $G$ is cyclic with a
generator $\sigma\in G$ and
$(\rho(\sigma)^{\ast}x,
\rho(\sigma)^{\ast}y,\rho(\sigma)^{\ast}z)=(\zeta^ax,\zeta^by,\zeta^cz)$,
where $a$,
$b$, $c\in {\mbox{\boldmath $Z$}}$ and $\zeta$ is a primitive $r$-th
root of unity for some coordinate $x$, $y$ and $z$ of $\mbox{\boldmath
$C$}^3$ at $0$, we shall use the
notation
$V_1(r;a,b,c)$ ( resp. $V_2(r;a,b,c)$ ) instead of $V_1(G,\rho)$ ( resp.
$V_2(G,\rho)$ ).
\end{em}\end{df}

\begin{rk}\begin{em} We note that if $(X,D)$ has singularity of type
$V_i(G,\rho)$ at $p$, the local fundamental group at $p$ of the singularity of $X$ is isomorphic to $G$ by its definition.
\end{em}\end{rk}

Let $f:X\rightarrow {\cal D}$ be a log minimal degeneration and let
$\Theta=\sum_i\Theta_i$ be the irreducible decomposition and put 
$\Delta_i:=\mbox{Diff}_{\Theta_i}(\Theta-\Theta_i)$ for any $i$. For $p\in X$, let $d(p)$ be
the  number of irreducible components of $\Theta$ passing through $p\in X$. Then the followings
hold.
\begin{description}
\item[ (a) ] For any $i$, $\Theta_i$ is normal, $\Delta_i$ is a standard boundary (see Definition~\ref{df:standard}) and $(\Theta_i,\Delta_i)$ is log
terminal (see \cite{shokurov}, Lemma 3.6, (3.2.3) and Corollary 3.10).   
\item[ (b) ] $d(p)\leq 3$.
\item[ (c) ] If $d(p)=2$, $(X,\Theta)$ has singularity of type $V_2(r;a,b,1)$ at $p$,
where
$r\in \mbox{\boldmath $N$}$, $a,b\in \mbox{\boldmath $Z$}$ and
$(r,a,b)=1$ (see \cite{corti/utah}, Theorem 16.15.2).
\item[ (d) ] If $d(p)=3$, $p\in \Theta\subset X$ is analytically isomorphic to the germ of 
the origin 
$0\in \{(x,y,z);xyz=0\}\subset {\mbox{\boldmath $C$}}^3$ (see \cite{corti/utah}, Theorem 16.15.1).
\item[ (e) ] For any $i$ and $p\in \Theta_i\setminus \mbox{\rm Supp }\Delta_i$, if $\Theta_i$
is smooth at $p$, then $X$ is smooth at $p$ (see \cite{shokurov}, Corollary 3.7).
\end{description}

One of the aims of this paper is to give the following Euler characteristic formula for log minimal
degenerations with 
 $K_X+\Theta$ being Cartier. We note here that the study of log minimal degenerations of surfaces with Kodaira dimension zero reduces to this case by taking the log canonical cover
with respect to $K_X+\Theta$ globally (see \S 6). 
\begin{th}\label{th:euler characteristic formula} 
Let $f:X\rightarrow {\cal D}$ be a projective log minimal degeneration of surfaces with Kodaira dimension zero such that 
$K_X+\Theta$ is Cartier. Let $f^{\ast}(0)=\sum_{i}m_i\Theta_i$ be the irreducible decomposition. 
Then for 
$t\in {\cal D}^{\ast}:={\cal D}\setminus \{0\}$, the following formula holds.
$$
e_{\mbox{\em top}}(X_t)=\sum m_i(e_{\mbox{\em orb}}(\Theta_i\setminus \Delta_i)+\sum_{p\in
\Theta_i\setminus \Delta_i}\delta_p(X,\Theta_i)),
$$
where $X_t:=f^{\ast}(t)$, $e_{\mbox{\em orb}}(\Theta_i\setminus \Delta_i)$ is the orbifold Euler number of
$\Theta_i\setminus \Delta_i$ and 
$\delta_p(X,\Theta_i)$ is the invariant of the singularity of the pair $(X,\Theta_i)$ at $p\in
\Theta_i\setminus \Delta_i$ 
which is well defined and can be calculated explicitly as explained in the next section.
\end{th}
The above formula turns out to be quite useful for further study of degenerations. In fact, we
apply the following corollary to the study on non-semistable degenerations of abelian or hyperelliptic
surfaces.

\begin{co}\label{co:chi0}
Let notation and assumptions be as in Theorem~\ref{th:euler characteristic formula}. Assume that
$e_{\mbox{\em top}}(X_t)=0$ for 
$t\in {\cal D}^{\ast}$. Then, for any $i$, we have $e_{\mbox{\em orb}}(\Theta_i\setminus \Delta_i)=0$ and for any $p\in \Theta_i\setminus\Delta_i$, $(X,\Theta)$ has only singularity of
type
$V_1(r;a,-a,1)$ at $p$, where $(r,a)=1$.
\end{co}

Based on the result of Corollary~\ref{co:chi0}, we shall prove the following theorem.

\begin{th}\label{th:main th.ab} Let $f:X\rightarrow {\cal D}$ be a
projective log minimal degeneration of abelian or hyperelliptic
surfaces, not neccesarily assuming that $K_X+\Theta$ is Cartier. Then the possible singularities of $(X,\Theta)$ at $p\in X$ are
the following three types $:$
\begin{description}
\item[(0)] $X$ is smooth at $p\in X$ and $\Theta$ has only normal
crossing singularity at $p$,
\item[(1)] $(X,\Theta)$ has singularity of type $V_2(r;a,b,1)$ at $p$,
where
$r\in \mbox{\boldmath $N$}$, $a,b\in \mbox{\boldmath $Z$}$ and
$(r,a,b)=1$.
\item[(2)] $(X,\Theta)$ has singularity of type $V_1(G,\rho)$ at $p$.
\end{description} 
More precisely, if $f$ is of type {\rm II}, we have $r=2,3,4$ or $6$ in
$(1)$, and $G\simeq \mbox{\boldmath $Z$}/n\mbox{\boldmath $Z$}$ or
$\mbox{\boldmath $Z$}/2\mbox{\boldmath $Z$}\oplus \mbox{\boldmath
$Z$}/n\mbox{\boldmath $Z$}$, where $n=2,3,4$ or $6$ in $(2)$. The dual graph of $\Theta$ is a linear chain or a cycle. Moreover,
there exists a  projective bimeromorphic morphism
$\psi:X\rightarrow X^s$ over
${\cal D}$ such that for the induced  projective degeneration
$f^s:X^s\rightarrow {\cal D}$, we have  
$K_{X^s}\sim_{\mbox{\boldmath
$Q$}}0$ and
$f^{s \ast}(0)=m\Theta ^s$ for some $m\in \mbox{\boldmath $N$}$,
where 
$\Theta ^s:=\psi_{\ast}\Theta$. The possible types of singularity of
$(X^s,\Theta ^s)$ and the dual graph of the support of the singular fibre are the same as ones of $(X,\Theta)$ $($ but the
components of the singular fiber may become non-normal $)$. If
$f$ is of type {\rm III}, we have $r=2$ in $(1)$, and $(2)$ is reduced to
the following three types. 
\begin{description}
\item[(III-2.1)] $(X,\Theta)$ has singularity of type $V_1(r;a,-a,1)$ at
$p$, where
$r=2,3,4$ or $6$, $a\in \mbox{\boldmath $Z$}$ and $(r,a)=1$,
\item[(III-2.2)] $(X,\Theta)$ has singularity of type $V_1(2;1,0,1)$ at
$p$,
\item[(III-2.3)] $(X,\Theta)$ has singularity of type $V_1(G,\rho)$ at
$p$, where $G\simeq \mbox{\boldmath $Z$}/2\mbox{\boldmath $Z$}\oplus
\mbox{\boldmath
$Z$}/2\mbox{\boldmath $Z$}$ and letting $\{\sigma, \tau\}$ denote a set of
generators,
$$
\rho(\sigma)=\left(
\begin{array}{ccc}
-1 & 0 & 0\\
0 & -1 & 0\\
0 & 0 & -1
\end{array}
\right),
\quad
\rho(\tau)=
\left(
\begin{array}{ccc}
0 & 1 & 0\\
1 & 0 & 0\\
0 & 0 & -1
\end{array}
\right).
$$
In particular, if $f$ is of type {\rm III}, then $X$ has only canonical
quotient singularity.
\end{description}
\end{th}

For the definition of types I, II and III, see Definition~\ref{df:type.deg}.

\begin{prob}\label{prob:sIII}\begin{em} Let $f:X\rightarrow {\cal D}$ a
projective log minimal degeneration of abelian or hyperelliptic
surfaces of type III. Applying the log minimal program to $f$ with
respect to $K_X$, we see that there exists a  projective bimeromorphic map
$\psi:X-\rightarrow X^s$ over
${\cal D}$ such that for the induced  projective degeneration
$f^s:X^s\rightarrow {\cal D}$, we have  
$K_{X^s}\sim_{\mbox{\boldmath
$Q$}}0$ and
$f^{s \ast}(0)=m\Theta ^s$ for some $m\in \mbox{\boldmath $N$}$,
where 
$\Theta ^s:=\psi_{\ast}\Theta$ and that $X^s$ has only canonical
singularity but the possible types of singularity of
$(X^s,\Theta ^s)$ may differ from the ones of $(X,\Theta)$. So
determination of the types of singularity of
$(X^s,\Theta ^s)$ remains to be done. 
\end{em}\end{prob}

Refining a canonical bundle formula in \cite{fujino-mori} by the Log Minimal Model Program, we obtain the following theorem as an
application of Theorem~\ref{th:main th.ab}.
\cite{fujino-mori}.

\begin{th}[Theorem~\ref{th:bouding the number of singular fibres}] 
The number of singular fibres of abelian fibred Calabi Yau threefolds over a smooth rational curves are bounded from above not depending on relative polarizations.

\end{th}

\begin{rk}

\begin{em}(1) The number of singular fibres and the bounding problem of Euler numbers are closely related to each other. \cite{hunt} asserts the boundedness of Euler numbers of fibred
Calabi-Yau threefolds, but unfortunately,
\cite{hunt} contains a several crucial gaps. In fact, one of the aims of this paper is to remedy the results in \cite{hunt}. For example, the crucial Lemma 4 in \cite{hunt} has a
counter-example as follows. Let 
$f:X\rightarrow {\cal D}$ be a projective connected morphism from a complex manifold $X$ onto a 1-dimensional complex disk ${\cal D}$. Assume that $f_{\circ}:=f|_{f^{-1}({\cal
D}^{\ast})}$ is a smooth family of abelian varieties, where ${\cal D}^{\ast}:={\cal D}\setminus\{0\}$. Assume moreover $f$ does not admit any sections. According to \cite{ueno3},
there exists a projective morphism $f^b:X^b\rightarrow {\cal D}$ from a complex manifold $X^b$ onto ${\cal D}$ such that $f_{\circ}^b:=f^b|_{f^{b -1}({\cal
D}^{\ast})}$ is a basic polarized bundle associated with $f_{\circ}$, that is, the pairs of period maps and the monodromies associated with $f_{\circ}$ and $f_{\circ}^b$ are
equivalent. But $f$ and $f^b$ are not bimeromorphically equivalent because $f^b$ admits a section while $f$ does not. By the same reason, 
the assertion in \S 3, Step 1 in \cite{hunt} saying that $\pi^{\prime}:X^{\prime}\rightarrow Y^{\prime}$ is birational to the pull back of 
$g^{\prime}:{\cal F}_{\Gamma^{\prime}}\rightarrow \Gamma^{\prime}\setminus {\cal D}^{\wedge}$ is incorrect. \\

(2) Moreover, fixing the degrees of direct image sheaves of relative dualizing sheaves as in \cite{hunt} gives no condition on the number of singular
fibres. For example, for any given integer
$g$, there exists a elliptic surface
$f_g:X_g\rightarrow B$ over a projective line
$B$ such that
$\deg f_{g\ast}{\cal O}(K_{X_g}/B)=0$ and
$f_g$ has
$2g+2$ singular fibres. These can be constructed by taking quotients of products of elliptic curves and hyperelliptic curves by involutions which are products of translations by
torsion point of order two and canonical involutions of hyperelliptic curves. These are not the counter-example of boundedness of Euler numbers of fibred Calabi-Yau
threefolds, but at least one has to care about the number of singular fibres which contribute to the Euler numbers. It seems that there is no argument like that in
\cite{hunt}. \\

(3) As for independence of boudedness on the relative polarizations, \cite{hunt} seems to be using the fact that abelian varieties defined over an algebraically closed field are
isogeneous to a principally polarized abelian varieties. The problem is to bound the degree of the neccesay base change but the argument \cite{hunt}, pp.150, Corollary does not seem
to be successful (It seems that $p_{\Gamma}$  in the proof is just the isomorphism). Instead, we used the Zarhin's trick in our argument.

\end{em}\end{rk}

In \S 2, we define the invariant $\delta_p$, and prove Riemann-Roch formula for divisors on singular 3-fold under some assumptions ( Proposition~\ref{pr:generalized
cm-formula} ) to prove Theorem~\ref{th:euler characteristic formula} and Corollary~\ref{co:chi0}. In \S 3, we classify type II and III log surfaces ( for the definition, see 
Definition~\ref{df:123} ) under certain typical assumption which are supposed to appear canonically as the components of the singular fibres of the degeneration of surfaces with Kodaira
dimension zero. In \S 4, we give a theory to calculate local fundamental groups seeing differents, which will be used to determine the singularity of the total spaces from the
information of log surfaces obtained in the previous section. In \S 6, we prove Theorem~\ref{th:main th.ab} by using the results in 
the previous sections. In \S 5, \S 7.1 and  \S 7.2, we give a systematic treatment of degenerations or fibrations of surfaces with Kodaira dimension zero including
Kodaira-Mori's canonical bundle formulae and in \S 7.3, we shall argue about the abelian case.

\begin{center}
Notation and Conventions
\end{center}
Let $X$ be a normal variety defined over an algebraically closed field $k$ (if the characteristic of $k$ is not zero, we assume the existence of a embedded resolution). 
An elements of $\mbox{Weil X}\otimes \mbox{\boldmath $Q$}$ is called a {\it {\bf Q}-divisor}. {\bf Q}-divisor $D$ has the unique irreducible decomposition 
$D=\sum_{\Gamma}(\mbox{mult}_{\Gamma}D)\Gamma$, where $\mbox{mult}_{\Gamma}D\in \mbox{\boldmath $Q$}$ and the summation is taken over all the prime divisors $\Gamma$ on $X$. 
{\bf Q}-divisor $\Delta$ is called a {\it {\bf Q}-boundary} if 
$\mbox{mult}_{\Gamma}\Delta\in [0,1]\cap \mbox{\boldmath $Q$}$ for any prime divisor $\Gamma$. 
{\bf Q}-divisor
$D$ is said to be {\it {\bf Q}-Cartier} if $rD\in \mbox{Div }X$ for some $r\in \mbox{\boldmath $Q$}$. $X$ is said to be {\bf Q}-Gorenstein if a canonical divisor $K_X$ is {\bf
Q}-Cartier. $X$ is said to be {\it {\bf Q}-factorial} if any Weil divisor on $X$ is {\bf Q}-Cartier. A pair
$(X,\Delta)$ which consists of a normal variety
$X$ and {\bf Q}-boundary
$\Delta$ on $X$ is called a normal {\it log variety}. For a normal log variety $(X,\Delta)$, a resolution $\mu:Y\rightarrow X$ is called a {\it log resolution} of $(X,\Delta)$ if each
component of the support of 
$\mu^{-1}_{\ast}\Delta+\sum_{i\in I}E_i$ are smooth and cross normally, where $\{E_i\}_{i\in I}$ is a set of all the exceptional divisors of $\mu$. Assume that $K_X+\Delta$ is {\bf
Q}-Cartier. The {\it log discrepancy} $a_l(E_i;X,\Delta)\in \mbox{\boldmath $Q$}$ of $E_i$ with respect to $(X,\Delta)$ is defined by 
$$a_l(E_i;X,\Delta):=\mbox{mult}_{E_i}(K_Y+\mu^{-1}_{\ast}\Delta+\sum_{i\in I}E_i-\mu^{\ast}(K_X+\Delta))\in \mbox{\boldmath $Q$}$$ and the {\it discrepancy} $a(E_i;X,\Delta)\in
\mbox{\boldmath $Q$}$ is defined by 
$
a(E_i;X,\Delta):=a_l(E_i;X,\Delta)+1.
$ 
The closure of $\mu(E_i)\subset X$ is called a {\it center} of $E_i$ at $X$ which is denoted by $\mbox{Center}_{X}(E_i)$.
The above definitions of discrepancies, a log discrepancies and centers are known to be well-defined and depend only on the rank one discrete valuation of the function field of $X$
associated with $E_i$'s. $E_i$'s are called a {\it exceptional divisors of the function field of $X$} and it has its meaning saying discrepancies, a log discrepancies and centers of
exceptional divisors of the function field of $X$. A normal log variety $(X,\Delta)$ is said to be {\it terminal} (resp. {\it canonical}, resp. {\it purely log terminal}) 
if $a(E_i;X,\Delta)>0$ (resp. $a(E_i;X,\Delta)\geq 0$, resp. $a_l(E_i;X,\Delta)>0$) for any log resolution $\mu$ and any $i\in I$. For some log resolution $\mu$, if
$a_l(E_i;X,\Delta)>0$ for any $i\in I$, $(X,\Delta)$ is said to be {\it log terminal}, moreover if the exceptional loci of $\mu$ is purely one
codimensional, $(X,\Delta)$ is said to be {\it divisorially log terminal}. We shall say that $X$ has only terminal (resp. canonical, resp. log terminal) singularity if $(X,0)$ is
terminal (resp. canonical, resp. log terminal) as usual (see
\cite{kollar/utah} or
\cite{shokurov} and see also
\cite{kawamata crep.},
\S 1 for the treatment in the complex analytic case ).

\vskip 5mm
In this paper, 
we shall use the following notation:
\begin{description}
\item[ $\nu:X^{\nu}\rightarrow X$ ]: The normalization of a scheme $X$.

\item[ $\mbox{\rm Diff}_{\Gamma^{\nu}}(\Delta)$ ]: $\mbox{\boldmath
$Q$}$-divisor which is called Shokurov's different satisfying 
$$
\nu^{\ast}(K_X+\Gamma+\Delta)=K_{\Gamma^{\nu}}+\mbox{Diff}_{\Gamma^{\nu}}(\Delta),
$$ where $\Gamma$ is a reduced divisor on a normal variety $X$ and $\Gamma+\Delta$ is a
$\mbox{\boldmath
$Q$}$-boundary on $X$ such that $K_X+\Gamma+\Delta$ is $\mbox{\boldmath $Q$}$-Cartier. (see
\cite{shokurov}, \S $3$, \cite{corti/utah}, \S $16$). 

\item[ $\Delta^Y$ ]: $\mbox{\boldmath $Q$}$-divisor on $Y$ satisfying 
$K_Y+\Delta^{Y}=f^{\ast}(K_X+\Delta)$, where $f:Y\rightarrow X$ is a birational morphism between normal
varieties and 
$\Delta$ is a $\mbox{\boldmath $Q$}$-boundary on $X$ such that $K_X+\Delta$ is $\mbox{\boldmath
$Q$}$-Cartier.

\item[ $\mbox{\rm ind}_p(D)$ ]: The smallest positive integer $r$ such that $rD$
is Cartier on the germ of $X$ at $p$, where $D$ is a $\mbox{\boldmath $Q$}$-Cartier $\mbox{\boldmath $Q$}$-divisor on a normal variety or a normal complex analytic space $X$.

\item[ $\mbox{\rm Ind}(D)$ ]: The smallest positive integer $r$ such that $rD\sim 0$, where $D$ is a $\mbox{\boldmath $Q$}$-Cartier $\mbox{\boldmath $Q$}$-divisor on a normal variety
$X$ such that $D\sim_{\mbox{\boldmath $Q$}}0$.

\item[ $\mbox{\rm Exc}f$ ]: Exceptional loci of a birational morphism $f:X\rightarrow Y$ between varieties $X$ and $Y$, that is, loci of points in $X$ in a neighbourhood of which $f$
is not isomorphic.

\item[ $\sim$ ]: Linear equivalence.

\item[ $\sim_{\mbox{\boldmath $Q$}}$ ]: {\boldmath $Q$}-linear equivalence.

\item[ $\lceil\Delta\rceil$ ]: Round up of a $\mbox{\boldmath $Q$}$-divisor $\Delta$.

\item[ $\lfloor \Delta \rfloor$ ]: Round down of a $\mbox{\boldmath $Q$}$-divisor $\Delta$.

\item[ $\{\Delta\}$ ]: Fractional part of the boundary $\Delta$.

\item[ $e_{\mbox{\rm top}}$ ]: Topological Euler characteristic.

\item[ $\rho(X/Y)$ ]: Relative Picard number of a normal $\mbox{\boldmath
$Q$}$-factorial variety $X$ over a variety $Y$.

\item[ $\Sigma_d$ ]: Hirzebruch surface of degree $d$.

\item[ $\mbox{\rm Card }{\cal S}$ ]: Cardinality of a set ${\cal S}$.

\end{description}

For a normal complete surface $S$ with at worst Du Val
singularities,  we shall write 
$$
\mbox{Sing S}=\sum_{{\cal T}}\nu({\cal T}){\cal T},
$$ where
$\nu({\cal T})$ denotes the number of singular points on $S$ of type ${\cal T}$. For a quasi projective complex surface $S$ with only quotient singularity, recall that {\it the orbifold
euler number} $e_{\mbox{\rm orb}}(S)\in \mbox{\boldmath $Q$}$ of $S$ is defined by 
$$
e_{\mbox{\rm orb}}(S):=e_{\mbox{\rm top}}(S)-\sum_{p\in S}(1-\frac{1}{ \mbox{\rm Card}\pi_{S,p} }),
$$
where $\pi_{S,p}$ denotes the local fundamental group of $S$ at $p\in S$ (see \cite{kawamata ab.}, page 233 or \cite{megyesi}, Definition 10.7).

\vskip 5mm
\noindent{\bf Acknowledgment.} The author would like to express his deep gratitude to his
thesis advisor, Prof. Yujiro Kawamata by whom he has been greatly influenced, to Prof. Shigefumi Mori for
calling his attention to the importance of the calculation of the Euler characteristic of the structure shaves of
singular fibres and kindly showing him the preliminary version of the preprint \cite{fujino-mori}, also to Dr. Osamu Fujino for showing him the expanded version, to Prof. Masayoshi
Miyanishi for pointing out to him that the existence of the case
$(1)$ in Proposition~\ref{pr:cltypeIII-2.2} is known, to Prof. Akira Fujiki for
providing Lemma~\ref{lm:universal} and helping him to work on the analytic spaces, to Prof. R.V.Gurjar for letting him know the existence of
\cite{dimca},  to Prof. Makoto Namba for kindly showing him the beautiful book \cite{namba}, which enabled him to improve the first draft of this paper and to Profs. Keiji Oguiso and
Yoshinori Namikawa for stimulating discussions.

\section{The Euler characteristic formula}
Firstly, let us recall the following result due to Crauder and Morrison.
\begin{pr}[\cite{crauder}, Proposition (A.1)]\label{pr:crauder-morrison's formula} Let $X$ be a smooth $3$-fold and let $D$ 
be a complete effective divisor on $X$ whose support has only simple normal crossing singularities.
Then the following holds.
$$
\chi({\cal O}_D)=\sum_{i}m_i\chi({\cal O}_{D_i})+{\frac {1}{6}}(D^3-\sum_{i}m_iD^3_i)+{\frac {1}{4}}(D^2-\sum_im_iD^2_i)K_X,
$$
where $D=\sum_{i}m_iD_i$ is the irreducible decomposition.
\end{pr}

Let $(X,p)$ be a germ of 3-dimensional terminal singularity at $p$ whose index $r$ is equal to or
greater than $2$. Take a Du Val element $S\in |-K_X|$ passing through $p$, where we say that $S\in
|-K_X|$ is a Du Val element, if $S$ is a reduced normal {\bf Q}-Cartier divisor on $X$ passing
through $p$ such that $S$ has a  Du Val singularity at $p$. The canonical cover $\pi:\tilde
X\rightarrow X$ with respect to $K_X$ induces a covering of Du Val singularities 
$\pi:\tilde S:=\pi^{-1}(S)\rightarrow S$. There is a coordinate system $x$, $y$ and $z$ of ${\mbox{\boldmath $C$}}^3$ 
which are semi-invariant under the action of the Galois group $\mbox{Gal }(\tilde S/S)$ such
that 
$\tilde p:=\pi^{-1}(p)\in \tilde S$ is analytically isomorphic to the germ of the origin of the hypersurface defined by a 
equation $f(x,y,z)=0$. Let $\sigma$ be a generator of $\mbox{Gal }(\tilde S/S)$ and let
$\zeta$ be a  primitive r-th root of unity. 
The actions of $\sigma$ are completely classified into the following $6$ types (see \cite{reid2}).
\begin{description}
\item[(1)]\ $\tilde p\in \tilde S$ is of type $A_{n-1}$ and $p\in S$ is of type $A_{rn-1}$ ($n\geq 1$). $f=xy+z^n$,$\sigma^{\ast}x=\zeta^ax$, 
$\sigma^{\ast}y=\zeta^{-a}y$ and $\sigma^{\ast}z=z$, where $(r,a)=1$.
\item[(2)]\ $\tilde p\in \tilde S$ is of type $A_{2n-2}$ and $p\in S$ is of type $D_{2n+1}$ ($n\geq2$). $r=4$, $f=x^2+y^2+z^{2n-1}$,
$\sigma^{\ast}x=\zeta x$, $\sigma^{\ast}y=\zeta^{3}y$ and $\sigma^{\ast}z=\zeta^{2}z$.
\item[(3)]\ $\tilde p\in \tilde S$ is of type $A_{2n-1}$ and $p\in S$ is of type $D_{n+2}$ ($n\geq2$). $r=2$, $f=x^2+y^2+z^{2n}$,
$\sigma^{\ast}x=x$, $\sigma^{\ast}y=-y$ and $\sigma^{\ast}z=-z$.
\item[(4)]\ $\tilde p\in \tilde S$ is of type $D_{4}$ and $p\in S$ is of type $E_{6}$. $r=3$, $f=x^2+y^3+z^{3}$,
$\sigma^{\ast}x=x$, $\sigma^{\ast}y=\zeta y$ and $\sigma^{\ast}z=\zeta^{2}z$.
\item[(5)]\ $\tilde p\in \tilde S$ is of type $D_{n+1}$ and $p\in S$ is of type $D_{2n}$. $r=2$, $f=x^2+y^2z+z^{n}$,
$\sigma^{\ast}x=-x$, $\sigma^{\ast}y=-y$ and $\sigma^{\ast}z=z$.
\item[(6)]\ $\tilde p\in \tilde S$ is of type $E_{6}$ and $p\in S$ is of type $E_{7}$. $r=2$, $f=x^2+y^3+z^{4}$,
$\sigma^{\ast}x=-x$, $\sigma^{\ast}y=y$ and $\sigma^{\ast}z=-z$.
\end{description}

\begin{df}\begin{em}
For $p\in S\subset X$ as above, we define the rational number $c_p(X,S)\in \mbox{\boldmath $Q$}$ as
follows: 
$$
c_p(X,S):=\left\{
\begin{array}{ll}
0 & \mbox{ Case $p\in X$ Gorenstein, }\\
n\{r- (1/r)\} & \mbox{ Case (1),}\\
3(2n+3)/4 & \mbox{ Case (2),}\\
3 & \mbox{ Case (3),}\\
16/3 & \mbox{ Case (4),}\\
3n/2 & \mbox{ Case (5),}\\
9/2 & \mbox{ Case (6).}
\end{array}
\right.
$$
\end{em}\end{df}

\begin{df}\begin{em}
Let $p\in S\subset X$ be as above. we define the rational number $\delta_p(X,S)\in \mbox{\boldmath
$Q$}$ as follows:
$$
\delta_p(X,S):=e_p(S)-{\frac {1}{o_p(S)}}-c_p(X,S)\in \mbox{\boldmath
$Q$},
$$
where $e_p(S)$ is the Euler number of the inverse image of $p$ by the morphism induced by the minimal resolution and 
$o_p(S)$ is the order of the local fundamental group of $S$ at $p$.
\end{em}\end{df}

If the index of $X$ at $p$ is equal to or greater than $2$, we obtain the following table.
\begin{center}
Table I\\
\begin{tabular}{l|l|c|c|c}\hline
 & $e_p(S)$ & $o_p(S)$ & $c_p(X,S)$ & $\delta_p(X,S)$ \\ \hline \hline
(1) & $rn$ & $rn$ & $n\{r- (1/r)\}$ & $(n^2-1)/rn$\\ \hline
(2)  & $2n+2$ & $8n-4$ & $3(2n+3)/4$ & $n(n-1)/(2n-1)$\\ \hline
(3)  & $n+3$ & $4n$  & $3$ & $(4n^2-1)/4n$\\ \hline
(4)  & $7$ & $24$ & $16/3$ &  $13/8$\\ \hline
(5)  & $2n+1$ & $8(n-1)$ & $3n/2$ & $(4n^2+4n-9)/8(n-1)$\\ \hline
(6)  & $8$ & $48$ & $9/2$ & $167/48$\\
\hline
\end{tabular}
\end{center}
\vskip 0.3in

\begin{pr}\label{pr:delta} $\delta_p(X,S)\geq 0$. $\delta_p(X,S)=0$ if and only if $(X,S)$ has only singularity of
type
$V_1(r;a,-a,1)$ at $p$, where $(r,a)=1$.
\end{pr}
{\it Proof.}\ If $p\in X$ is Gorenstein, it is easy to see that $\delta_p(X,S)=e_p(S)-1/o_p(S)\geq 0$ 
and that $\delta_p(X,S)=0$ 
if and only if $X$ and $S$ is smooth at $p$. Assume that the index of $p\in X$ is equal to or greater than $2$. 
If we have $\delta_p(X,S)=0$, we infer that $\tilde p\in \tilde S$ is smooth (hence $\tilde p\in \tilde X$) from Table I. 
Thus we get the assertion.
\hfill \bsquare
\vskip 5mm
We give a proof of the following Reid's Riemann-Roch formula, which seems to be more clear than the one in \cite{reid2}, (9.2) to see the last statement in the theorem which is a
crucial point for our subsequent argument.

\begin{th}[\cite{reid2}, Theorem 9.1 (I)]\label{th:rrr} Let $X$ be a projective surface with 
at worst Du Val singularities and 
let $D$ be a Weil divisor on $X$. Then 
$$
\chi({\cal O}_X(D))=\chi({\cal O}_X)+{\frac {1}{2}}D(D-K_X)+\sum_{p\in X}c_p(D),
$$
where $c_p(D)$ is the rational number which depends only on the local analytic type of $p\in X$ and ${\cal O}_X(D)$.
\end{th}

{\it Proof.} Let $\mu:Y\rightarrow X$ be the minimal resolution of $X$. Put $\Gamma:=\lceil \mu^{\ast}D\rceil-\lfloor\mu^{\ast}D\rfloor$. 
Then there exists the following exact sequence:
$$
0\rightarrow {\cal O}_Y(\lfloor\mu^{\ast}D\rfloor)\rightarrow{\cal O}_Y(\lceil \mu^{\ast}D\rceil)\rightarrow
{\cal O}_{\Gamma}(\lceil \mu^{\ast}D\rceil)\rightarrow 0.
$$
Since we have $\mu_{\ast}{\cal O}_Y(\lfloor\mu^{\ast}D\rfloor)\simeq {\cal O}_Y(D)$ by \cite{sakai}, Theorem 2.1 and $R^i\mu_{\ast}{\cal O}_Y(\lceil \mu^{\ast}D\rceil)=0$ 
for $i>0$ by \cite{sakai}, Theorem 2.2, we obtain the following exact sequence:
$$
0\rightarrow {\cal O}_X(D)\rightarrow \mu_{\ast}{\cal O}_Y(\lceil \mu^{\ast}D\rceil)\rightarrow
\mu_{\ast}{\cal O}_{\Gamma}(\lceil \mu^{\ast}D\rceil)\rightarrow R^1\mu_{\ast}{\cal O}_Y(\lfloor\mu^{\ast}D\rfloor)\rightarrow 0
$$
to get  
\begin{eqnarray*}
\chi({\cal O}_X(D))&=&\chi(\mu_{\ast}{\cal O}_Y(\lceil \mu^{\ast}D\rceil))\\
&&-{\mbox {length Ker}}\{\mu_{\ast}{\cal O}_{\Gamma}(\lceil \mu^{\ast}D\rceil)\rightarrow R^1\mu_{\ast}{\cal O}_Y(\lfloor\mu^{\ast}D\rfloor)\}.
\end{eqnarray*}
Put $\Delta:=\lceil \mu^{\ast}D\rceil-\mu^{\ast}D$. Then $\chi(\mu_{\ast}{\cal O}_Y(\lceil \mu^{\ast}D\rceil))$ can be written as follows:
\begin{eqnarray*}
\chi(\mu_{\ast}{\cal O}_Y(\lceil \mu^{\ast}D\rceil))&=&\chi({\cal O}_Y(\lceil \mu^{\ast}D\rceil))\\
&=&\chi({\cal O}_Y)+{\frac {1}{2}}\lceil\mu^{\ast}D\rceil(\lceil\mu^{\ast}D\rceil-K_Y)\\
&=&\chi({\cal O}_X)+{\frac {1}{2}}(\mu^{\ast}D+\Delta)(\mu^{\ast}D+\Delta-\mu^{\ast}K_X)\\
&=&\chi({\cal O}_X)+{\frac {1}{2}}(D^2-DK_X)+{\frac {1}{2}}\Delta^2.
\end{eqnarray*}
Putting $c_p(D):=(1/2)(\Delta|_{\mu^{-1}(p)})^2-\mbox{length }{\cal T}(D)_p$, where 
$${\cal T}(D):={\mbox {Ker}}\{\mu_{\ast}{\cal O}_{\Gamma}(\lceil \mu^{\ast}D\rceil)\rightarrow R^1\mu_{\ast}{\cal O}_Y(\lfloor\mu^{\ast}D\rfloor)\},
$$
we get
$$
\chi({\cal O}_X(D))=\chi({\cal O}_X)+{\frac {1}{2}}D(D-K_X)+\sum_{p\in X}c_p(D).
$$
As for the last assertion, 
for any two Weil divisors $D_1$ and $D_2$ such
that $D_1-D_2$ is Cartier at $p\in X$, since $\Delta$ for $D=D_1$ and $D=D_2$ is the same and we have ${\cal T}(D_1)_p={\cal T}(D_2)_p\otimes_{{\cal O}_{X,p}}{\cal O}_X(D_1-D_2)_p$, 
we infer that 
$c_p(D_1)=c_p(D_2)$.  Thus we get the assertion.\hfill\bsquare

\begin{lm}\label{lm:canonical} Let $X$ be a germ of reduced irreducible normal $\mbox{\boldmath
$Q$}$-Gorenstein analytic spaces and $D$ be a non zero $\mbox{\boldmath $Q$}$-boundary on $X$.
Assume that
$(X,D)$ is canonical and that the center on $X$ of any divisor $E_j$ with discrepancy zero is
contained in the support of $D$, then $X$ has only terminal singularities.   
\end{lm}

{\it Proof.} Let $\mu:Y\rightarrow X$ be a Hironaka resolution of $X$ and write
$K_Y=\mu^{\ast}K_X+\sum_{j\in J}a_jE_j$ and $\mu^{\ast}D=\mu^{-1}_{\ast}D+\sum_{j\in J}\nu_jE_j$,
where $\{E_j|j\in J\}$ is the set of all the $\mu$-exceptional divisors and $a_j$, $\nu_j$ are
non-negative rational numbers for any $j\in J$. By the assumption, we have
$a_j\geq \nu_j$ for any $j\in J$. By the choice of our resolution, $\nu_j=0$ implies that $E_j$ is
obtained by blowing up the center which is not contained in the support of the weak transform of
some multiple of $D$. Thus we get the assertion. \hfill \bsquare

\begin{pr}\label{pr:generalized cm-formula}
Let $X$ be a normal {\bf Q}-Gorenstein $3$-fold and $D$ be an effective complete Cartier divisor on 
$X$ such that the log $3$-fold 
 $(X,D_{{\rm red}})$ is divisorially log terminal and that $X$ is smooth outside the support of $D$. Assume
that
$K_X+D_{{\rm red}}$ is Cartier and each irreducible component of $D$ is algebraic and {\bf Q}-Cartier. 
Then the following formula holds $:$
\begin{eqnarray*}
\chi({\cal O}_D)&=&\sum_{i}m_i\chi({\cal O}_{D_i})+{\frac {1}{6}}(D^3-\sum_{i}m_iD^3_i)+{\frac {1}{4}}(D^2-\sum_im_iD^2_i)K_X\\
&&-{\frac {1}{12}}\sum_im_i\sum_{p\in D_i^o}c_p(X,D_i),
\end{eqnarray*}
where $D=\sum_{i}m_iD_i$ is the irreducible decomposition and $D_i^o:=D_i\setminus \cup_{j\neq i}D_j$.
\end{pr}
{\it Proof.}\  We calculate the contribution of singularities to the formula
in Proposition~\ref{pr:crauder-morrison's formula} using Reid's Riemann-Roch for surfaces with Du
Val singularities and examining the original proof of
Proposition~\ref{pr:crauder-morrison's formula}. By the assumptions, we infer that  
if $X$ is not smooth or $X$ is smooth but $D_{{\rm red}}$ is not a normal crossing divisor at $p\in X$, then $p\in D^o_i$ for some $i$, $X$ has only terminal singularity by
Lemma~\ref{lm:canonical} and $D_i$ has at worst Du Val singularity at $p\in D_i\subset X$ since $(X,D_i)$ is
canonical for any
$i$. We note that if
a singularity $p\in X$ is Gorenstein then this singularity does not contribute to the
Riemann-Roch. Assume that $p\in D^o_i$ for some $i$ and that $r:=\mbox{ind }_pK_X\geq 2$. Take the canonical cover 
$\pi:\tilde X\rightarrow X$ locally at $p$ (here we used the same notation $X$ for an open neighbourhood of $p\in X$). 
Putting $\tilde p:=\pi^{-1}(p)$ and $\tilde D_i:=\pi^{-1}(D_i)$, let $\vartheta_i\in {\cal O}_{\tilde X}$ be a 
defining equation of $\tilde D_i$ and let $\sigma$ be the generator of the covering action. For any $l\in {\mbox{\boldmath N}}$, 
there is the following exact sequence locally at $\tilde p$ 
which is compatible with the action of $\sigma$. 
$$
0\rightarrow \vartheta_i^l{\cal O}_{\tilde D_i}\rightarrow {\cal O}_{(l+1)\tilde D_i}\rightarrow {\cal O}_{l\tilde D_i}\rightarrow 0.
$$
Since $H^1(<\sigma >,\vartheta_i^l{\cal O}_{\tilde D_i})=0$, taking the invariant part of the above exact sequence, 
we obtain the following exact sequence locally at $p$: 
$$
0\rightarrow (\vartheta_i^l{\cal O}_{\tilde D_i})^{<\sigma>}\rightarrow {\cal O}_{(l+1)D_i}\rightarrow {\cal O}_{l D_i}\rightarrow 0,
$$
where $(\vartheta_i^l{\cal O}_{\tilde D_i})^{<\sigma>}$ is the $\sigma$-invariant part of $\vartheta_i^l{\cal O}_{\tilde D_i}$. We note that in fact, 
$(\vartheta_i^l{\cal O}_{\tilde D_i})^{<\sigma>}$ is a restriction of the divisorial sheaf 
${\cal F}_l:=\mbox{Ker }\{{\cal O}_{(l+1)D_i}\rightarrow {\cal O}_{l D_i}\}$ on $D_i$.
Define $a\in {\mbox{\boldmath $N$}}$  so that $\sigma^{\ast}\vartheta_i=\zeta^{-a}\vartheta_i$
where $\zeta$ is a primitive $r$-th root of unity and note that $(a,r)=1$ since 
$K_X+D_{\mbox{red}}$ is Cartier. we note that for any $\varphi\in {\cal O}_{\tilde D_i}$, 
$\vartheta_i^l\varphi\in \vartheta_i^l{\cal O}_{\tilde D_i}$ is $\sigma$-invariant if and only if $\sigma^{\ast}\varphi=\zeta^{al}\varphi$, so we have 
${\cal F}_l\simeq \{\varphi\in \pi_{\ast}{\cal O}_{\tilde D_i}; \sigma^{\ast}\varphi=\zeta^{al}\varphi\}$ locally at $p\in D_i^o$. 
Assuming that $\tilde p\in \tilde D_i$ is smooth (the case (1), $n=1$) for simplicity, we calculate the summation of 
the contribution to the Riemann-Roch for ${\cal F}_l$  ($1\leq l\leq m_i-1$), that is, $\sum_{l=1}^{m_i-1}c_p({\cal F}_l)$ as 
follows. Note that there exists a natural number $d_i\in {\mbox{\boldmath $N$}}$ such that 
$m_i=rd_i$ since $m_iD_i$ is Cartier at $p$. 
\begin{eqnarray*}
\sum_{l=1}^{m_i-1}c_p({\cal F}_l)&=&\sum_{l=1}^{m_i-1}-{\frac {\overline{al}(r-\overline{al})}{2r}}=
-d_i\sum_{l=1}^{r-1}{\frac {\overline{al}(r-\overline{al})}{2r}}=-d_i\sum_{l=1}^{r-1}{\frac {l(r-l)}{2r}}\\
&=&-m_i{\frac {r^2-1}{12r}},
\end{eqnarray*}
where $\overline{al}\in {\mbox{\boldmath $Z$}}$ is the unique integer such that 
$\overline{al}\equiv al$ (mod $r$) and $0\leq \overline{al}\leq r-1$. 
The other cases can be treated similarly using the {\bf Q}-smoothing method as explained in \cite{reid2}. 
\hfill\bsquare
\vskip 5mm
{\it Proof of Theorem~\ref{th:euler characteristic formula}.}\ Since $f$ is flat, we have 
$\chi({\cal O}_{X_t})=\chi({\cal O}_{X_0})$ for $t\in {\cal D}^{\ast}$, where $X_0:=f^{\ast}(0)$. 
Using Proposition~\ref{pr:generalized cm-formula} and 
the assumption that $K_X+\Theta\sim_{ \mbox{\boldmath $Q$} }0$, we obtain
\begin{eqnarray}
e_{\mbox{top}}(X_t)&=&12\chi({\cal O}_{X_t})\nonumber \\
&=&\sum_im_i(12\chi( {\cal O}_{\Theta_i})-2\Theta_i^3+3\sum_j\Theta_i^2\Theta_j-
\sum_{p\in \Theta_i\setminus \Delta_i}c_p(X,\Theta_i))\label{eqn:cm2}.
\end{eqnarray}
On the other hand, since for any $i$, $\Delta_i$ is either $0$, disjoint union of smooth elliptic curves 
or a cycle of rational curves ( see Lemma~\ref{lm:ind 1 ls} ), we have
\begin{eqnarray*}
K_{\Theta_i}^2+e_{\mbox{\rm top}}(\Delta_i)=\Delta_i^2+e_{\mbox{\rm top}}(\Delta_i)
&=&\sum_{j\neq i}\Theta_i\Theta_j^2+3\sum_{j,k\neq i {\rm{ and }}
j<k}\Theta_i\Theta_j\Theta_k\\ 
&=&\sum_{j\neq i}\Theta_i\Theta_j^2+{\frac
{3}{2}}\sum_{i\neq j,k{\rm{ and }}j\neq k}\Theta_i\Theta_j\Theta_k \\
&=&2\Theta_i^3-{\frac {1}{2}}\sum_j\Theta_i\Theta_j^2-3\sum_j\Theta_i^2\Theta_j+{\frac {3}{2}}\sum_{j,k}\Theta_i\Theta_j\Theta_k.
\end{eqnarray*}
Therefore, we have
\begin{equation}
\sum_im_i(K_{\Theta_i}^2+e_{\mbox{\rm
 top}}(\Delta_i))=\sum_im_i(2\Theta_i^3-3\sum_j\Theta_i^2\Theta_j)\label{eqn:cm8},
\end{equation}
since $\sum_im_i\Theta_i\sim 0$. From (\ref{eqn:cm2}) and (\ref{eqn:cm8}), we obtain 
$$
e_{\mbox{\rm top}}(X_t)=
\sum_im_i(12\chi({\cal O}_{\Theta_i})-K_{\Theta_i}^2-e_{\mbox{\rm top}}(\Delta_i)-\sum_{p\in
\Theta_i\setminus \Delta_i}c_p(X,\Theta_i)).
$$
For any $i$, let $\Theta_i^{\prime}\rightarrow \Theta_i$ be the minimal resolution of the 
singularities of $\Theta_i$. Since we have
$$
e_{\mbox{\rm
top}}(\Theta_i^{\prime})=e_{\mbox{\rm top}}(\Theta_i)+\sum_{p\in\Theta_i}(e_p(\Theta_i)-1)
$$
and $\Theta_i$ has only Du Val singularities, we have
$$
e_{\mbox{\rm top}}(\Theta_i)=12\chi({\cal
O}_{\Theta_i})-K_{\Theta_i}^2-\sum_{p\in\Theta_i}(e_p(\Theta_i)-1)
$$
by Noether's equality.
Thus we obtain
\begin{eqnarray*}
e_{\mbox{\rm top}}(X_t)&=&\sum_im_i\{e_{\mbox{\rm top}}(\Theta_i\setminus \Delta_i)+
\sum_{p\in \Theta_i\setminus \Delta_i}(e_p(\Theta_i)-1-c_p(X,\Theta_i))\}\\
&=&\sum m_i(e_{\mbox{\rm orb}}(\Theta_i\setminus \Delta_i)+\sum_{p\in \Theta_i\setminus
\Delta_i}\delta_p(X,\Theta_i)).
\end{eqnarray*}
\hfill \bsquare
\vskip 5mm
{\it Proof of Corollary~\ref{co:chi0}.}\ Since we have $e_{\mbox{\rm orb}}(\Theta_i\setminus
\Delta_i)
\geq 0$ for any $i$ by Miyaoka's inequality (see \cite{miyaoka}, Theorem 1.1 or \cite{megyesi}, Theorem 10.14), 
we obtain the assertion from Proposition~\ref{pr:delta}.
\hfill \bsquare

\section{Structures of log surfaces with a standard boundary whose log canonical
divisors are numerically trivial}
\subsection{${\cal S}$-extractions and ${\cal S}$-elementary transformations}
\begin{df}[\cite{shokurov-cp}]\label{df:standard}\begin{em} A normal log variety $(X,\Delta)$ defined over an algebraically 
closed field $k$ is said to be {\it a log variety with a standard boundary}, if 
${\mbox{mult}}_{\Gamma}\Delta\in {\cal S}:=\{(b-1)/b|b\in \mbox{\boldmath$N$}\cup \{\infty\}\}$ 
for any prime divisor $\Gamma$.
\end{em}\end{df}
\begin{df}\begin{em} Let $(S,\Gamma+\Delta)$ be a normal log surface with a $\mbox{\boldmath $Q$}$-boundary
$\Gamma+\Delta$, where $\Gamma$ is a prime divisor on
$S$. For a point $p\in \Gamma$, we define a
rational number $m_p(\Gamma^{\nu};\Delta)\in \mbox{\boldmath $Q$}$ as follows.
$$
m_p(\Gamma^{\nu};\Delta):=\mbox{mult}_p\mbox{Diff}_{\Gamma^{\nu}}(\Delta).
$$
\end{em}\end{df}

Let $(S,\Delta)$ be a normal log surface with a standard boundary defined over an 
algebraically closed field $k$. Assume that $(S,\Delta)$ is 
log canonical and $\lfloor\Delta\rfloor\neq 0$. 
Take $p\in \lfloor\Delta\rfloor$ and let $\Gamma$ be an irreducible component of the germ of 
$\lfloor\Delta\rfloor$ passing through $p$. Then 
$m_p(\Gamma;\Delta-\Gamma)$ can be written as follows: 
$$
m_p(\Gamma;\Delta-\Gamma)=
\frac{n-1}{n}+\sum_{b=2}^{\infty}\frac{b-1}{b} \frac{k_b}{n},
$$  
where $n\in \mbox{\boldmath $N$}$ and $k_b$ is a non-negative integer for all $b$ and $k_b=0$
except for finite number of $b$. We note that this is well known in the case of characteristic
$0$ but in fact this is just a conclusion of the intersection theory, so we don't have to worry
about the characteristic of $k$ about this matter. Since $(S,\Delta)$ is log canonical, we have
$$
\frac{n-1}{n}+\sum_{b=2}^{\infty}\frac{b-1}{b} \frac{k_b}{n}\leq 1.
$$  
Noting that $(b-1)/b\geq 1/2$, we get $\sum_{b=2}^{\infty}k_b\leq 2$. Therefore, 
one of the following three cases ($1$), ($2$) and ($3$) occurs. 
\begin{description}
\item[($1$)] $k_b=0$ for all $b$,\par
\item[($2$)] There exists $b^{\prime}$ such that $k_{b^{\prime}}=1$ and $k_b=0$ if  
$b\neq b^{\prime}$,\par
\item[($3$)] $k_2=2$ and $k_b=0$ if $b\neq 2$.
\end{description}
So we have
$$
m_p(\Gamma;\Delta-\Gamma)=\left\{
\begin{array}{ll}
(n-1)/n & \mbox{Case (1),}\\
(bn-1)/bn & \mbox{Case (2),}\\
1 & \mbox{Case (3)}.
\end{array}
\right.
$$
Let $\rho:T\rightarrow S$ be a birational morphism from a normal surface $T$ which is the blow up 
of $p$ if $S$ is smooth at $p$ or is obtained from the minimal resolution of $p\in S$ by contracting 
all the exceptional divisors that do not intersect the strict transform of $\Gamma$ 
if $S$ is singular at $p$. We call this $\rho$ an {\it ${\cal S}$-extraction}. 
\begin{lm}\label{lm:multE}
$(T,\Delta^T)$ is a log surface with a standard boundary. 
\end{lm}
{\it Proof.}\
Let $E$ be the exceptional prime divisor for $\rho$. We only have to show 
$m:=\mbox{mult}_{E}\Delta^T\in {\cal S}$. By definitions, we have 
$$
K_{\Gamma^{\prime}}+(\Delta^{\prime}-\Gamma^{\prime})|_{\Gamma^{\prime}}+mE|_{\Gamma^{\prime}}
=\rho^{\ast}(K_{\Gamma}+{\mbox{Diff}}_{\Gamma}(\Delta-\Gamma)),
$$
where $\Delta^{\prime}$ and $\Gamma^{\prime}$ are 
the strict transform of $\Delta$ and $\Gamma$ respectively, hence 
$$
(\Delta^{\prime}-\Gamma^{\prime},\Gamma^{\prime})_{p^{\prime}}+m=m_p(\Gamma;\Delta-\Gamma),
$$
where $p^{\prime}:=\rho^{-1}(p)\cap \Gamma^{\prime}$ and 
$(\Delta^{\prime}-\Gamma^{\prime},\Gamma^{\prime})_{p^{\prime}}$ denotes the local intersection 
number of $\Delta^{\prime}-\Gamma^{\prime}$ and $\Gamma^{\prime}$ at $p^{\prime}$. 
If we are in the cases ($1$) or ($2$), we see that 
$(\Delta^{\prime}-\Gamma^{\prime},\Gamma^{\prime})_{p^{\prime}}=0$, hence 
$m=m_p(\Gamma;\Delta-\Gamma)\in {\cal S}$. Assume that we are in the case ($3$). Since 
$(\Delta^{\prime}-\Gamma^{\prime},\Gamma^{\prime})_{p^{\prime}}\in (1/2)\mbox{\boldmath$Z$}$, 
 $m\geq 0$ and $m=1-(\Delta^{\prime}-\Gamma^{\prime},\Gamma^{\prime})_{p^{\prime}}$, 
we have $m=0$, $1/2$ or $1$. Thus we get the assertion. 
\hfill \bsquare

Assume that there exists a proper surjective connected $k$-morphism 
$\varphi:S\rightarrow B$ to a smooth curve $B$ defined over $k$ such that $\rho(S/B)=1$ and 
that the support of $\lfloor\Delta\rfloor$ is not entirely contained in a fibre of 
$\varphi$ and that $K_S+\Delta$ is relatively numerically trivial over $B$. 
Let $\hat \rho:T\rightarrow \hat S$ be a birational morphism to 
a normal surface $\hat S$ obtained by contracting the other component of the fibre of 
$\varphi\circ\rho$ which contains $E$ and put $\hat\Delta:=\hat\rho_{\ast}\Delta^T$. Then 
$(\hat S,\hat\Delta)$ with the induced morphism $\hat\varphi:\hat S\rightarrow B$ satisfies the 
same conditions as $(S,\Delta)$ and $\varphi:S\rightarrow B$. We call this transformation a 
{\it ${\cal S}$-elementary transformation}.
\begin{df}\begin{em} For $t\in B$, we define the $\mbox{\boldmath $Q$}$-boundaries 
$\Delta^{+}_{\varphi}(t)$ and $\Delta^{-}_{\varphi}(t)$ on $S$ as follows. 
$$
\begin{array}{rcl}
\Delta^{+}_{\varphi}(t)&:=&\Delta+(1-{\mbox{mult}}_{C_{\varphi}(t)}\Delta)C_{\varphi}(t),\\
\Delta^{-}_{\varphi}(t)&:=&\Delta-({\mbox{mult}}_{C_{\varphi}(t)}\Delta)C_{\varphi}(t),
\end{array}
$$
where $C_{\varphi}(t):=\varphi^{\ast}(t)_{\mbox{red}}$.
\end{em}\end{df}
\begin{lm}\label{lm:makelc}
Assume that $(S,\Delta^{+}_{\varphi}(t))$ is not log canonical over $t\in B$. 
Then one of the followings $(a)$ or 
$(b)$ holds.
\begin{description}
\item[($a$)] $\varphi$ is smooth over $t\in B$ and $\Delta^{-}_{\varphi}(t)=\lfloor\Delta^{-}_{\varphi}(t)\rfloor$.  
$\Delta^{-}_{\varphi}(t)$ is smooth in the neighbourhood of $C_{\varphi}(t)$ and $\Delta^{-}_{\varphi}(t)\cdot C_{\varphi}(t)=2p$, 
where $\Delta^{-}_{\varphi}(t)\cdot C_{\varphi}(t)$ denotes the intersection cycle. \par
\item[($b$)] There exists a sequence of ${\cal S}$-elementary transformations over $t\in B$ 
such that for the resulting 
log surface $(\hat S,\hat\Delta)$, $(\hat S,\Delta^+(\hat\varphi;t))$ is log canonical 
in the neighbourhood of the fibre over $t\in B$. 
\end{description}
\end{lm}
{\it Proof.}\ Assume that $(S,\Delta^{+}_{\varphi}(t))$ is not log canonical at $p_0\in C_{\varphi}(t)$. 
From 
$$\sum_{p\in C_{\varphi}(t)}m_{p}(C_{\varphi}(t);\Delta^{-}_{\varphi}(t))=2
\mbox{ and }
m_{p_0}(C_{\varphi}(t);\Delta^{-}_{\varphi}(t))>1,
$$ we deduce that
$m_{p}(C_{\varphi}(t);\Delta^{-}_{\varphi}(t))<1$ for any $p\neq p_0$, that is,
$(S,\Delta^{+}_{\varphi}(t))$ is purely log terminal at $p\in C_{\varphi}(t)$ if $p\neq p_0$. In
particular,  we see that
$\lfloor\Delta^{-}_{\varphi}(t)\rfloor\cap C_{\varphi}(t)=\{p_0\}$. Let
$\rho:T\rightarrow S$ be the ${\cal S}$-extraction at 
$p_0$ 
and let $E$ be the exceptional divisor for $\rho$. Put $\psi:=\varphi\circ \rho$. Then we have
$$
2=(\Delta^{-}_{\varphi}(t),\psi^{\ast}(t))\geq (\Delta^{-}_{\varphi}(t)^{\prime},C_{\varphi}(t)^{\prime})+(\Delta^{-}_{\varphi}(t)^{\prime},E)\geq 
(\Delta^{\prime}_-,C_{\varphi}(t)^{\prime})+1,
$$
where $\Delta^{\prime}_-$ and $C_{\varphi}(t)^{\prime}$ are the strict 
transforms of $\Delta^{-}_{\varphi}(t)$ and $C_{\varphi}(t)$ respectively. Let $p_0^{\prime}\in T$ be as in the proof of 
Lemma~\ref{lm:multE}. Suppose that $p_0^{\prime}\in C_{\varphi}(t)^{\prime}$. Since 
$(\Delta^{-}_{\varphi}(t)^{\prime},C_{\varphi}(t)^{\prime})_{p_0^{\prime}}\geq 1$, we have 
$(\Delta^{-}_{\varphi}(t)^{\prime},C_{\varphi}(t)^{\prime})=(\Delta^{-}_{\varphi}(t)^{\prime},E)=1$ and $\psi^{\ast}(t)=C_{\varphi}(t)^{\prime}+E$. 
Noting that $C_{\varphi}(t)^{\prime}\cap E=\{p_0^{\prime}\}$ and $T$ is smooth at $p_0^{\prime}$, 
we see that $T$ is smooth in the neighbourhood of the fibre over $t\in B$ and that 
$C_{\varphi}(t)^{\prime}$ and $E$ are ($-1$)-curves with $(C_{\varphi}(t)^{\prime},E)=1$, hence we are in the case ($a$). 
In the case that $p_0^{\prime}\notin C_{\varphi}(t)^{\prime}$, by a ${\cal S}$-elementary transformation, 
we may assume that $S$ is smooth at $p_0$ and $(C_{\varphi}(t), \lfloor\Delta^{-}_{\varphi}(t)\rfloor)_{p_0}=1$ at first. 
Moreover we may assume that $(\lfloor\Delta^{-}_{\varphi}(t)\rfloor, \{\Delta^{-}_{\varphi}(t)\})_{p_0}=0$ after operating 
${\cal S}$-elementary transformations. Thus we are in the case ($b$).
\hfill \bsquare
\begin{lm}\label{lm:dynkin}
Assume that $(S,\Delta^{+}_{\varphi}(t))$ is log canonical over $t\in B$. Then one of the followings $(${\em I}$)$ 
or $(${\em II}$)$ holds.
\begin{description}
\item[ (I) ] $\deg\lfloor \mbox{\em Diff}_{C_{\varphi}(t)}(\Delta^{-}_{\varphi}(t))\rfloor=1$, then
the dual graph of 
$\mbox{\em Supp }\mu^{\ast}\varphi^{\ast}(t)\cup \mbox{\em Supp }\mu^{-1}_{\ast} \Delta^{-}_{\varphi}(t)$ is one of the 
three types $(${\em I}-$1)_b$, $(${\em I}-$2)_b$ or $(${\em I}-$3)_b$ as below, 
where $\mu:M\rightarrow S$ is 
the minimal resolution of the singularities of $S$ over $t\in B$.

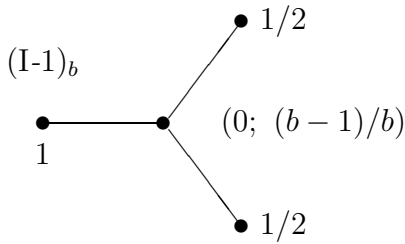
\unitlength 0.8mm
\begin{picture}(100,30)(-10,0)
\put(0,10){\makebox(0,0){$(${\em I}-$1)_b$}}
\put(0,0){\circle*{2}}\put(0,-5){\makebox(0,0){$1$}}
\put(0,0){\line(1,0){20}}
\put(20,0){\circle*{2}}
\put(40,0){\makebox(10,0){$(0;\ (b-1)/b)$}}
\put(21,1){\line(3,4){11.5}}
\put(21,-1){\line(3,-4){11.5}}
\put(33,17){\circle*{2}}\put(40,17){\makebox(0,0){$1/2$}}
\put(33,-17){\circle*{2}}\put(40,-17){\makebox(0,0){$1/2$}}
\end{picture}

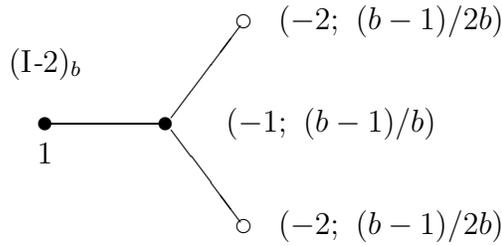
\begin{picture}(100,30)(-90,-35)
\put(0,10){\makebox(0,0){$(${\em I}-$2)_b$}}
\put(0,0){\circle*{2}}\put(0,-5){\makebox(0,0){$1$}}
\put(0,0){\line(1,0){20}}
\put(20,0){\circle*{2}}
\put(40,0){\makebox(15,0){$(-1;\ (b-1)/b)$}}
\put(21,1){\line(3,4){11.5}}
\put(21,-1){\line(3,-4){11.5}}
\put(33,17){\circle{2}}\put(50,17){\makebox(15,0){$(-2;\ (b-1)/2b)$}}
\put(33,-17){\circle{2}}\put(50,-17){\makebox(15,0){$(-2;\ (b-1)/2b)$}}
\end{picture}

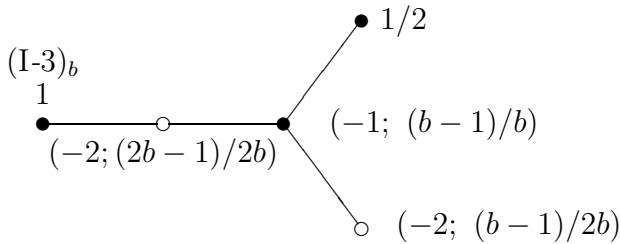
\begin{picture}(100,30)(-10,0)
\put(0,10){\makebox(0,0){$(${\em I}-$3)_b$}}
\put(0,0){\circle*{2}}
\put(0,5){\makebox(0,0){$1$}}
\put(0,0){\line(1,0){19}}
\put(20,0){\circle{2}}
\put(20,-5){\makebox(0,0){$(-2; (2b-1)/2b)$}}
\put(21,0){\line(1,0){19}}
\put(40,0){\circle*{2}}
\put(60,0){\makebox(10,0){$(-1;\ (b-1)/b)$}}
\put(40,0){\line(3,4){13}}
\put(40,0){\line(3,-4){12.5}}
\put(53,17){\circle*{2}}\put(60,17){\makebox(0,0){$1/2$}}
\put(53,-17.5){\circle{2}}\put(70,-17){\makebox(15,0){$(-2;\ (b-1)/2b)$}}
\end{picture}
\vskip 20mm

\item[ (II) ] If $\deg\lfloor \mbox{\em Diff}_{C_{\varphi}(t)}(\Delta^{-}_{\varphi}(t))\rfloor=2$, then there exists a 
${\cal S}$-elementary transformation such that for the resulting log surface $(\hat S,\hat\Delta)$, 
the dual graph of 
$\mbox{\em Supp }\hat\mu^{\ast}\hat\varphi^{\ast}(t)\cup \mbox{\em Supp }\hat\mu^{-1}_{\ast} \hat\Delta^{-}_{\varphi}(t)$ is one of the 
three types $(${\em II}-$1)_b$, $(${\em II}-$2)_b$ or $(${\em II}-$3)_{b,k}$ $(k\geq 1)$ as below, where 
$\hat\mu:\hat M\rightarrow \hat S$ is 
the minimal resolution of the singularities of $\hat S$ over $t\in B$. 
\end{description}

\unitlength 0.8mm
\begin{picture}(100,30)(-10,0)
\put(0,10){\makebox(0,0){$(${\em II}-$1)_b$}}
\put(0,0){\circle*{2}}\put(0,-5){\makebox(0,0){$1$}}
\put(0,0){\line(1,0){20}}
\put(20,0){\circle*{2}}
\put(20,-5){\makebox(0,0){$(0;\ (b-1)/b)$}}
\put(20,0){\line(1,0){20}}
\put(40,0){\circle*{2}}\put(45,0){\makebox(0,0){$1$}}
\end{picture}

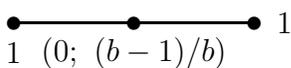
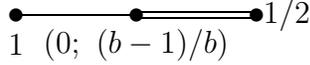
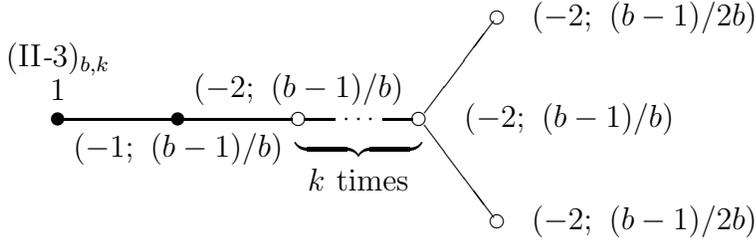
\begin{picture}(100,30)(-80,-30)
\put(0,10){\makebox(0,0){$(${\em II}-$2)_b$}}
\put(0,0){\circle*{2}}
\put(0,-5){\makebox(0,0){$1$}}
\put(0,0){\line(1,0){20}}
\put(20,0){\circle*{2}}
\put(20,-5){\makebox(0,0){$(0;\ (b-1)/b)$}}
\put(20,0.5){\line(1,0){20}}
\put(20,-0.5){\line(1,0){20}}
\put(40,0){\circle*{2}}\put(45,0){\makebox(0,0){$1/2$}}
\end{picture}

\begin{picture}(100,30)(-10,0)
\put(0,10){\makebox(0,0){$(${\em II}-$3)_{b,k}$}}
\put(0,0){\circle*{2}}
\put(0,5){\makebox(0,0){$1$}}
\put(0,0){\line(1,0){20}}
\put(20,0){\circle*{2}}
\put(20,-5){\makebox(0,0){$(-1;\ (b-1)/b)$}}
\put(20,0){\line(1,0){19}}
\put(40,0){\circle{2}}
\put(40,5){\makebox(0,0){$(-2;\ (b-1)/b)$}}
\put(41,0){\line(1,0){5}}
\put(48,0){\circle*{0.5}}
\put(50,0){\circle*{0.5}}
\put(50,-5){\makebox(0,0){$\underbrace{\qquad \qquad}$}}
\put(50,-10){\makebox(0,0){$k$ \mbox{\rm times}}}
\put(52,0){\circle*{0.5}}
\put(54,0){\line(1,0){5}}
\put(60,0){\circle{2}}\put(80,0){\makebox(10,0){$(-2;\ (b-1)/b)$}}
\put(61,1){\line(3,4){11.5}}
\put(61,-1){\line(3,-4){11.5}}
\put(73,17){\circle{2}}\put(90,17){\makebox(15,0){$(-2;\ (b-1)/2b)$}}
\put(73,-17){\circle{2}}\put(90,-17){\makebox(15,0){$(-2;\ (b-1)/2b)$}}
\end{picture}
\vskip 20mm
In the above dual graphs, $\bullet$ denotes a germ of curves or a smooth rational curve 
which is not 
$\mu$-, $($resp., $\hat\mu$ $)$-exceptional and $\circ$ denotes a $\mu$-, 
$($resp., $\hat\mu$ $)$-exceptional curve. 
The numbers attached on $\bullet$ are the multiplicities in $\Delta^M$, 
$($resp., $\hat\Delta^{\hat M}$ $)$. 
$(\ast;\ast)$ denotes $($ the self-intersection number; the multiplicities in 
$\Delta^M$, $($resp., $\hat\Delta^{\hat M}$ $)$ $)$. $b$ is a natural number or $\infty$ such that 
$\mbox{\em mult}_{C_{\varphi}(t)}\Delta$, $($resp., $\mbox{\em
mult}_{C_{\hat\varphi}(t)}\hat\Delta$ $)$ $=(b-1)/b$. 
\end{lm}
{\it Proof.}\ Since we have $\sum_{p\in
C_{\varphi}(t)}m_p(C_{\varphi}(t);\Delta^{-}_{\varphi}(t))=2$, we have $\deg\lfloor
\mbox{Diff}_{C_{\varphi}(t)}(\Delta^{-}_{\varphi}(t))\rfloor=1$ or $2$. Firstly assume that $\deg\lfloor
\mbox{Diff}_{C_{\varphi}(t)}(\Delta^{-}_{\varphi}(t))\rfloor=1$. Since we have 
$m_p(C_{\varphi}(t);\Delta^{-}_{\varphi}(t))\in{\cal S}$ for any $p\in C_{\varphi}(t)$, there exists three points $p_1$, $p_2$ and $p_3\in C_{\varphi}(t)$ such that
$$
m_p(C_{\varphi}(t);\Delta^{-}_{\varphi}(t))=\left\{
\begin{array}{ll}
1 & \mbox{if $p=p_0$},\\
1/2 & \mbox{if $p=p_1$ or $p_2$},\\
0& \mbox{otherwise}
\end{array}
\right.
$$
and one of the following three case occurs. $(1$-$1)$ $S$ is smooth at $p_1$ and $p_2$, 
$(1$-$2)$ $S$ has Du val singularities of type $A_1$ at $p_1$ and $p_2$ or 
$(1$-$3)$ $S$ is smooth at $p_1$ and has a Du val singularity of type $A_1$ at $p_2$. 
Assume that $S$ is not smooth at $p_0$ and let $\mu:M\rightarrow S$ be the minimal resolution of 
$\{p_i|i=0, 1, 2\}$. Write
$
\mu^{\ast}C_{\varphi}(t)=\mu^{-1}_{\ast}C_{\varphi}(t)+\sum_{i=1}^{s}(1/2)E_i+\sum_{j=1}^nl_jF_j,
$
where $l_j>0$ for any $j$ and $\{E_i|1\leq i \leq s\}$, $\{F_j|1\leq j \leq n\}$ are $\mu$-exceptional divisors over 
$p_i$ ($i=1$, $2$) and $p_0$ respectively. Say $(F_1,\mu^{-1}_{\ast}C_{\varphi}(t))=1$, then we have
$$
0=(\mu^{\ast}C_{\varphi}(t),\mu^{-1}_{\ast}C_{\varphi}(t))=(\mu^{-1}_{\ast}C_{\varphi}(t))^2+(1/2)s+l_1 
$$
and $l_1+(1/2)s=1$, since $(\mu^{-1}_{\ast}C_{\varphi}(t))^2=-1$. We know that $0<l_1<1$, so we see that 
$S$ must be smooth at $p_0$ in the cases $(1$-$1)$ or $(1$-$2)$ and that we are 
in the cases (I-$1)$ or (I-$2)$ respectively. Assume that we are 
in the case $(1$-$3)$. Let $\rho:T\rightarrow S$ be the ${\cal S}$-extraction at $p_0$ and 
$\bar F_1$ be the $\rho$-exceptional divisor. Since we have 
$
\rho^{\ast}C_{\varphi}(t):=\rho^{-1}_{\ast}C_{\varphi}(t)+(1/2)\bar F_1,
$
we see that $\bar F_1^2=-2$ and that $T$ is smooth in the neighbourhood of $\bar F_1$, 
hence  we are in the case (I-$3)$.  Secondarily, assume that $\deg\lfloor
\mbox{Diff}_{C_{\varphi}(t)}(\Delta^{-}_{\varphi}(t))\rfloor=2$. Take two points
$p_0$ and $p_1\in C_{\varphi}(t)$ such that 
$$
m_p(C_{\varphi}(t);\Delta^{-}_{\varphi}(t))=\left\{
\begin{array}{ll}
1 & \mbox{if $p=p_0$ or $p_1$},\\
0& \mbox{otherwise}.
\end{array}
\right.
$$
By operating a ${\cal S}$-elementary transformation at $p_0\in \lfloor\Delta^{-}_{\varphi}(t)\rfloor$, 
we may assume that $S$ is smooth at $p_0$. If $S$ is smooth at $p_1$, then 
we see that we are in the cases (II-$1)$ or (II-$2)$. Assume $S$ is not smooth at $p_1$ and 
let $\mu:M\rightarrow S$ be the minimal resolution of 
$p_1$. Write
$
\mu^{\ast}C_{\varphi}(t)=\mu^{-1}_{\ast}C_{\varphi}(t)+\sum_{j=1}^nl_jF_j
$
and
$
\mu^{\ast}K_S=K_M+\sum_{j=1}^na_jF_j
$
where $a_j\geq 0$, $l_j>0$ for any $j$ and $\{F_j|1\leq j \leq n\}$ are $\mu$-exceptional divisors over 
$p_1$. Say $(F_1,\mu^{-1}_{\ast}C_{\varphi}(t))=1$, then we have
$
0=(\mu^{\ast}C_{\varphi}(t),\mu^{-1}_{\ast}C_{\varphi}(t))=(\mu^{-1}_{\ast}C_{\varphi}(t))^2+l_1, 
$
hence $l_1=1$ and $a_1=0$. So we deduce that $(S,C_{\varphi}(t))$ is not purely log terminal but 
$S$ has a Du Val singularity at $p_1$ and $(C_{\varphi}(t),\Delta^{-}_{\varphi}(t))_{p_1}=0$. So we are in the case (II-$3)_{b,k}$.
\hfill \bsquare
\vskip 5mm
\par
Using the technique of ${\cal S}$-elementary transformations as above, we can obtain partial
generalization of Proposition~\ref{pr:bound} to the positive characteristic case. We note that
Proposition~\ref{pr:bound} is proved by the covering trick and the Hodge theory both of which does
not work well in the positive characteristic case.

\begin{pr}\label{pr:12} Let $(S,\Delta)$ be a normal projective log surface with a standard boundary 
defined over an algebraically  closed field $k$ such that $(S,\Delta)$ is log
canonical and
$K_S+\Delta$ is numerically trivial. Assume that there exists a proper surjective  connected
$k$-morphism $\varphi:S\rightarrow B$ onto a smooth projective curve $B$ defined over $k$  and that
there exists an irreducible component of $\lfloor\Delta\rfloor$ which is not contained in a  fibre
of $\varphi$. Then we have $8(K_S+\Delta)\sim 0$ or $12(K_S+\Delta)\sim 0$ and in particular if 
$\mbox{\em char }k\neq 2$, we have $4(K_S+\Delta)\sim 0$ or $6(K_S+\Delta)\sim 0$, 
where $D\sim 0$ for $\mbox{\boldmath$Q$}$-divisor $D$ means that 
$D$ is integral and is linearly equivalent to $0$. 
\end{pr}
{\it Proof.}\  By contracting all the components of singular fibres of $\varphi$ except one in 
each fibre and operating ${\cal S}$-elementary transformations, we may assume that $\rho(S/B)=1$ and all 
of the fibres of $\varphi$ are of type ($a$) as in Lemma~\ref{lm:makelc} or of types (I-$i)_b$ or (II-$i)_b$ ($i=1$, $2$ or $3$) as in 
Lemma~\ref{lm:dynkin}. Let $\Gamma$ be an irreducible component of $\lfloor\Delta\rfloor$. 
Since $m_p(\Gamma;\Delta-\Gamma)\in {\cal S}$ and $\sum_{p\in \Gamma}m_p(\Gamma;\Delta-\Gamma)=2$, 
the possible values of $m_p(\Gamma;\Delta-\Gamma)$ are $0$, $1$, $1/2$, $2/3$, $3/4$ or $5/6$. 
So we deduce that $4(K_M+\Delta^M)$ or $6(K_M+\Delta^M)$is not integral if and only if $\varphi$ has a fibre of types 
(I-$2)_b$ or (II-$3)_b$ ($b=4$, $6$). Assume that $\mbox{char }k\neq 2$ and that $\varphi$ has a fibre of type 
(I-$2)_4$ or (II-$3)_4$. Then there exists distinct three points $p_0$, $p_1$ and $p_2$ in $\Gamma$ 
such that 
$$
m_p(\Gamma;\Delta-\Gamma)=\left\{
\begin{array}{ll}
1/2 & \mbox{if $p=p_0$},\\
3/4 & \mbox{if $p=p_1$ or $p_2$},\\
0& \mbox{otherwise}
\end{array}
\right.
$$
and the induced morphism $\varphi:\Gamma\simeq \mbox{\boldmath$P$}^1\rightarrow B$ has degree $2$ 
which is separable by the assumption and branches at $\{p_i|0\leq i\leq2\}$ but which is 
absurd by the Hurwitz's formula. Assume that $\varphi$ has a fibre of type 
(I-$2)_6$ or (II-$3)_6$. Then there exists distinct three points $p_0$, $p_1$ and $p_2$ in $\Gamma$ 
such that 
$$
m_p(\Gamma;\Delta-\Gamma)=\left\{
\begin{array}{ll}
1/2 & \mbox{if $p=p_0$},\\
2/3 & \mbox{if $p=p_1$},\\
5/6 & \mbox{if $p=p_2$},\\
0& \mbox{otherwise}
\end{array}
\right.
$$
and we can derive the absurdity as in the same way as above. Therefore, we only have to check that
if $r(K_M+\Delta^M)$ is an integral divisor on $M$ for some $r\in \mbox{\boldmath $N$}$, 
then we have $r(K_M+\Delta^M)\sim 0$ or $(r-1)(K_M+\Delta^M)\sim 0$. Put $D:=r(K_M+\Delta^M)$. 
If the genus of $B$ is zero, 
then $M$ is a rational surface, hence $D\sim 0$. Assume that the genus of $B$ is positive. 
Then $\lfloor\Delta^M\rfloor$ contains an smooth elliptic curve $\Gamma$ 
and $M$ is birationally elliptic ruled. From $\chi({\cal O}_M(D))=\chi({\cal O}_M)=0$, we have 
$h^0({\cal O}_M(D))+h^2({\cal O}_M(D))=h^1({\cal O}_M(D))$. If $h^2({\cal O}_M(D))\neq 0$, then 
$h^0({\cal O}_M(K_M-D))\neq 0$ by the Serre duality and for an ample divisor $H$, we have 
$0\leq (K_M-D,H)=(K_M,H)=-(\Delta^M,H)<0$ which is absurd. 
Thus we get $h^2({\cal O}_M(D))=0$ and $h^0({\cal O}_M(D))=h^1({\cal O}_M(D))$. 
Assume that $h^2({\cal O}_M(D-\Gamma))\neq 0$. Then $h^0({\cal O}_M(K_M-D+\Gamma))\neq 0$ 
by the Serre duality again and $(K_M+\Gamma, H)\geq 0$, hence $\Delta^M=\Gamma$ and 
$(r-1)(K_M+\Delta^M)\sim 0$. So we may assume that $h^2({\cal O}_M(D-\Gamma))=0$. 
By the assumption, there exists a surjection: 
$
H^1({\cal O}_M(D))\rightarrow H^1({\cal O}_{\Gamma}(D))\simeq H^1({\cal O}_{\Gamma}),
$
hence $h^0({\cal O}_M(D))=h^1({\cal O}_M(D))>0$. Thus we get the assertion.
\hfill \bsquare
\vskip 5mm

\subsection{Case $e_{\mbox{\rm orb}}(\tilde S\setminus \tilde\Delta)=0$ and 
$\lfloor\Delta\rfloor\neq 0$}
Let $(S,\Delta)$ be a projective log surface with a standard boundary defined over 
an algebraically closed field $k$. Assume that $(S,\Delta)$ is log terminal and 
$K_S+\Delta$ is numerically trivial. $(S,\Delta)$ can be roughly classified into the following three 
types.  
\begin{description}
\item[I]: $\lfloor \Delta\rfloor=0$,\par
\item[II]: $\lfloor \Delta\rfloor \neq 0$ and 
$\lfloor \mbox{Diff}_{\lfloor \Delta\rfloor^{\nu}}(\Delta-\lfloor \Delta\rfloor)\rfloor=0$,\par
\item[III]: $\lfloor \Delta\rfloor \neq 0$ and 
$\lfloor \mbox{Diff}_{\lfloor \Delta\rfloor^{\nu}}(\Delta-\lfloor \Delta\rfloor)\rfloor\neq 0$,
\end{description}
where $\nu$ denotes the normalization map 
$\nu:\lfloor \Delta\rfloor^{\nu}\rightarrow \lfloor \Delta\rfloor$.

\begin{df}\label{df:123}\begin{em} Log surfaces $(S,\Delta)$ with the above conditions are said to
 be {\it of type {\rm I}, {\rm II} and {\rm III}} respectively.
\end{em}\end{df}
The classification of the log surfaces as above in the case 
that $K_S+\Delta$ is Cartier is well known.

\begin{lm}\label{lm:ind 1 ls} Assume that $K_S+\Delta$ is Cartier, then one of the following hold. 
\begin{description}
\item[ $(1)$ ] $S$ is either an abelian surface, a hyperelliptic surface or a normal surface with at worst Du Val singularities
whose minimal resolution is a $K3$ surface and $\Delta=0$. \par
\item[ $(2)$ ] $S$ is a rational or birationally elliptic ruled surface with Du Val singularities and $\Delta$ is either 
a smooth elliptic curve or a disjoint union of two smooth elliptic curves and the support of 
$\Delta$ is disjoint from the singular loci of $S$. \par
\item[ $(3)$ ] $S$ is a rational surface with Du Val singularities and $\Delta$ is a connected cycle of smooth rational curves which is 
disjoint from the singular loci of $S$.
\end{description}
\end{lm}
{\it Proof.} By the definition and assumptions, the proof reduces to the well known results 
(see \cite{kulikov}). \hfill \bsquare
\vskip 5mm
In what follows we assume that $\mbox{char }k=0$.
\begin{lm}\label{lm:local cover} Let $(S,\Delta)$ be a germ of 
log surfaces with a standard boundary and let 
$(\tilde S, \tilde\Delta)$ be the log canonical cover of $(S,\Delta)$.  
If $(S,\Delta)$ is log terminal, then $(\tilde S, \tilde \Delta)$ 
is also log terminal.
\end{lm}
{\it Proof.} If the number of irreducible components of $\lfloor \Delta\rfloor$
 is less than $2$, then $(S,\Delta)$ is purely log terminal, hence so is 
$(\tilde S, \tilde \Delta)$ (see \cite{shokurov}, Corollary 2.2). If not, we get the assertion since the index of 
$K_S+\Delta$ is $1$ in this case. 
\hfill \bsquare
\vskip 5mm
By the log abundance theorem for surfaces, there exists a 
global log canonical cover $(\tilde S,\tilde\Delta)$ of $(S,\Delta)$.
\begin{lm} If $(S,\Delta)$ is of type {\em I} $($ resp., {\em II}, resp., {\em III} $)$, then 
$(\tilde S,\tilde\Delta)$ is also of type {\em I} $($ resp., {\em II}, resp., {\em III} $)$.
\end{lm}
{\it Proof.} Firstly, we prove that $(\tilde S,\tilde\Delta)$ is log terminal. By 
Lemma~\ref{lm:local cover}, we only have to check that any irreducible components of 
$\lfloor \Delta\rfloor$ does not have self-intersections. 
Assume that there exists an irreducible component $\tilde \Gamma$ which has self-intersections. 
Let $\Gamma$ be the image of $\tilde \Gamma$ on $S$. Since we have
$$
K_{\tilde S}+\pi^{-1}(\Gamma)=\pi^{\ast}(K_S+\Gamma+\{\Delta \}), 
$$
where $\pi$ is the covering morphism $\pi:\tilde S\rightarrow S$, $(\tilde S, \pi^{-1}(\Gamma))$ is purely log terminal, 
hence so is $(\tilde S,\tilde \Gamma)$, but which is absurd. Noting that $\tilde \Delta=\pi^{-1}(\Delta)$ 
and that 
$
\lfloor \mbox{Diff}_{\lfloor \tilde\Delta\rfloor^{\nu}}(0)\rfloor=
\pi^{-1}\lfloor \mbox{Diff}_{\lfloor \Delta\rfloor^{\nu}}(\Delta-\lfloor \Delta\rfloor)\rfloor,
$
we get the assertion.
\hfill \bsquare
\vskip 5mm
The following result is essentially due to S. Tsunoda.
\begin{pr}[c.f. \cite{tsunoda1}, Theorem 2.1]\label{pr:bound} \ 

\begin{description}
\item[ $(1)$ ] If $(S,\Delta)$ is of type {\em II}, then $4(K_S+\Delta)\sim 0$ or $6(K_S+\Delta)\sim 0$. \par 
\item[ $(2)$ ] If $(S, \Delta)$ is of type {\em III}, then $2(K_S+\Delta)\sim 0$.
\end{description}
\end{pr}
{\it Proof.} The argument in the proof of \cite{tsunoda1}, Theorem 2.1 can also be applied in the  
case in which $S$ is rational. If $S$ is not rational, we only have to apply the argument in
the last part of the proof of Proposition\ref{pr:12}\hfill
\bsquare
\vskip 5mm
\par
Let $(S,\Delta)$ be a projective log surface with a standard boundary defined over the complex number
field 
$\mbox{\boldmath$C$}$. 
Assume that $(S,\Delta)$ is log terminal and that $K_S+\Delta$ is numerically trivial. 
Let $(\tilde S,\tilde\Delta)$ be the log canonical cover. 
To classify log surfaces as above, we need some conditions for 
$(\tilde S,\tilde\Delta)$ especially in the case $\lfloor\Delta\rfloor\neq 0$. 
In this section, we shall classify $(S,\Delta)$ with $\lfloor\Delta\rfloor\neq 0$ under the 
condition that $e_{\mbox{\rm orb}}(\tilde S\setminus \tilde\Delta)=0$ which seems to be the most fundamental
 case. In general, applying the log minimal program for $S$ with respect $K_S+\{\Delta\}$, we get a birational morphism 
$\tau:S\rightarrow S^{\flat}$ to a projective normal surface $S^{\flat}$ such that 
(1) $K_{S^{\flat}}+\{\Delta^{\flat}\}$ is nef, where $\Delta^{\flat}:=\tau_{\ast}\Delta$, 
(2) $-(K_{S^{\flat}}+\{\Delta^{\flat}\})$ is ample and $\rho(S^{\flat})=1$ or (3) 
there exists a projective surjective morphism $\varphi^{\flat}:S^{\flat}\rightarrow B$ onto a 
smooth projective curve $B$ with $\rho(S^{\flat}/B)=1$ and 
$-(K_{S^{\flat}}+\{\Delta^{\flat}\})$ is relatively ample with respect to $\varphi^{\flat}$. 
Assume that we are in the case (1). 
Then we have 
$\lfloor\Delta^{\flat}\rfloor=0$ and $K_S+\Delta=\tau^{\ast}(K_{S^{\flat}}+\{\Delta^{\flat}\})$. On the other hand, we have 
$K_S+\{\Delta\}-\tau^{\ast}(K_{S^{\flat}}+\{\Delta^{\flat}\})\geq 0$ by our construction. Thus we see that the condition (1) 
is equivalent to $\lfloor\Delta\rfloor=0$. The following is the key lemma for our purpose.
\begin{lm}\label{lm:K-cont} Assume that $e_{\mbox{\rm{orb}}}(\tilde S\setminus \tilde\Delta)=0$. 
 Then $\mbox{\em{Exc} }\tau\subset \lfloor\Delta\rfloor$. 
\end{lm}
{\it Proof.} We may assume that $\lfloor\Delta\rfloor\neq 0$. Let $\tau^{\natural}:S\rightarrow 
S^{\natural}$  be a birational extremal contraction with respect to $K_S+\{\Delta\}$. Let $E$
denote the  exceptional divisor for $\tau^{\natural}$ and put
$\Delta^{\natural}:=\tau^{\natural}_{\ast}\Delta$. Assume that 
$K_S+\Delta$ is Cartier firstly. By \cite{lipman-sommese}, Theorem 0.1, we see that $E$ is a
$(-1)$-curve and that
$(E,\Delta)=1$. Moreover,  there exists a point $p\in E$ such that $E\cap\mbox{Sing
}S=\emptyset$ or $\{p\}$ and
$E$ contracts  to a smooth point of $S^{\natural}$. Therefore, if $E$ is not contained in
$\mbox{Supp }\Delta$, we have
$$
e_{\mbox{\rm orb}}(S\setminus \Delta)-e_{\mbox{\rm orb}}(S^{\natural}\setminus
\Delta^{\natural})=
\rho(S)-\rho(S^{\natural})-(1-\frac{1}{o_p(S)})=\frac{1}{o_p(S)}. 
$$
On the other hand, we have $e_{\mbox{\rm orb}}(S\setminus \Delta)=0$ and 
$e_{\mbox{\rm orb}}(S^{\natural}\setminus \Delta^{\natural})\geq 0$, which is absurd. 
If $K_S+\Delta$ is not Cartier, we take log canonical covers $(\tilde S,\tilde\Delta)$ and 
$(\tilde S^{\natural},\tilde\Delta^{\natural})$ of $(S,\Delta)$ and 
$(S^{\natural},\Delta^{\natural})$ respectively such that there exists a
birational morphism 
$\tilde\tau^{\natural}:\tilde S\rightarrow \tilde S^{\natural}$ induced by $\tau^{\natural}$. 
Let $\pi$ denote the covering morphism $\tilde S\rightarrow S$. Since we have
$
(K_{\tilde S}, \pi^{\ast}E)=\deg \pi(K_S+\{\Delta\},E)<0,
$
we see that $\pi^{-1}E\subset \mbox{Supp }\tilde \Delta=\mbox{Supp }\pi^{-1}\lfloor\Delta\rfloor$ 
by the previous argument, hence we conclude that $E\subset \lfloor\Delta\rfloor$.
\hfill\bsquare

\subsubsection{Case Type II}
Let $(S,\Delta)$ be a projective log surface with a standard boundary and assume that there exists a
structure of a conic fibration $\varphi:S\rightarrow B$ with $\rho(S/B)=1$, where $B$ is a
smooth projective connected curve such that
$\varphi$ has only fibres of types listed in Lemma~\ref{lm:dynkin}. We denote the number of fibres
 of type ${\cal T}$ by $\nu({\cal T})$ if it is finite. 
\begin{lm}\label{lm:typeII} Assume that $(S,\Delta)$ is of type {\em II} and that 
$e_{\mbox{\rm{orb}}}(\tilde S\setminus \tilde\Delta)=0$. Then the followings holds.
\begin{description}
\item[($1$)] There exists a projective surjective morphism $\varphi:S\rightarrow B$ onto a 
smooth projective curve $B$ with $\rho(S/B)=1$,\par
\item[($2$)] $(S,\Delta^{+}_{\varphi}(t))$ is log canonical for any $t\in B$ and\par
\item[($3$)] $\mbox{\em{Supp }}\lfloor\mbox{\em{Diff}}_{C_{\varphi}(t)}(\Delta^{-}_{\varphi}(t))\rfloor
\subset\lfloor\Delta^{-}_{\varphi}(t)\rfloor$ for any $t\in B$.
\end{description}  
\end{lm}
{\it Proof.} Firstly, assume that $K_S+\Delta$ is Cartier. 
If ($1$) does not hold, $S$ is a rank one Gorenstein log del Pezzo surface by
Lemma~\ref{lm:K-cont}. Therefore we have 
$
0=e_{\mbox{\rm{orb}}}(S\setminus \Delta)=e_{\mbox{\rm{orb}}}(S)=3-\sum_{p\in
S}\{1-(1/o_p(S))\}, 
$
which contradicts the Table IV.  
Hence ($1$) holds under the assumption that $K_S+\Delta$ is Cartier. 
Assume that $h^1({\cal O}_S)=1$. Then we have
$
0=e_{\mbox{\rm orb}}(S\setminus \Delta)=-\sum_{p\in S}\{1-(1/o_p(S))\}, 
$
hence 
$\varphi$ gives a structure of relatively minimal elliptic ruled surface on $S$. 
Noting that the induced morphism $\varphi:\Delta\rightarrow B$ is \'etale, 
we see that ($2$) and ($3$) hold. If $h^1({\cal O}_S)=0$, we have 
$$
4=\sum_{t\in B}\sum_{p\in C(t)}(1-\frac{1}{o_p(S)})=
\nu(\mbox{(I-2)}_1)+\sum_{k\geq 1}\frac{4k-1}{4k}\nu(\mbox{(II-3)}_{1,k}). 
$$
We note here that we do not have to operate any ${\cal S}$-elementary transformation 
by examining the proof of Lemma~\ref{lm:dynkin}. By the Hurwitz's formula we have 
$\nu(\mbox{(I-2)}_1)+\sum_{k\geq 1}\nu(\mbox{(II-3)}_{1,k})\leq 4$, hence we deduce that 
$\nu(\mbox{(I-2)}_1)=4$ and $\nu(\mbox{(II-3)}_{1,k})=0$ $(k\geq 1)$. 
Thus we get the assertion ($2$) and ($3$) under the assumption that $K_S+\Delta$ is Cartier. 
We note that $\Delta^2=K_S^2=0$ in each of the above cases. 
Assume that $K_S+\Delta$ is not necessarily Cartier and that there exists a birational extremal 
contraction 
$\tau^{\natural}:S\rightarrow  S^{\natural}$ with respect to 
$K_S+\{\Delta\}$ and let $E$ be as in Lemma~\ref{lm:K-cont}. 
Then by the assumption and  Lemma~\ref{lm:K-cont}, we have
$
(K_S+\Delta,E)=(K_S+\{\Delta\},E)+E^2+(\lfloor\Delta\rfloor-E,E)<0,
$
since $(K_S+\{\Delta\},E)<0$, $E^2<0$ and $(\lfloor\Delta\rfloor-E,E)=0$, but which is absurd. 
Thus if ($1$) does not hold, $-(K_{S}+\{\Delta\})$ is ample and $\rho(S)=1$, which implies that
 $\tilde \Delta=\pi^{\ast}\lfloor\Delta\rfloor$ is ample. Thus we get the absurdity. 
By taking the Stein factorization, 
there exists a finite morphism 
$\varpi:\tilde B\rightarrow B$ from a smooth projective curve $\tilde B$ 
and a proper surjective connected morphism 
$\tilde \varphi:\tilde S\rightarrow \tilde B$ 
such that $\varphi\circ \pi=\varpi\circ\tilde \varphi$. Since we have 
$
K_{\tilde S}+\sum_{\varpi(\tilde t)=t}\Delta^{+}(\tilde \varphi;\tilde t)=
\pi^{\ast}(K_S+\Delta^{+}_{\varphi}(t))
$
and we know that $(\tilde S,\sum_{\varpi(\tilde t)=t}\Delta^{+}(\tilde \varphi;\tilde t))$ is 
log canonical, we conclude that ($2$) holds. As for ($3$), we get the assertion by the following diagram:
\begin{eqnarray*}
\mbox{Supp}\sum_{\varpi(\tilde t)=t}
\mbox{Diff}_{C(\tilde \varphi;\tilde t)}(\Delta^{-}(\tilde \varphi;\tilde t))&=&
\mbox{Supp }\pi^{-1}\lfloor\mbox{Diff}_{C_{\varphi}(t)}(\Delta^{-}_{\varphi}(t))\rfloor \\
\cap \qquad  & &\qquad \cap\\
\sum_{\varpi(\tilde t)=t}\Delta^{-}(\tilde \varphi;\tilde t)&=& 
\pi^{-1}\lfloor\Delta^{-}_{\varphi}(t)\rfloor.
\end{eqnarray*}
\hfill\bsquare

Let $(S,\Delta)$ be a projective log surface with a standard boundary and assume that there exists a
structure of a conic fibration $\varphi:S\rightarrow B$ with $\rho(S/B)=1$, where $B$ is a
smooth projective connected curve such that
$\varphi$ has only fibres of types listed in Lemma~\ref{lm:dynkin}. If $\varphi$ has fibres of type
${\cal T}_1$, $\dots$, ${\cal T}_r$ and generically of type ${\cal T}$, then we shall write 
$\mbox{Typ}(S,\Delta;\varphi)=({\cal T}_1+ \cdots +{\cal T}_r;{\cal T})$.

\begin{pr}\label{pr:cltypeII} Type $\mbox{\em II}$ log
surfaces $(S,\Delta)$ with 
$e_{\mbox{\rm{orb}}}(\tilde S\setminus \tilde\Delta)=0$ are 
classified into the following $22$ types modulo ${\cal S}$-elementary transformations.
\begin{description}
\item[ $(a)$ ] $S\simeq \mbox{\boldmath $P$}_E({\cal E})$, where ${\cal E}$ is a rank two vector
bundle on an elliptic curve $E$ and 
$\Delta=\Gamma$ $(a\mbox{-}1)$, $\Gamma_1+\Gamma_2$ $(a\mbox{-}2)$, $\Gamma+(1/2)\Xi$ 
$(a\mbox{-}3)$ or 
$\Gamma+(1/2)\Xi_1+(1/2)\Xi_2$ $(a\mbox{-}4)$, where $\Gamma$, $\Gamma_i$, $\Xi$ and $\Xi_j$ are 
smooth elliptic curves which are disjoint from each other.
\item[ $(b)$ ] $S\simeq \mbox{\boldmath $P$}^1\times \mbox{\boldmath $P$}^1$ and 
$$
\Delta=\Gamma+\sum_{i=1}^{2}(1/2)\Xi_{1,i}+\left\{
\begin{array}{l}
\sum_{j=1}^4(1/2)\Xi_{2,j},\\
\sum_{j=1}^3(2/3)\Xi_{2,j},\\
(1/2)\Xi_{2,1}+\sum_{j=2}^3(3/4)\Xi_{2,j} \mbox{ or }\\
(1/2)\Xi_{2,1}+(2/3)\Xi_{2,2}+(5/6)\Xi_{2,3},
\end{array}
\right.
$$
where $\Gamma$ and $\Xi_{1,i}$ are fibres of the first projection for any $i$ and 
$\Xi_{2,j}$ is a fibre of the second projection for any $j$.
\item[ $(c)$ ] There exists a sequence of ${\cal S}$-elementary transformations 
such that for the resulting log surface $(\hat S,\hat \Delta)$, 
$\hat S\simeq \mbox{\boldmath $P$}^1\times \mbox{\boldmath $P$}^1$ and 
$$
\hat\Delta=\hat\Gamma_1+\hat\Gamma_2+\left\{
\begin{array}{l}
\sum_{j=1}^4(1/2)\hat\Xi_j,\\
\sum_{j=1}^3(2/3)\hat\Xi_j,\\
(1/2)\hat\Xi_1+\sum_{j=2}^3(3/4)\hat\Xi_j \mbox{ or }\\
(1/2)\hat\Xi_1+(2/3)\hat\Xi_2+(5/6)\hat\Xi_3,
\end{array}
\right.
$$ 
where $\hat\Gamma_i$ is a fibre of the first projection for any $i$ and 
$\hat\Xi_j$ is a fibre of the second projection for any $j$.
\item[ $(d)$ ] There exists a structure of $\mbox{\boldmath $P$}^1$-fibration 
$\varphi:S\rightarrow \mbox{\boldmath $P$}^1$ with $\rho(S/\mbox{\boldmath $P$}^1)=1$ such that one
of the followings in 
 {\em Table II} holds.
\begin{center}
 {\em Table II}\\
\begin{tabular}{l|c|c}\hline
 & $\Delta$ & $\mbox{\em Typ}(S,\Delta;\varphi)$  \\ 
\hline \hline
$(d\mbox{-}1)$ & $\Gamma$  & $(4(\mbox{\em I}$-$2)_1$; $(\mbox{\em II}$-$1)_1)$ \\ 
\hline
$(d\mbox{-}2)$ & $\Gamma+(1/2)\Xi_h$ &$(4(\mbox{\em I}$-$3)_1$; $(\mbox{\em I}$-$1)_1)$ \\
\hline
$(d\mbox{-}3)$ & $\Gamma+(1/2)\Xi_{h}+(1/2)\Xi_{v,1}+(1/2)\Xi_{v,2}$ &$(2(\mbox{\em I}$-$1)_2+2(\mbox{\em
I}$-$3)_1$; $(\mbox{\em I}$-$1)_1)$ \\
\hline
$(d\mbox{-}4)$ & $\Gamma+(1/2)\Xi_h+(1/2)\Xi_{v,1}+(3/4)\Xi_{v,2}$ & $((\mbox{\em
I}$-$1)_4+(\mbox{\em I}$-$3)_2+(\mbox{\em I}$-$3)_1$; $(\mbox{\em I}$-$1)_1)$ \\
\hline
$(d\mbox{-}5)$ & $\Gamma+(1/2)\Xi_h+\sum_{j=1}^{3}(1/2)\Xi_{v,j}$ & $(2(\mbox{\em
I}$-$3)_2+(\mbox{\em I}$-$1)_2$; $(\mbox{\em I}$-$1)_1)$ \\
\hline
$(d\mbox{-}6)$ & $\Gamma+(1/2)\Xi_h+(2/3)\Xi_{v,1}+(2/3)\Xi_{v,2}$ & $((\mbox{\em
I}$-$1)_3+(\mbox{\em I}$-$3)_1+(\mbox{\em I}$-$3)_3$; $(\mbox{\em I}$-$1)_1)$ \\
\hline
\end{tabular}
\end{center}
$\Gamma$ in $(d\mbox{-}1)$ and $\Xi_h$ in $(d\mbox{-}2)$ denote a smooth elliptic curve with
$\Gamma^2=\Xi_h^2=0$. 
$\Gamma$ in $(d\mbox{-}2)$, $(d\mbox{-}3)$, $(d\mbox{-}4)$, $(d\mbox{-}5)$ and $(d\mbox{-}6)$, $\Xi_h$ in $(d\mbox{-}3)$, $(d\mbox{-}4)$, $(d\mbox{-}5)$
and
$(d\mbox{-}6)$ denote smooth rational curves with $\Gamma^2=\Xi_h^2=0$ which are horizontal with
respect to $\varphi$.
$\Xi_{v,j}$ in $(d\mbox{-}3)$, $(d\mbox{-}4)$, $(d\mbox{-}5)$ and
$(d\mbox{-}6)$ denote smooth rational curves which are vertical with respect to $\varphi$ for any $j$.

\item[ $(e)$ ] There exists a structure of $\mbox{\boldmath $P$}^1$-fibration 
$\varphi:S\rightarrow \mbox{\boldmath $P$}^1$ with $\rho(S/\mbox{\boldmath $P$}^1)=1$ and 
there exists a sequence of ${\cal S}$-elementary transformations 
such that for the resulting log surface $(\hat S,\hat \Delta)$, one of the followings in 
 {\em Table III} holds.
\begin{center}
 {\em Table III}\\
\begin{tabular}{l|c|c}\hline
 & $\hat\Delta$ & $\mbox{\em Typ}(\hat S,\hat\Delta;\hat\varphi)$  \\ 
\hline \hline
$(e\mbox{-}1)$ & $\hat\Gamma+\sum_{j=1}^{3}(1/2)\hat\Xi_j$ & $(2(\mbox{\em I}$-$2)_2+(\mbox{\em
II}$-$1)_2$;
$(\mbox{\em II}$-$1)_1)$ \\
\hline
$(e\mbox{-}2)$ & $\hat\Gamma+\sum_{j=1}^{2}(1/2)\hat\Xi_j$ & $(2(\mbox{\em I}$-$2)_1+2(\mbox{\em
II}$-$1)_2$;
$(\mbox{\em II}$-$1)_1)$ \\
\hline
$(e\mbox{-}3)$ & $\hat\Gamma+\sum_{j=1}^{2}(2/3)\hat\Xi_j$ & $((\mbox{\em I}$-$2)_1+(\mbox{\em
I}$-$2)_3+(\mbox{\em II}$-$1)_3$; $(\mbox{\em II}$-$1)_1)$ \\
\hline
$(e\mbox{-}4)$ & $\hat\Gamma+(1/2)\hat\Xi_1+(3/4)\hat\Xi_2$ & $((\mbox{\em I}$-$2)_1+(\mbox{\em
I}$-$2)_2+(\mbox{\em II}$-$1)_4$; $(\mbox{\em II}$-$1)_1)$ \\
\hline
\end{tabular}
\end{center}
$\hat\Gamma$ denote a smooth rational curve with $(\hat\Gamma,\hat\varphi^{\ast}(t))=2$ for $t\in
\mbox{\boldmath $P$}^1$ and $\hat\Xi_i$ denote a smooth rational curve which are vertical with
respect to $\hat\varphi$ for any $i$.
\end{description}
\end{pr}

{\it Proof.} Let $\Gamma_1$ be an irreducible component of $\lfloor\Delta\rfloor$. Then we have
$2g(\Gamma_1)-2+\deg\mbox{Diff}_{\Gamma_1}(\Delta-\Gamma_1)=0$, where $g(\Gamma_1)$ denote the
genus of $\Gamma_1$. (1) Assume that $g(\Gamma_1)=1$. Then
$\deg\mbox{Diff}_{\Gamma_1}(\Delta-\Gamma_1)=0$ and $\varphi$ has only fibres of types (I-1)$_1$,
(I-2)$_1$ or (II-1)$_1$ by Lemma~\ref{lm:dynkin} and Lemma~\ref{lm:typeII}. (1-1) Assume that 
$g(B)=1$. Since the induced morphism $\Gamma_1\rightarrow B$ is \'etale, type (I-2)$_1$ fibre does
not exist, hence we are in the case $(a)$. (1-2) Assume that $g(B)=0$. Applying the Hurwitz's
formula to the induced double covering $\Gamma_1\rightarrow B$, we obtain $\nu((\mbox{I-}2)_1)=4$,
hence we are in the case $(d\mbox{-}1)$. (2) Assume that $g(\Gamma_1)=0$. Then we have $\sum_{p\in
\Gamma_1}m_p(\Gamma_1;\Delta-\Gamma_1)=2$. (2-1) Assume that $\Gamma_1$ is a section of $\varphi$.
Then $\varphi$ has only fibres of types (I-1)$_b$, (I-3)$_b$ or (II-1)$_b$. (2-1-1) Assume that
$\Gamma_1=\lfloor\Delta\rfloor$. Then $\varphi$ has only fibres of types (I-1)$_b$ or (I-3)$_b$ and
we have
$\sum_{b\geq 1}\{(b-1)/b\}\nu((\mbox{I-}1)_b)+\sum_{b^{\prime}\geq
1}\{(2b^{\prime}-1)/2b^{\prime}\}\nu((\mbox{I-}3)_{b^{\prime}})$. Let $\Xi_h$ be the sum
of all the horizontal components of $\mbox{Supp}\{\Delta\}$. If $\Xi_h$ is reducible, then
$\varphi$ has only fibres of type (I-1)$_b$, hence we are in one of the cases $(b)$. Assume that
$\Xi_h$ is irreducible. Then we have
$e_{\mbox{top}}(\Xi_h)=4-\sum_{b^{\prime}\geq 1}\nu((\mbox{I-}3)_{b^{\prime}})$ by the Hurwitz's
formula. On the other hand, since we have $K_S+\Gamma_1+\Xi_h+\{\Delta\}_v\sim_{\mbox{\boldmath
$Q$}}(1/2)\Xi_h$, where $\{\Delta\}_v$ denotes the vertical part of $\{\Delta\}$, we have
$
e_{\mbox{top}}(\Xi_h)=\deg\mbox{Diff}_{\Xi_h}(\{\Delta\}_v)-(1/2)\Xi_h^2
=2\sum_{b\geq1}\{(b-1)/b\}\nu((\mbox{I-}1)_b)
+\sum_{b^{\prime}\geq1}\{(b^{\prime}-1)/b^{\prime}\}\nu((\mbox{I-}3)_{b^{\prime}})-(1/2)\Xi_h^2$,
hence $\Xi_h^2=0$. Since $e_{\mbox{top}}(\Xi_h)\leq 2$ and $e_{\mbox{top}}(\Xi_h)\in
2\mbox{\boldmath $Z$}$, we have
$\sum_{b^{\prime}\geq 1}\nu((\mbox{I-}3)_{b^{\prime}})=2$ or $4$. If $\sum_{b^{\prime}\geq
1}\nu((\mbox{I-}3)_{b^{\prime}})=4$, then we have
$e_{\mbox{top}}(\Xi_h)=0$ and $\nu((\mbox{I-}3)_{1})=4$, hence we are in the case $(d\mbox{-}2)$.
Assume that $\sum_{b^{\prime}\geq 1}\nu((\mbox{I-}3)_{b^{\prime}})=2$. Then we have
$e_{\mbox{top}}(\Xi_h)=2$ and we see that we are in the cases $(d\mbox{-}3)$, $(d\mbox{-}4)$,
$(d\mbox{-}5)$ and $(d\mbox{-}6)$. (2-1-2) Assume that $\lfloor\Delta\rfloor$ is decomposed into
two sections $\Gamma_1$ and $\Gamma_2$. Then $\varphi$ has only fibres of type (II-1)$_b$ after
${\cal S}$-elementary transformations, hence we are in the cases $(c)$. (2-2) Assume that
$(\Gamma_1,\varphi^{\ast}(t))=2$ for $t\in \mbox{\boldmath $P$}^1$. Then $\varphi$ has only fibres 
of types (I-2)$_b$ or (II-1)$_b$ after ${\cal S}$-elementary transformations. We note
that $\sum_{b\geq 1}\{(b-1)/b\}\nu((\mbox{I-}2)_b)+2\sum_{b^{\prime}\geq
1}\{(b^{\prime}-1)/b^{\prime}\}\nu((\mbox{II-}1)_{b^{\prime}})=2$. By the Hurwitz's
formula, we have $\sum_{b\geq 1}\nu((\mbox{I-}2)_b)=2$. Thus we see that we are in the cases 
$(e\mbox{-}1)$, $(e\mbox{-}2)$, $(e\mbox{-}3)$ and $(e\mbox{-}4)$.
\hfill\bsquare
\vskip 5mm
We prove the lemma needed later which was taught to the author by Prof. A. Fujiki.

\begin{lm}\label{lm:universal} Let $X$ be a complex manifold and let $G$ be a finite subgroup of
the holomorphic automorphism group $\mbox{\em Aut }X$. Let $\Pi:{\cal U}\rightarrow X$ be the
universal cover of $X$. Then there exists a discrete subgroup $\tilde G$ of $\mbox{\em Aut }{\cal
U}$ which acts on ${\cal U}$ properly discontinuously such that ${\cal U}/\tilde G\simeq X/G$
\end{lm}
{\it Proof.} Let $\mbox{Aut}_{\Pi}{ \cal U}$ be the subgroup of $\mbox{Aut }{\cal U}$ which
consists of all the elements $\tilde\gamma\in \mbox{Aut }{\cal U}$ such that there exists an
element $\gamma\in \mbox{Aut }X$ which satisfies $\gamma\circ \Pi=\Pi\circ\tilde\gamma$. Then
there exists the following exact sequence of groups:
$
1\longrightarrow \pi_1(X)\longrightarrow \mbox{Aut}_{\Pi}{\cal U}\longrightarrow \mbox{Aut
}X\longrightarrow 1.
$
Let $\tilde G$ be the inverse image of $G$ by the third map in the above exact sequence. We note
that $\tilde G$ is a discrete subgroup of $\mbox{Aut}_{\Pi}{ \cal U}$ and that we have the
following exact sequence:
$
1\longrightarrow \pi_1(X)\longrightarrow \tilde G\longrightarrow G\longrightarrow 1.
$
Therefore, we see that $\tilde G$ acts on ${\cal U}$ properly discontinuously and that ${\cal
U}/\tilde G\simeq ({\cal U}/\pi_1(X))/G\simeq X/G$.
\hfill\bsquare
\vskip 5mm

Let $(X,D)$ be a normal log variety with reduced boundary $D$. Assume that a finite group
$G$ acts on $X$ faithfully preserving $D$. Let $f:X\rightarrow X/G$ be the quotient
morphism and assume that any component of $D$ is not contained in the ramification divisor of $f$.
For a prime divisor
$\Gamma$ on
$X/G$, take a prime divisor
$\tilde\Gamma$ on
$X$ contained in
$f^{-1}(\Gamma)$ and let $G_{\tilde\Gamma}$ denote the subgroup of $G$ consisting of all the
element of $G$ which acts on $\tilde\Gamma$ trivially. The order $|G_{\tilde\Gamma}|$
of $G_{\tilde\Gamma}$, which is nothing but the ramification index of $f$ at $\tilde\Gamma$, does
not depend on the choice of $\tilde\Gamma$. Put $e_{\Gamma}:=|G_{\tilde\Gamma}|$. We
define the $\mbox{\boldmath $Q$}$-boundary
$D_G$ on $X/G$ as follows.
$$
D_G:=D/G+\sum_{\Gamma;prime}\frac{e_{\Gamma}-1}{e_{\Gamma}}\Gamma,
$$
where the summation is taken over all the prime divisors $\Gamma$ on $X/G$. We note
that by the definition, we have
$K_X+D=f^{\ast}(K_{X/G}+D_G)$ if $K_{X/G}+D_G$ is $\mbox{\boldmath $Q$}$-Cartier.
\vskip 5mm

Here is a corollary of Proposition~\ref{pr:cltypeII}.

\begin{co} Let $(S,\Delta)$ be a type $\mbox{\em II}$ log surface with 
$e_{\mbox{\rm{orb}}}(\tilde S\setminus \tilde\Delta)=0$. Then there exists a relatively
minimal elliptic ruled surface $\pi:P\rightarrow E$ over an elliptic curve $E$ which admits a $\pi$-equivariant action of
a finite abelian group $G$ with $|G|=l\mbox{\em Ind}(K_S+\Delta)$, where $l=1$, $2$ or $4$ and there exists a
$G$-invariant reduced divisor $D$ on $P$ which is a smooth elliptic curve or a disjoint union of two smooth elliptic
curves such that $(S,\Delta)\simeq (P/G, D_G)$. In particular, $(S,\Delta)$ is uniformizable to $\mbox{\boldmath
$C$}\times \mbox{\boldmath $P$}^1$ in the sense of R. Kobayashi, S. Nakamura and F. Sakai
$($see {\em \cite{kns}, \cite{kobayashi}}$)$.
\end{co}

{\it Proof.} We only have to consider the case where $B\simeq \mbox{\boldmath $P$}^1$. We define a
$\mbox{\boldmath $Q$}$-divisor 
$\delta$ on $B$ as $\delta:=\sum_{t\in B}\{(m_tb_t-1)/m_tb_t\}t$, where 
$m_t:=\mbox{mult}_{C_{\varphi}(t)}\varphi^{\ast}(t)$ and
$b_t:=(1-\mbox{mult}_{C_{\varphi}(t)}\Delta)^{-1}$ for $t\in B$. Then by checking the list in
Proposition\ref{pr:cltypeII}, we see that $\deg \delta=2$. Put
$r:=\mbox{Ind}(K_{B}+\delta)$. We can also check that $r=\mbox{Ind}(K_S+\Delta)$.
Let $\mbox{\boldmath
$C$}(x)\subset K:=\mbox{\boldmath $C$}(x,y)$ be the field extension induced by $\varphi$. We
define a finite Galois field extension $K\subset L$ as follows. 
Take a rational function $\alpha(x)$ on $B$ such that
$\mbox{div }\alpha(x)=r(K_{B}+\delta)$. If 
$\Xi_h=0$, we put
$L:=K(\root{r}\of{\alpha(x)})$. If $\Xi_h\neq 0$. Let $\Xi_{h,\eta}$ be the
restriction of $\Xi_h$ to the generic fibre $S_{\eta}\simeq \mbox{\boldmath
$P$}^1(\mbox{\boldmath $C$}(x))$ of $\varphi$. Then the defining equation of $\Xi_{h,\eta}$ can be
written as $y^2=\beta(x)$. Let $\Xi_{h,\bar\eta}$ be the pull back of $\Xi_{h,\eta}$ on
$\mbox{\boldmath
$P$}^1(\mbox{\boldmath $C$}(x,\sqrt{\beta(x)}))$. By substituting $y$ for
$y^{\prime}=(y-\sqrt{\beta(x)})/(y+\sqrt{\beta(x)})$, we may assume that $\Xi_{h,\bar\eta}$
corresponds to $y=0$, $\infty$. We put $L:=K(\root{r}\of{\alpha(x)},\sqrt{\beta(x)},\sqrt{y})$. We
note that $L$ is a Kummer extension of $K$ by the construction. Let
$f:P\rightarrow S$ be the normalization of $S$ in $L$. We show that
$P$ is a desired relatively minimal elliptic ruled surface. Let
$f_1:P_1\rightarrow S$ be the normalization of $S$ in
$K(\root{r}\of{\alpha(x)})$ and $\varphi_1:P_1\rightarrow B_1$ be the induced
morphism from $\varphi$ by taking the Stein factorization. By our construction, $B_1$ is nothing
but the log canonical cover of $(B,\delta)$, hence $B_1$ is an elliptic curve. Since locally on
$S$, we can write $\mbox{div }\varphi^{\ast}\alpha(x)=\{r(b_t-1)/b_t\}C_{\varphi}(t)$,
$f_1^{\ast}C_{\varphi}(t)$ is Cartier, hence $\varphi_1$ is smooth and the ramification index at
the generic point of prime divisors in $f_1^{-1}(C_{\varphi}(t))$ is $b_t$. Thus in the case where 
$\Xi_h=0$, we get the assertion. Assume that $\Xi_h\neq 0$. Then we see that
$f_1^{-1}(\Xi_h)$ is smooth and $e_{\mbox{top}}(f_1^{-1}(\Xi_h))=0$ from
$K_{P_1}+f_1^{-1}(\Gamma)+(1/2)f_1^{-1}(\Xi_h)=f_1^{\ast}(K_S+\Delta)\sim_
{\mbox{\boldmath $Q$}}0$ and $(f_1^{-1}(\Xi_h))^2=0$. Let
$f_2:P_2\rightarrow P_1$ be the normalization of
$P_1$ in
$K(\root{r}\of{\alpha(x)},\sqrt{\beta(x)})$ and $\varphi_2:P_2\rightarrow B_2$
be the induced morphism from $\varphi_1$. If $K(\root{r}\of{\alpha(x)},\sqrt{\beta(x)})\neq
K(\root{r}\of{\alpha(x)})$, then $f_1^{-1}(\Xi_h)$ is an elliptic curve and
$B_2\simeq f_1^{-1}(\Xi_h)$ by the construction. Therefore, $P_2$ is also
a relatively minimal elliptic ruled surface and $f_2^{-1}f_1^{-1}(\Xi_h)$ is a disjoint
union of two sections over which $f_3:P\rightarrow P_2$ ramifies with the index two, where $f_3$ is
the normalization map of $P_2$ in $L$. Thus we conclude that
$P$ is a relatively minimal elliptic ruled surface with a section
$f^{-1}(\lfloor\Delta\rfloor)$. The last assertion follows from
Lemma~\ref{lm:universal}.\hfill\bsquare

\subsubsection{Case Type III}

\begin{lm}[c.f. \cite{an}, \cite{catanese} and \cite{kobayashi}, \S
5. Appendix]\label{lm:i2} Let
$(S,\Delta;p)$ be a germ of normal log surface with $\mbox{\boldmath
$Q$}$-boundary $\Delta$ at $p$ such that $\lfloor\Delta\rfloor=0$. 
Assume that $(S,\Delta)$ is purely log terminal and $\mbox{\em ind}_{p}(K_S+\Delta)=2$. 
Then there exists a resolution $f:T\rightarrow S$ such that 
$f^{-1}(p)\cup \mbox{\em Supp }f^{-1}_{\ast}\Delta$ has only normal crossing singularities 
whose dual graph is one of the following types. 
\begin{description}
\item[$(1)$] Case $p\in S$ smooth.
\vskip 10mm
\unitlength 0.7mm
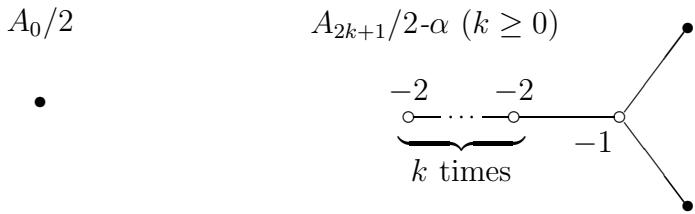
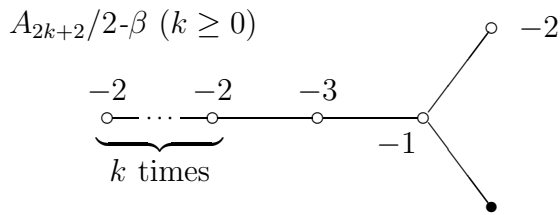
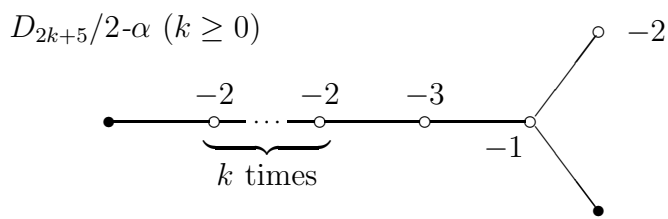
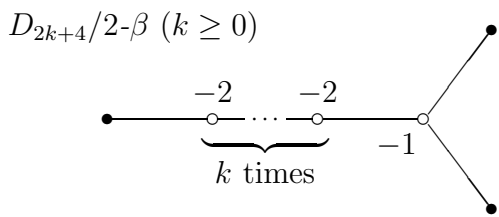
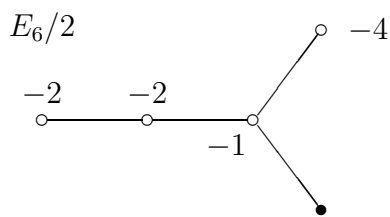
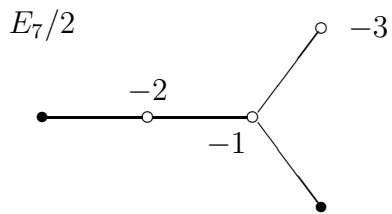
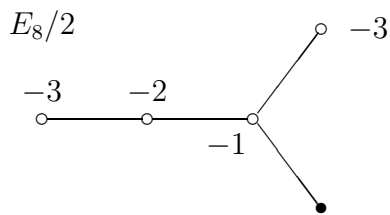
\begin{picture}(5,5)(-10,0)
\put(0,15){\makebox(0,0){$A_0/2$}}
\put(0,0){\circle*{2}}
\end{picture}

\begin{picture}(60,30)(-80,-30)
\put(0,15){\makebox(10,5){$A_{2k+1}/2$-$\alpha$ $(k\geq 0)$}}
\put(0,0){\circle{2}}
\put(0,5){\makebox(0,0){$-2$}}
\put(1,0){\line(1,0){5}}
\put(8,0){\circle*{0.5}}
\put(10,0){\circle*{0.5}}
\put(10,-5){\makebox(0,0){$\underbrace{\qquad \qquad}$}}
 \put(10,-10){\makebox(0,0){$k$ \mbox{\rm times}}}
\put(12,0){\circle*{0.5}}
\put(14,0){\line(1,0){5}}
\put(20,0){\circle{2}}\put(20,5){\makebox(0,0){$-2$}}
\put(21,0){\line(1,0){18}}
\put(40,0){\circle{2}}\put(35,-5){\makebox(0,0){$-1$}}
\put(41,1){\line(3,4){11.5}}
\put(41,-1){\line(3,-4){11.5}}
\put(53,17){\circle*{2}}
\put(53,-17){\circle*{2}}
\end{picture}

\unitlength 0.7mm
\begin{picture}(100,30)(-10,0)
\put(0,15){\makebox(10,5){$A_{2k+2}/2$-$\beta$ $(k\geq 0)$}}
\put(0,0){\circle{2}}
\put(0,5){\makebox(0,0){$-2$}}
\put(1,0){\line(1,0){5}}
\put(8,0){\circle*{0.5}}
\put(10,0){\circle*{0.5}}
\put(10,-5){\makebox(0,0){$\underbrace{\qquad \qquad}$}}
\put(10,-10){\makebox(0,0){$k$ \mbox{\rm times}}}
\put(12,0){\circle*{0.5}}
\put(14,0){\line(1,0){5}}
\put(20,0){\circle{2}}\put(20,5){\makebox(0,0){$-2$}}
\put(21,0){\line(1,0){18}}
\put(40,0){\circle{2}}\put(40,5){\makebox(0,0){$-3$}}
\put(41,0){\line(1,0){18}}
\put(60,0){\circle{2}}\put(55,-5){\makebox(0,0){$-1$}}
\put(61,1){\line(3,4){11.5}}
\put(61,-1){\line(3,-4){11.5}}
\put(73,17){\circle{2}}\put(82,17){\makebox(0,0){$-2$}}
\put(73,-17){\circle*{2}}
\end{picture}

\unitlength 0.7mm
\begin{picture}(100,30)(-130,-30)
\put(0,15){\makebox(10,5){$D_{2k+5}/2$-$\alpha$ $(k\geq 0)$}}
\put(0,0){\circle*{2}}
\put(0,0){\line(1,0){19}}
\put(20,0){\circle{2}}
\put(20,5){\makebox(0,0){$-2$}}
\put(21,0){\line(1,0){5}}
\put(28,0){\circle*{0.5}}
\put(30,0){\circle*{0.5}}
\put(30,-5){\makebox(0,0){$\underbrace{\qquad \qquad}$}}
\put(30,-10){\makebox(0,0){$k$ \mbox{\rm times}}}
\put(32,0){\circle*{0.5}}
\put(34,0){\line(1,0){5}}
\put(40,0){\circle{2}}\put(40,5){\makebox(0,0){$-2$}}
\put(41,0){\line(1,0){18}}
\put(60,0){\circle{2}}\put(60,5){\makebox(0,0){$-3$}}
\put(61,0){\line(1,0){18}}
\put(80,0){\circle{2}}\put(75,-5){\makebox(0,0){$-1$}}
\put(81,1){\line(3,4){11.5}}
\put(81,-1){\line(3,-4){11.5}}
\put(93,17){\circle{2}}\put(102,17){\makebox(0,0){$-2$}}
\put(93,-17){\circle*{2}}
\end{picture}

\unitlength 0.7mm
\begin{picture}(100,30)(-10,0)
\put(0,15){\makebox(10,5){$D_{2k+4}/2$-$\beta$ $(k\geq 0)$}}
\put(0,0){\circle*{2}}
\put(0,0){\line(1,0){19}}
\put(20,0){\circle{2}}\put(20,5){\makebox(0,0){$-2$}}
\put(21,0){\line(1,0){5}}
\put(28,0){\circle*{0.5}}
\put(30,0){\circle*{0.5}}
\put(30,-5){\makebox(0,0){$\underbrace{\qquad \qquad}$}}
\put(30,-10){\makebox(0,0){$k$ \mbox{\rm times}}}
\put(32,0){\circle*{0.5}}
\put(34,0){\line(1,0){5}}
\put(40,0){\circle{2}}\put(40,5){\makebox(0,0){$-2$}}
\put(41,0){\line(1,0){18}}
\put(60,0){\circle{2}}\put(55,-5){\makebox(0,0){$-1$}}
\put(61,1){\line(3,4){11.5}}
\put(61,-1){\line(3,-4){11.5}}
\put(73,17){\circle*{2}}
\put(73,-17){\circle*{2}}
\end{picture}

\unitlength 0.7mm
\begin{picture}(100,30)(-130,-30)
\put(0,15){\makebox(0,5){$E_6/2$}}
\put(0,0){\circle{2}}
\put(0,5){\makebox(0,0){$-2$}}
\put(1,0){\line(1,0){18}}
\put(20,0){\circle{2}}
\put(20,5){\makebox(0,0){$-2$}}
\put(21,0){\line(1,0){18}}
\put(40,0){\circle{2}}\put(35,-5){\makebox(0,0){$-1$}}
\put(41,1){\line(3,4){11.5}}
\put(41,-1){\line(3,-4){11.5}}
\put(53,17){\circle{2}}\put(62,17){\makebox(0,0){$-4$}}
\put(53,-17){\circle*{2}}
\end{picture}

\unitlength 0.7mm
\begin{picture}(60,30)(-10,0)
\put(0,15){\makebox(0,5){$E_7/2$}}
\put(0,0){\circle*{2}}
\put(1,0){\line(1,0){18}}
\put(20,0){\circle{2}}
\put(20,5){\makebox(0,0){$-2$}}
\put(21,0){\line(1,0){18}}
\put(40,0){\circle{2}}\put(35,-5){\makebox(0,0){$-1$}}
\put(41,1){\line(3,4){11.5}}
\put(41,-1){\line(3,-4){11.5}}
\put(53,17){\circle{2}}\put(62,17){\makebox(0,0){$-3$}}
\put(53,-17){\circle*{2}}
\end{picture}

\unitlength 0.7mm
\begin{picture}(60,30)(-130,-40)
\put(0,15){\makebox(0,5){$E_8/2$}}
\put(0,0){\circle{2}}
\put(0,5){\makebox(0,0){$-3$}}
\put(1,0){\line(1,0){18}}
\put(20,0){\circle{2}}
\put(20,5){\makebox(0,0){$-2$}}
\put(21,0){\line(1,0){18}}
\put(40,0){\circle{2}}\put(35,-5){\makebox(0,0){$-1$}}
\put(41,1){\line(3,4){11.5}}
\put(41,-1){\line(3,-4){11.5}}
\put(53,17){\circle{2}}\put(62,17){\makebox(0,0){$-3$}}
\put(53,-17){\circle*{2}}
\end{picture}

\item[$(2)$] Case $p\in S$ singular.
\vskip 10mm
\unitlength 0.7mm
\begin{picture}(5,5)(-10,0)
\put(0,15){\makebox(0,5){$A_1/2$-$\gamma$}}
\put(0,0){\circle{2}}
\put(0,5){\makebox(0,0){$-4$}}
\end{picture}

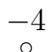
\begin{picture}(100,30)(-80,-35)
\put(0,15){\makebox(10,5){$A_{2k+3}/2$-$\delta$ $(k\geq 0)$}}
\put(0,0){\circle{2}}
\put(0,5){\makebox(0,0){$-3$}}
\put(1,0){\line(1,0){18}}
\put(20,0){\circle{2}}
\put(20,5){\makebox(0,0){$-2$}}
\put(21,0){\line(1,0){5}}
\put(28,0){\circle*{0.5}}
\put(30,0){\circle*{0.5}}
\put(30,-5){\makebox(0,0){$\underbrace{\qquad \qquad}$}}
\put(30,-10){\makebox(0,0){$k$ \mbox{\rm times}}}
\put(32,0){\circle*{0.5}}
\put(34,0){\line(1,0){5}}
\put(40,0){\circle{2}}\put(40,5){\makebox(0,0){$-2$}}
\put(41,0){\line(1,0){18}}
\put(60,0){\circle{2}}\put(60,5){\makebox(0,0){$-3$}}
\end{picture}

\unitlength 0.7mm
\begin{picture}(100,30)(-10,0)
\put(0,15){\makebox(10,5){$A_{2k+2}/2$-$\epsilon$ $(k\geq 0)$}}
\put(0,0){\circle*{2}}
\put(1,0){\line(1,0){18}}
\put(20,0){\circle{2}}
\put(20,5){\makebox(0,0){$-2$}}
\put(21,0){\line(1,0){5}}
\put(28,0){\circle*{0.5}}
\put(30,0){\circle*{0.5}}
\put(30,-5){\makebox(0,0){$\underbrace{\qquad \qquad}$}}
\put(30,-10){\makebox(0,0){$k$ \mbox{\rm times}}}
\put(32,0){\circle*{0.5}}
\put(34,0){\line(1,0){5}}
\put(40,0){\circle{2}}\put(40,5){\makebox(0,0){$-2$}}
\put(41,0){\line(1,0){18}}
\put(60,0){\circle{2}}\put(60,5){\makebox(0,0){$-3$}}
\end{picture}

\begin{picture}(100,30)(-130,-35)
\put(0,15){\makebox(10,5){$A_{2k+1}/2$-$\zeta$ $(k\geq 1)$}}
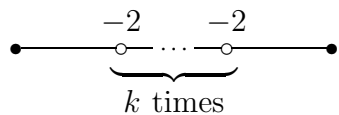
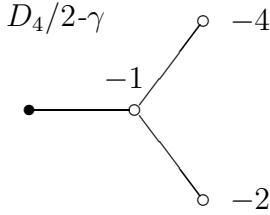
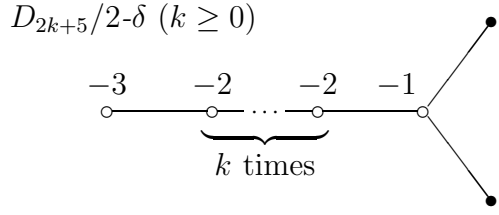
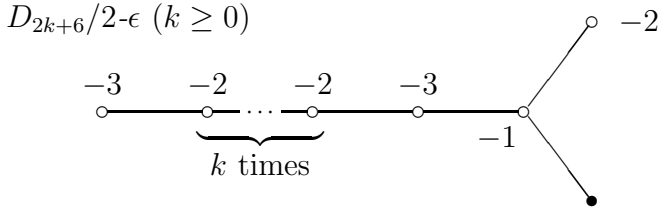
\put(0,0){\circle*{2}}
\put(1,0){\line(1,0){18}}
\put(20,0){\circle{2}}
\put(20,5){\makebox(0,0){$-2$}}
\put(21,0){\line(1,0){5}}
\put(28,0){\circle*{0.5}}
\put(30,0){\circle*{0.5}}
\put(30,-5){\makebox(0,0){$\underbrace{\qquad \qquad}$}}
\put(30,-10){\makebox(0,0){$k$ \mbox{\rm times}}}
\put(32,0){\circle*{0.5}}
\put(34,0){\line(1,0){5}}
\put(40,0){\circle{2}}\put(40,5){\makebox(0,0){$-2$}}
\put(41,0){\line(1,0){18}}
\put(60,0){\circle*{2}}
\end{picture}

\unitlength 0.7mm
\begin{picture}(100,40)(-10,0)
\put(0,15){\makebox(10,5){$D_4/2$-$\gamma$}}
\put(0,0){\circle*{2}}
\put(0,0){\line(1,0){19}}
\put(20,0){\circle{2}}\put(18,6){\makebox(0,0){$-1$}}
\put(21,1){\line(3,4){11.5}}
\put(21,-1){\line(3,-4){11.5}}
\put(33,17){\circle{2}}\put(42,17){\makebox(0,0){$-4$}}
\put(33,-17){\circle{2}}\put(42,-17){\makebox(0,0){$-2$}}
\end{picture}

\unitlength 0.7mm
\begin{picture}(100,40)(-130,-40)
\put(0,15){\makebox(10,5){$D_{2k+5}/2$-$\delta$ $(k\geq 0)$}}
\put(0,0){\circle{2}}
\put(0,5){\makebox(0,0){$-3$}}
\put(1,0){\line(1,0){18}}
\put(20,0){\circle{2}}
\put(20,5){\makebox(0,0){$-2$}}
\put(21,0){\line(1,0){5}}
\put(28,0){\circle*{0.5}}
\put(30,0){\circle*{0.5}}
\put(30,-5){\makebox(0,0){$\underbrace{\qquad \qquad}$}}
\put(30,-10){\makebox(0,0){$k$ \mbox{\rm times}}}
\put(32,0){\circle*{0.5}}
\put(34,0){\line(1,0){5}}
\put(40,0){\circle{2}}\put(40,5){\makebox(0,0){$-2$}}
\put(41,0){\line(1,0){18}}
\put(60,0){\circle{2}}\put(55,5){\makebox(0,0){$-1$}}
\put(61,1){\line(3,4){11.5}}
\put(61,-1){\line(3,-4){11.5}}
\put(73,17){\circle*{2}}
\put(73,-17){\circle*{2}}
\end{picture}

\unitlength 0.7mm
\begin{picture}(100,40)(-10,0)
\put(0,15){\makebox(10,5){$D_{2k+6}/2$-$\epsilon$ $(k\geq 0)$}}
\put(0,0){\circle{2}}
\put(0,5){\makebox(0,0){$-3$}}
\put(1,0){\line(1,0){18}}
\put(20,0){\circle{2}}
\put(20,5){\makebox(0,0){$-2$}}
\put(21,0){\line(1,0){5}}
\put(28,0){\circle*{0.5}}
\put(30,0){\circle*{0.5}}
\put(30,-5){\makebox(0,0){$\underbrace{\qquad \qquad}$}}
\put(30,-10){\makebox(0,0){$k$ \mbox{\rm times}}}
\put(32,0){\circle*{0.5}}
\put(34,0){\line(1,0){5}}
\put(40,0){\circle{2}}\put(40,5){\makebox(0,0){$-2$}}
\put(41,0){\line(1,0){18}}
\put(60,0){\circle{2}}\put(60,5){\makebox(0,0){$-3$}}
\put(61,0){\line(1,0){18}}
\put(80,0){\circle{2}}\put(75,-5){\makebox(0,0){$-1$}}
\put(81,1){\line(3,4){11.5}}
\put(81,-1){\line(3,-4){11.5}}
\put(93,17){\circle{2}}\put(102,17){\makebox(0,0){$-2$}}
\put(93,-17){\circle*{2}}
\end{picture}

\vskip 20mm

\end{description}

In the above dual graphs $\bullet$ denotes the strict transform of $\mbox{\em Supp }\{\Delta\}$
by $f$ and $\circ$ denotes an $f$-exceptional rational curve. The number attached on $\circ$ is
its self-intersection number.

\end{lm}

{\it Proof. }We give a proof for the convenience of the readers. In the case where $S$ is smooth at $p$, the proof reduces to the well know
fact about the resolution of plane curves, so we omit the proof in this case. Assume that $S$ is
singular at $p$ and let
$\mu:M\rightarrow S$ be the minimal resolution of
$p\in S$. Put $\Xi:=\mbox{Supp }\Delta$ and $\Xi^{\prime}:=\mu_{\ast}^{-1}\Xi$. Let
$\mu^{-1}(p)=\cup_{i=1}^lE_i$ be the irreducible decomposition. Then we can write
$K_M+(1/2)\Xi^{\prime}+\sum_{i=1}^l a_iE_i=\mu^{\ast}(K_S+\Delta)$, where $a_i=0$ or $1/2$ for any
$i$ by the assumption. Since we have $(K_M+(1/2)\Xi^{\prime}+\sum_{i=1}^l a_iE_i,E_j)=0$,
we have 
$(\Xi^{\prime}+\sum_{i\neq j}2a_iE_i,E_j)=4+2(1-a_j)E_j^2$ for any $j$. Assume that $a_{i_0}=0$ for
some $i_0$, say for $i_0=1$. Then we see that $(\Xi^{\prime}+\sum_{i\neq
1}2a_iE_i,E_1)=0$ since $E_1^2\leq -2$, which implies that
$\mu^{-1}(p)=\cup_{a_i=0}E\amalg\cup_{a_i=1/2}E_i$. Since $\mu^{-1}(p)$ is connected, we infer that 
$a_i=0$ and $(\Xi^{\prime},E_i)=0$ for any $i$, hence $\mbox{ind}_p(K_S+\Delta)=1$, which is
absurd. Thus we obtain $a_i=1/2$ for any $i$ and $(\Xi^{\prime}+\sum_{i\neq
j}E_i,E_j)=4+E_j^2$ for any $j$, which implies that $E_i^2=-2$, $-3$ or $-4$ for any $i$. Assume
that there exists $E_{i_0}$, say $E_1$, such that $E_{1}^2=-4$. Then since we have $(\Xi^{\prime}+\sum_{i\neq
j}E_i,E_1)=0$, we see that
$\mu^{-1}(p)=E_1$ and $(\Xi^{\prime},E_1)=0$, which implies that we are in the case
$A_1/2\mbox{-}\gamma$. Thus we may assume that $E_i^2=-2$, $-3$. We note that we have
$(\Xi^{\prime}+\sum_{i\neq j}E_i,E_j)=1$ if $E_j^2=-3$ and $2$ if $E_j^2=-2$. Assume that there
exists $E_{i_0}$, say $E_1$, such that $E_{1}^2=-3$. Then it is easily seen that we are in the
cases $A_{2k+3}/2\mbox{-}\delta$ $(k\geq 0)$ or $A_{2k+2}/2\mbox{-}\epsilon$ $(k\geq 0)$. Thus we
may assume that
$E_i^2=-2$ for any $i$. Assume that $(\Xi^{\prime},E_{i_0})=1$ for some $i_0$. Then by the
inductive argument, we see that we are in the cases $A_{2k+1}/2\mbox{-}\zeta$ $(k\geq 2)$. Assume
that $(\Xi^{\prime},E_{i_0})=2$ for some $i_0$. Then we have $\mu^{-1}(p)=E_{i_0}$ since 
$(\sum_{i\neq i_0}E_i,E_{i_0})=0$, from which we deduce that we are in the cases
$A_3/2\mbox{-}\zeta$, $D_4/2\mbox{-}\gamma$, $D_{2k+5}/2\mbox{-}\delta$ $(k\geq 0)$ and
$D_{2k+6}/2\mbox{-}\epsilon$ $(k\geq 0)$. \hfill\bsquare

\begin{df}{\em A normal complex projective surface $S$ is called a {\it rank one log del Pezzo surface} if 
$S$ has at worst quotient singularities, $-K_S$ is ample and the Picard number one and moreover, 
if $S$ is Gorenstein, $S$ is called a {\it rank one Gorenstein log del Pezzo surface}.
}\end{df}

Singular rank one Gorenstein log del Pezzo surfaces are classified into the following $27$ types listed as follows
(see \cite{furushima}, \cite{miyanishi-zhang}). We need the list below to study type III log surfaces.
\begin{center}
Table IV\\
\begin{tabular}{l|c|c|c|l|c|c|c}\hline
 & $(-K_S)^2$ & $\mbox{ Sing }S$ & $e_{\mbox{\rm orb}}(S)$  & &$(-K_S)^2$  &$\mbox{ Sing }S$ &$e_{\mbox{\rm orb}}(S)$\\ \hline \hline
$(1)$ & $8$ & $A_1$ & $5/2$  & $(15)$ & $1$ & $E_8$ & $241/120$\\ \hline
$(2)$ & $6$ & $A_2+A_1$ & $11/6$ & $(16)$ & $1$ & $E_7+A_1$ & $73/48$\\ \hline
$(3)$ & $5$ & $A_4$ & $11/5$ & $(17)$ & $1$ & $E_7+A_2$ & $65/48$\\ \hline
$(4)$ & $4$ & $D_5$ & $25/12$ & $(18)$ & $1$ & $A_8$ & $19/9$\\ \hline
$(5)$ & $4$ & $A_3+2A_1$ & $5/4$ & $(19)$ & $1$ & $A_7+A_1$ & $13/8$\\ \hline
$(6)$ & $3$ & $E_6$ & $49/24$ & $(20)$ & $1$ & $A_5+A_2+A_1$ & $1$\\ \hline
$(7)$ & $3$ & $A_5+A_1$ & $5/3$ & $(21)$ & $1$ & $D_8$ & $49/24$\\ \hline
$(8)$ & $3$ & $3A_2$ & $1$ & $(22)$ & $1$ & $D_6+2A_1$ & $17/16$\\ \hline
$(9)$ & $2$ & $E_7$ & $97/48$ & $(23)$ & $1$ & $D_5+A_3$ & $4/3$\\ \hline
$(10)$ & $2$ & $D_6+A_1$ & $25/16$ & $(24)$ & $1$ & $2D_4$ & $5/4$\\ \hline
$(11)$ & $2$ & $A_7$ & $17/8$ & $(25)$ & $1$ & $4A_2$ & $1/3$\\ \hline
$(12)$ & $2$ & $D_4+A_3$ & $11/8$ & $(26)$ & $1$ & $2A_3+2A_1$ & $1/2$\\ \hline
$(13)$ & $2$ & $A_5+A_2$ & $3/2$ & $(27)$ & $1$ & $2A_4$ & $7/5$\\ \hline
$(14)$ & $2$ & $2A_3+A_1$ & $1$\\ \cline{1-4}
\end{tabular}
\end{center}

\begin{pr}\label{pr:cltypeIII-1} Let $(S,\Delta)$ be of type $\mbox{\em III}$. 
Assume that $K_S+\Delta$ is Cartier and that $e_{\mbox{\rm{orb}}}(S\setminus \Delta)=0$. 
Then there exists a birational morphism $\tau:S\rightarrow S^{\flat}$ 
which is composed of contractions of $(-1)$-curves with $\mbox{\em Exc }\tau\subset \Delta$ to a normal projective surface 
$S^{\flat}$ and $(S^{\flat},\Delta^{\flat})$ satisfies one of the followings. 
\begin{description}
\item[ $(a)$ ]  $S^{\flat}\simeq \Sigma_d$, where $\Sigma_d$ is a Hirzebruch surface and 
$\Delta^{\flat}=\sum_{i=1}^{4}\Gamma^{\flat}_i$ where $\Gamma^{\flat}_1$ and $\Gamma^{\flat}_2$ are 
two disjoint sections, $\Gamma^{\flat}_3$ and $\Gamma^{\flat}_4$ are two fibres.

\item[ $(b)$ ]  $S^{\flat}$ is a rank one Gorenstein log del Pezzo surface with 
$\mbox{\em Sing }S^{\flat}=3A_2$, $A_1+2A_3$, or $A_1+A_2+A_5$ and $\Delta^{\flat}$ 
is a rational curve with only one node as its singularities.

\item[ $(c)$ ] $S^{\flat}$ has a structure of a conic fibration
$\varphi^{\flat}:S^{\flat}\rightarrow \mbox{\boldmath $P$}^1$ with
$\rho(S^{\flat}/\mbox{\boldmath $P$}^1)=1$ and 
$\mbox{\em Typ}(S^{\flat},\Delta^{\flat};\varphi^{\flat})=
(2(\mbox{\em I}\mbox{-}2)_1+(\mbox{\em II}\mbox{-}1)_{\infty};(\mbox{\em II}\mbox{-}1)_1)$. 
$\Delta^{\flat}=\Gamma^{\flat}_1+\Gamma^{\flat}_2$, where 
$\Gamma^{\flat}_1$ is a smooth rational curve such that $(\Gamma^{\flat}_1,\varphi^{\flat
\ast}(t))=2$ for $t\in \mbox{\boldmath $P$}^1$ and that $(\Gamma^{\flat}_1)^2=0$ and 
$\Gamma^{\flat}_2$ is a fibre of $\varphi^{\flat}$. 
\end{description}
\end{pr}
{\it Proof.} As we have already seen in Lemma~\ref{lm:K-cont}, by applying the log minimal program on $S$ with respect to $K_S$, we get a birational morphism 
$\tau:S\rightarrow S^{\flat}$ with $\mbox{Exc }\tau\subset \Delta$ to a normal projective surface $S^{\flat}$ ($1$) 
which is a rank one Gorenstein log del Pezzo surface or ($2$) which admits a structure of conic fibration 
$\varphi^{\flat}:S^{\flat}\rightarrow \mbox{\boldmath $P$}^1$ 
with $\rho(S^{\flat}/\mbox{\boldmath $P$}^1)=1$. Assume that we are in the case ($1$). 
Since we have $e_{\mbox{\rm orb}}(S^{\flat})=e_{\mbox{\rm top}}(\Delta^{\flat})\in \mbox{\boldmath $N$}$, 
we see that $S^{\flat}\simeq \mbox{\boldmath $P$}^2$ and $e_{\mbox{\rm top}}(\Delta^{\flat})=3$ or 
$\mbox{Sing }S^{\flat}=3A_2$, $A_1+2A_3$, or $A_1+A_2+A_5$ and $e_{\mbox{\rm top}}(\Delta^{\flat})=1$ 
by Table IV. Thus we are in the cases ($a$) and $(b)$. 
Assume that we are in the case ($2$). By the assumption, we have
$
4=e_{\mbox{\rm top}}(\Delta^{\flat})+\nu((\mbox{I}\mbox{-}2)_1)+\sum_{k\geq
1}\{(4k-1)/(4k)\}\nu((\mbox{II}\mbox{-}3)_{1,k}).
$
where $\nu((\mbox{I}\mbox{-}2)_1)$ and $\nu((\mbox{II}\mbox{-}3)_{1,k})$ denote the number of fibres of type
$(\mbox{I}\mbox{-}2)_1$ and 
$(\mbox{II}\mbox{-}3)_{1,k}$ respectively.
 Assume that $\Delta^{\flat}$ contains an irreducible component $\Gamma^{\flat}_0$ with 
$(\Gamma^{\flat}_0,\varphi^{\flat\ast}(t))=2$ for $t\in \mbox{\boldmath $P$}^1$. Since
$\Delta^{\flat}-\Gamma^{\flat}_0$  is composed of fibres of $\varphi^{\flat}$, we can write
$\Delta^{\flat}=\Gamma^{\flat}_0+\Gamma^{\flat}_1$ where $\Gamma^{\flat}_1$ is a fibre of
$\varphi^{\flat}$ and we see that $e_{\mbox{\rm top}}(\Delta^{\flat})=2$ and that
$
\nu((\mbox{I}\mbox{-}2)_1)+\sum_{k\geq 1}\{(4k-1)/(4k)\}\nu((\mbox{II}\mbox{-}3)_{1,k})=2.
$
On the other hand, since $\varphi^{\flat}$ induces a double cover 
$\varphi^{\flat}:\Gamma^{\flat}_0\simeq \mbox{\boldmath $P$}^1\rightarrow \mbox{\boldmath $P$}^1$,
we have
$
\nu((\mbox{I}\mbox{-}2)_1)+\sum_{k\geq 1}\nu((\mbox{II}\mbox{-}3)_{1,k})\leq 2
$
by the Hurwitz's formula. So we have $\nu((\mbox{I}\mbox{-}2)_1)=2$ and $\nu((\mbox{II}\mbox{-}3)_{1,k})=0$ for any $k$, 
hence we are in the case $(c)$. Assume that there exist two sections of $\varphi^{\flat}$,
$\Gamma^{\flat}_1$ and $\Gamma^{\flat}_2$.  Then we see that
$\nu((\mbox{I}\mbox{-}2)_1)=\nu((\mbox{II}\mbox{-}3)_{1,k})=0$ for any $k$ and that 
$\varphi^{\flat}:S^{\flat}\rightarrow \mbox{\boldmath $P$}^1$ is smooth and $e_{\mbox{\rm
top}}(\Delta^{\flat})=4$. Thus we are in the case $(a)$.
\hfill\bsquare
\vskip 5mm
\par
Let $(S,\Delta)$ be of type \mbox{III} and assume that $K_S+\Delta$ is not Cartier and  
let $\lfloor\Delta\rfloor=\sum_{i=0}^{b}\Gamma_i$ be the irreducible decomposition. Then $b\geq
1$ and $\Gamma_i$'s form a linear chain of smooth rational curves, that is,    
$(\Gamma_i,\Gamma_j)=0$, if $|i-j|>1$, $1$, if $|i-j|=1$ and 
$\{p\in \lfloor\Delta\rfloor|\mbox{ ind}_p(K_S+\Delta)=2\}$, namely the set of points 
$p\in \lfloor\Delta\rfloor$ at which $K_S+\Delta$ is not Cartier consists of four points $p_i$ $(1\leq i\leq 4)$, 
where $p_i\in \Gamma_0$ for $i=1$, $2$ and $p_i\in \Gamma_b$ for $i=3$, $4$ by \cite{ohno}, Lemma 3.2 (1). 
Comparing the Euler numbers, we have 
\begin{equation}
e_{\mbox{top}}(\tilde \Delta)=2e_{\mbox{top}}(\lfloor\Delta\rfloor)-4.\label{eqn:topofb1}
\end{equation}
Let $(S^{\flat},\Delta^{\flat})$ be a log surface obtained from $(S,\Delta)$ by applying log
minimal program with respect to $K_{S}+\{\Delta\}$. Write
$\Delta^{\flat}=\Gamma^{\flat}+(1/2)\Xi^{\flat}$, where 
$\Gamma^{\flat}:=\lfloor\Delta^{\flat}\rfloor$ and $\Xi^{\flat}:=2\{\Delta^{\flat}\}$. Put
$d:=(-K_{\tilde S^{\flat}})^2$. 
We note that we have 
\begin{equation}
e_{\mbox{top}}(\Xi^{\flat \nu})=
\deg\mbox{Diff}_{\Xi^{\flat \nu}}(\Gamma^{\flat})-\frac{1}{2}(\Xi^{\flat})^2\label{eqn:topofx}
\end{equation}
from $K_{S^{\flat}}+\Gamma^{\flat}+(1/2)\Xi^{\flat}\sim_{\mbox{\boldmath $Q$}}0$. Let
$\mu^{\flat}:M^{\flat}\rightarrow S^{\flat}$ be the minimal resolution of $S^{\flat}$ and put 
$\delta K^2:=K_{S^{\flat}}^2-K_{M^{\flat}}^2$, $\delta \rho:=\rho(M^{\flat})-\rho(S^{\flat})$.
Then we have 
$$
K_{S^{\flat}}^2=K_{M^{\flat}}^2+\delta K^2=10-\rho(M^{\flat})+\delta K^2=10-\rho(S^{\flat})-\delta
\rho+\delta K^2. $$
On the other hand, we have 
$$
K_{S^{\flat}}^2=(\Gamma^{\flat}+\frac{1}{2}\Xi^{\flat})^2=\frac{d}{2}+
(\Gamma^{\flat},\Xi^{\flat})+\frac{1}{4}(\Xi^{\flat})^2,
$$
hence
\begin{equation}
\delta K^2-\delta \rho=
\frac{d}{2}+(\Gamma^{\flat},\Xi^{\flat})+\frac{1}{4}(\Xi^{\flat})^2+\rho(S^{\flat})-10.\label{eqn:k-r}
\end{equation}
If we have $\rho(S^{\flat})=1$, 
for any divisor $D$ on $S^{\flat}$, there exists $a\in \mbox{\boldmath $Q$}$ such that 
$D$ is numerically equivalent to $a\Gamma^{\flat}$. 
Noting that $D^2=a^2(\Gamma^{\flat})^2=(1/2)da^2$ by 
$(\Gamma^{\flat})^2=(1/2)(\tilde \Delta^{\flat})^{2}=(1/2)(-K_{\tilde S^{\flat}})^2=d/2$ and that 
$(\Gamma^{\flat},D)=a(\Gamma^{\flat})^2=(1/2)da$, we deduce that 
\begin{equation}
D^2=(2/d)(\Gamma^{\flat},D)^2.\label{eqn:r1}
\end{equation} 
Assume that 
$e_{\mbox{\rm{orb}}}(\tilde S\setminus \tilde\Delta)=0$. Then from Lemma~\ref{lm:K-cont}, we have
\begin{equation}
\rho(S)-\rho(S^{\flat})=
e_{\mbox{top}}(\lfloor\Delta\rfloor)-e_{\mbox{top}}(\lfloor\Delta^{\flat}\rfloor).\label{eqn:dr=de}
\end{equation}
Let
$\mbox{LC }(S^{\flat},\Delta^{\flat})$ be the set which consists of all point $p\in S^{\flat}$ such
that
$(S^{\flat},\Delta^{\flat})$ is not log terminal at $p$. We see that
$\mbox{LC}(S^{\flat},\Delta^{\flat})$ consists of at most two points and for $p\in
\mbox{LC}(S^{\flat},\Delta^{\flat})$, $\Gamma^{\flat}$ and $\Xi^{\flat}$ are smooth at $p$ and 
$(\Gamma^{\flat}, \Xi^{\flat})_p=2$. Let $\pi^{\flat}:\tilde S^{\flat}\rightarrow S^{\flat}$ be the log
canonical cover of
$(S^{\flat},\Delta^{\flat})$. Put $\tilde \Delta^{\flat}:=\pi^{\flat
-1}\Delta^{\flat}$ and $s:=\mbox{Card }\mbox{LC}(S^{\flat},\Delta^{\flat})$. Comparing the Euler
numbers again, we have
\begin{equation}
e_{\mbox{top}}(\tilde
\Delta^{\flat})=2e_{\mbox{top}}(\lfloor\Delta^{\flat}\rfloor)-4+s.\label{eqn:topofb2}
\end{equation}
If $\mbox{ind}(K_{S^{\flat}}+\Delta^{\flat})=1$ for $p\in S^{\flat}\setminus
\lfloor\Delta^{\flat}\rfloor$,  then $p\notin  \Xi^{\flat}$ and $S^{\flat}$ has at worst Du Val
singularity at $p$. In this situation, we shall say that $(S^{\flat},\Delta^{\flat})$ has 
singularity of type $A_n$ at $p$, if $S^{\flat}$ has Du Val singularity of type $A_n$ at $p$ for
instance.  Let $\psi:U\rightarrow B$ be a complex elliptic surface over a smooth complete curve
$B$ with a section such that $U$ is minimal over $B$.  We shall write 
\begin{eqnarray*}
\mbox{Typ}(U;\psi)&=&\sum_{k\geq 1}\nu(\mbox{I}_k)\mbox{I}_k+\sum_{l\geq 0}\nu(\mbox{I}_l^{\ast})\mbox{I}_l^{\ast}+
\nu(\mbox{II})\mbox{II}+\nu(\mbox{II}^{\ast})\mbox{II}^{\ast}+\nu(\mbox{III})\mbox{III}+\nu(\mbox{III}^{\ast})\mbox{III}^{\ast}\\
&&+\nu(\mbox{IV})\mbox{IV}+\nu(\mbox{IV}^{\ast})\mbox{IV}^{\ast},
\end{eqnarray*}
where $\nu({\cal T})$ denotes the number of singular fibres of type ${\cal T}$ in the Kodaira's 
notation (\cite{kodaira}). It turns out that the following lemma works better than the fixed point
formula in our cases.
\begin{lm}\label{lm:pic} Let $G$ denote a finite subgroup of $k$-automorphism group of normal
variety $X$ defined over an algebraically closed field $k$ and let $f:X\rightarrow
X/G$ be the quotient morphism. Assume that
$X/G$ is $\mbox{\boldmath
$Q$}$-factorial. Then $f$ induces the isomorphism 
$\mbox{\em Pic }(X/G)\otimes \mbox{\boldmath $Q$}\simeq (\mbox{\em Pic }X\otimes \mbox{\boldmath
$Q$})^G$.
\end{lm}
{\it Proof.} The injectivity is trivial. As for the surjectivity, we only have to note that
$f^{\ast}f_{\ast}D=\sum_{g\in G}g^{\ast}D$ for any $\mbox{\boldmath $Q$}$-divisor $D$ on $X$,
where $f_{\ast}D$ denotes the pushforward of $D$ by $f$ as a cycle.
\hfill\bsquare
\vskip 3mm
Assume that $\Xi^{\flat}\neq 0$ and that there exists a birational morphism $\eta:W\rightarrow
S^{\flat}$ from a normal surface
$W$ with $\Delta^W\geq 0$ such that $\eta^{-1}_{\ast}\Xi^{\flat}$ is nef and
$(\eta^{-1}_{\ast}\Xi^{\flat})^2=0$. By applying the log abundance theorem (see \cite{fujita} or
\cite{fong-mc/utah}, Theorem 11.1.3) to 
$K_{W}+\Delta^W+\varepsilon\Xi^{\flat}\sim_{\mbox{\boldmath $Q$}}\varepsilon\Xi^{\flat}$ for
sufficiently small positive rational number $\varepsilon$, we obtain a proper surjective connected
morphism $\varphi:W\rightarrow \mbox{\boldmath $P$}^1$. In what follows, we shall frequently use
this technique. 
In Propositions~\ref{pr:cltypeIII-2.1}, \ref{pr:cltypeIII-2.2}, \ref{pr:cltypeIII-2.3},
\ref{pr:cltypeIII-2.4} and \ref{pr:cltypeIII-2.5}, a log surface $(S,\Delta)$ is assumed to be of
type III such that $K_S+\Delta$ is not Cartier and $e_{\mbox{\rm{orb}}}(\tilde S\setminus
\tilde\Delta)=0$.
\begin{pr}\label{pr:cltypeIII-2.1} Assume that $\tilde S$
is smooth. Then there exists a birational morphism 
$\tau:S\rightarrow S^{\flat}$ with $\mbox{\em Exc }\tau\subset
\lfloor\Delta\rfloor$ such that there exists a conic fibration $\varphi^{\flat}:S^{\flat}\rightarrow
\mbox{\boldmath $P$}^1$ with $\rho(S^{\flat}/\mbox{\boldmath $P$}^1)=1$ and that $(S^{\flat},\Delta^{\flat})$ is log terminal, where
$\Delta^{\flat}:=\tau_{\ast}\Delta$. Moreover, there exists a
sequence of ${\cal S}$-elementary transformation such that for the resulting log surface 
$(\hat S^{\flat},\hat \Delta^{\flat})$, 
$\hat S^{\flat}\simeq \mbox{\boldmath $P$}^1\times \mbox{\boldmath $P$}^1$ with 
$\hat\varphi^{\flat}$ being the second projection and 
$\hat \Delta^{\flat}=\sum_{i=0}^2\hat \Gamma^{\flat}_i+(1/2)\sum_{j=1}^2\hat \Xi^{\flat}_j$, 
where
$\hat \Gamma^{\flat}_i$ is a fibre of the first projection for $i=1$, $2$ and $\hat
\Gamma^{\flat}_0$,
$\hat \Xi^{\flat}_j$ are fibres of the second projection for $j=1$, $2$. 
\end{pr}
{\it Proof.} From the exact sequence:
$
0\rightarrow \Omega_{\tilde S}^1\rightarrow \Omega_{\tilde S}^1(\log \tilde \Delta)\rightarrow
{\cal O}_{\tilde \Delta^{\flat \nu}}\rightarrow 0,
$ we have
the following exact sequence: 
$$
0\rightarrow H^0(\Omega_{\tilde S}^1(\log \tilde \Delta))\rightarrow H^0({\cal O}_{\tilde 
\Delta^{\flat \nu}})\rightarrow \mbox{Pic }\tilde S\otimes \mbox{\boldmath $C$}\rightarrow 0,
$$
since $\mbox{Pic }\tilde S\otimes \mbox{\boldmath $Q$}$ is generated by all the irreducible
components of $\tilde \Delta$. By taking the invariant subspaces under the induced
action of
$G:=\mbox{Gal}(\tilde S/S)$, we get the the following exact sequence:
$$
0\rightarrow H^0(\Omega_{\tilde S}^1(\log \tilde \Delta))^G\rightarrow H^0({\cal
O}_{\lfloor\Delta\rfloor ^{\nu}})\rightarrow \mbox{Pic }S\otimes \mbox{\boldmath $C$}\rightarrow 0,
$$
where $H^0(\Omega_{\tilde S}^1(\log \tilde \Delta))^G$ is the invariant subspace under
the $G$-action. Since $\tilde S\setminus
\tilde
\Delta$ is a two dimensional algebraic torus, 
we see that
$h^0(\Omega_{\tilde S}^1(\log\tilde
\Delta))=2$ and that the natural morphism $\bigwedge ^2 H^0(\Omega_{\tilde S}^1(\log \tilde
\Delta))\rightarrow H^0({\cal O}_{\tilde S}(K_{\tilde S}+\tilde \Delta))$ is an isomorphism. 
Since the eigenvalue of the induced action
of a generator of $G$ on
$H^0({\cal O}_{\tilde S}(K_{\tilde S}+\tilde \Delta))$ is $-1$, $H^0(\Omega_{\tilde S}^1(\log
\tilde \Delta))^G$ is one dimensional. Thus we
get $\rho(S)=e_{\mbox{top}}(\lfloor\Delta\rfloor)-2$. Firstly, we show that there exists a $\mbox{\boldmath
$P$}^1$-fibration $\varphi:S\rightarrow \mbox{\boldmath $P$}^1$. Assume that
$\rho(S^{\flat})=1$. From (\ref{eqn:dr=de}), we get 
$e_{\mbox{top}}(\lfloor\Delta^{\flat}\rfloor)=3$ which implies $e_{\mbox{top}}(\tilde
\Delta^{\flat})=2+s$ by (\ref{eqn:topofb2}), hence $\rho(\tilde S^{\flat})=e_{\mbox{top}}(\tilde
\Delta^{\flat})-2=s$. We note that $\Xi^{\flat}\neq 0$ since $s>0$. By Noether's
equality, we have $K_{\tilde S^{\flat}}^2=10-s$. Since $S^{\flat}$ has only Du Val singularities of
type $A_1$, we have $\delta K^2=0$. Combining this with $\delta
\rho=4-(\Gamma^{\flat},\Xi^{\flat})$, we obtain $(\Xi^{\flat})^2-2s=0$ from (\ref{eqn:k-r}). Let
$\iota:S^{\flat \prime}\rightarrow S^{\flat}$ be the extraction of all the divisors over
$\mbox{LC}(S^{\flat},\Delta^{\flat})$ whose log discrepancy with respect to
$(S^{\flat},\Delta^{\flat})$ are $0$. Noting that 
$(\iota^{-1}_{\ast}\Xi^{\flat})^2=(\Xi^{\flat})^2-2s=0$ by our construction and that $\Xi^{\flat}$
is irreducible by our assumption that $\rho(S^{\flat})=1$ and $\tilde S^{\flat}$ is smooth, some
multiple of $\iota^{-1}_{\ast}\Xi^{\flat}$ defines a $\mbox{\boldmath $P$}^1$-fibration
$\varphi^{\flat \prime}:S^{\flat \prime}\rightarrow \mbox{\boldmath $P$}^1$. Since $\tau$ factors into
$\iota\circ \tau^{\prime}$, where $\tau^{\prime}:S\rightarrow S^{\flat \prime}$ is a birational
morphism, we conclude that there exists a $\mbox{\boldmath $P$}^1$-fibration $\varphi:S\rightarrow
\mbox{\boldmath $P$}^1$ and we may assume that $S^{\flat}$ has a structure of conic fibration
$\varphi^{\flat}:S^{\flat}\rightarrow \mbox{\boldmath $P$}^1$ with $\rho(S^{\flat}/\mbox{\boldmath
$P$}^1)=1$. Under this assumption, we have
$e_{\mbox{top}}(\lfloor\Delta^{\flat}\rfloor)=4$, $e_{\mbox{top}}(\tilde \Delta^{\flat})=4+s$,
$\rho(\tilde S^{\flat})=2+s$, $K_{\tilde S^{\flat}}^2=8-s$, $(\Gamma^{\flat})^2=4-(1/2)s$ and 
$(\Xi^{\flat})^2=2s$ in the same way as above. Let
$\Gamma^{\flat}=\sum_{i=0}^2\Gamma^{\flat}_i$ be the irreducible decomposition and assume that
$\Gamma^{\flat}_0$ is horizontal with respect to $\varphi^{\flat}$. We may assume that
$\Gamma^{\flat}_1$ is contained in a fibre of $\varphi^{\flat}$ since
$(\Gamma^{\flat},\varphi^{\ast}(t))\leq 2$ for $t\in \mbox{\boldmath $P$}^1$. Firstly, assume that 
$(\Gamma^{\flat}_0,\varphi^{\ast}(t))=2$ for $t\in \mbox{\boldmath $P$}^1$. Then
$\Gamma^{\flat}_1$, $\Gamma^{\flat}_2$ and $\Xi^{\flat}$ are contained in fibres of
$\varphi^{\flat}$. Noting that $(\Gamma^{\flat}_0,\Gamma^{\flat}_i)=1$ for $i=1$, $2$, we have
$(\Gamma^{\flat}_0,\Xi^{\flat})=0$. Thus we conclude that $\Xi^{\flat}=0$ and that $s=0$, hence
$(\Gamma^{\flat}_0)^2=0$. Since
$(-K_{S^{\flat}},C)=(\Gamma^{\flat},C)>0$ for any
irreducible curve $C$ on $S^{\flat}$, we see that
$\overline{\mbox{NE}}_{K_{S^{\flat}}}(S^{\flat})=0$ and that $\overline{\mbox{NE}}(S^{\flat})$ is
spanned by exactly two extremal rays with respect to $K_{S^{\flat}}$ one of which corresponds to
$\varphi^{\flat}$. Let $\check\varphi^{\flat}$ be the extremal contraction of 
the other extremal ray. If $\check\varphi^{\flat}$ is birational, 
then, by Lemma~\ref{lm:K-cont}, $\mbox{Exc }\check\varphi^{\flat}\subset \Gamma^{\flat}$, which is
absurd. Thus  $\check\varphi^{\flat}$ defines another conic fibration 
$\check\varphi^{\flat}\rightarrow \mbox{\boldmath $P$}^1$ with 
$\rho(S^{\flat}/\mbox{\boldmath $P$}^1)=1$ such that $\Gamma^{\flat}_1$ and $\Gamma^{\flat}_2$ is
horizontal with respect to $\check\varphi^{\flat}$. Thus we may assume that $\Gamma_0$ is a section
of $\varphi$ at first. Assume that $\Gamma^{\flat}_2$ is contained in a fibre of
$\varphi^{\flat}$. Then we have $(\Gamma^{\flat}_0)^2=-(1/2)s$. Since $(\Gamma^{\flat}_0)^2\in
\mbox{\boldmath $Z$}$, we have $s=0$ or $2$. We may assume that $s=0$ for if we assume that $s=2$,
then we have $(\Gamma^{\flat}_0)^2=-1$ and we can use the argument in the case
$\rho(S^{\flat})=1$. Under the assumption that $s=0$, there exists another another conic fibration 
$\check\varphi^{\flat}\rightarrow \mbox{\boldmath $P$}^1$ with 
$\rho(S^{\flat}/\mbox{\boldmath $P$}^1)=1$ such that $\Gamma^{\flat}_1$ and $\Gamma^{\flat}_2$ is
horizontal with respect to $\check\varphi^{\flat}$, by the same argument as
above. Thus we may assume that $\Gamma_0$ and $\Gamma_2$ are sections of $\varphi^{\flat}$ at
first. Then we have $(\Xi^{\flat},\varphi^{\flat\ast}(t))=0$ for $t\in \mbox{\boldmath $P$}^1$ and
$s=0$. Since we have $\rho(\tilde S^{\flat})=2$, $(S^{\flat},\Delta^{\flat
+}_{\varphi^{\flat}}(t))$  is 
log canonical for any $t\in \mbox{\boldmath $P$}^1$ and $\mbox{Supp
}\lfloor\mbox{Diff}_{C_{\varphi^{\flat}}(t)}(\Delta^{\flat -}_{\varphi^{\flat}}(t))\rfloor
\subset\lfloor\Delta^{\flat -}_{\varphi^{\flat}}(t)\rfloor$ for any $t\in \mbox{\boldmath $P$}^1$
by the argument in the proof of Lemma~\ref{lm:typeII}. Therefore, possibly after ${\cal S}$-elementary transformations, $\varphi^{\flat}$ has only fibres of type $(\mbox{II-}1)_b$ as in
Lemma~\ref{lm:dynkin}. Thus we get the assertion.
\hfill\bsquare
\vskip 5mm
Assume that $\tilde S$ is not smooth. By Proposition~\ref{pr:cltypeIII-1}, we see that  
$\rho(\tilde S)=e_{\mbox{top}}(\tilde \Delta)$ and all the irreducible components of
$\tilde\Delta$ give a complete basis of $\mbox{Pic }\tilde S\otimes \mbox{\boldmath $Q$}$, hence
we get  
$$
\rho(S)=e_{\mbox{top}}(\tilde \Delta)/2+1=e_{\mbox{top}}(\lfloor\Delta\rfloor)-1.
$$
by Lemma~\ref{lm:pic}. Combining this with (\ref{eqn:dr=de}), we obtain
$e_{\mbox{top}}(\lfloor\Delta^{\flat}\rfloor)=1+\rho(S^{\flat})$. From (\ref{eqn:topofb2}), we
have $\rho(\tilde S^{\flat})=e_{\mbox{top}}(\tilde\Delta^{\flat})=2\rho(S^{\flat})-2+s$. Let
$\tilde \mu^{\flat}:\tilde M^{\flat}\rightarrow \tilde S^{\flat}$ be the minimal resolution of
$\tilde S^{\flat}$ and put $\delta\tilde\rho:=\rho(\tilde M^{\flat})-\rho(\tilde S^{\flat})$. Then
by Noether's equality, we have
\begin{eqnarray}
d&=&K_{\tilde M^{\flat}}^2=10-\rho(\tilde M^{\flat})\nonumber \\ 
&=&10-2\rho(S^{\flat})-s-\delta\tilde\rho \nonumber \\ 
&=&-2\rho(S^{\flat})-s+\left\{
\begin{array}{ll}
8, & \mbox{if Sing }\tilde S=4A_1,\\ 
6, & \mbox{if Sing }\tilde S=3A_2,\\ 
5, & \mbox{if Sing }\tilde S=A_1+2A_3,\\ 
4 & \mbox{if Sing }\tilde S=A_1+A_2+A_5. \label{eqn:d}
\end{array}
\right.
\end{eqnarray}
\begin{pr}\label{pr:cltypeIII-2.2} Assume that $\mbox{\em Sing }\tilde S=3A_2$. 
Then one of the followings holds.
\begin{description}
\item[ (1) ] $S^{\flat}\simeq \mbox{\boldmath $P$}^2$ and
$\Delta^{\flat}=\Gamma^{\flat}+(1/2)\Xi^{\flat}$, where $\Xi^{\flat}$ is an irreducible quartic
curve with three cusps and $\Gamma^{\flat}$ is a double tangent.
\item[ (2) ] $S^{\flat}$ is a rank one Gorenstein log del Pezzo surface with
$\mbox{\em Sing }S^{\flat}=A_1+A_2$ and
 there exists a birational morphism $\lambda:U\rightarrow S^{\flat}$ from a smooth projective surface $U$ 
such that $U$ admits a structure of elliptic surface with a section 
$\psi:U\rightarrow \mbox{\boldmath $P$}^1$ which is minimal over $\mbox{\boldmath $P$}^1$ with  
$\mbox{\em Typ}(U;\psi)=\mbox{\em I}^{\ast}_1+\mbox{\em II}+\mbox{\em I}_3$ and that 
$\Delta^{\flat \ U}=(1/2)\psi^{\ast}(t+u)$ where $\psi^{\ast}(t)$ and $\psi^{\ast}(u)$ are the singular fibres of type 
$\mbox{\em I}^{\ast}_1$ and $\mbox{\em II}$ respectively.
\end{description}
\end{pr}
{\it Proof.} The possible singular types of $(S^{\flat}, \Delta^{\flat})$ on 
$S^{\flat}\setminus \lfloor\Delta^{\flat}\rfloor$ are types $A_2$, $A_0/2$, $A_2/2$-$\beta$ and 
$A_2/2$-$\epsilon$ and we see that $\Xi^{\flat}$ is irreducible if $\rho(S^{\flat})=1$. In what
follows,
$\nu({\cal T})$ denotes the number of points of type
${\cal T}$. We note that we have 
$
2\nu(A_2)+\nu(A_2/2\mbox{-}\beta)+\nu(A_2/2\mbox{-}\epsilon)=3.
$
Firstly, we
consider the case in which $S^{\flat}$ is a rank one log del Pezzo surface. Assume that
$s=2$. Then we have $\rho(\tilde S^{\flat})=2$, $d=2$, $(\Gamma^{\flat})^2=1$,  
$(\Gamma^{\flat},\Xi^{\flat})=4$ and   
$(\Xi^{\flat})^2=16$. Thus we have $\delta K^2-\delta\rho=0$. On the other hand, since we have
$\delta K^2=(1/3)\nu(A_2/2\mbox{-}\epsilon)$ and $\delta\rho=2\nu(A_2)+\nu(A_2/2\mbox{-}\epsilon)$,
we have $\delta K^2-\delta\rho=-2\nu(A_2)-(2/3)\nu(A_2/2\mbox{-}\epsilon)$, hence
$\nu(A_2)=\nu(A_2/2\mbox{-}\epsilon)=0$ and $\nu(A_2/2\mbox{-}\beta)=3$, which implies that we are
in the case $(1)$. Assume that $s=1$. Then we have $d=3$ and $\rho(\tilde S^{\flat})=1$. We note
that we have
$(\Gamma^{\flat},\Xi^{\flat})\equiv 1$
$(\mbox{mod }2)$ from
$(\Gamma^{\flat})^2=3/2$, so we have 
$(\Gamma^{\flat},\Xi^{\flat})=3$, $(\Xi^{\flat})^2=6$ and $\delta K^2-\delta\rho=-3$. On the other hand, since we have 
$\delta K^2=(1/3)\nu(A_2/2\mbox{-}\epsilon)$ 
and 
$\delta \rho=1+2\nu(A_2)+\nu(A_2/2\mbox{-}\epsilon)$, 
we get
$
\delta K^2-\delta\rho=-1-2\nu(A_2)-(2/3)\nu(A_2/2\mbox{-}\epsilon).
$
Therefore, we obtain $3\nu(A_2)+2\nu(A_2/2\mbox{-}\epsilon)=3$. Hence $(\nu(A_2), \nu(A_2/2\mbox{-}\epsilon))=(1,0)$ or $(0,3)$. 
Since we have $\deg\mbox{Diff}_{\Xi^{\flat \nu}}(\Gamma^{\flat})=3+2\nu(A_2/2\mbox{-}\beta)+(2/3)\nu(A_2/2\mbox{-}\epsilon)$, we have 
$e_{\mbox{top}}(\Xi^{\flat \nu})=2\nu(A_2/2\mbox{-}\beta)+(2/3)\nu(A_2/2\mbox{-}\epsilon)$. 
 Assume that 
$(\nu(A_2), \nu(A_2/2\mbox{-}\epsilon))=(1,0)$. Then we have $\nu(A_2/2\mbox{-}\beta)=1$ and 
$e_{\mbox{top}}(\Xi^{\flat \nu})=2$, i.e., $\Xi^{\flat \nu}\simeq \mbox{\boldmath $P$}^1$. We can see that there exists a 
birational morphism $\eta:W\rightarrow S^{\flat}$ from a smooth rational surface 
with $\rho(W)=13$ such that 
$\mbox{Supp }\eta^{-1}_{\ast}\Delta^{\flat}\cup \mbox{Exc }\eta$ has only simple normal crossing singularities whose dual 
graph is as follows.

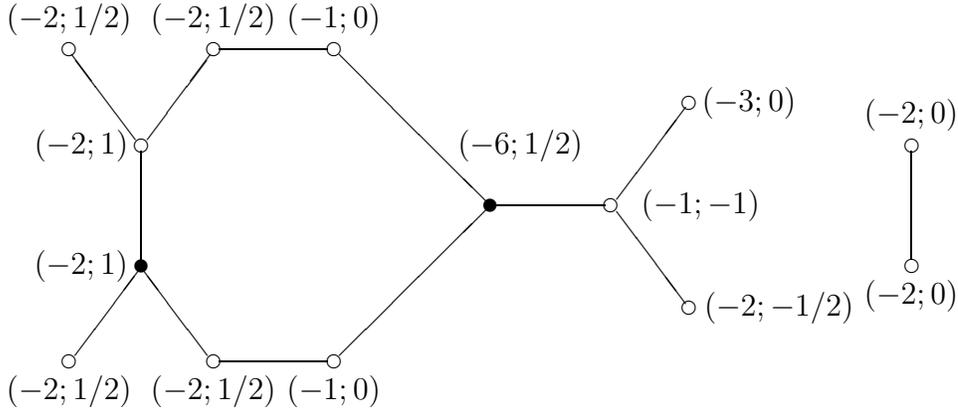
\unitlength 0.8mm
\begin{picture}(150,75)(-30,0)
\put(0,0){\circle{2}}\put(0,0){\makebox(0,-10){$(-2;1/2)$}}
\put(1,1){\line(3,4){11.5}}
\put(12,16){\circle*{2}}\put(12,16){\makebox(-20,0){$(-2;1)$}}
\put(24,0){\circle{2}}\put(24,0){\makebox(0,-10){$(-2;1/2)$}}
\put(23,1){\line(-3,4){11.5}}
\put(12,16){\line(0,1){19}}
\put(12,36){\circle{2}}\put(12,36){\makebox(-20,0){$(-2;1)$}}
\put(11,37){\line(-3,4){10.5}}
\put(0,52){\circle{2}}\put(0,52){\makebox(0,10){$(-2;1/2)$}}
\put(13,37){\line(3,4){10.5}}
\put(24,52){\circle{2}}\put(24,52){\makebox(0,10){$(-2;1/2)$}}

\put(25,52){\line(1,0){18}}
\put(44,52){\circle{2}}\put(44,52){\makebox(0,10){$(-1;0)$}}
\put(45,51){\line(1,-1){25}}
\put(25,0){\line(1,0){18}}
\put(44,0){\circle{2}}\put(44,0){\makebox(0,-10){$(-1;0)$}}
\put(45,1){\line(1,1){25}}
\put(70,26){\circle*{2}}\put(70,26){\makebox(10,20){$(-6;1/2)$}}

\put(71,26){\line(1,0){18}}
\put(90,26){\circle{2}}\put(90,26){\makebox(30,0){$(-1;-1)$}}
\put(91,27){\line(3,4){11.5}}
\put(91,25){\line(3,-4){11.5}}
\put(103,43){\circle{2}}\put(103,43){\makebox(20,0){$(-3;0)$}}
\put(103,9){\circle{2}}\put(103,9){\makebox(30,0){$(-2;-1/2)$}}

\put(140,36){\circle{2}}\put(140,36){\makebox(0,10){$(-2;0)$}}
\put(140,17){\line(0,1){18}}
\put(140,16){\circle{2}}\put(140,16){\makebox(0,-10){$(-2;0)$}}
\end{picture}
\vskip 20mm
In the above dual graph, $\bullet$ denotes a curves 
which is not 
$\eta$-exceptional and $\circ$ denotes a $\eta$-exceptional curve. 
$(\ast;\ast)$ denotes $($ the self-intersection number; the multiplicities in 
$\Delta^{\flat W}$. We shall use the same notations in the subsequent dual graphs.

We see that there exists a birational morphism $\upsilon:W\rightarrow U$ 
to a smooth rational surface $U$ with $\rho(U)=10$ such that $\eta=\lambda\circ\upsilon$, where 
$\lambda:U\rightarrow S^{\flat}$ is a birational morphism and that
$\lambda_{\ast}^{-1}\Xi^{\flat}$ is nef with 
$(\lambda_{\ast}^{-1}\Xi^{\flat})^2=0$ and 
$(\Delta^{\flat U}-(1/2)\lambda_{\ast}^{-1}\Xi^{\flat},\lambda_{\ast}^{-1}\Xi^{\flat})=0$. We note
that $\Delta^{\flat U}\geq 0$. Applying the log abundance theorem to 
$K_U+\Delta^{\flat U}+\varepsilon\lambda_{\ast}^{-1}\Xi^{\flat}$ for sufficiently small 
positive rational number $\varepsilon$, we obtain a connected surjective morphism 
$\psi:U\rightarrow \mbox{\boldmath $P$}^1$, 
which turns out to be a minimal elliptic fibration since 
$K_U+\Delta^{\flat U}\sim_{\mbox{\boldmath $Q$}}0$ and $\rho(U)=10$. Put 
$t:=\psi(\lambda_{\ast}^{-1}\Gamma^{\flat})$ and $u:=\psi(\lambda_{\ast}^{-1}\Xi^{\flat})$. 
Then we see that $\psi^{\ast}(t)$ and $\psi^{\ast}(u)$ are singular fibre of type 
$\mbox{I}^{\ast}_1$ and $\mbox{II}$ respectively and that $\Delta^{\flat \ U}=(1/2)\psi^{\ast}(t+u)$. 
Since we have 
$
\sum_{v\neq t,u}e_{\mbox{top}}(\psi^{-1}(v))=3
$
and there exists $v\in \mbox{\boldmath $P$}^1$ such that the singular fibre $\psi^{\ast}(v)$ 
has the dual graph containing $A_2$ as its dual subgraph, we deduce that $\psi^{\ast}(v)$ is of type $\mbox{I}_3$ 
and that $\psi$ is smooth over $\mbox{\boldmath $P$}^1\setminus \{t,u,v\}$. Thus we are in the case
$(2)$. We shall show that the case $(\nu(A_2),
\nu(A_2/2\mbox{-}\epsilon))=(0,3)$ does not occurs.  Assume that $(\nu(A_2),
\nu(A_2/2\mbox{-}\epsilon))=(0,3)$. Then we have 
$\nu(A_2/2\mbox{-}\beta)=0$ and $\Xi^{\flat \nu}\simeq \mbox{\boldmath $P$}^1$.  
We see that that there exists a birational morphism $\eta:W\rightarrow S^{\flat}$, 
where $W$ is a rational surface with $\rho(W)=11$, such that 
$\mbox{Supp }\eta^{-1}_{\ast}\Delta^{\flat}\cup \mbox{Exc }\eta$ has only simple normal crossing singularities whose dual graph is 
 as follows.

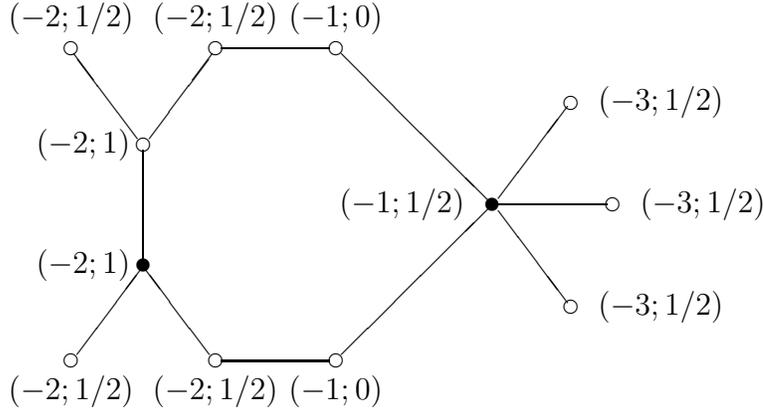
\unitlength 0.8mm
\begin{picture}(150,75)(-30,0)
\put(0,0){\circle{2}}\put(0,0){\makebox(0,-10){$(-2;1/2)$}}
\put(1,1){\line(3,4){11.5}}
\put(12,16){\circle*{2}}\put(12,16){\makebox(-20,0){$(-2;1)$}}
\put(24,0){\circle{2}}\put(24,0){\makebox(0,-10){$(-2;1/2)$}}
\put(23,1){\line(-3,4){11.5}}
\put(12,16){\line(0,1){19}}
\put(12,36){\circle{2}}\put(12,36){\makebox(-20,0){$(-2;1)$}}
\put(11,37){\line(-3,4){10.5}}
\put(0,52){\circle{2}}\put(0,52){\makebox(0,10){$(-2;1/2)$}}
\put(13,37){\line(3,4){10.5}}
\put(24,52){\circle{2}}\put(24,52){\makebox(0,10){$(-2;1/2)$}}

\put(25,52){\line(1,0){18}}
\put(44,52){\circle{2}}\put(44,52){\makebox(0,10){$(-1;0)$}}
\put(45,51){\line(1,-1){25}}
\put(25,0){\line(1,0){18}}
\put(44,0){\circle{2}}\put(44,0){\makebox(0,-10){$(-1;0)$}}
\put(45,1){\line(1,1){25}}
\put(70,26){\circle*{2}}\put(70,26){\makebox(-30,0){$(-1;1/2)$}}

\put(71,26){\line(1,0){18}}
\put(90,26){\circle{2}}\put(90,26){\makebox(30,0){$(-3;1/2)$}}
\put(71,27){\line(3,4){11.5}}
\put(71,25){\line(3,-4){11.5}}
\put(83,43){\circle{2}}\put(83,43){\makebox(30,0){$(-3;1/2)$}}
\put(83,9){\circle{2}}\put(83,9){\makebox(30,0){$(-3;1/2)$}}
\end{picture}
\vskip 20mm
Let $\upsilon:W\rightarrow U$ be a contraction of $\eta^{-1}_{\ast}\Xi^{\flat}$. 
Then we see that $U$ is smooth rational surface which admits a structure of an elliptic fibration $\psi:U\rightarrow \mbox{\boldmath $P$}^1$ by the same 
argument as above, such that there exists two points $t$, $u\in \mbox{\boldmath $P$}^1$ 
with $\psi^{\ast}(t)$ and $\psi^{\ast}(u)$ being singular fibres of type $\mbox{I}^{\ast}_1$ and $\mbox{IV}$ respectively. 
From the above dual graph, we can see that there exists a smooth rational curve $E$ on $U$ such that $(\psi^{\ast}(t),E)=1$ and 
$(\psi^{\ast}(u),E)=3$ but which is absurd. Secondarily, we consider the case in which there exists
a conic fibration $\varphi^{\flat}:S^{\flat}\rightarrow \mbox{\boldmath $P$}^1$ with
$\rho(S^{\flat}/\mbox{\boldmath $P$}^1)=1$. In this case, we note that we have
$e_{\mbox{top}}(\lfloor\Delta^{\flat}\rfloor)=3$. Let
$\Gamma^{\flat}=\Gamma^{\flat}_0+\Gamma^{\flat}_1$ be the irreducible decomposition. Assume that
$s=2$. Then we have $d=0$ and $\rho(\tilde S^{\flat})=4$, hence the self-intersection numbers of all
the irreducible components of $\tilde \Delta^{\flat}$ are $(a)$ $-1$, $-2$, $-2$, $-3$ or $(b)$
$-1$,
$-1$, $-3$, $-3$. The case $(a)$ is excluded for if we are in the case $(a)$, the $-1$-curve and
the
$(-3)$-curve is invariant under the action of the covering transformation group, but which is
absurd. Thus we may assume that $(\Gamma^{\flat}_0)^2=0$, $(\Gamma^{\flat}_1)^2=-2$. Since
$\mbox{Pic }S^{\flat}\otimes \mbox{\boldmath $Q$}$ is generated by $\Gamma^{\flat}_0$ and
$\Gamma^{\flat}_1$, we have $\Xi^{\flat}\sim _{\mbox{\boldmath
$Q$}}6\Gamma^{\flat}_0+2\Gamma^{\flat}_1$ from $(\Gamma^{\flat}_i,\Xi^{\flat})=2$ for $i=0$, $1$,
which implies that $(\Xi^{\flat})^2=16$. Thus we have $0=\delta K^2-\delta
\rho=-2\nu(A_1)-(2/3)\nu(A_2/2\mbox{-}\epsilon)$, hence $\nu(A_2)=\nu(A_2/2\mbox{-}\epsilon)=0$ and
$\nu(A_2/2\mbox{-}\beta)=3$. Moreover, since we have $\deg\mbox{Diff}_{\Xi^{\flat
\nu}}(\Gamma^{\flat})=4+2\nu(A_2/2\mbox{-}\beta)=10$, we have $e_{\mbox{top}}(\Xi^{\flat
\nu})=2$. Since $\varphi^{\flat}$ is smooth, It is easily seen that $\Gamma^{\flat}_0$ is a fibre,
$\Gamma^{\flat}_0$ is a section of $\varphi^{\flat}$ and $\Xi^{\flat}$ is
irreducible, hence $\Xi^{\flat\nu}\simeq \mbox{\boldmath $P$}^1$. $\varphi^{\flat}$ induces a double
cover $\varphi^{\flat}:\Xi^{\flat\nu}\rightarrow \mbox{\boldmath $P$}^1$ with at least four
branching point, which contradicts the Hurwitz's formula. Assume that $s=1$. Then we have $d=1$
and $\rho(\tilde S^{\flat})=3$, hence the self-intersection numbers of all the irreducible
components of $\tilde \Delta^{\flat}$ are $-1$, $-1$, $-2$. Assume that
$\mbox{LC}(S^{\flat},\Delta^{\flat})\cap \Gamma^{\flat}_0\neq \emptyset$. Then we have 
$(\Gamma^{\flat}_0)^2=-1$ and $(\Gamma^{\flat}_1)^2=-1/2$ and we see that $\Gamma^{\flat}_1$
supports an extremal ray with respect to $K_{S^{\flat}}+\Delta^{\flat}$. Contracting
$\Gamma^{\flat}_1$, this case reduces to the rank one log del Pezzo case.  Assume that $s=0$. 
Then we have $d=2$
and $\rho(\tilde S^{\flat})=2$, hence the self-intersection numbers of all the irreducible
components of $\tilde \Delta^{\flat}$ are $-1$, $-1$, hence $(\Gamma^{\flat}_i)^2=-1/2$ for
$i=0$, $1$. By the same way as in the previous argument, we conclude that this case also reduces to
the rank one log del Pezzo case. Thus we get the assertion. \hfill\bsquare

\begin{ex}\label{ex:ex1}{\rm The example of the case $(1)$ is well known and goes back to
\cite{zariski}. To show the existence, we only have to take a dual curve of a nodal cubic as
$\Xi^{\flat}$. There exists exactly one double tangent by the Pl\"ucker's formulae. It is also well
known that for a suitable choice of homogeneous coordinates, a defining equation
$f(X,Y,Z)$ of a nodal cubic curve can be written as  $f(X,Y,Z)=Y^2Z-X^2(X+Z)$ and the defining
equation $\hat f(X,Y,Z)$ of its dual curve is calculated to be $\hat
f(X,Y,Z)=4X^3Z+4X^4-36XY^2Z-8X^2Y^2-27Y^2Z^2+4Y^4$. (see, for example,
\cite{brieskorn}, p.585, \cite{dimca}, p.131, Exercise (4.7)\footnote{Unfortunately, there is a
minor misprint in \cite{dimca}, p131 which says \lq\lq \ldots cuspidal cubic \ldots ''. \lq\lq
cuspidal '' should be \lq\lq nodal ''.} or
\cite{miranda-persson}, Table 6.8). In particular, we see that such a log surface
$(S^{\flat},\Delta^{\flat})$ is unique up to isomorphisms. As for the case $(2)$, the existence of
a minimal rational elliptic surface with a section 
$\psi:U\rightarrow \mbox{\boldmath $P$}^1$ with $\mbox{Typ}(U;\psi)=\mbox{I}^{\ast}_1+\mbox{II}+\mbox{I}_3$ is known (see
\cite{persson}).  By the list in \cite{oguiso-shioda}, we see that $\mbox{MW}(U_{\eta})=\mbox{\boldmath $Z$}P$ for some $P\in
\mbox{MW}(U_{\eta})$  with $<P,P>=1/12$, where $\mbox{MW}(U_{\eta})$ denotes the Mordell-Weil group of the generic fibre $U_{\eta}$
and 
$<\ast,\ast>$ denotes the height paring in the Shioda's sense. Put $Q:=3P$. Then from the formula
$(8.12)$ in \cite{shioda}, we have
$$
\frac{3}{4}=<Q,Q>=2+2(QO)-\sum_{v\in R}\mbox{contr}_v(Q),
$$
where we followed the notations in \cite{shioda}.
Let $\psi^{\ast}(t)=\Theta_{t,0}+\Theta_{t,1}+\Theta_{t,2}+\Theta_{t,3}+2\Theta_{t,4}+2\Theta_{t,5}$ be the 
type $\mbox{I}^{\ast}_1$ singular fibre 
and $\psi^{\ast}(v)=\Theta_{v,0}+\Theta_{v,1}+\Theta_{v,2}$ be the type $\mbox{I}_3$ fibre, 
where $\Theta_{t,0}$ and $\Theta_{v,0}$ are the components which intersect the section $(O)$.
From $(8.16)$ in \cite{shioda},  we have 
$$
\mbox{contr}_t(Q)=\left\{
\begin{array}{ll}
0 & \mbox{ if }(Q\Theta_{t,0})=1,\\
1 & \mbox{ if }(Q\Theta_{t,1})=1,\\
5/4 & \mbox{ if }(Q\Theta_{t,2})=1 \mbox{ or } (Q\Theta_{t,3})=1
\end{array}
\right.
$$
and 
$$
\mbox{contr}_v(Q)=\left\{
\begin{array}{ll}
0 & \mbox{ if }(Q\Theta_{v,0})=1,\\
2/3 & \mbox{ if }(Q\Theta_{v,1})=1 \mbox{ or } (Q\Theta_{v,2})=1.
\end{array}
\right.
$$
Noting that $(QO)\in \mbox{\boldmath $Z$}$, we see that 
$(Q\Theta_{t,2})=1$ or $(Q\Theta_{t,3})=1$ and that $(Q\Theta_{v,0})=1$ and $(QO)=0$. 
Since $\omega_U=\psi^{\ast}{\cal O}_{\mbox{\boldmath $P$}^1}(-1)$, we have $2K_U+\psi^{\ast}(t+u)\sim 0$, 
where $\psi^{\ast}(u)$ is the type $\mbox{II}$ singular fiber. 
Let $\lambda:U\rightarrow S^{\flat}$ be the contraction of all the curves $(O)$, $(Q)$, 
$\Theta_{t,i}$ $(0\leq i\leq 4)$ and $\Theta_{v,j}$ $(1\leq j\leq 2)$ and put 
$\Delta^{\flat}:=(1/2)\lambda_{\ast} \psi^{\ast}(t+u)$. Then $(S^{\flat},\Delta^{\flat})$ gives an 
example of the log surfaces of type $(2)$ as in 
Proposition~\ref{pr:cltypeIII-2.2}. }
\end{ex}

\begin{pr}\label{pr:cltypeIII-2.3} Assume that $\mbox{\em Sing }\tilde S=A_1+2A_3$. 
Then there exists a birational morphism $\lambda:U\rightarrow S^{\flat}$ from a smooth projective surface $U$ which 
admits a structure of an elliptic surface with a section $\psi:U\rightarrow \mbox{\boldmath $P$}^1$ 
which 
is minimal over $\mbox{\boldmath $P$}^1$. And one of the followings holds.
\begin{description}
\item[ $(1)$ ] $S^{\flat}\simeq \mbox{\boldmath $P$}^2$ and 
$\mbox{\em Typ}(U;\psi)=\mbox{\em I}^{\ast}_1+\mbox{\em III}+\mbox{\em I}_2$,

\item[ $(2)$ ] $S^{\flat}$ is a rank one Gorenstein log del Pezzo surface with 
$\mbox{\em Sing }S^{\flat}=2A_1+A_3$ and 
$\mbox{\em Typ}(U;\psi)=\mbox{\em I}^{\ast}_1+\mbox{\em I}_1+\mbox{\em I}_4$. 

\end{description}
Moreover, $\Delta^{\flat \ U}=(1/2)\psi^{\ast}(t+u)$ where $\psi^{\ast}(t)$ is the singular fibres
of type $\mbox{\em I}^{\ast}_1$ and $\psi^{\ast}(u)$ is the singular fibres of type
$\mbox{\em III}$ in the case $(1)$, $\mbox{\em I}_1$ in the case $(2)$.
\end{pr}
{\it Proof.} The possible singular types of $(S^{\flat}, \Delta^{\flat})$ on 
$S^{\flat}\setminus \lfloor\Delta^{\flat}\rfloor$ are types 
$A_3$, $A_0/2$, $A_1/2$-$\alpha$, $A_1/2$-$\gamma$, $A_3/2$-$\alpha$, $A_3/2$-$\delta$ and
$A_3/2$-$\zeta$. We note that we have
$
\nu(A_1/2\mbox{-}\alpha)+\nu(A_1/2\mbox{-}\gamma)=1 
$
and
$
2\nu(A_3)+\nu(A_3/2\mbox{-}\alpha)+\nu(A_3/2\mbox{-}\delta)+\nu(A_3/2\mbox{-}\zeta)=2.
$
From 
$
\delta K^2=\nu(A_1/2\mbox{-}\gamma)+\nu(A_3/2\mbox{-}\delta)
$
and
$
\delta\rho=4-(\Gamma^{\flat},\Xi^{\flat})+3\nu(A_3)+\nu(A_1/2\mbox{-}\gamma)+2\nu(A_3/2\mbox{-}\delta)+\nu(A_3/2\mbox{-}\zeta),
$
we obtain 
$
\delta K^2-\delta\rho=(\Gamma^{\flat},\Xi^{\flat})-4-3\nu(A_3)-\nu(A_3/2\mbox{-}\delta)-\nu(A_3/2\mbox{-}\zeta),
$
hence we have 
$
3\nu(A_3)+\nu(A_3/2\mbox{-}\delta)+\nu(A_3/2\mbox{-}\zeta)=4-(1/4)(\Xi^{\flat})^2.
$
Moreover, since we have 
$
\deg\mbox{Diff}_{\Xi^{\flat \nu}}(\Gamma^{\flat})=
(\Gamma^{\flat},\Xi^{\flat})+2\nu(A_1/2\mbox{-}\alpha)+4\nu(A_3/2\mbox{-}\alpha)+2\nu(A_3/2\mbox{-}\zeta),
$
we have 
$
e_{\mbox{top}}(\Xi^{\flat \nu})=
(\Gamma^{\flat},\Xi^{\flat})+2\nu(A_1/2\mbox{-}\alpha)+4\nu(A_3/2\mbox{-}\alpha)+2\nu(A_3/2\mbox{-}\zeta)-(1/2)(\Xi^{\flat})^2.
$
Firstly, we consider the case in which $S^{\flat}$ is a rank one log del Pezzo surface. In this
case, we have $d=3-s$, hence $(\Gamma^{\flat})^2=(3-s)/2$. We note that we have $s=1$ or $s=2$
since $\rho(\tilde S^{\flat})=s$. Assume that $s=2$. Then we have $(\Gamma^{\flat})^2=1/2$, which
is absurd since $\Gamma^{\flat}\cap\mbox{Sing }S^{\flat}=\emptyset$ by our assumption. Thus we have
$d=2$, hence $(\Gamma^{\flat})^2=1$ and
$(\Xi^{\flat})^2=(\Gamma^{\flat},\Xi^{\flat})^2$. We note that $(\Gamma^{\flat},\Xi^{\flat})=2$
or $4$ since $(\Gamma^{\flat})^2\in \mbox{\boldmath $Z$}$. Assume that
$(\Gamma^{\flat},\Xi^{\flat})=4$.  Then we have
$3\nu(A_3)+\nu(A_3/2\mbox{-}\delta)+\nu(A_3/2\mbox{-}\zeta)=0$, hence we have 
$\nu(A_3)=\nu(A_3/2\mbox{-}\delta)=\nu(A_3/2\mbox{-}\zeta)=0$ and $\nu(A_3/2\mbox{-}\alpha)=2$. We note that we have 
$e_{\mbox{top}}(\Xi^{\flat \nu})=4+2\nu(A_1/2\mbox{-}\alpha)\geq 4$, 
so we see that $\Xi^{\flat }$ is reducible and consists of two or three irreducible components 
since $\rho(S^{\flat})=1$. 
In fact, we show that $\Xi^{\flat}$ consists of exactly three irreducible components. 
Assume that $\Xi^{\flat }$ consists of two irreducible components $\Xi^{\flat}_1$ and $\Xi^{\flat}_2$. 
Since we have $e_{\mbox{top}}(\Xi^{\flat \nu})\leq 4$, we have $\nu(A_1/2\mbox{-}\alpha)=0$, $\nu(A_1/2\mbox{-}\gamma)=1$ and 
$\Xi^{\flat \nu}_1$, $\Xi^{\flat \nu}_1\simeq \mbox{\boldmath $P$}^1$. We note that we have
$(\Xi^{\flat}_1,\Xi^{\flat}_2)=4$ since $(\Xi^{\flat}_i)^2=(\Gamma^{\flat},\Xi^{\flat}_i)^2=4$ for
$i=1$, $2$. Thus we see that there exists a birational morphism
$\eta^{\prime}:W^{\prime}\rightarrow S^{\flat}$,  where $W^{\prime}$ is a rational surface with
$\rho(W^{\prime})=14$ such that 
$\mbox{Supp }\eta^{\prime -1}_{\ast}\Delta^{\flat}\cup \mbox{Exc }\eta^{\prime}$ has only simple normal crossing 
singularities whose dual graph is as follows.

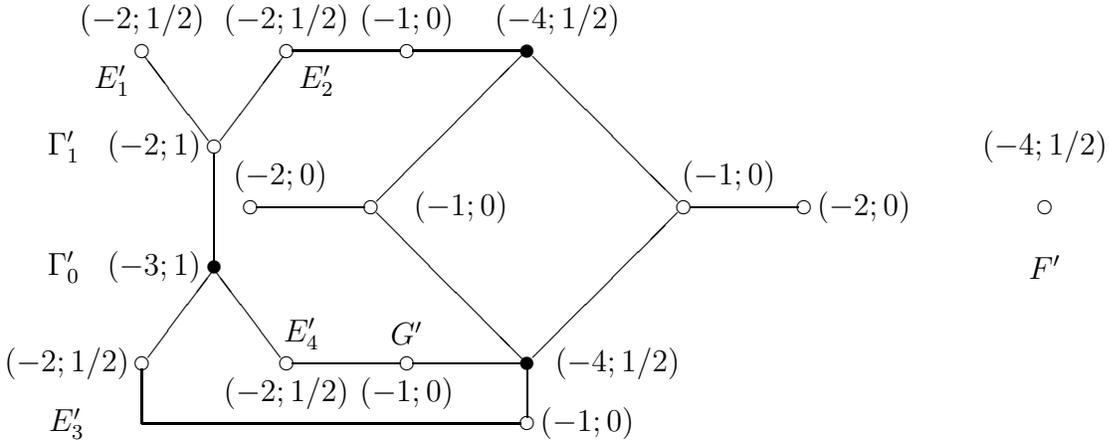
\unitlength 0.8mm
\begin{picture}(150,75)(-30,0)
\put(0,0){\circle{2}}\put(0,0){\makebox(-25,0){$(-2;1/2)$}}\put(0,0){\makebox(-25,-20){$E_3^{\prime}$}}
\put(1,1){\line(3,4){11.5}}
\put(12,16){\circle*{2}}\put(12,16){\makebox(-20,0){$(-3;1)$}}\put(12,16){\makebox(-50,0){$\Gamma_0^{\prime}$}}
\put(24,0){\circle{2}}\put(24,0){\makebox(0,-10){$(-2;1/2)$}}\put(24,0){\makebox(5,10){$E_4^{\prime}$}}

\put(23,1){\line(-3,4){11.5}}
\put(12,16){\line(0,1){19}}
\put(12,36){\circle{2}}\put(12,36){\makebox(-20,0){$(-2;1)$}}\put(12,36){\makebox(-50,0){$\Gamma_1^{\prime}$}}
\put(11,37){\line(-3,4){10.5}}
\put(0,52){\circle{2}}\put(0,52){\makebox(0,10){$(-2;1/2)$}}\put(0,52){\makebox(-10,-10){$E_1^{\prime}$}}
\put(13,37){\line(3,4){10.5}}
\put(24,52){\circle{2}}\put(24,52){\makebox(0,10){$(-2;1/2)$}}\put(24,52){\makebox(10,-10){$E_2^{\prime}$}}

\put(25,52){\line(1,0){18}}
\put(44,52){\circle{2}}\put(44,52){\makebox(0,10){$(-1;0)$}}
\put(45,52){\line(1,0){18}}
\put(64,52){\circle*{2}}\put(64,52){\makebox(10,10){$(-4;1/2)$}}
\put(65,51){\line(1,-1){24}}

\put(25,0){\line(1,0){18}}
\put(44,0){\circle{2}}\put(44,0){\makebox(0,-10){$(-1;0)$}}\put(44,0){\makebox(0,10){$G^{\prime}$}}
\put(45,0){\line(1,0){18}}
\put(64,0){\circle*{2}}\put(64,0){\makebox(30,0){$(-4;1/2)$}}
\put(65,1){\line(1,1){24}}

\put(90,26){\circle{2}}\put(90,26){\makebox(15,10){$(-1;0)$}}
\put(91,26){\line(1,0){18}}
\put(110,26){\circle{2}}\put(110,26){\makebox(20,0){$(-2;0)$}}

\put(18,26){\circle{2}}\put(18,26){\makebox(10,10){$(-2;0)$}}
\put(19,26){\line(1,0){18}}
\put(38,26){\circle{2}}\put(48,26){\makebox(10,0){$(-1;0)$}}
\put(39,27){\line(1,1){25}}
\put(39,25){\line(1,-1){25}}

\put(0,-10){\line(0,1){9}}
\put(0,-10){\line(1,0){63}}
\put(64,-10){\circle{2}}\put(64,-10){\makebox(20,0){$(-1;0)$}}
\put(64,-9){\line(0,1){8}}

\put(150,26){\circle{2}}\put(150,26){\makebox(0,20){$(-4;1/2)$}}\put(150,26){\makebox(0,-20){$F^{\prime}$}}
\end{picture}
\vskip 20mm
 From the above dual graph, we see that there exists a birational morphism 
$\upsilon:W^{\prime}\rightarrow W$ to a rational surface with $\rho(W)=12$ such that 
$\eta^{\prime}$ factors into $\upsilon\circ\eta$ where 
$\eta:W\rightarrow S^{\flat}$ is a birational morphism and 
that $\eta^{-1}_{\ast}\Xi^{\flat}$ is nef and $(\eta^{-1}_{\ast}\Xi^{\flat})^2=0$. 
By the log abundance theorem, some multiple of $\eta^{-1}_{\ast}\Xi^{\flat}$ 
determines the structure of elliptic fibration with a section 
$\psi_W:W\rightarrow \mbox{\boldmath $P$}^1$. We see that 
$\psi^{\ast}_W(u)=\eta^{-1}_{\ast}\Xi^{\flat}$ is a singular fibre of type $\mbox{III}$ where 
$u:=\psi_W(\eta^{-1}_{\ast}\Xi^{\flat})$. Write 
$
\Delta^{\flat
W}=\Gamma_0+\Gamma_1+(1/2)(\sum_{i=1}^4E_i+\eta^{-1}_{\ast}\Xi^{\flat}+F),
$
where $\Gamma_0:=\upsilon_{\ast}\Gamma_0^{\prime}$, $\Gamma_1^{\prime}:=\upsilon_{\ast}\Gamma_1^{\prime}$, 
$E_i:=\upsilon_{\ast}E_i^{\prime}$ 
$(1\leq i\leq 4)$ and $F:=\upsilon_{\ast}F^{\prime}$. Let $e_0$ be a $(-1)$-curve on $W$ which is contained in a fibre of 
$\psi_W$ and let $\upsilon_0:W:=W_0\rightarrow W_1$ be the contraction of $e_0$.
Assume that $(e_0,F)=0$. We can see that $(e_0,\sum_{i=1}^4E_i)<2$ by the semi-negativity
of fiber components.  Since we have 
$(e_0,\Delta^{\flat W})=(e_0,-K_W)=1$, we have $(e_0,\sum_{i=1}^4E_i)\in 2\mbox{\boldmath $Z$}$, 
hence we obtain $(e_0,\sum_{i=1}^4E_i)=0$ and $(e_0,\Gamma_0+\Gamma_1)=1$. Noting that $\rho(W)=12$, 
we have $(e_0,\Gamma_0)=1$ and $(e_0,\Gamma_1)=0$. 
Let $\psi_{W_1}:W_1\rightarrow \mbox{\boldmath $P$}^1$ be the induced morphism from $\psi$ and put 
$t:=\psi_W(\Gamma_0)$. Then we see that $\psi_{W_1}^{\ast}(t)$ is a singular fibre of type
$\mbox{I}^{\ast}_1$. Let $e_1$ be a $(-1)$-curve on $W_1$ which is contained in a fibre of 
$\psi_{W_1}$, $\upsilon_1:W_1\rightarrow W_2$ be the contraction of $e_1$ and put
$F^{(1)}:=\upsilon_{0\ast}F$. We note that $(e_1,F^{(1)})=2$. Let $\psi_{W_2}:W_2\rightarrow
\mbox{\boldmath $P$}^1$  be the induced morphism from $\psi_{W_1}$ and put $v:=\psi_{W_1}(F^{(1)})$.
Then we see that $\psi_{W_2}^{\ast}(v)$ is a  singular fibre of type $\mbox{II}$ or of type
$\mbox{I}_1$ and that
$2K_{W_2}+\psi_{W_2}^{\ast}(t+u+v)\sim 0$. On the other hand, since $\psi_{W_2}:W_2\rightarrow
\mbox{\boldmath $P$}^1$ is a rational elliptic  surface with a section which is minimal over
$\mbox{\boldmath $P$}^1$, we have 
$\omega_{W_2}=\psi_{W_2}^{\ast}{\cal O}_{\mbox{\boldmath $P$}^1}(-1)$, which is a contradiction. 
Thus we may assume that $(e_0,F)>0$. In the same way as the above argument, we see that $(e_0,F)=1$, hence 
$(e_0,\sum_{i=1}^4E_i)=1$. Since $\rho(W)=12$, we see that $(e_0,E_3)=1$ or $(e_0,E_4)=1$. Say $(e_0,E_4)=1$. By contracting $e_0$ and $E_4$, 
we obtain a minimal elliptic fibration such that the singular fibre over $t$ is of type $\mbox{I}^{\ast}_1$. Thus we see that 
$
\psi^{\ast}_W(t)=2\Gamma_0+2\Gamma_1+\sum_{i=1}^3E_i+3E_4+4e_0+F.
$
Put $G:=\upsilon_{\ast}G^{\prime}$. The above decomposition implies that $(\psi^{\ast}_W(t), G)\geq
3(E_4,G)=3$.  On the other hand, we have $(\psi^{\ast}_W(u), G)=1$, but which is absurd. Thus we
conclude that
$\Xi^{\flat}$ has three  irreducible components assuming $(\Gamma^{\flat},\Xi^{\flat})=4$. Let 
$\Xi^{\flat}=\sum_{i=1}^{3}\Xi^{\flat}_i$ be the irreducible decomposition such that 
$(\Gamma^{\flat},\Xi^{\flat}_1)=2$ and 
$(\Gamma^{\flat},\Xi^{\flat}_i)=1$ for $i=2$, $3$. For $i=2$, $3$, we have 
$
(\Xi^{\flat}_1+\Xi^{\flat}_i)^2=(\Gamma^{\flat},\Xi^{\flat}_1+\Xi^{\flat}_i)^2=9,
$
hence $(\Xi^{\flat}_1,\Xi^{\flat}_i)=2$. In the same way, we have 
$
(\Xi^{\flat}_2+\Xi^{\flat}_3)^2=(\Gamma^{\flat},\Xi^{\flat}_2+\Xi^{\flat}_3)^2=4,
$
hence $(\Xi^{\flat}_2,\Xi^{\flat}_3)=1$. Thus we infer that $\nu(A_1/2\mbox{-}\alpha)=1$ and $\nu(A_1/2\mbox{-}\gamma)=0$, that is, 
$S^{\flat}\simeq \mbox{\boldmath $P$}^2$, $\Xi^{\flat}_1$ is a conic and $\Xi^{\flat}_i$ is a line for $i=2$, $3$. 
Moreover there exists a birational morphism $\eta:W\rightarrow S^{\flat}$ from a smooth rational surface $W$ with $\rho(W)=12$ 
such that $\mbox{Supp }\eta^{-1}_{\ast}\Delta^{\flat}\cup \mbox{Exc }\eta$ has only simple normal crossing singularities whose dual 
graph is as follows.

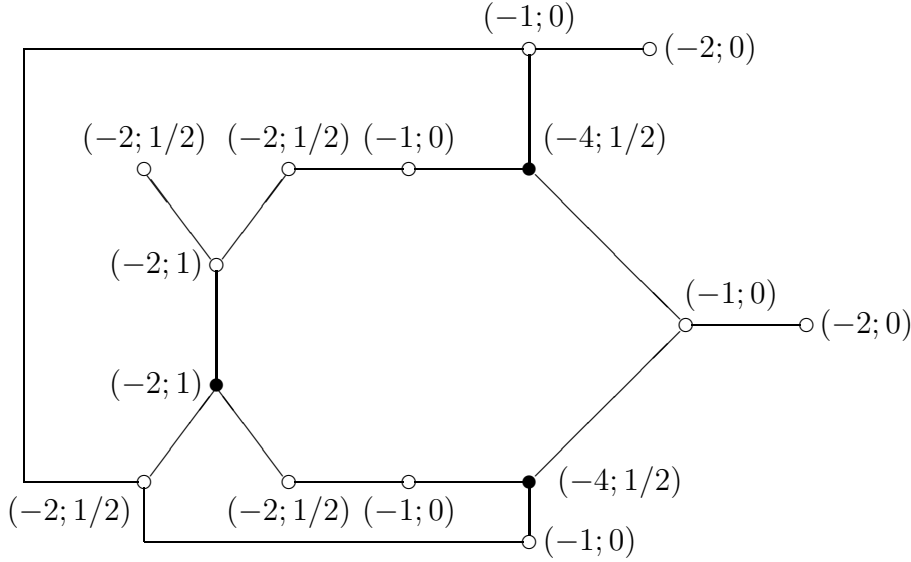
\unitlength 0.8mm
\begin{picture}(150,85)(-50,0)
\put(0,0){\circle{2}}\put(0,0){\makebox(-25,-10){$(-2;1/2)$}}
\put(1,1){\line(3,4){11.5}}
\put(12,16){\circle*{2}}\put(12,16){\makebox(-20,0){$(-2;1)$}}
\put(24,0){\circle{2}}\put(24,0){\makebox(0,-10){$(-2;1/2)$}}
\put(23,1){\line(-3,4){11.5}}
\put(12,16){\line(0,1){19}}
\put(12,36){\circle{2}}\put(12,36){\makebox(-20,0){$(-2;1)$}}
\put(11,37){\line(-3,4){10.5}}
\put(0,52){\circle{2}}\put(0,52){\makebox(0,10){$(-2;1/2)$}}
\put(13,37){\line(3,4){10.5}}
\put(24,52){\circle{2}}\put(24,52){\makebox(0,10){$(-2;1/2)$}}

\put(25,52){\line(1,0){18}}
\put(44,52){\circle{2}}\put(44,52){\makebox(0,10){$(-1;0)$}}
\put(45,52){\line(1,0){18}}
\put(64,52){\circle*{2}}\put(64,52){\makebox(25,10){$(-4;1/2)$}}
\put(65,51){\line(1,-1){24}}

\put(25,0){\line(1,0){18}}
\put(44,0){\circle{2}}\put(44,0){\makebox(0,-10){$(-1;0)$}}
\put(45,0){\line(1,0){18}}
\put(64,0){\circle*{2}}\put(64,0){\makebox(30,0){$(-4;1/2)$}}
\put(65,1){\line(1,1){24}}

\put(90,26){\circle{2}}\put(90,26){\makebox(15,10){$(-1;0)$}}
\put(91,26){\line(1,0){18}}
\put(110,26){\circle{2}}\put(110,26){\makebox(20,0){$(-2;0)$}}

\put(0,-10){\line(0,1){9}}
\put(0,-10){\line(1,0){63}}
\put(64,-10){\circle{2}}\put(64,-10){\makebox(20,0){$(-1;0)$}}
\put(64,-9){\line(0,1){8}}

\put(-20,0){\line(1,0){19}}
\put(-20,0){\line(0,1){72}}
\put(64,52){\line(0,1){19}}
\put(64,72){\circle{2}}\put(64,72){\makebox(0,10){$(-1;0)$}}
\put(-20,72){\line(1,0){83}}

\put(65,72){\line(1,0){18}}
\put(84,72){\circle{2}}\put(84,72){\makebox(20,0){$(-2;0)$}}
\end{picture}
\vskip 20mm

We see that there exists a birational morphism $\lambda:W\rightarrow U$ 
to a smooth rational surface $U$ with $\rho(U)=10$ such that $\lambda^{-1}_{\ast}\Xi^{\flat}$ is nef and 
$(\lambda^{-1}_{\ast}\Xi^{\flat})^2=0$. Thus we get a rational elliptic surface with a section 
$\psi:U\rightarrow \mbox{\boldmath $P$}^1$ which is minimal over $\mbox{\boldmath $P$}^1$ 
with singular fibres $\psi^{\ast}(t)$ of type $\mbox{I}^{\ast}_1$ and $\psi^{\ast}(u)$ of 
type $\mbox{III}$. Noting that we have $\sum_{v\neq t,u}e_{\mbox{top}}(\psi^{-1}(v))=2$ and that
there exists a singular fibre
$\psi^{\ast}(v)$ whose dual graph contains a subgraph of type $A_1$, we see that 
$\psi^{\ast}(v)$ is of type $\mbox{I}_2$ and $\psi$ is smooth over $\mbox{\boldmath $P$}^1\setminus
\{t,$ $u$, $v\}$.  Thus we see that we are in the case $(1)$. Assume that
$(\Gamma^{\flat},\Xi^{\flat})=2$. Then we have 
$(\Xi^{\flat})^2=(\Gamma^{\flat},\Xi^{\flat})^2=4$, which implies that $\Xi^{\flat}$ is irreducible, and 
$
3\nu(A_3)+\nu(A_3/2\mbox{-}\delta)+\nu(A_3/2\mbox{-}\zeta)=3.
$
Since we have $\nu(A_3/2\mbox{-}\delta)+\nu(A_3/2\mbox{-}\zeta)\leq 2$, we see that $\nu(A_3)=1$, 
$\nu(A_3/2\mbox{-}\alpha)=\nu(A_3/2\mbox{-}\delta)=\nu(A_3/2\mbox{-}\zeta)=0$ and we get
$e_{\mbox{top}}(\Xi^{\flat \nu})=2\nu(A_1/2\mbox{-}\alpha)$. 
We note that $\nu(A_1/2\mbox{-}\alpha)=0$ or $1$ but in fact we can show that $\nu(A_1/2\mbox{-}\alpha)=1$ as follows. 
Assume that $\nu(A_1/2\mbox{-}\alpha)=0$. Then we have $\nu(A_1/2\mbox{-}\gamma)=1$ and $e_{\mbox{top}}(\Xi^{\flat \nu})=0$, 
hence $\Xi^{\flat}$ is a smooth elliptic curve. We see that there exists a birational morphism 
$\eta:W\rightarrow S^{\flat}$ from a smooth rational surface $W$ with $\rho(W)=11$ such that 
$\mbox{Supp }\eta^{-1}_{\ast}\Delta^{\flat}\cup \mbox{Exc }\eta$ has only simple normal crossing singularities whose dual 
graph is as follows.

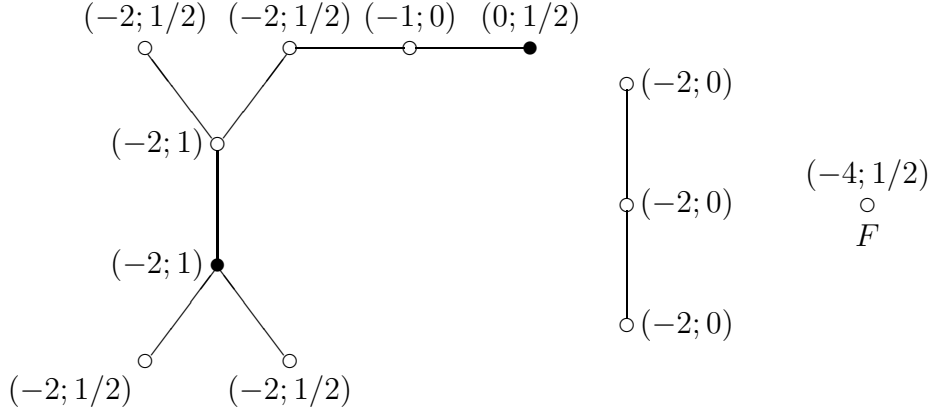
\unitlength 0.8mm
\begin{picture}(150,80)(-50,0)
\put(0,0){\circle{2}}\put(0,0){\makebox(-25,-10){$(-2;1/2)$}}
\put(1,1){\line(3,4){11.5}}
\put(12,16){\circle*{2}}\put(12,16){\makebox(-20,0){$(-2;1)$}}
\put(24,0){\circle{2}}\put(24,0){\makebox(0,-10){$(-2;1/2)$}}
\put(23,1){\line(-3,4){11.5}}
\put(12,16){\line(0,1){19}}
\put(12,36){\circle{2}}\put(12,36){\makebox(-20,0){$(-2;1)$}}
\put(11,37){\line(-3,4){10.5}}
\put(0,52){\circle{2}}\put(0,52){\makebox(0,10){$(-2;1/2)$}}
\put(13,37){\line(3,4){10.5}}
\put(24,52){\circle{2}}\put(24,52){\makebox(0,10){$(-2;1/2)$}}

\put(25,52){\line(1,0){18}}
\put(44,52){\circle{2}}\put(44,52){\makebox(0,10){$(-1;0)$}}
\put(45,52){\line(1,0){18}}
\put(64,52){\circle*{2}}\put(64,52){\makebox(0,10){$(0;1/2)$}}

\put(80,26){\circle{2}}\put(80,26){\makebox(20,0){$(-2;0)$}}
\put(80,27){\line(0,1){18}}
\put(80,46){\circle{2}}\put(80,46){\makebox(20,0){$(-2;0)$}}
\put(80,7){\line(0,1){18}}
\put(80,6){\circle{2}}\put(80,6){\makebox(20,0){$(-2;0)$}}

\put(120,26){\circle{2}}\put(120,26){\makebox(0,10){$(-4;1/2)$}}\put(120,26){\makebox(0,-10){$F$}}
\end{picture}
\vskip 20mm
Some multiple of $\eta^{-1}_{\ast}\Xi^{\flat}$ determines an elliptic fibration $\psi_W:W\rightarrow
\mbox{\boldmath $P$}^1$ with a  section. 
Let $e_0$ be a $(-1)$-curve on $W$ which is contained in a fibre of $\psi_W$. Then we have
$(e_0,F)=2$ from 
$2(K_W+\Delta^{\flat W})\sim 0$. Let $\upsilon:W\rightarrow U$ be the contraction $e_0$ and 
$\psi:U\rightarrow \mbox{\boldmath $P$}^1$ be the induced morphism $\psi_W$. Then $\psi$ is
minimal since $\rho(U)=10$ and 
$\upsilon_{\ast}F$ supports a singular fibre of type $\mbox{I}_1$ or $\mbox{II}$. 
We see that $\upsilon_{\ast}\Delta^{\flat W}=(1/2)\psi^{\ast}(t+u+v)$, where $\psi^{\ast}(t)$ is
a singular fibre of type 
$\mbox{I}^{\ast}_{1}$, $\psi^{\ast}(u)=\upsilon_{\ast}\eta^{-1}_{\ast}\Xi^{\flat}$ and
$\psi^{\ast}(v)=\upsilon_{\ast}F$, but which is absurd since $\omega_U=\psi^{\ast}{\cal
O}_{\mbox{\boldmath $P$}^1}(-1)$.  Thus we conclude that $\nu(A_1/2\mbox{-}\alpha)=1$,
$\nu(A_1/2\mbox{-}\gamma)=0$, hence $e_{\mbox{top}}(\Xi^{\flat})=2$, that is, 
$\Xi^{\flat}\simeq \mbox{\boldmath $P$}^1$. We see that there exists a birational morphism 
$\lambda:U\rightarrow S^{\flat}$ from a smooth rational surface $U$ with $\rho(U)=10$ such that 
$\mbox{Supp }\lambda^{-1}_{\ast}\Delta^{\flat}\cup \mbox{Exc }\lambda$ has only normal crossing singularities whose dual 
graph is as follows.

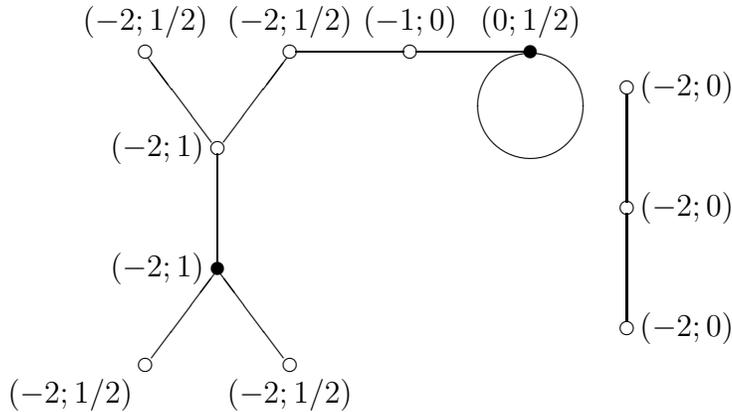
\unitlength 0.8mm
\begin{picture}(150,70)(-50,0)
\put(0,0){\circle{2}}\put(0,0){\makebox(-25,-10){$(-2;1/2)$}}
\put(1,1){\line(3,4){11.5}}
\put(12,16){\circle*{2}}\put(12,16){\makebox(-20,0){$(-2;1)$}}
\put(24,0){\circle{2}}\put(24,0){\makebox(0,-10){$(-2;1/2)$}}
\put(23,1){\line(-3,4){11.5}}
\put(12,16){\line(0,1){19}}
\put(12,36){\circle{2}}\put(12,36){\makebox(-20,0){$(-2;1)$}}
\put(11,37){\line(-3,4){10.5}}
\put(0,52){\circle{2}}\put(0,52){\makebox(0,10){$(-2;1/2)$}}
\put(13,37){\line(3,4){10.5}}
\put(24,52){\circle{2}}\put(24,52){\makebox(0,10){$(-2;1/2)$}}

\put(25,52){\line(1,0){18}}
\put(44,52){\circle{2}}\put(44,52){\makebox(0,10){$(-1;0)$}}
\put(45,52){\line(1,0){18}}
\put(64,52){\circle*{2}}\put(64,52){\makebox(0,10){$(0;1/2)$}}

\put(64,43){\circle{20}}

\put(80,26){\circle{2}}\put(80,26){\makebox(20,0){$(-2;0)$}}
\put(80,27){\line(0,1){18}}
\put(80,46){\circle{2}}\put(80,46){\makebox(20,0){$(-2;0)$}}
\put(80,7){\line(0,1){18}}
\put(80,6){\circle{2}}\put(80,6){\makebox(20,0){$(-2;0)$}}
\end{picture}
\vskip 20mm
Some multiple of $\lambda^{-1}_{\ast}\Xi^{\flat}$ determines a 
minimal elliptic fibration with a section $\psi:U\rightarrow \mbox{\boldmath $P$}^1$ 
with singular fibres $\psi^{\ast}(t)$ of type $\mbox{I}^{\ast}_{1}$ and $\psi^{\ast}(u)$ of type $\mbox{I}_1$. 
Since we have $\sum_{v\neq t,u}e_{\mbox{top}}(\psi^{-1}(v))=4$ and there exists a singular fibre
$\psi^{\ast}(v)$  whose dual graph has a subgraph of type $A_3$, we see that $\psi^{\ast}(v)$ is of
type
$\mbox{I}_4$ and $\psi$ is smooth except over $t$, $u$ and $v$. Thus we are in the case $(2)$.
Secondarily, we consider the case in which $S^{\flat}$ has a structure of a conic fibration
$\varphi^{\flat}:S^{\flat}\rightarrow \mbox{\boldmath $P$}^1$. In this case, we have $d=1-s$, hence
$(\Gamma^{\flat})^2=(1-s)/2$. Assume that $s=2$. Then we have $(\Gamma^{\flat})^2=-1/2$, which is
absurd since $\Gamma^{\flat}\cap \mbox{Sing }S^{\flat}$ by our assumption. Assume that $s=1$. Then
we have $\rho(\tilde S^{\flat})=3$, hence the self-intersection numbers of all the irreducible
components of $\tilde \Delta^{\flat}$ are $-1$, $-2$ and $-3$, which implies that all the irreducible
components of $\tilde \Delta^{\flat}$ are invariant under the action of the covering
transformation group. This contradicts the assumption $s=1$. Thus we conclude that $s=0$. 
Since
$\rho(\tilde S^{\flat})=2$, the self-intersection numbers of all the irreducible
components of $\tilde \Delta^{\flat}$ are $-1$ and $-2$. Thus we infer that there exists an
irreducible component $\Gamma^{\flat}_0\subset \Gamma^{\flat}$ such that
$(\Gamma^{\flat}_0)^2=-1/2$. Contracting $\Gamma^{\flat}_0$, our case reduces to the rank one log
del Pezzo case.\hfill\bsquare
\begin{ex}\label{ex:ex2}{\rm By \cite{persson}, there exists minimal rational elliptic surface with a section 
$\psi:U\rightarrow \mbox{\boldmath $P$}^1$ with $\mbox{Typ}(U;\psi)=\mbox{I}^{\ast}_1+\mbox{I}_1+\mbox{I}_4$ 
from which we can easily construct an example of $(S^{\flat},\Delta^{\flat})$ in $(2)$. 
We can also construct an example of $(S^{\flat},\Delta^{\flat})$ in $(1)$ from a pair 
$(\mbox{\boldmath $P$}^1, \mbox{line}+(1/2)\mbox{line}+(1/2)\mbox{line}+(1/2)\mbox{conic})$ with properly chosen alignment. 
}
\end{ex}

\begin{pr}\label{pr:cltypeIII-2.4} Assume that $\mbox{\em Sing }\tilde S=A_5+A_2+A_1$. 
Then $S^{\flat}$ is a rank one Gorenstein log del Pezzo surface with $\mbox{\em Sing }S^{\flat}=A_1$ and
 there exists a birational morphism $\lambda:U\rightarrow S^{\flat}$ from a smooth projective surface $U$ 
such that $U$ admits a structure of elliptic surface with a section 
$\psi:U\rightarrow \mbox{\boldmath $P$}^1$ which is minimal over $\mbox{\boldmath $P$}^1$ with  
$\mbox{\em Typ}(U;\psi)=\mbox{\em I}^{\ast}_1+\mbox{\em II}+\mbox{\em I}_3$ and that 
$\Delta^{\flat \ U}=(1/2)\psi^{\ast}(t+u)$ where $\psi^{\ast}(t)$ and $\psi^{\ast}(u)$ are the singular fibres of type 
$\mbox{\em I}^{\ast}_1$ and $\mbox{\em II}$ respectively.
\end{pr}
{\it Proof.} The possible singular types of $(S^{\flat},\Delta^{\flat})$ on 
$S^{\flat}\setminus \lfloor\Delta\rfloor$ are types $A_0/2$, $A_1/2\mbox{-}\alpha$, $A_1/2\mbox{-}\gamma$, 
$A_2/2\mbox{-}\beta$, $A_2/2\mbox{-}\epsilon$, $A_5/2\mbox{-}\alpha$,
$A_5/2\mbox{-}\delta$ and $A_5/2\mbox{-}\zeta$. We note that 
$\nu(A_1/2\mbox{-}\alpha)+\nu(A_1/2\mbox{-}\gamma)=1$, $\nu(A_2/2\mbox{-}\beta)+\nu(A_2/2\mbox{-}\epsilon)=1$ and 
$\nu(A_5/2\mbox{-}\alpha)+\nu(A_5/2\mbox{-}\delta)+\nu(A_5/2\mbox{-}\zeta)=1$. Firstly we consider
the case in which $S^{\flat}$ is a rank one log del Pezzo surface. Assume that $s=2$. Then we have
$\rho(\tilde S^{\flat})=2$ and the self-intersection numbers of all the irreducible components of
$\tilde \Delta^{\flat}$ are $-1$ and $-3$, which is absurd. Thus we have $s=1$, hence $\rho(\tilde
S^{\flat})=1$. Since
$d=1$, we have
$(\Gamma^{\flat})^2=1/2$, hence $(\Gamma^{\flat},\Xi^{\flat})=3$, $(\Xi^{\flat})^2=18$ and $\delta
K^2-\delta\rho=-1$. From
$
\delta  K^2=\nu(A_1/2\mbox{-}\gamma)+(1/3)\nu(A_2/2\mbox{-}\epsilon)+\nu(A_5/2\mbox{-}\delta)
$
and
$
\delta\rho=1+\nu(A_1/2\mbox{-}\gamma)+\nu(A_2/2\mbox{-}\epsilon)+3\nu(A_5/2\mbox{-}\delta)+3\nu(A_5/2\mbox{-}\zeta)
$
we obtain
$
\delta
K^2-\delta\rho=-1-(2/3)\nu(A_2/2\mbox{-}\epsilon)-2\nu(A_5/2\mbox{-}\delta)-3\nu(A_5/2\mbox{-}\zeta).
$
Thus we get
$
(2/3)\nu(A_2/2\mbox{-}\epsilon)+2\nu(A_5/2\mbox{-}\delta)+3\nu(A_5/2\mbox{-}\zeta)=0,
$
hence $\nu(A_2/2\mbox{-}\epsilon)=\nu(A_5/2\mbox{-}\delta)=\nu(A_5/2\mbox{-}\zeta)=0$ and 
$\nu(A_2/2\mbox{-}\beta)=\nu(A_5/2\mbox{-}\alpha)=1$. Moreover, since we have

\begin{eqnarray*}
\deg\mbox{Diff}_{\Xi^{\flat \nu}}(\Gamma^{\flat})&=&3+2\nu(A_1/2\mbox{-}\alpha)+2\nu(A_2/2\mbox{-}\beta)+6\nu(A_5/2\mbox{-}\alpha)\\
&=&11+2\nu(A_1/2\mbox{-}\alpha),
\end{eqnarray*}
we have $e_{\mbox{top}}(\Xi^{\flat \nu})=\deg\mbox{Diff}_{\Xi^{\flat \nu}}(\Gamma^{\flat})-9=2+2\nu(A_1/2\mbox{-}\alpha)$. 
The number of irreducible components of $\Xi^{\flat}$ is one or two. Firstly, 
we shall show that $\Xi^{\flat}$ consists of two irreducible components. 
Assume that $\Xi^{\flat}$ is irreducible. Then since we have $e_{\mbox{top}}(\Xi^{\flat \nu})\leq 2$, 
we have $\nu(A_1/2\mbox{-}\alpha)=0$ and $\nu(A_1/2\mbox{-}\gamma)=1$. Thus we see that there exists a birational morphism 
$\eta^{\prime}:W^{\prime}\rightarrow S^{\flat}$ from a smooth rational surface $W^{\prime}$ with $\rho(W^{\prime})=15$ such that 
$\mbox{Supp }\eta^{\prime -1}_{\ast}\Delta^{\flat}\cup \mbox{Exc }\eta^{\prime}$ has only simple normal crossing singularities whose dual 
graph is as follows.

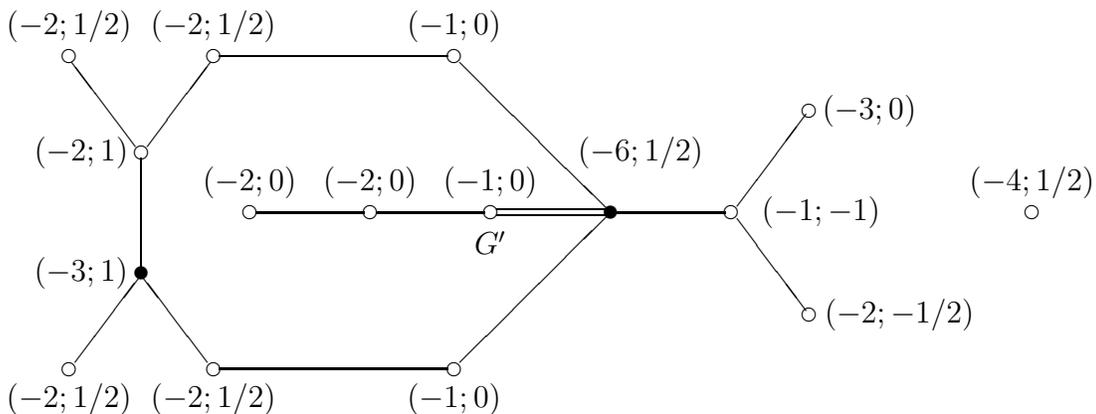
\unitlength 0.8mm
\begin{picture}(150,75)(-30,0)
\put(0,0){\circle{2}}\put(0,0){\makebox(0,-10){$(-2;1/2)$}}
\put(1,1){\line(3,4){11.5}}
\put(12,16){\circle*{2}}\put(12,16){\makebox(-20,0){$(-3;1)$}}
\put(24,0){\circle{2}}\put(24,0){\makebox(0,-10){$(-2;1/2)$}}
\put(23,1){\line(-3,4){11.5}}
\put(12,16){\line(0,1){19}}
\put(12,36){\circle{2}}\put(12,36){\makebox(-20,0){$(-2;1)$}}
\put(11,37){\line(-3,4){10.5}}
\put(0,52){\circle{2}}\put(0,52){\makebox(0,10){$(-2;1/2)$}}
\put(13,37){\line(3,4){10.5}}
\put(24,52){\circle{2}}\put(24,52){\makebox(0,10){$(-2;1/2)$}}

\put(25,52){\line(1,0){38}}
\put(64,52){\circle{2}}\put(64,52){\makebox(0,10){$(-1;0)$}}
\put(65,51){\line(1,-1){25}}
\put(25,0){\line(1,0){38}}
\put(64,0){\circle{2}}\put(64,0){\makebox(0,-10){$(-1;0)$}}
\put(65,1){\line(1,1){25}}
\put(90,26){\circle*{2}}\put(90,26){\makebox(10,20){$(-6;1/2)$}}

\put(89,26.5){\line(-1,0){18}}
\put(89,25.5){\line(-1,0){18}}
\put(70,26){\circle{2}}\put(70,26){\makebox(0,10){$(-1;0)$}}\put(70,26){\makebox(0,-10){$G^{\prime}$}}
\put(69,26){\line(-1,0){18}}
\put(50,26){\circle{2}}\put(50,26){\makebox(0,10){$(-2;0)$}}
\put(49,26){\line(-1,0){18}}
\put(30,26){\circle{2}}\put(30,26){\makebox(0,10){$(-2;0)$}}

\put(91,26){\line(1,0){18}}
\put(110,26){\circle{2}}\put(110,26){\makebox(30,0){$(-1;-1)$}}
\put(111,27){\line(3,4){11.5}}
\put(111,25){\line(3,-4){11.5}}
\put(123,43){\circle{2}}\put(123,43){\makebox(20,0){$(-3;0)$}}
\put(123,9){\circle{2}}\put(123,9){\makebox(30,0){$(-2;-1/2)$}}

\put(160,26){\circle{2}}\put(160,26){\makebox(0,10){$(-4;1/2)$}}
\end{picture}
\vskip 20mm

From the above dual graph, we see that there exists a birational 
morphism $\upsilon:W^{\prime}\rightarrow W$ to a rational surface $W$ with $\rho(W)=12$ such that 
$\eta^{\prime}$ factors into $\eta\circ\upsilon$ where $\eta:W\rightarrow S^{\flat}$ 
is a birational morphism and that $(\eta^{-1}_{\ast}\Xi^{\flat})^2=0$. We note that
$\Delta^{\flat W}\geq 0$. Some multiple of
$\eta^{-1}_{\ast}\Xi^{\flat}$  determines an elliptic fibration $\psi_W:W\rightarrow \mbox{\boldmath
$P$}^1$ with a section such that 
$\psi^{\ast}_W(u)=\eta^{-1}_{\ast}\Xi^{\flat}$ is a singular fibre of type $\mbox{II}$. 
Put $G:=\upsilon_{\ast}G^{\prime}$. Since $(G,\psi^{\ast}_W(u))=(G,\eta^{-1}_{\ast}\Xi^{\flat})=2$
we have $(G,\psi^{\ast}(t))=2$,  where $t:=\psi_W(\eta^{-1}_{\ast}\Gamma^{\flat})$. On the other
hand, by the previous argument in the  proof of Proposition~\ref{pr:cltypeIII-2.3}, we have
$(G,\psi^{\ast}_W(t))\geq 4$, which is a contradiction.  Thus we conclude that $\Xi^{\flat}$
consists of two irreducible components.  Let $\Xi^{\flat}=\Xi^{\flat}_1+\Xi^{\flat}_2$ be the
irreducible decomposition such that 
$(\Gamma^{\flat},\Xi^{\flat}_1)=2$ and $(\Gamma^{\flat},\Xi^{\flat}_2)=1$. We note that since we have
$(\Xi^{\flat}_1)^2=2(\Gamma^{\flat},\Xi^{\flat}_1)^2=8$ and $(\Xi^{\flat}_2)^2=2(\Gamma^{\flat},\Xi^{\flat}_2)^2=2$, 
we have $18=(\Xi^{\flat})^2=10+2(\Xi^{\flat}_1,\Xi^{\flat}_2)$, 
hence $(\Xi^{\flat}_1,\Xi^{\flat}_2)=4$,  
which implies that $\nu(A_1/2\mbox{-}\alpha)=1$ and $\nu(A_1/2\mbox{-}\gamma)=0$,
hence $e_{\mbox{top}}(\Xi^{\flat \nu})=4$, that is, 
$\Xi^{\flat \nu}_i\simeq \mbox{\boldmath $P$}^1$ $(i=1$, $2)$. From 
$K_{S^{\flat}}+\Xi^{\flat}_2+\Gamma^{\flat}+(1/2)\Xi^{\flat}_1\sim_{\mbox{\boldmath $Q$}}(1/2)\Xi^{\flat}_2$, 
we obtain 
$$
\deg\mbox{Diff}_{\Xi^{\flat \nu}_2}(\Gamma^{\flat}+\frac{1}{2}\Xi^{\flat}_1)=
e_{\mbox{top}}(\Xi^{\flat \nu}_2)+\frac{1}{2}(\Xi^{\flat}_2)^2=3.
$$ 
On the other hand, we have 
\begin{eqnarray*}
\deg\mbox{Diff}_{\Xi^{\flat \nu}_2}(\Gamma^{\flat}+\frac{1}{2}\Xi^{\flat}_1) & = & (\Gamma^{\flat},\Xi^{\flat}_2)+\frac{1}{2}(\Xi^{\flat}_1,\Xi^{\flat}_2)+\deg\mbox{Diff}_{\Xi^{\flat \nu}_2}(0)\\
& = & 3+\deg\mbox{Diff}_{\Xi^{\flat \nu}_2}(0).
\end{eqnarray*}
Thus we get $\deg\mbox{Diff}_{\Xi^{\flat \nu}_2}(0)=0$ and consequently, we infer that the
type 
$A_2/2\mbox{-}\beta$ point lies on $\Xi^{\flat}_1$. Thus we see that there exists a birational morphism 
$\eta:W\rightarrow S^{\flat}$ from a smooth rational surface $W$ with $\rho(W)=13$ such that 
$\mbox{Supp }\eta^{-1}_{\ast}\Delta^{\flat}\cup \mbox{Exc }\eta$ has only simple normal crossing singularities whose dual 
graph is as follows.

\unitlength 0.8mm
\begin{picture}(150,75)(-30,0)
\put(0,0){\circle{2}}\put(0,0){\makebox(0,-10){$(-2;1/2)$}}
\put(1,1){\line(3,4){11.5}}
\put(12,16){\circle*{2}}\put(12,16){\makebox(-20,0){$(-2;1)$}}
\put(24,0){\circle*{2}}\put(24,0){\makebox(20,10){$(-2;1/2)$}}
\put(23,1){\line(-3,4){11.5}}
\put(12,16){\line(0,1){19}}
\put(12,36){\circle{2}}\put(12,36){\makebox(-20,0){$(-2;1)$}}
\put(11,37){\line(-3,4){10.5}}
\put(0,52){\circle{2}}\put(0,52){\makebox(0,10){$(-2;1/2)$}}
\put(13,37){\line(3,4){10.5}}
\put(24,52){\circle{2}}\put(24,52){\makebox(0,10){$(-2;1/2)$}}

\put(25,52){\line(1,0){18}}
\put(44,52){\circle{2}}\put(44,52){\makebox(0,10){$(-1;0)$}}
\put(45,51){\line(1,-1){25}}
\put(25,0){\line(1,0){18}}
\put(44,0){\circle{2}}\put(44,0){\makebox(0,-10){$(-1;0)$}}
\put(45,1){\line(1,1){25}}
\put(70,26){\circle*{2}}\put(70,26){\makebox(10,20){$(-6;1/2)$}}

\put(71,26){\line(1,0){18}}
\put(90,26){\circle{2}}\put(90,26){\makebox(30,0){$(-1;-1)$}}
\put(91,27){\line(3,4){11.5}}
\put(91,25){\line(3,-4){11.5}}
\put(103,43){\circle{2}}\put(103,43){\makebox(20,0){$(-3;0)$}}
\put(103,9){\circle{2}}\put(103,9){\makebox(30,0){$(-2;-1/2)$}}

\put(24,0){\line(0,-1){10}}
\put(24,-10){\line(1,0){45}}
\put(70,-10){\circle{2}}\put(70,-10){\makebox(0,-10){$(-1;0)$}}
\put(70,-9){\line(0,1){36}}
\put(71,-10){\line(1,0){18}}
\put(90,-10){\circle{2}}\put(90,-10){\makebox(0,-10){$(-2;0)$}}
\put(91,-10){\line(1,0){18}}
\put(110,-10){\circle{2}}\put(110,-10){\makebox(0,-10){$(-2;0)$}}

\end{picture}
\vskip 20mm
We see that there exists a birational morphism $\upsilon:W\rightarrow U$ to a rational surface with 
$\rho(U)=10$ and $\Delta^{\flat U}\geq 0$ such that $\eta$ 
factors into $\lambda\circ\upsilon$, where $\lambda:U\rightarrow S^{\flat}$ is a birational 
morphism and that $\lambda^{-1}_{\ast}\Xi^{\flat}$ is nef
with $(\lambda^{-1}_{\ast}\Xi^{\flat})^2=0$. Some multiple $\lambda^{-1}_{\ast}\Xi^{\flat}$
determines an elliptic fibration $\psi:U\rightarrow \mbox{\boldmath $P$}^1$ and we see that 
$\mbox{Typ}(U;\psi)=\mbox{I}^{\ast}_1+\mbox{II}+\mbox{I}_3$ as in the same way in the proof of 
Proposition~\ref{pr:cltypeIII-2.2}. Secondarily, we consider the case in which there exists a
structure of a conic fibration $\varphi^{\flat}:S^{\flat}\rightarrow \mbox{\boldmath $P$}^1$.
Assume that $s=2$. Then we have $\rho(\tilde S^{\flat})=4$ and the self-intersection numbers of all
the irreducible components are $(a)$ $-1$, $-2$, $-2$ and $-5$, $(b)$ $-1$, $-2$, $-3$ and $-4$ or
$(c)$
$-1$, $-1$, $-3$ and $-5$, which is absurd with the assumption $s=2$. Assume that $s=1$. 
Then we have $\rho(\tilde S^{\flat})=3$ and the self-intersection numbers of all
the irreducible components are $-1$, $-2$ and $-4$, which is absurd with the assumption $s=1$.
Assume that $s=0$. 
Then we have $\rho(\tilde S^{\flat})=2$ and the self-intersection numbers of all
the irreducible components are $-1$, $-3$. Thus we see that there there exists an irreducible
component $\Gamma^{\flat}_0\subset \Gamma^{\flat}$ such that $(\Gamma^{\flat}_0)^2=-1/2$.
Contracting $\Gamma^{\flat}_0$, our case reduces to the rank one log del Pezzo case. Thus we get the
assertion.
\hfill\bsquare
\begin{ex}\label{ex:ex3}{\rm  We shall follow the notations in Example~\ref{ex:ex1}. 
By
\cite{shioda}, Theorem $8.6$, $(8.12)$, we have 
$$
\frac{1}{12}=<-P,-P>=2+2((-P),(O))-\mbox{contr}_{t}(-P)-\mbox{contr}_{v}(-P)
$$
and $((-P),(O))\in \mbox{\boldmath $Z$}$, we see that $((-P),\Theta_{t,i})=1$ for $i=2$ or $3$ and
that 
$((-P),\Theta_{v,j})=1$ for $j=1$ or $2$, hence $((-P),(O))=0$. From \cite{shioda}, Theorem $8.6$,
$(8.11)$, we have
$$
-\frac{1}{4}=<-P,Q>=1-((-P),(Q))-\mbox{contr}_{t}(-P,Q).
$$
Noting that by \cite{shioda}, $(8.16)$,
$$
\mbox{contr}_{t}(-P,Q)=\left\{
\begin{array}{ll}
5/4 & \mbox{if $(-P)$ and $(Q)$ intersects the same components of $\psi^{\ast}(t)$}\\
3/4& \mbox{otherwise}
\end{array}
\right.
$$
and $((-P),(Q))\in \mbox{\boldmath $Z$}$, we infer that $((-P),(Q))=0$ 
and $((-P),\Theta_{t,i})=((Q),\Theta_{t,i})=1$ for $i=2$ or $3$. 
Let $\lambda:U\rightarrow S^{\flat}$ be the contraction of all the curves $(O)$, $(-P)$, $(Q)$, 
$\{\Theta_{t,j}|j\neq i,\ 5\}$ and $\{\Theta_{v,k}|k=1,\ 2\}$ and put 
$\Delta^{\flat}:=(1/2)\lambda_{\ast} \psi^{\ast}(t+u)$. Then $(S^{\flat},\Delta^{\flat})$ gives an 
example of the log surfaces as in Proposition~\ref{pr:cltypeIII-2.4}.
}\end{ex}
\begin{pr}\label{pr:cltypeIII-2.5} Assume that $\mbox{\em Sing }\tilde
S=4A_1$. Then one of the following holds.
\begin{description}
\item[ $(1)$ ] $S^{\flat}\simeq \mbox{\boldmath $P$}^1\times\mbox{\boldmath $P$}^1$ and
$\Delta^{\flat}=\sum_{i=0}^1\Gamma^{\flat}_i+(1/2)\sum_{j=1}^4\Xi^{\flat}_j$, where
$\Gamma^{\flat}_0$, $\Xi^{\flat}_j$ are fibres of the first
projection for $j=1$, $2$ and $\Gamma^{\flat}_1$, $\Xi^{\flat}_j$ are fibres of the second 
projection for $j=3$, $4$.
\item[ $(2)$ ] $(S^{\flat},\Delta^{\flat})$ is log terminal and $S^{\flat}$ has a structure of
conic fibration
$\varphi^{\flat}:S^{\flat}\rightarrow
\mbox{\boldmath $P$}^1$ with $\rho(S^{\flat}/\mbox{\boldmath $P$}^1)=1$ and $\mbox{\em
Typ}(S^{\flat},\Delta^{\flat};\varphi^{\flat})=((\mbox{\em
I}\mbox{-}2)_{\infty}+(\mbox{\em I}\mbox{-}2)_1+(\mbox{\em
II}\mbox{-}1)_2;(\mbox{\em II}\mbox{-}1)_1)$ possibly after operating ${\cal S}$-elementary transformations.
$\Delta^{\flat}=\Gamma^{\flat}_0+\Gamma^{\flat}_1+(1/2)\Xi^{\flat}$, where $\Gamma^{\flat}_0$ is a
smooth rational curve with $(\Gamma^{\flat}_0, \varphi^{\flat \ast}(t))=2$ for $t\in
\mbox{\boldmath
$P$}^1$ and $(\Gamma^{\flat}_0)^2=0$, $\Gamma^{\flat}_1$ and $\Xi^{\flat}$ are fibres of
$\varphi^{\flat}$ with reduced structure.
\item[ $(3)$ ] $S^{\flat}$ is a rank one Gorenstein log del Pezzo surface with $\mbox{\em Sing
}S^{\flat}=A_1$ and
 there exists a birational morphism $\lambda:U\rightarrow S^{\flat}$ from a smooth projective surface $U$ 
such that $U$ admits a structure of elliptic surface with a section 
$\psi:U\rightarrow \mbox{\boldmath $P$}^1$ which is minimal over $\mbox{\boldmath $P$}^1$ with  
$\mbox{\em Typ}(U;\psi)=\mbox{\em I}^{\ast}_2+2\mbox{\em I}_1$ and that 
$\Delta^{\flat \ U}=(1/2)\psi^{\ast}(t+u)$ where $\psi^{\ast}(t)$ and $\psi^{\ast}(u)$ are the singular fibres of type 
$\mbox{\em I}^{\ast}_2$ and $\mbox{\em I}_1$ respectively. 
\end{description}
\end{pr}
{\it Proof.} Firstly, we consider the case in which $S^{\flat}$ is a rank one log del Pezzo surface.
In this case, since we have
$s=\rho(\tilde S^{\flat})\geq 2$, we have $s=2$ and $(\Gamma^{\flat},\Xi^{\flat})=4$, hence $d=4$
by (\ref{eqn:d}),
$(\Gamma^{\flat})^2=2$ and
$(\Xi^{\flat})^2=8$ by (\ref{eqn:r1}). The possible of singular types on $S^{\flat}\setminus
\lfloor\Delta^{\flat}\rfloor$ are types
$A_1$, $A_0/2$,
$A_1/2\mbox{-}\alpha$,
$A_1/2\mbox{-}\gamma$. We note that we have
$2\nu(A_1)+\nu(A_1/2\mbox{-}\alpha)+\nu(A_1/2\mbox{-}\gamma)=4$. Since we have $\delta
K^2=\nu(A_1/2\mbox{-}\gamma)$ and $\delta \rho=\nu(A_1)+\nu(A_1/2\mbox{-}\gamma)$, we have $\delta
K^2-\delta\rho=-\nu(A_1)$. On the other hand, we have $\delta
K^2-\delta\rho=-1$ by (\ref{eqn:k-r}). Thus we obtain $\nu(A_1)=1$ and
$\nu(A_1/2\mbox{-}\alpha)+\nu(A_1/2\mbox{-}\gamma)=2$. Noting that we have
$\deg\mbox{Diff}_{\Xi^{\flat\nu}}(\Gamma^{\flat})=4+2\nu(A_1/2\mbox{-}\alpha)$, we obtain  
$e_{\mbox{top}}(\Xi^{\flat\nu})=2\nu(A_1/2\mbox{-}\alpha)$. We shall show that $\Xi^{\flat}$ is
reducible. Assume the contrary. Since we have
$e_{\mbox{top}}(\Xi^{\flat\nu})=2\nu(A_1/2\mbox{-}\alpha)\leq 2$, we have
$\nu(A_1/2\mbox{-}\alpha)\geq 1$, hence $\nu(A_1/2\mbox{-}\gamma)\geq 1$. Therefore, we see that
there exists a birational morphism $\eta:W\rightarrow S^{\flat}$ from a smooth rational surface with
$\rho(W)=10+\nu(A_1/2\mbox{-}\gamma)$ such that $\Delta^{\flat W}\geq 0$, 
$(\Delta^{\flat W}-(1/2)\eta^{-1}_{\ast}\Xi^{\flat},\eta^{-1}_{\ast}\Xi^{\flat})=0$ and
$(\eta^{-1}_{\ast}\Xi^{\flat})^2=0$. Some multiple of $\eta^{-1}_{\ast}\Xi^{\flat}$ defines an
elliptic fibration with a section 
$\psi_W:W\rightarrow \mbox{\boldmath $P$}^1$. Put $t:=\psi_W(\eta^{-1}_{\ast}\Gamma^{\flat})$ and
$u:=\psi_W(\eta^{-1}_{\ast}\Xi^{\flat})$. We see that $\psi_W^{\ast}(t)$ is a singular fibre of type
$\mbox{I}^{\ast}_2$ and  
$\psi_W^{\ast}(u)$ is a singular fibre of type $\mbox{I}_0$, if
$\nu(A_1/2\mbox{-}\gamma)=2$, of type $\mbox{I}_1$, if $\nu(A_1/2\mbox{-}\gamma)=1$, in
particular, $\psi_W$ is minimal over $t$ and $u$. Consequently, we can write
$K_W+(1/2)\psi^{\ast}_W(t+u)+(1/2)F\sim_{\mbox{\boldmath $Q$}}0$, where $F$ is a sum of disjoint
$(-4)$-curves which comes from the minimal resolution of singularities of type
$A_1/2\mbox{-}\gamma$. Let
$\psi:U\rightarrow \mbox{\boldmath $P$}^1$ be the minimal model of $\psi_W$ and
$\upsilon:W\rightarrow U$ be the induced morphism. Then we have
$K_U+(1/2)\psi^{\ast}(t+u)+(1/2)\upsilon_{\ast}F\sim_{\mbox{\boldmath $Q$}}0$ which implies that
$\upsilon_{\ast}F=0$ since
$\omega_U\simeq \psi^{\ast}{\cal O}_{\mbox{\boldmath $P$}^1}(-1)$. Thus we have
$K_W+(1/2)\psi^{\ast}_W(t+u)+(1/2)F=\upsilon^{\ast}(K_U+(1/2)\psi^{\ast}(t+u))$,
hence $K_W+(1/2)F=\upsilon^{\ast}K_U$, but which is absurd. Thus we conclude that $\Xi^{\flat}$ is
reducible. We note that we can write $\Xi^{\flat}=\Xi^{\flat}_1+\Xi^{\flat}_2$, where
$\Xi^{\flat}_i$ are irreducible curves with $(\Gamma^{\flat},\Xi^{\flat}_i)=2$ and hence 
$(\Xi^{\flat}_i)^2=2$ by (\ref{eqn:r1}) for $i=1$, $2$. From
$8=(\Xi^{\flat})^2=4+2(\Xi^{\flat}_1, \Xi^{\flat}_2)$, we have $(\Xi^{\flat}_1, \Xi^{\flat}_2)=2$,
hence $\nu(A_1/2\mbox{-}\alpha)=2$ and $\nu(A_1/2\mbox{-}\gamma)=0$. Thus we see that there exists
a birational morphism $\lambda:U\rightarrow S^{\flat}$ from a smooth rational surface
with $\rho(U)=10$ such that $\Delta^{\flat U}\geq 0$ and $(\Delta^{\flat
U}-\lambda^{-1}_{\ast}\Xi^{\flat},\lambda^{-1}_{\ast}\Xi^{\flat})=0$,
$\lambda^{-1}_{\ast}\Xi^{\flat}$ is nef and
$(\lambda^{-1}_{\ast}\Xi^{\flat})^2=0$. Some multiple defines an elliptic fibration with a
section $\psi:U\rightarrow \mbox{\boldmath $P$}^1$ which is necessarily minimal over
$\mbox{\boldmath $P$}^1$ such that
$\psi_W^{\ast}(t)$ is a singular fibre of type
$\mbox{I}^{\ast}_2$, 
$\psi_W^{\ast}(u)$ is a singular fibre of type
$\mbox{I}_2$, where $t:=\psi_W(\lambda^{-1}_{\ast}\Gamma^{\flat})$
and $u:=\psi_W(\lambda^{-1}_{\ast}\Xi^{\flat})$. Noting that there exists a singular fibre
$\psi^{\ast}(v)$ whose dual graph contains a Dynkin diagram of type $A_1$ as a subgraph and that we
have $\sum_{v\neq t,u}e_{\mbox{top}}(\psi^{-1}(v))=2$, we infer that $\psi^{\ast}(v)$ is of type
$\mbox{I}_2$ and $\psi$ is smooth over $\mbox{\boldmath $P$}^1\setminus\{t,u,v\}$. Thus we
conclude that we are in the case $(2)$. Secondarily, we consider the case in which there
exists a structure of a conic fibration $\varphi^{\flat}:S^{\flat}\rightarrow \mbox{\boldmath
$P$}^1$ with $\rho(S^{\flat}/\mbox{\boldmath $P$}^1)=1$. Let
$\Gamma^{\flat}=\Gamma^{\flat}_0+\Gamma^{\flat}_1$ be the irreducible decomposition. We note that
we have
$d=4-s$ by (\ref{eqn:d}), hence we have
$(\Gamma^{\flat})^2=2-(1/2)s$ and $(\Gamma^{\flat}_0)^2+(\Gamma^{\flat}_1)^2=-(1/2)s$. Since we 
have $\delta
K^2=\nu(A_1/2\mbox{-}\gamma)$ and $\delta
\rho=\nu(A_1)+\nu(A_1/2\mbox{-}\gamma)+4-(\Gamma^{\flat},\Xi^{\flat})$, we have
$\delta K^2-\delta\rho=-\nu(A_1)-4+(\Gamma^{\flat},\Xi^{\flat})$. On the other hand, we have 
$\delta K^2-\delta\rho=(\Gamma^{\flat},\Xi^{\flat})+(1/4)(\Xi^{\flat})^2-6-(1/2)s$ by
(\ref{eqn:k-r}), hence we obtain $(\Xi^{\flat})^2=8+2s-4\nu(A_1)$. Combining this with 
(\ref{eqn:topofx}), we obtain $e_{\mbox{top}}(\Xi^{\flat
\nu})=(\Gamma^{\flat},\Xi^{\flat})-4-s+2\nu(A_1)+2\nu(A_1/2\mbox{-}\alpha)$. 
We note that if $\rho(\tilde S^{\flat})=2$, 
 $(S^{\flat},
\Delta^{+}(\varphi^{\flat};t))\subset\lfloor\Delta^{\flat -}_{\varphi^{\flat}}(t)\rfloor$ is log
canonical and
$\mbox{Supp}\lfloor\mbox{Diff}_{C_{\varphi^{\flat}}(t)}(\Delta^{\flat
-}_{\varphi^{\flat}}(t))\rfloor\subset \lfloor\Delta^{\flat
-}_{\varphi^{\flat}}(t)\rfloor$ for any $t\in \mbox{\boldmath $P$}^1$ by the same argument as in
the proof of Lemma~\ref{lm:typeII}. Assume that $(\Gamma^{\flat}_0,\varphi^{\flat\ast}(t))=2$ for
$t\in \mbox{\boldmath $P$}^1$. Then
$\Gamma^{\flat}_1$ and $\Xi^{\flat}$ are contained in fibres of $\varphi^{\flat}$, hence $s=0$ and
$\rho(\tilde S^{\flat})=2$. Therefore, $\varphi$ has only fibres of type $(\mbox{I-}2)_b$ and
$(\mbox{II-}1)_b$ by Lemma~\ref{lm:dynkin}. We note that
$\varphi^{\flat\ast}(\varphi^{\flat}(\Gamma^{\flat}_1))$ is a fibre of type
$(\mbox{I-}2)_{\infty}$. Since
$\varphi^{\flat}$ induces a double cover $\varphi^{\flat}:\Gamma^{\flat}_0\rightarrow
\mbox{\boldmath $P$}^1$, $\varphi^{\flat}$ has exactly one fibre of type $(\mbox{I-}2)_1$ by
Hurwitz's formula. Thus we conclude that $\varphi^{\flat}$ has exactly one fibre of type
$(\mbox{II-}1)_2$ since $\sum_{p\in
\Gamma^{\flat}_0}m_p(\Gamma_0^{\flat},\Delta^{\flat}-\Gamma^{\flat}_0)=2$ and the other fibres are
of type
$(\mbox{II-}1)_1$, hence we are in the case
$(2)$. 
Assume that $(\Gamma^{\flat}_0,\varphi^{\flat\ast}(t))=1$ for
$t\in \mbox{\boldmath $P$}^1$. We shall show that
$\Gamma^{\flat}_1$ is contained in a fibre. Assume the contrary. Then $(\Xi^{\flat})^2=0$, hence 
$s=0$, which implies $\nu(A_1)=2$. Take type $A_1$ point
$p\in S^{\flat}\setminus\lfloor\Delta^{\flat}\rfloor$ and take $C_{\varphi^{\flat}}(t)$ passing
through $p$. Since $S^{\flat}$ has only Du Val singularities of type $A_1$, we can write
$\varphi^{\flat}(t)=2C_{\varphi^{\flat}}(t)$, hence $(\Gamma^{\flat}_i,C_{\varphi^{\flat}}(t))=1/2$
for $i=0$, $1$. From $(K_{S^{\flat}}+\Gamma^{\flat},C_{\varphi^{\flat}}(t))=0$, we have
$\deg\mbox{Diff}_{C_{\varphi^{\flat}}(t)}(\Gamma^{\flat})=2$, which implies
$m_p(C_{\varphi^{\flat}}(t);\Gamma^{\flat})=0$. Thus we get absurdity and we conclude that
$\Gamma^{\flat}_1$ is contained in a fibre. We note that $(\Gamma^{\flat}_0)^2=-(1/2)s$, hence the
case $s=1$ reduces to the case in which $S^{\flat}$ is rank one log del Pezzo surface by contracting
$\Gamma^{\flat}_0$. Let $\Xi^{\flat}_h$ be the horizontal part of $\Xi^{\flat}$ with respect to
$\varphi^{\flat}$. Assume that $s=0$. In this case, we may assume that $\Xi^{\flat}_h$ is reducible
for if 
$\Xi^{\flat}_h$ is irreducible, since
$\varphi^{\flat}$ has only fibres of type $(\mbox{I-}1)_1$,
$(\mbox{I-}1)_2$ and $(\mbox{I-}3)_1$, we have $\nu((\mbox{I-}1)_2)=0$ and $\nu((\mbox{I-}3)_1)=2$,
hence $(\Xi^{\flat}_h)^2=0$ by applying Hurwitz's formula to the double cover
$\varphi^{\flat}:\Xi^{\flat}_h\rightarrow
\mbox{\boldmath $P$}^1$, thus some multiple defines another conic fibration
$\check\varphi^{\flat}:S^{\flat}\rightarrow \mbox{\boldmath $P$}^1$ with
$(\Gamma^{\flat}_1,\check\varphi^{\flat\ast}(t))=2$, which case was already considered.
Assuming that $\Xi^{\flat}_h$ is reducible, we see that $\varphi^{\flat}$ has only fibres of type
$(\mbox{I-}2)_b$ $(b=1$, $2$, $\infty)$. Thus we see that we are in the case $(1)$. Assume that
$s=2$. Since $\Gamma^{\flat}_0$ is a section of $\varphi^{\flat}$ which does not pass through any
singular point of $S^{\flat}$, $\varphi^{\flat}$ is smooth, hence
$\nu(A_1)=\nu(A_1/2\mbox{-}\gamma)=0$ and $\nu(A_1/2\mbox{-}\alpha)=4$, which implies that 
$(\Xi^{\flat})^2=12$ and $e_{\mbox{top}}(\Xi^{\flat\nu})=6$. We see that $\Xi^{\flat}$ is
reducible. Let $\Xi^{\flat}=\Xi^{\flat}_1+\Xi^{\flat}_2$ be the irreducible decomposition of
$(\Gamma^{\flat}_0,\Xi^{\flat}_1)=2$ and $(\Gamma^{\flat}_1,\Xi^{\flat}_2)=2$. Since
$\Xi^{\flat}_1$ is not a fibre of $\varphi^{\flat}$, $(\Gamma^{\flat}_1,\Xi^{\flat}_1)>0$, which
contradicts $(K_{S^{\flat}}+\Delta^{\flat},\Gamma^{\flat}_1)=0$.
\hfill\bsquare

\begin{ex}\label{ex:ex3}{\rm  The existence of minimal rational elliptic surface with a section 
$\psi:U\rightarrow \mbox{\boldmath $P$}^1$ with
$\mbox{Typ}(U;\psi)=\mbox{I}^{\ast}_2+2\mbox{I}_2$ is known (see
\cite{persson}). From the list in \cite{oguiso-shioda}, we see that
$\mbox{MW}(U_{\eta})=(\mbox{\boldmath $Z$}/2\mbox{\boldmath $Z$})^{\oplus 2}$. Take $P\in
\mbox{MW}(U_{\eta})\setminus \{0\}$.  From 
\cite{shioda}, Theorem $8.6$, $(8.12)$, we have 
$$
0=<P,P>=2+2(PO)-\mbox{contr}_{t}(P)-\mbox{contr}_{u}(P)-\mbox{contr}_{v}(P),
$$
hence 
$$
(PO)=\frac{1}{2}(\mbox{contr}_{t}(P)+\mbox{contr}_{u}(P)+\mbox{contr}_{v}(P))-1.
$$
Let $\psi^{\ast}(t)=\sum_{i=0}^3\Theta_{t,i}+2\sum_{i=4}^6\Theta_{t,i}$ be the type
$\mbox{I}_2^{\ast}$ singular fibre with $(\Theta_{t,5}, \Theta_{t,i})=0$ for any $i$ such that
$0\leq i\leq 3$,
$\psi^{\ast}(u)=\Theta_{u,1}+\Theta_{u,i}$,
$\psi^{\ast}(t)=\Theta_{u,1}+\Theta_{u,i}$
be the type $\mbox{I}_2$ singular fibre and assume that the $0$-section intersects $\Theta_{t,0}$,
$\Theta_{u,0}$ and $\Theta_{v,0}$. Then from
\cite{shioda}, 
$(8.16)$, we have 
$$
\mbox{contr}_{t}(P)=\left\{
\begin{array}{ll}
0 & \mbox{ if $(P\Theta_{t,0})=1$, }\\
1 & \mbox{ if $(P\Theta_{t,1})=1$, }\\
3/2& \mbox{otherwise}
\end{array}
\right.
\mbox{ and }
\mbox{contr}_{v^{\prime}}(P)=\left\{
\begin{array}{ll}
0 & \mbox{ if $(P\Theta_{v^{\prime},0})=1$, }\\
1/2 & \mbox{ if $(P\Theta_{v^{\prime},1})=1$ }
\end{array}
\right.
$$
for $v^{\prime}=u$, $v$. Since $\mbox{contr}_{t}(P)+\mbox{contr}_{u}(P)+\mbox{contr}_{v}(P)\in
2\mbox{\boldmath $N$}$, there exists exactly two types of $P$, $(1)$ $(P\Theta_{t,1})=1$,
$(P\Theta_{u,1})=1$ and $(P\Theta_{v,1})=1$ or $(2)$ $(P\Theta_{t,i})=1$ for some $i>1$,
$(P\Theta_{u,1})=1$ and $(P\Theta_{v,0})=1$. Take two elements $P$, $Q\in
\mbox{MW}(U_{\eta})\setminus \{0\}$ which is distinct from each other. Assume that both of $P$ and
$Q$ are of type $(1)$. Then from \cite{shioda}, Theorem $8.6$,
$(8.11)$, we have
$$
0=<P,Q>=1-(PQ)-\mbox{contr}_{t}(P,Q)-\mbox{contr}_{u}(P,Q)-\mbox{contr}_{v}(P,Q),
$$
hence $(PQ)=-1$ from \cite{shioda}, $(8.16)$, which is absurd. Consequently, there exists $P\in
\mbox{MW}(U_{\eta})$ of type $(2)$. We may assume that $(P\Theta_{v,0})=1$. 
Let
$\lambda:U\rightarrow S^{\flat}$ be the contraction of all the curves $(O)$, $(P)$, 
$\{\Theta_{t,i}|i\neq  5\}$ and $\{\Theta_{v,j}|j=1,\ 2\}$ and put 
$\Delta^{\flat}:=(1/2)\lambda_{\ast} \psi^{\ast}(t+u)$. Then $(S^{\flat},\Delta^{\flat})$ gives an 
example of the log surfaces as in Proposition~\ref{pr:cltypeIII-2.5}.
}\end{ex}

\section{Generalized local fundamental groups for analytic singularities with $\mbox{\boldmath $Q$}$-divisors} 
In this section, we give a theory to calculate local fundamental groups from differents. Let us review here the theory due to Prill.  
For a germ of normal complex analytic spaces $(X,p)$, put 
$\mbox{Reg }X:=\mbox{projlim}_{p\in {\cal U};{\rm open }}\mbox{Reg }{\cal U}$, where $\mbox{Reg }{\cal U}$ is the smooth loci of ${\cal U}$. 
$\pi_1^{{\rm loc}}(\mbox{Reg
}X):=\mbox{projlim}_{p\in {\cal U};{\rm open }}\pi_1(\mbox{Reg }{\cal U})$ called the {\it local fundamental group} for $(X,p)$. We denote by $\hat\pi_1^{{\rm loc}}(\mbox{Reg }X)$ its
profinite completion which is called the {\it local algebraic fundamental group} of $(X,p)$. Let $\Sigma$ be an analytically closed proper subset of $X$.  According to Prill
(\cite{prill}, \S IIB), there exists a contractible open neighbourhood $U$ of $p$ such that there exists a neighbourhood basis
$\{U_{\lambda}\}_{\lambda\in
\Lambda}$ of $p$ satisfying the condition that
$U_{\lambda}\setminus\Sigma$ is a deformation retract of $U\setminus \Sigma$ for any $\lambda\in
\Lambda$. By the definition, we have $\pi_1(U\setminus\Sigma)=\mbox{projlim}_{p\in {\cal U};{\rm open }}\pi_1({\cal U}\setminus\Sigma)$. 
We call such $U$ as above a {\it Prill's good neighbourhood with regard to
$\Sigma$} and we say that 
$\{U_{\lambda}\}_{\lambda\in
\Lambda}$ is a {\it neighbourhood basis associated with $U$}. Recall that $U_{\lambda}$ is also a Prill's good neighbourhood with
regard to $\Sigma$ and for any two Prill's good neighbourhood $U$ and $U^{\prime}$, $U\setminus\Sigma$ and $U^{\prime}\setminus\Sigma$ have the same homotopy type.
In particular, we have $\pi_1^{{\rm loc}}(\mbox{Reg }X)\simeq\pi_1(\mbox{Reg }U)$ for a Prill's good neighbourhood $U$ with
regard to $\mbox{Sing }X$.

\subsection{B(D)-local fundamental groups}
In this section, we introduce a generalized local fundamental group for a pair consisting of a germ of a normal complex analytic space and a Weil divisor on it, which turns out to
appear canonically  when we calculate local fundamental groups from differents. To introduce the generalized notion of local fundamental groups, let us briefly review here the theory
of universal ramified coverings due to M. Kato (\cite{kato}), M. Namba (\cite{namba}) and J.P. Serre (\cite{serre}, Appendix 6.4) according to M. Namba.  Let $B$ be an integral
effective divisor on a connected complex manifold $M$ and let
$B:=\sum_{i\in I}b_iB_i$ be the irreducible decomposition of
$B$.  Fix a base point $x\in M\setminus\mbox{Supp }B$ and let $\gamma_i$ be a loop which starts from $x$ and goes around $B_i$ once in a counterclockwise direction with the center
being a smooth point of
$\mbox{Supp }B$ on $B_i$. Let ${\cal N}(M,B,x)\subset \pi_1(M\setminus \mbox{Supp }B,x)$ denote the normal subgroup generated by all the conjugates of the loops
$\{\gamma_i^{b_i}\}_{i\in I}$. Recall that ${\cal N}(M,B,x)$ is known to be independent from the choice of such loops. We define a {\it $B$-fundamental group} of $M$ by putting  
$$
\overline\pi_1^B(M,x):=\pi_1(M\setminus \mbox{Supp }B,x)/{\cal N}(M,B,x).
$$

Here, let us fix our terminology from the category theory. By a {\it projective system}, we mean a category ${\cal
I}$ such that $\mbox{Hom }_{{\cal I}}(\lambda,\mu)$ is empty or consists of exactly one element
$f_{\lambda,\mu}$ satisfying $f_{\lambda,\mu}\circ f_{\mu,\nu}=f_{\lambda,\nu}$ for any
$\lambda$, $\mu$,  $\nu\in \mbox{Ob }{\cal I}$. An object $\alpha\in \mbox{Ob }{\cal I}$ (
resp. $\omega\in \mbox{Ob }{\cal I}$ ) is called an {\it initial object} ( resp. a {\it final
object} ) if $\mbox{Card Hom }_{\cal I}(\alpha, \lambda)=1$ ( resp. $\mbox{Card Hom }_{\cal
I}(\lambda,\omega)=1$) for any $\lambda\in \mbox{Ob }{\cal I}$. A projective system ${\cal I}$
is said to be {\it cofilterd} if, for any given two objects $\lambda$, $\mu\in \mbox{Ob }{\cal I}$,
there exists $\nu\in \mbox{Ob }{\cal I}$ with $\mbox{Card Hom }_{{\cal I}}(\nu,\lambda)=\mbox{Card Hom }_{{\cal I}}(\nu,\mu)=1$.  A covariant functor $\Phi:{\cal I}^{\circ}\rightarrow
{\cal I}^{\prime \circ}$ between injective systems ${\cal I}^{\circ}$ and ${\cal I}^{\prime \circ}$ is said to be {\it cofinal}, if, for any given
$\lambda^{\prime}\in \mbox{Ob }{\cal I}^{\prime \circ}$, there exists $\lambda\in \mbox{Ob }{\cal I}^{\circ}$
such that $\mbox{Card Hom }_{{\cal I}^{\prime \circ}}(\lambda^{\prime},\Phi(\lambda))=1$. 
We shall also say that a projective subsystem ${\cal I}^{\prime}$ in a projective system ${\cal I}$ is cofinal in ${\cal I}$ 
if the dual embedding functor from ${\cal I}^{\prime \circ}$
to ${\cal I}^{\circ}$ is cofinal.  (see \cite{artin-mazur}, Appendix (1.5), \cite{sga4}, Expos\'e I, Definition 2.7 and Definition
8.1.1). A finite covering $f:N\rightarrow M$ from a connected normal complex analytic space $N$ which is \'etale over $M\setminus \mbox{Supp }B$ is said to be branching at most (resp.
branching ) at
$B$, if the ramification index $e_{\tilde B_{i,j}}(f)$ of $f$ at any prime divisor $\tilde B_{i,j}$ such that $f(\tilde B_{i,j})=B_i$ divides (resp. is equals to ) $b_i$ for any $i\in
I$. Let
$FC^{\leq B}(M)$ (resp. $FC^B(M)$ ) denote the category of finite coverings over $M$  branching at most (resp. branching ) at $B$.  Let $FGC^{\leq B}(M)$ (resp. $FGC^B(M)$ denote the
full subcategory of $FC^{\leq B}(M)$ whose objects consists of Galois covers over $M$. Triplet $(N, f, y)$, where $(N,f)\in \mbox{Ob }FC^{\leq B}(M)$ and $y\in f^{-1}(x)$ are called
{\it pointed finite coverings branching at most at $B$}. Pointed finite coverings branching at most at $B$ and morphisms $f_{\lambda,\mu}\in
\mbox{Hom}_{FC^{\leq B}(M)}((N_{\mu},f_{\mu}),(N_{\lambda},f_{\lambda}))$ such that
$f_{\lambda,\mu}(y_{\mu})=y_{\lambda}$, where $(N_{\mu},f_{\mu}, y_{\mu})$ and $(N_{\lambda},f_{\lambda}, y_{\lambda})$ are two pointed finite coverings branching at most at $B$ 
form a
projective system denoted by $FC^{\leq B}(M)^p$. We also define the projective subsystems $F(G)C^{(\leq) B}(M)^p$ in the same way. 
From \cite{namba}, Lemma 1.3.1, Theorem 1.3.8 and 
Theorem 1.3.9, we see that there exists a canonical functor $\Psi$ from $FC^{\leq B}(M)^p$ to the projective system of subgroups of finite indices in $\overline\pi_1^B(M,x)$ such that 
$$
\Psi((N,f,y))=f_{\ast}\pi_1(N\setminus \mbox{Supp }f^{-1}B,y)/{\cal N}(M,B,x)\subset \overline\pi_1^B(M,x)
$$
for 
$(N,f,y)\in \mbox{Ob}FC^{\leq B}(M)^p$ and that the functor $\Psi$ defines an equivalence between the above two projective systems. Thus by using the basic group theory, we obtain the
following lemma. 

\begin{lm}\label{lm:fgc}
$FGC^{\leq B}(M)^p$ $($ resp. $FGC^B(M)^p$ $)$ is cofilterd and cofinal in $FC^{\leq B}(M)^p$ $($ resp. $FC^B(M)^p$ $)$ and hence in particular, $FGC^B(M)^p$ is cofinal in
$FGC^{\leq B}(M)^p$ if
$FGC^B(M)^p$ is not empty. 
\end{lm}

\begin{rk}{\em Let $\overline\pi_1^B(M,x)^{\wedge}$ denote the profinite completion of $\overline\pi_1^B(M,x)$ called the {\it $B$-algebraic fundamental group} of $M$.
Assume that $FGC^B(M)^p$ is not empty. Then by Lemma~\ref{lm:fgc}, we have
$$
\overline\pi_1^B(M,x)^{\wedge}\simeq \mbox{projlim}_{(N,f,y)\in {\rm Ob}FGC^B(M)^p}\mbox{Gal }(N/M), 
$$
where $\mbox{Gal }(N/M):=\overline\pi_1^B(M,x)/f_{\ast}\overline\pi_1^B(N,y)$.
}\end{rk}

In what follows, we shall use the following notation. Let $X$ be a normal Stein space or a germ of normal complex analytic spaces with a point $p\in X$. 
$\mbox{Weil }X$ is the free abelian group generated by prime divisors on $X$ and 
$\mbox{Div }X$ is the subgroup of $\mbox{Weil }X$ generated by Cartier divisors.  $\mbox{Div}_{\mbox{\boldmath $Q$}}X$ is the 
{\bf Q}-submodule of $\mbox{Weil }X\otimes \mbox{\boldmath $Q$}$ generated by $\mbox{Div }X$. Let $f:Y\rightarrow X$ be a finite morphism between 
normal Stein spaces or germs of normal complex analytic spaces. The pull-back homomorphism 
$f^{\ast}:\mbox{Weil }X\rightarrow \mbox{Weil }Y$ canonically  extends to a homomorphism 
$f^{\ast}:\mbox{Weil }X\otimes{\mbox{\boldmath $Q$}}\rightarrow \mbox{Weil }Y\otimes{\mbox{\boldmath $Q$}}$.

\begin{df}{\em 
For a germ of normal complex analytic spaces $(X,p)$ and $B\in \mbox{Weil }X$, we 
define a {\it $B$-local fundamental group} of $X$ with respect to $B$ as follows:
$$
\overline\pi_{1,{\rm loc}}^B(\mbox{Reg }X):=\mbox{projlim}_{p\in {\cal U};{\rm open}}\overline\pi_1^B(\mbox{Reg }{\cal U}). 
$$
Moreover, by $\overline\pi_{1,{\rm loc}}^B(\mbox{Reg }X)^{\wedge}$, we
mean the profinite completion of $\overline\pi_{1,{\rm loc}}^B(\mbox{Reg }X)$.
}\end{df}

\begin{rk}\label{rk:relation with usual ones}{\em We note that apparently, we have $\overline\pi_{1,{\rm loc}}^B(\mbox{Reg }X)=\pi_1^{{\rm loc}}(\mbox{Reg }X)$ and if $B$ is reduced.
}\end{rk}

\begin{df}{\em  For $D\in {\mbox{Weil }}X\otimes {\mbox{\boldmath $Q$}}$, a finite
surjective morphism $f:Y\rightarrow X$, where $Y$ is a germ of irreducible normal complex analytic
spaces such that 
$f^{\ast}D\in \mbox{Weil }Y$ is called a {\it integral cover
with respect to }$D$. A integral cover $f:Y\rightarrow X$ with respect to $D$ is called a {\it strict integral
cover},  if $e_{\tilde \Gamma}(f)=e_{\Gamma}(D)$ for any prime divisors $\tilde\Gamma$ on $Y$ and $\Gamma$ on $X$ such that $f(\tilde\Gamma)=\Gamma$, where $e_{\tilde \Gamma}(f)$
denotes the ramification index of $f$ at $\tilde\Gamma$ and $e_{\Gamma}(D):=[\mbox{\boldmath $Z$}(\mbox{mult}_{\Gamma}D):\mbox{\boldmath
$Z$}(\mbox{mult}_{\Gamma}D)\cap \mbox{\boldmath $Z$}]\in \mbox{\boldmath $N$}$. By $\mbox{Int}^m(X;D)$ (
resp.
$\mbox{Int}^{\dag}(X;D)$ ), we mean a category of integral covers ( resp. a category of strict integral covers ) with respect to $D$. We shall also define categories
$\mbox{Int}^{m(\dag)}(G)(X;D)^{(p)}$ similarly as before. }\end{df}

Let ${\cal X}$ be an arcwise connected, locally arcwise
connected, Hausdorff topological space, A continuous map $f:{\cal Y}\rightarrow {\cal X}$ from a Hausdorff topological space ${\cal Y}$ with discrete finite fibres
is called a {\it finite topological covering} if for any $x\in{\cal X}$, there exists an arcwise connected open neighbourhood ${\cal U}$ of
$x\in{\cal X}$ such that the restriction of $\pi$ to each arcwise connected component of of $\pi^{-1}({\cal U})$ gives a homeomorphism
onto ${\cal U}$. The following lemma is nothing but a consequence from the first covering homotopy theorem (see \cite{steenrod}, 11.3).

\begin{lm}\label{lm:lifting of homotopy} Let $f:{\cal Y}\rightarrow {\cal X}$ be a connected finite topological covering. Assume that
${\cal X}$ is paracompact and let
${\cal Z}\subset {\cal X}$ be a topological subspace which is a deformation retract of ${\cal X}$. Then $\tilde
{\cal Z}:=\pi^{-1}({\cal Z})\subset {\cal Y}$ is a deformation retract of ${\cal Y}$. In particular, $\tilde {\cal Z}$ is also 
arcwise connected and
$\pi_1(\tilde {\cal Z})=\pi_1({\cal Y})$.
\end{lm}
For a germ of normal complex analytic space $(X,p)$, let $U$ be a Prill's good
neighbourhood with regard to a proper analytically closed subset $\Sigma\subset X$ and let $\{U_{\lambda}\}_{\lambda\in\Lambda}$ be a neighbourhood basis
associated with $U$. We put $U^{-}:=U\setminus \Sigma$ and $U^{-}_{\lambda}:=U_{\lambda}\setminus \Sigma$. Take any connected finite topological covering $f^{-}:V^{-}\rightarrow
U^{-}$.  Then $V^-$ has the unique analytic structure such
that $f^-:V^-\rightarrow U^{-}$ is \'etale. By the Grauert-Remmert's theorem, $f^-$
extends uniquely to a finite cover
$f:V\rightarrow U$, where $V$ is a normal complex analytic space such that $f^{-1}(U^{-})=V^-$ 
(see \cite{grauert-remmert}, \S 2, Satz 8 and \cite{sga1}, XII,
Theorem 5.4). Recall here that $f^{-1}(p)$ consists of exactly one point, for if $f^{-1}(p)=\{q_1,\dots ,q_n\}$ and $n\geq
2$, where
$q_i$ are distinct from each other, then by \cite{fischer}, 1.10, Lemma 2, there exists $\lambda\in \Lambda$ such that
$f^{-1}(U_{\lambda})=\coprod_{i=1}^n W_{\lambda,i}$, where $W_{\lambda,i}$ is an open neighbourhood of $q_i$ for $i=1,\dots, n$.
Since $f^{-1}(U_{\lambda}^{-})\subset f^{-1}(U_{\lambda})$ is connected by
Lemma~\ref{lm:lifting of homotopy}, there exists $i$, say $i_0$, such that 
$f^{-1}(U_{\lambda}^{-})\subset W_{\lambda,i_0}$ and $f^{-1}(U_{\lambda}^{-})\cap W_{\lambda, i}$ is
empty if $i\neq i_0$, but which is absurd for  
$f^{-1}(U_{\lambda}^{-})\cap W_{\lambda, i}=W_{i,\lambda}\setminus f^{-1}(\Sigma)$ is non-empty for any $i$. Thus we see that any connected finite topological covering
$f^{-}:V^{-}\rightarrow U^{-}$ determines a finite surjective morphism $f:Y:=(V,q)\rightarrow (X,p)$ from a germ of normal complex analytic spaces $Y$ uniquely up to isomorphisms, where
$f^{-1}(p)=\{q\}$. For two connected finite topological coverings $f^{-}_1:V^{-}_1\rightarrow U^{-}$ and $f^{-}_2:V^{-}_2\rightarrow U^{-}$, 
let $f_i:V_i\rightarrow U^{-}$ be the
extended finite covers and $f:Y_i\rightarrow X$ be the corresponding finite surjective morphisms as above for
$i=1$, $2$. By
\cite{sga1}, Expos\'e XII, Proposition 5.3, the restriction map 
$r:\mbox{Hom}_U(V_1,V_2)\rightarrow \mbox{Hom}_{U^{-}}( V^-_1, V^-_2 )$ is
bijective and composed with the canonical injection $\mbox{Hom}_U(V_1,V_2)\rightarrow \mbox{Hom}_X(Y_1,Y_2)$, $r^{-1}$ gives a canonical injection $\mbox{Hom}_{U^{-}}( V^-_1, V^-_2
)\rightarrow  \mbox{Hom}_{X}(Y_1,Y_2)$. From the above argument, we see that we have a canonical faithful functor ${\cal P}$ called a {\it Prill functor} from the category of
connected topological finite coverings over $U^{-}$ denoted by $FT(U^{-})$ to the category of germs of normal complex analytic spaces which is finite over $X$ and
\'etale outside over
$\Sigma$ denoted by 
$FC(X,\Sigma)$.

\begin{lm}\label{lm:prill-fun} A Prill functor defines an equivalence between the categories $FT(U^{-})$ and $FC(X,\Sigma)$.
\end{lm}

{\it Proof. } We note that the restriction functor ${\cal R}_{\lambda}:FT(U^{-})\rightarrow FT(U_{\lambda}^{-})$ defines an equivalence of categories between
$FT(U^{-})$ and
$FT(U_{\lambda}^{-})$ since these are known to be determined up to equivalences by the corresponding fundamental
groups. Put 
$(V_{i,\lambda}^-, f_{i,\lambda}^-):={\cal R}_{\lambda}((V_i^-,f_i^-))\in \mbox{Ob }FT(U_{\lambda})$ and
$V_{i,\lambda}:=f^{-1}_i(U_{\lambda})$ for $i=1,2$. Note also that the canonical map 
$$
\mbox{injlim}_{\lambda\in \Lambda} \mbox{Hom
}_{U_{\lambda}}(V_{1,\lambda},V_{2,\lambda})\rightarrow \mbox{Hom }_{FC(X,\Sigma)}((Y_1,f_1),(Y_2,f_2))
$$ is bijective and 
 that we have 
$$
\mbox{injlim}_{\lambda\in \Lambda} \mbox{Hom }_{FT(U^{-}_{\lambda})}((V_{1,\lambda}^-,f_{1,\lambda}^-),(V_{2,\lambda}^-,f_{2,\lambda}^-))=\mbox{Hom
}_{FT(U^{-})}((V_1^-,f_1^-),(V_2^-,f_2^-)).
$$ Therefore, we conclude that the canonical map 
$$\mbox{Hom }_{FT(U^{-})}( (V^-_1,f^-_1), (V^-_2,f^-_2) )\rightarrow  \mbox{Hom }_{FC(X,\Sigma)}((Y_1,f_1),(Y_2,f_2))
$$ is
bijective, which implies that the functor ${\cal P}$ is faithfully full. Take any $(Y,f)\in \mbox{Ob }FC(X,\Sigma)$. Then $f$ is
represented by a finite cover $f:V_{\lambda}\rightarrow U_{\lambda}$ for some $\lambda \in \Lambda$, where $V_{\lambda}$ is
connected. Since 
$f$ is
\'etale over
$U_{\lambda}^-$ and $V_{\lambda}^-:=f^{-1}(U_{\lambda}^-)$ 
is also connected, we obtain an object 
$(V_{\lambda}^-, f|_{V_{\lambda}^-})\in \mbox{Ob }FT(U^{-}_{\lambda})$ which goes to $(Y,f)\in \mbox{Ob }FC(X,\Sigma)$ via ${\cal P}\circ {\cal
R}_{\lambda}^{-1}$.  Thus we conclude that ${\cal P}$ is essentially surjective, and hence
${\cal P}$ defines an equivalence.
\hfill\bsquare

\begin{rk}\label{rk:prill-fun} {\em It is obvious that a Prill functor also defines an equivalence between the full
subcategory of Galois objects of $FT(U^{-})$ and $FC(X,\Sigma)$.  Note that giving a pointing to an object of $FT(U^{-})$ and $FC(X,\Sigma)$ has essentially the same meaning since the
number of pointings for $(V^-,f^-)\in \mbox{Ob }FT(U^{-})$  and the number of pointings for $(Y,f)\in \mbox{Ob }FC(X,\Sigma)$ are both $\deg f=\deg f^-$.  Thus we see
that ${\cal P}$ induces an equivalence between 
$FT(U^{-})^p$ and $FC(X,\Sigma)^p$. 
}\end{rk}

For {\bf Q}-divisor $D$ on $X$, let ${\cal B}_X(D)$ be the set of all the prime divisors on $X$ such that $e_{\Gamma}(D)>1$ and put $D^{\vee}:=\sum_{\Gamma\in {\cal
B}_X(D)}e_{\Gamma}(D)\Gamma\in
\mbox{Weil }X$. Combined with Lemma~\ref{lm:fgc}, Lemma~\ref{lm:prill-fun} yields the following proposition.

\begin{pr}\label{pr:top-int} There exists a canonical functor ${\cal P}:FC^{D^{\vee}}(\mbox{\em Reg }U)\rightarrow \mbox{\em Int}^{\dag}(X;D)$ which defines an equivalence between the
categories
$FC^{D^{\vee}}(\mbox{\em Reg }U)$ and $\mbox{\em Int}^{\dag}(X;D)$. In particular, $\mbox{\em Int}^{\dag}G(X;D)^p$
is cofilterd and cofinal in $\mbox{\em Int}^{\dag}(X;D)^p$. 
\end{pr}

\begin{rk}\label{rk:top-int}{\em From Proposition~\ref{pr:top-int}, we deduce that $\mbox{projlim}_{(Y,f,i_Y)\in {\rm Ob\ Int}^{\dag}G(X;D)^p}\mbox{Gal
}(Y/X)$ is isomorphic to $\overline\pi_{1,{\rm loc}}^{D^{\vee}}(\mbox{Reg }X)^{\wedge}$ in the category of profinite groups.
}\end{rk}

\begin{df}\label{df:log pi1}{\em  For {\bf Q}-divisor $D$ on $X$, we define a group $\pi_{1,X,p}^{{\rm loc}}[D]$
by 
$
\pi_{1,X,p}^{{\rm loc}}[D]:=\overline\pi_{1,{\rm loc}}^{D^{\vee}}(\mbox{Reg }X),
$ which is called {\it the $D$-local fundamental group} for $((X,p),D)$ and we denote by $\hat\pi_{1,X,p}^{{\rm loc}}[D]$ its profinite completion.
}\end{df}

\begin{rk}\label{rk:integral case}{\em We note that $\pi_{1,X,p}^{{\rm loc}}[D]$ depends only on the class $[D]\in\mbox{Weil }X\otimes {\mbox{\boldmath $Q$}}/\mbox{Weil }X$ and that if
$D\in
\mbox{Weil }X$,
$\pi_{1,X,p}^{{\rm loc}}[D]\simeq \pi_1^{{\rm loc}}(\mbox{Reg
}X)$.   }\end{rk}

\subsection{The category of Cartier covers and Comparison Theorem}
In this section, we introduce a category which is easier to handle with than the category of integral covers. Let
$(X,p)$ be a germ of irreducible normal complex analytic spaces and let ${\cal M}_X$ denote the field of
germs of meromorphic functions on $X$. In what follows, we fix an algebraic
closure $\overline{\cal M}_X$ of ${\cal M}_X$ and the inclusion $i_X:{\cal M}_X\rightarrow \overline{\cal M}_X$. 
Take any $D\in \mbox{Div}_{\mbox{\boldmath $Q$}}X$ and fix it. Recall that a holomorphic map between complex analytic spaces
 is said to be {\it finite}, if it is proper
with discrete finite fibres.

\begin{df}{\em A finite
surjective morphism $f:Y\rightarrow X$, where $Y$ is a germ of irreducible normal complex analytic
spaces such that 
$f^{\ast}D$ is integral and Cartier, is called an {\it Cartier cover
with respect to }$D$. A Cartier cover $f:Y\rightarrow X$ with respect to $D$ is called {\it Cartier Galois cover
with respect to }$D$ if $f$ is Galois.
}
\end{df}

\begin{df}{\em Cartier covers (resp. Cartier Galois covers) of $X$ with respect to $D$ form a full subcategory of complex analytic germs 
denoted by
${\rm Cart}^m(X;D)$ (resp. ${\rm Cart}^mG(X;D)$). For $(Y,f)\in \mbox{Ob }{\rm Cart}^m(X;D)$, an injective homomorphism $i_Y:{\cal M}_Y\rightarrow
\overline{\cal M}_X$, where ${\cal M}_Y$ is the meromorphic
function field of $Y$, such that $i_Y\circ f^{\ast}=i_X$ is called a {\it pointing}. Triplet
$(Y, f, i_Y)$ composed of  $(Y,f)\in \mbox{Ob }{\rm Cart}^m(X;D)$ and a pointing $i_Y$ are called {\it pointed Cartier
covers with respect to }$D$. Pointed Cartier covers (resp. pointed Cartier Galois covers) with respect to $D$ and morphisms $f_{\lambda,\mu}\in \mbox{Hom
}_{{\rm Cart}^m(X;D)}((Y_{\mu},f_{\mu}),(Y_{\lambda},f_{\lambda}))$ satisfying
$i_{Y_{\mu}}\circ f_{\lambda,\mu}^{\ast}=i_{Y_{\lambda}}$ form a projective system denoted by ${\rm Cart}^m(X;D)^p$ (resp. ${\rm Cart}^mG(X;D)^p$). 
}\end{df}

Non-zero $\mbox{\boldmath $C$}$-algebra ${\cal A}$ is called
{\it a complex analytic ring} if there exists a surjective
$\mbox{\boldmath $C$}$-algebra homomorphism ${\cal O}^{\rm an}_{E,0}
\rightarrow {\cal A}$, where $E\simeq\mbox{\boldmath $C$}^n$ for some $n$. Recall that the category of complex analytic rings ${\cal A}$ which are finite ${\cal O}_X$-modules and that
category of germs of complex analytic spaces which are finite over $X$ are dual to each other via the contravariant functor $\mbox{Specan}_X$ (see
\cite{fischer} and \cite{grothendieck}, VI). The
structure morphism $f^{\ast}:{\cal O}_X\rightarrow {\cal A}$ is injective if and only if $f:\mbox{Specan}_X{\cal A}\rightarrow X$
is surjective by the Remmert's proper mapping theorem (see, for example, \cite{fischer}, 1.18). Let $\varphi$ be a meromorphic function on $X$
such that ${\cal O}_X(-rD)=\varphi{\cal O}_X$ and let 
$\pi:\tilde X\rightarrow X$ be the index one cover with respect to $D$ obtained by
taking a
$r$-th root of $\varphi$, where $r:=\mbox{ind}_pD$ and fix a pointing $i_{\tilde X}$. Take any
$(Y,f,i_Y)\in
\mbox{Ob }{\rm Cart}^m(X;D)^p$. For simplicity, assume that 
${\cal M}_X\subset{\cal M}_{\tilde X}\subset\overline{\cal M}_X$ and ${\cal M}_X\subset{\cal M}_Y\subset\overline{\cal M}_X$. The assumption
on
$Y$ implies that there exists a meromorphic function $\psi$ on $Y$ such that
$\varphi{\cal O}_Y=\psi^{r}{\cal O}_Y$, which implies that there
exists a unit
$u\in {\cal O}_Y^{\times}$ such that $\psi^r=u\varphi$. Since ${\cal
O}_Y$ is a henselian local ring whose residue field is the complex number field (see, for example, \cite{abhyankar}, Ch. III, \S 20,
Proposition 20.6), we see that $\root{r}\of{u}\in {\cal O}_Y^{\times}$, hence
${\cal M}_X(\root{r}\of{\varphi})\subset {\cal M}_Y$. Consequently, there exists a ${\cal O}_X$-homomorphism 
$\pi_{\ast}{\cal O}_{\tilde
 X}\rightarrow f_{\ast}{\cal O}_{Y}$ which induces a finite surjective morphism
$\varpi_Y:Y=\mbox{Specan}_Xf_{\ast}{\cal O}_Y\rightarrow \tilde X=\mbox{Specan}_{X}\pi_{\ast}{\cal
O}_{\tilde X}$ satisfying $f=\pi\circ\varpi_Y$. The above argument implies that $\mbox{Card Hom }_{{\rm Cart}^m(X;D)^p}((Y,f,i_Y),(\tilde X,\pi,i_{\tilde X}))=1$ for any $(Y,f,i_Y)\in
\mbox{Ob }{\rm Cart}^m(X;D)^p$, that is, $(\tilde X,\pi,i_{\tilde X})$ is a final object, or equivalently, a colimit of ${\rm Cart}^m(X;D)^p$. Let $\varpi_Y(i_Y,i_{\tilde X})$ denote the element of $\mbox{Hom }_{{\rm Cart}^m(X;D)^p}((Y,f,i_Y),(\tilde X,\pi,i_{\tilde X}))$.

\begin{df}{\em A Cartier cover $f:Y\rightarrow X$ with respect to $D$ is called a {\it strict Cartier cover with respect to} $D$, if $\varpi_Y(i_Y,i_{\tilde X})$ is \'etale in codimension one for any pointings $i_Y$, $i_{\tilde X}$.
}\end{df}

\begin{rk}\label{rk:et1} {\em $\varpi_Y(i_Y,i_{\tilde X})$ is \'etale in codimension one if and only if $\varpi_Y(i_Y,i_{\tilde X})$ is \'etale over $\mbox{Reg }X$ by the
purity of branch loci (see 
\cite{abhyankar}, V, \S 39, (39.8) or \cite{fischer}, 4.2).
}\end{rk}

\begin{df}{\em A strict Cartier cover $f:Y\rightarrow X$ with respect to $D$ is called a {\it strict Cartier Galois cover} with respect to $D$, if $f$ is Galois.
}\end{df}

\begin{rk}{\em Take another pointings $i_Y^{\prime}$ and $i_{\tilde X}^{\prime}$ of
$(Y,f)$,
$(\tilde X,\pi)\in \mbox{Ob }{\rm Cart}^m(X;D)^p$ respectively and assume that $f:Y\rightarrow X$ is Galois. Then, by the Galois theory, 
there exist two 
isomorphisms $\alpha(i_Y,i_Y^{\prime})\in \mbox{Hom }_{{\rm Cart}^m(X;D)^p}((Y,f,i_Y), (Y,f,i_Y^{\prime}))$ and
$\beta(i_{\tilde X},i_{\tilde X}^{\prime})\in \mbox{Hom }_{{\rm Cart}^m(X;D)^p}((\tilde X,\pi,i_{\tilde X}),
(\tilde X,\pi,i_{\tilde X}^{\prime}))$ such that the following diagram in ${\rm Cart}^m(X;D)^p$ commutes.

\[\begin{array}{ccc}
(Y, f,i_Y)&\mapright{\alpha(i_Y,i_Y^{\prime})}&(Y, f,i_Y^{\prime})\\
\mapdown{\varpi_Y(i_Y,i_{\tilde X})}&&\mapdownr{\varpi_Y(i_Y^{\prime},i_{\tilde X}^{\prime})} \\
(\tilde X,\pi,i_{\tilde X})&\mapright{\beta(i_{\tilde X},i_{\tilde X}^{\prime})} & (\tilde X,\pi,i_{\tilde X}^{\prime})
\end{array}\]

Therefore, $(Y,f)\in \mbox{Ob }{\rm Cart}^m(X;D)^p$ is a strict Cartier Galois cover if one of $\varpi_Y(i_Y,i_{\tilde X})\in \mbox{Hom }_{{\rm Cart}^m(X;D)}((Y,f),(\tilde X,\pi))$ is \'etale in codimension one. One
can also check easily that the same holds even if $f$ is not Galois.
}
\end{rk}

Let ${\rm Cart}^{\dag}(X;D)$ (resp, ${\rm Cart}^{\dag}(X;D)^p$) denote the full subcategory of ${\rm Cart}^m(X;D)$ (resp.
projective subsystem ${\rm Cart}^m(X;D)^p$) whose objects are strict Cartier covers with respect to
$D$ (resp. pointed strict Cartier covers with respect to
$D$) and let ${\rm Cart}^{\dag}G(X;D)$ (resp, ${\rm Cart}^{\dag}G(X;D)^p$) denote the full subcategory of ${\rm Cart}^{\dag}(X;D)$ (resp.
projective subsystem ${\rm Cart}^{\dag}(X;D)^p$) whose objects are strict Cartier Galois covers with respect to
$D$ (resp. pointed strict Cartier Galois covers with respect to
$D$).  Let $f_{\lambda,\mu}:(Y_{\mu},f_{\mu}, i_{Y_{\mu}})\rightarrow (Y_{\lambda},f_{\lambda},i_{Y_{\lambda}})$ be a morphism
in
${\rm Cart}^{\dag}G(X;D)^p$ and assume that ${\cal M}_X\subset{\cal M}_{Y_{\lambda}}\subset {\cal M}_{Y_{\mu}}\subset \bar{\cal
M}_X$ for simplicity. Then by the Galois theory, there
exists a canonical surjective homomorphism
$g_{\lambda,\mu}:\mbox{Gal }({\cal M}_{Y_{\mu}}/{\cal M}_X)\rightarrow \mbox{Gal }({\cal M}_{Y_{\lambda}}/{\cal M}_X)$ which is
nothing but the restriction map. Thus Galois groups $\mbox{Gal }(Y/X):=\mbox{Gal }({\cal M}_Y/{\cal M}_X)$ for $(Y,f,i_Y)\in {\rm Cart}^{\dag}G(X;D)^p$ 
form a projective system with the induced morphisms from ${\rm Cart}^{\dag}G(X;D)^p$.  Here is our comparison theorem.

\begin{th}\label{th:comparison} There exists a canonical isomorphism $:$
$$
\hat\pi_{1,X,p}^{{\rm loc}}[D]\simeq\mbox{\rm projlim}_{(Y,f,i_Y)\in {\rm Ob\ }{\rm Cart}^{\dag}G(X;D)^p}\mbox{\rm Gal
}(Y/X)
$$ in the category of profinite groups for any $D\in \mbox{\em Div}_{\mbox{\boldmath $Q$}}X$.
\end{th}
To prove the above theorem, we need some lemmas and propositions as follows. Let ${\cal C}$ be a category and let 
$X\in \mbox{Ob }{\cal C}$ be an object of ${\cal C}$ and $G\subset \mbox{Aut }X$ be
a subgroup of the automorphism group of $X$.
\begin{df}{\em An epimorphism $f:X\rightarrow Y$ in ${\cal C}$ is said to be {\it Galois with the Galois group} $G$,
if
$G=\mbox{Aut }_Y X:=\{\sigma\in \mbox{Aut }X|f\circ\sigma=f \}$ and for any morphism $f^{\prime}:X\rightarrow Y^{\prime}$ such
that
$G\subset \mbox{Aut }_{Y^{\prime}} X$, there exists a unique morphism $\varphi:Y\rightarrow Y^{\prime}$ satisfying
$f^{\prime}=\varphi\circ f$.
}
\end{df}

\begin{rk}{\em Assume that two Galois morphisms $f:X\rightarrow Y$ and $f^{\prime}:X\rightarrow Y^{\prime}$ with the Galois group
$G$ are given. Then by the universal mapping property, there exists an isomorphism $\varphi:Y\rightarrow Y^{\prime}$ such that
$f^{\prime}=\varphi\circ f$, that is, Galois morphisms with the Galois group
$G$ is unique up to this equivalence.
}\end{rk}

\begin{ex}{\em Let ${\cal F}:=(\mbox{Fields})$ be a category of fields such that $\mbox{Hom}_{\cal F}(K_1,K_2)$ is empty or
consists of inclusions for any $K_1$, $K_2\in \mbox{Ob }{\cal F}$. For any finite extension $i:K_1\rightarrow K_2$, 
$i$ is a Galois extension if and only if its dual $i^{\circ}:K_2^{\circ}\rightarrow K_1^{\circ}$ in the dual category ${\cal
F}^{\circ}$ is Galois by the Galois theory.
}\end{ex}

\begin{df}{\em For any two morphisms $f\in \mbox{Hom}_{{\cal C}}(X,Y)$ and $g\in \mbox{Hom}_{{\cal C}}(Y,Z)$, we define a subgroup
$\mbox{Aut}^f_Z X\subset\mbox{Aut }X\times \mbox{Aut}_Z Y$ as $\mbox{Aut}^f_Z X:=\{(\tilde \sigma,\sigma)|f\circ\tilde
\sigma=\sigma\circ f\}$.
}\end{df}

\begin{lm}\label{lm:galois} Let $f\in \mbox{\em Hom}_{{\cal C}}(X,Y)$ and $g\in \mbox{\em Hom}_{{\cal C}}(Y,Z)$ be two Galois
morphisms and assume that the second projection $p_2:\mbox{\em Aut}^f_Z X\rightarrow \mbox{\em Aut}_Z Y$ is surjective. Then
$h:=g\circ f$ is also Galois.
\end{lm}

{\it Proof.} Take any $h^{\prime}\in \mbox{Hom }_{{\cal C}}(X,Z^{\prime})$ with $\mbox{Aut}_Z X\subset \mbox{Aut}_{Z^{\prime}}X$.
Since $\mbox{Aut}_Y X\subset \mbox{Aut}_Z X$ and $f$ is Galois, there exists a morphism $\psi:Y\rightarrow Z^{\prime}$ such that
$\psi\circ f=h^{\prime}$. Take any $\sigma\in \mbox{Aut}_Z Y$. Then there exists $\tilde \sigma\in \mbox{Aut }X$ such that
$f\circ\tilde\sigma=\sigma\circ f$ by the assumption. Since $\tilde\sigma\in \mbox{Aut}_Z X\subset \mbox{Aut}_{Z^{\prime}}X$, we
have
$h^{\prime}=h^{\prime}\circ\tilde\sigma=\psi\circ f\circ \tilde\sigma=\psi\circ\sigma \circ f$, hence 
$\psi\circ f=\psi\circ\sigma
\circ f$. Since $f$ is an epimorphism, we deduce that $\psi=\psi\circ \sigma$, that is, $\mbox{Aut}_Z Y\subset
\mbox{Aut}_{Z^{\prime}}Y$. Thus we conclude that there exists a morphism $\varphi:Z\rightarrow Z^{\prime}$ such that $\varphi\circ
g=\psi$.  Obviously $\varphi$ satisfies $\varphi\circ h=h^{\prime}$. As for the uniqueness of $\varphi$, Let
$\varphi^{\prime}:Z\rightarrow Z^{\prime}$ be another morphism satisfying $\varphi^{\prime}\circ h=h^{\prime}$. Then $\psi\circ
f=h^{\prime}=\varphi^{\prime}\circ h=\varphi^{\prime}\circ g\circ f$, hence $\varphi^{\prime}\circ g=\psi=\varphi\circ g$. Since $g$
is also an epimorphism, we obtain $\varphi^{\prime}=\varphi$.
\hfill \bsquare

\begin{rk}\label{rk:lift}{\em The assumption in Lemma~\ref{lm:galois} is satisfied in the following two theoretically important cases.

\par
(1) Let $f$ and $g$ are finite Galois covers between normal algebraic varieties over an algebraically closed field or normal connected complex
analytic spaces. 
Assume that there exits a Zariski closed subset or an analytic subset
$\Sigma$ on
$Y$ with
$\mbox{codim}_Y \Sigma\leq 2$ such that the restriction $f^-$ of $f$ to $X^-:=X\setminus
f^{-1}(\Sigma)$ gives the algebraic universal cover of $Y^-:=Y\setminus \Sigma$, that is,
$\hat\pi_1(X^-)=\{1\}$. Moreover assume that $Y^-$ is invariant under the action of
$\mbox{Gal }(Y/Z)$. Take
any
$\sigma\in
\mbox{Gal }(Y/Z)$. Since $\sigma$ acts on
$Y^-$, there exists an automorphism
$\tilde\sigma^-$ on
$X^-$ such that $f^-\circ\tilde \sigma^-=\sigma\circ f^-$ by the property of algebraic universal
cover.
$\tilde\sigma^-$ extends uniquely to an automorphism $\tilde\sigma$ on
$X$ satisfying $f\circ\tilde \sigma=\sigma\circ f$ by the normality (see also \cite{catanese}, \S1 and \cite{gabi}, Lemma
2.1). 
\par

(2) Let $f$ and $g$ are finite Galois covers between germs of normal complex
analytic spaces. Assume that $X$ is obtained from
$Y$ by taking a $r$-th root of a primitive principal divisor $P=\mbox{div }\varphi$ on $Y$ such that ${\cal 
O}_Y(P)\subset{\cal M}_Y$ is invariant under the action of $\mbox{Gal }(Y/Z)$ (for the definition
of primitive principal divisors, see \cite{shokurov}, 2.3). Take any $\sigma\in \mbox{Gal }(Y/Z)$. Then by the assumption,
we have $\sigma ^{\ast}\varphi=u\varphi$ for some unit $u\in {\cal O}_Y^{\times}$. 
As in the previous argument, there exists a unit $v\in {\cal O}_Y^{\times}$ such that $v^r=u$. Since we can write ${\cal M}_X={\cal M}_Y[T]/(T^r-\varphi)$, it is obvious that 
$\sigma
^{\ast}$ lifts to an automorphism $\tilde\sigma^{\ast}$ on ${\cal M}_X$ by putting $\tilde\sigma^{\ast}T=vT$.  Thus we see that any elements of $\mbox{Gal }(Y/Z)$ lift to 
elements of $\mbox{Gal }(X/Z)$.  
}\end{rk}

\begin{lm}[c.f., \cite{seidenberg2}]\label{lm:normalization} Let ${\cal A}$ be an integral complex analytic ring and ${\cal M}$ be its quotient field. Let
${\cal A}_{{\cal L}}$ be the normalization of ${\cal A}$ in a finite extension field ${\cal L}$ of ${\cal M}$. Then ${\cal
A}_{{\cal L}}$ is also an integral complex analytic ring which is a finite ${\cal A}$-module. 
\end{lm}
{\it Proof.} Recall that ${\cal A}$ is noetherian (\cite{grothendieck}, II, Proposition 2.3). \cite{seidenberg2}, Theorem 4
says that
${\cal A}$ is N-1, hence N-2 by \cite{matsumura}, Ch. 12, Corollary 1, that is, ${\cal
A}_{{\cal L}}$ is a finite ${\cal A}$-module. By \cite{seidenberg2}, Theorem 1, 
${\cal A}$ is a finite ${\cal O}^{\rm an}_{C^n,0}$-module for some $\mbox{\boldmath $C$}^n$, hence so is ${\cal
A}_{{\cal L}}$. Thus by \cite{seidenberg2}, Theorem 3, we conclude that ${\cal A}_{{\cal L}}$ is a complex analytic ring.
\hfill\bsquare

\begin{rk}\label{rk:normalization} {\em Lemma~\ref{lm:normalization} implies that if we are given a finite extension field ${\cal L}$ of the
meromorphic function field ${\cal M}_X$ of an irreducible germ of complex analytic spaces $X$, there exists a germ of normal
complex analytic spaces
$Y$ with a finite surjective morphism
$f:Y\rightarrow X$ such that ${\cal M}_Y$ is isomorphic to ${\cal L}$ over ${\cal M}_X$ and such $Y$ as above is uniquely determined up to isomorphisms over $X$.
}\end{rk}

The following proposition can be also derived from Proposition~\ref{pr:top-int}, but we shall give an algebraic proof for further research such as extending our theory to the positive
characteristic case.

\begin{pr}\label{pr:cofilterd} ${\rm Cart}^{\dag}G(X;D)^p$ is cofilterd and is cofinal in ${\rm Cart}^{\dag}(X;D)^p$ for any $D\in \mbox{\em Div}_{\mbox{\boldmath $Q$}}X$.
\end{pr}

{\it Proof. }Firstly, we prove the first statement. Take any two objects $(Y_{\lambda},f_{\lambda},i_{Y_{\lambda}})$, $(Y_{\mu},f_{\mu},i_{Y_{\mu}})\in {\rm Cart}^{\dag}G(X;D)^p$. Let
${\cal L}:=i_{Y_{\lambda}}({\cal M}_{Y_{\lambda}})\lor i_{Y_{\mu}}({\cal M}_{Y_{\mu}})$ be the minimal subfield of $\overline{\cal M}_X$ containing $i_{Y_{\lambda}}({\cal
M}_{Y_{\lambda}})$ and
$i_{Y_{\mu}}({\cal M}_{Y_{\mu}})$. We note that
${\cal L}$ is a finite Galois extension of
${\cal M}_X$ by its definition. Let $g:Z\rightarrow X$ be the normalization of $X$ in ${\cal L}$ as explained in Remark~\ref{rk:normalization}. By the construction, we get an
object $(Z, g, i_Z)\in \mbox{Ob }{\rm Cart}^m(X;D)^p$ dominating both of $(Y_{\lambda},f_{\lambda},i_{Y_{\lambda}})$ and $(Y_{\mu},f_{\mu},i_{Y_{\mu}})$ in ${\rm Cart}^m(X;D)^p$. Let
$(\tilde X,\pi,i_{\tilde X})\in \mbox{Ob }{\rm Cart}^{\dag}G(X;D)^p$ be a final object of ${\rm Cart}^{\dag}G(X;D)^p$. From the equality $i_{\tilde X}=i_Y\circ \varpi_Y(i_Y,i_{\tilde
X})^{\ast}$, we see that there exists a canonical embedding: 
\begin{eqnarray}\label{eq:embed}
\Phi_{i_{\tilde X}}:{\rm Cart}^{\dag}G(X;D)^p\rightarrow {\rm Cart}^{\dag}G(\tilde X;0)^p, 
\end{eqnarray} 
depending on the choice of pointings $i_{\tilde X}$ for $(\tilde X,\pi)\in \mbox{Ob }{\rm Cart}^{\dag}G(X;D)$, such that $\Phi_{i_{\tilde X}}((Y,f,i_Y))=(Y,\varpi_Y(i_Y,i_{\tilde X}),
i_Y)\in \mbox{Ob }{\rm Cart}^{\dag}G(\tilde X;0)^p$ for $(Y,f,i_Y)\in \mbox{Ob }{\rm Cart}^{\dag}G(X;D)^p$.
Since ${\rm Cart}^{\dag}G(\tilde X;0)^p$ is cofilterd as explained in Remark~\ref{rk:prill-fun}, there exists $(W,h,i_W)\in \mbox{Ob }{\rm Cart}^{\dag}G(\tilde X;0)^p$ which dominates
both of $(Y_{\lambda},\varpi_{Y_{\lambda}}(i_{Y_{\lambda}},i_{\tilde X}), i_{Y_{\lambda}})$ and $(Y_{\mu},\varpi_{Y_{\mu}}(i_{Y_{\mu}},i_{\tilde X}), i_{Y_{\mu}})$ in ${\rm Cart}^{\dag}G(\tilde X;0)^p$. By the construction of
$(Z,g,i_Z)\in \mbox{Ob }{\rm Cart}^m(X;D)^p$, $i_Z$ factors into $i_W\circ \tau^{\ast}$, where $\tau^{\ast}:{\cal M }_Z\rightarrow {\cal M}_W$ is an injective homomorphism.
Let
$\tau$ be the induced morphism $\tau:W\rightarrow Z$. Then we see that 
$h$ factors into $\varpi_Z(i_Z,i_{\tilde X})\circ \tau$. Since $h$ is \'etale in codimension one, hence so is $\varpi_Z(i_Z,i_{\tilde X})$. Thus we conclude that $(Z, g,
i_Z)\in \mbox{Ob }{\rm Cart}^{\dag}G(X;D)^p$ and consequently, ${\rm Cart}^{\dag}G(X;D)^p$ is cofilterd. As for second statement, take any $(Y,f,i_Y)\in \mbox{Ob }{\rm Cart}^{\dag}(X;D)^p$ and let $\{i_Y^{(k)}|k=1,2,\dots n\}$ be all the pointings for 
$(Y,f)\in \mbox{Ob }{\rm Cart}^{\dag}(X;D)$. Let 
${\cal L}$ be the minimal subfield of $\overline{\cal M}_X$ containing all the subfields $i_Y^{(1)}({\cal M}_Y)$, \dots, $i_Y^{(1)}({\cal M}_Y)$. Since ${\cal L}$ is a finite Galois
extension of ${\cal M}_X$ by its construction, we have an object $(Z,g,i_Z)\in \mbox{Ob }{\rm Cart}^mG(X;D)^p$ dominating all the $(Y,f,i_Y^{(1)})$, \dots, $(Y,f,i_Y^{(n)})\in \mbox{Ob
}{\rm Cart}^{\dag}(X;D)^p$, where $g:Z\rightarrow X$ is the normalization of $X$ in
${\cal L}$. In the same way as in the previous argument, we conclude that
$(Z,g,i_Z)\in \mbox{Ob }{\rm Cart}^{\dag}G(X;D)^p$.
\hfill\bsquare
\vskip 5mm

{\it Proof of Theorem~\ref{th:comparison}}. By Proposition~\ref{pr:top-int} and Remark~\ref{rk:top-int}, we only have to show that the full subcategory ${\rm Cart}^{\dag}G(X;D)^p$ is
cofinal in
$\mbox{Int}^{\dag}G(X;D)^p$  (see \cite{sga4}, Expos\'e I, Proposition 8.1.3 or \cite{artin-mazur}, Appendix, Corollary (2.5)). Choose any object $(Y,f,i_Y)\in \mbox{Ob
}\mbox{Int}^{\dag}G(X;D)^p$ and let
$\pi_Y:\tilde Y\rightarrow Y$ be the index one cover with respect to $f^{\ast}D$. We can choose a pointing
$i_{\tilde Y}$ so that a triple $(\tilde Y,\tilde f,i_{\tilde Y})$ becomes an object in $\mbox{Int}^{\dag}(X;D)^p$ dominating $(Y,f,i_Y)$. From Remark~\ref{rk:lift}, we deduce that
$(\tilde Y,\tilde f,i_{\tilde Y})\in \mbox{Ob }{\rm Cart}^{\dag}G(X;D)^p$ by its construction.
\hfill\bsquare

\subsection{Universal Cartier covers}
Let $U$ be a Prill's good neighbourhood with regard to $\mbox{Sing }X$ and $\{U_{\lambda}\}_{\lambda\in \Lambda}$ its associated neighbourhood basis. Take any $(Y,f, i_Y)\in \mbox{Ob
}{\rm Cart}^{\dag}(X;0)^p$ and put
$(V^{-}, f^-, y):={\cal P}^{-1}(Y,f, i_Y)\in FT(U^{-})$, where ${\cal P}$ is a Prill functor. Let $f:V\rightarrow U$ be the extended finite cover of $f^-$.
By Lemma~\ref{lm:lifting of homotopy}, $V$ is a Prill's good neighbourhood with regard to $f^{-1}(\mbox{Sing }X)$ with $\{V_{\lambda}\}_{\lambda\in \Lambda}$ being its associated
neighbourhood basis. Thus we have $\pi_1(V^-)=\mbox{projlim}_{q\in {\cal V};{\rm open}}\pi_1({\cal V}\setminus f^{-1}(\mbox{Sing }X))=\mbox{projlim}_{q\in {\cal V};{\rm
open}}\pi_1(\mbox{Reg }{\cal V})=\pi_1^{{\rm loc}}(\mbox{Reg }Y)$ since $\mbox{Reg }{\cal V}\cap f^{-1}(\mbox{Sing }X)$ is a closed analytic subspace of codimension at least two in
$\mbox{Reg }{\cal V}$. (see, for example, 
\cite{prill}, III, Corollary 2.). In particular, we see that $(Y,f, i_Y)\in \mbox{Ob
}{\rm Cart}^{\dag}(X;0)^p$ is an initial object of ${\rm Cart}^{\dag}(X;0)^p$ if and only if $\hat\pi_1^{{\rm loc}}(\mbox{Reg }Y)=\{1\}$.

\begin{df}\label{df:u-index one cover}{\em For $D\in \mbox{Div}_{\mbox{\boldmath $Q$}}X$, a strict Cartier Galois cover
$\pi^{\dag}:X^{\dag}\rightarrow X$ with respect to $D$ is called an {\it algebraic 
universal Cartier cover}, or abbreviated, a {\it universal Cartier cover} with respect to 
$D$ if $\hat\pi_1^{{\rm loc}}(\mbox{Reg }X^{\dag})=\{1\}$.
}\end{df}

\begin{rk}{\em Singularity with trivial local algebraic fundamental group is quite
restrictive one. For example, $\hat\pi_1^{{\rm loc}}(\mbox{Reg}X)=\{1\}$ implies $\mbox{Div}_{\mbox{\boldmath $Q$}}X\cap \mbox{Weil }X=\mbox{Div }X$. Moreover, if we assume, in
addition, that $(X,p)$ is analytically $\mbox{\boldmath
$Q$}$-factorial, then ${\cal O}_X$ is factorial (see, for example, \cite{brieskorn}, Satz 1.4).}
\end{rk}

\begin{pr}\label{pr:fcc}  For $D\in \mbox{\em Div}_{\mbox{\boldmath $Q$}}X$, take the index one cover $\pi:\tilde X\rightarrow X$ with respect to $D$, the there exists the universal
Cartier cover of
$X$ with respect to
$D$ if and only if
$\hat\pi_1^{{\rm loc}}(\mbox{\em Reg }\tilde X)$ is finite.
\end{pr}

{\it Proof.} Assume that $\hat\pi_1^{{\rm loc}}(\mbox{Reg }\tilde X)$ is finite and take a final object $(\tilde X, \pi,i_{\tilde X})\in \mbox{Ob }{\rm Cart}^{\dag}G(X;D)^p$.
By the assumption, ${\rm Cart}^{\dag}(\tilde X;0)^p$ has an initial object
$(Y,f,i_Y)\in
\mbox{Ob }{\rm Cart}^{\dag}G(\tilde X;0)^p$ such that 
$\hat\pi_1^{{\rm loc}}(\mbox{Reg }Y)=\{1\}$. We note that $i_Y$ is also a pointing for $(Y,\pi\circ f)\in \mbox{Ob }{\rm Cart}^m(X;D)$ since we have $i_Y\circ
f^{\ast}\circ\pi^{\ast}=i_{\tilde X}\circ \pi^{\ast}=i_X$. Consider an object
$(Y,\pi\circ f,i_Y)\in
\mbox{Ob }{\rm Cart}^m(X;D)^p$. By the argument in Remark~\ref{rk:lift}, (1), we see that
$\pi\circ f$ is Galois. Since $\varpi_Y(i_Y,i_{\tilde X})=f$ is \'etale in codimension one, we conclude that $(Y,\pi\circ f,i_Y)\in \mbox{Ob }{\rm Cart}^{\dag}G(X;D)^p$. Conversely,
assume that there exists a pointed universal Cartier cover $(X^{\dag},\pi^{\dag},i_{X^{\dag}})\in \mbox{Ob }{\rm Cart}^{\dag}G(X;D)^p$ with respect to $D$. Then 
$\Phi_{i_{\tilde X}}((X^{\dag},\pi^{\dag},i_{X^{\dag}}))\in \mbox{Ob }{\rm Cart}^{\dag}G(\tilde X;0)^p$ is an initial object of $\mbox{Ob }{\rm Cart}^{\dag}(\tilde X;0)^p$, hence
$\hat\pi_1^{{\rm loc}}(\mbox{Reg }\tilde X)$ is finite.
\hfill\bsquare

\begin{rk}\label{rk:s-w}{\em Assume that $(X,\Delta)$ is purely log terminal, where $\Delta$ is a standard {\bf Q}-boundary. Then $(\tilde X,\Delta_{\tilde X})$ is known to be
canonical, hence $\tilde X$ has only canonical singularity if we assume that $\lfloor\Delta\rfloor=0$ or
$\tilde X$ is {\bf Q}-Gorenstein. Thus if $\mbox{dim } X\leq 3$, then $\hat\pi_1^{{\rm loc}}(\mbox{Reg }\tilde X)$ is finite by
\cite{shepherd-wilson}, Theorem 3.6. }
\end{rk}

\begin{pr}[c.f., \cite{sga1}, Expos\'e IX, Remark 5.8]\label{pr:fund.ext seq} For $D\in \mbox{\em Div}_{\mbox{\boldmath $Q$}}X$, Let $\pi:\tilde X\rightarrow X$ be the index
one cover with respect to $D$. Then there exists the following exact sequence in the category of profinite groups $:$
\begin{eqnarray}\label{eq:g-extseq}
\{1\}\longrightarrow\hat\pi_1^{{\rm loc}}(\mbox{\em Reg }\tilde X)\longrightarrow \hat\pi_{1,X,p}^{{\rm loc}}[D]\longrightarrow \mbox{\em Gal }(\tilde X/X)\simeq
\mbox{\boldmath $Z$}/r\mbox{\boldmath $Z$}\longrightarrow
\{1\},
\end{eqnarray}
where $r:=\mbox{\em ind}_pD$.
\end{pr}

{\it Proof. } Recall that we have a canonical embedding
$\Phi_{i_{\tilde X}}:{\rm Cart}^{\dag}(X;D)^p\rightarrow {\rm Cart}^{\dag}(\tilde X;0)^p$ as in (\ref{eq:embed}). Since we have the exact sequence:
$$
\{1\}\longrightarrow\mbox{projlim}_{(Y,f,i_Y)\in {\rm Ob }{\rm Cart}^{\dag}G(X;D)^p}\mbox{Gal }(Y/\tilde X)\longrightarrow \hat\pi_{1,X,p}^{{\rm loc}}[D]\longrightarrow
\mbox{Gal }(\tilde X/X)\longrightarrow
\{1\},
$$
we only have to show that ${\rm Cart}^{\dag}G(X;D)^p$ is cofinal in ${\rm Cart}^{\dag}G(\tilde X;0)^p$ via the functor $\Phi_{i_{\tilde X}}$. Choose any object $(Y,f,i_Y)\in \mbox{Ob
}{\rm Cart}^{\dag}G(\tilde X;0)^p$. Then we see that $(Y,\pi\circ f,i_Y)\in \mbox{Ob }\mbox{Int}^{\dag}(X;D)^p$ since $\pi^{-1}(\mbox{Reg }X\setminus\mbox{Supp }B)\subset \mbox{Reg
}\tilde X$ and $\pi\circ f$ is \'etale over $\mbox{Reg }X\setminus\mbox{Supp }B$. By Proposition~\ref{pr:top-int}, There exists an object $(Z,g,i_Z)\in \mbox{Ob
}\mbox{Int}^{\dag}G(X;D)^p$ dominating the object $(Y,\pi\circ f,i_Y)$.  Since ${\rm Cart}^{\dag}G(X;D)^p$ is cofinal in $\mbox{Int}^{\dag}G(X;D)^p$ (see the proof of
Theorem~\ref{th:comparison}), $(Z,g,i_Z)$ is dominated by some object in ${\rm Cart}^{\dag}G(X;D)^p$. Thus we get the assertion.
\hfill\bsquare

\begin{co}\label{co:cor to fund.ext seq} A pointed universal Cartier cover $(X^{\dag},\pi^{\dag},i_{X^{\dag}})\in \mbox{\em Ob }{\rm Cart}^{\dag}G(X;D)^p$
 is an initial object, or equivalently, a limit
of ${\rm Cart}^{\dag}(X;D)^p$ and vice versa. In particular, a universal Cartier cover with respect to 
$D$ is unique up to isomorphisms over $X$ if it exists. 
\end{co}

{\it Proof. }Let $(\tilde X,\pi,i_{\tilde X})\in \mbox{Ob }{\rm Cart}^{\dag}G(X;D)^p$ be a final object of ${\rm Cart}^{\dag}(X;D)^p$. As we noted firstly in this
section, $(X^{\dag},\varpi_{X^{\dag}}(i_{X^{\dag}},i_{\tilde X}),i_{X^{\dag}})\in \mbox{Ob }{\rm Cart}^{\dag}(\tilde X;0)^p$ is an initial object of ${\rm Cart}^{\dag}(\tilde X;0)^p$,
 hence 
$(X^{\dag},\pi^{\dag},i_{X^{\dag}})\in \mbox{Ob }{\rm Cart}^{\dag}G(X;D)^p$ is also an initial object of
${\rm Cart}^{\dag}(X;D)^p$. On the contrary, assume that there exists an initial object $(X^{\dag\prime},\pi^{\dag\prime},i_{X^{\dag\prime}})$ of ${\rm Cart}^{\dag}(X;D)^p$. By
Proposition\ref{pr:cofilterd}, we see that $(X^{\dag\prime},\pi^{\dag\prime},i_{X^{\dag\prime}})\in \mbox{Ob }{\rm Cart}^{\dag}G(X;D)^p$. Since $\hat\pi_{1,X,p}^{{\rm loc}}[D]$ is
finite, $\hat\pi_1^{{\rm loc}}(\mbox{Reg }\tilde X)$ is also finite by Proposition~\ref{pr:fund.ext seq}. Therefore there exists a pointed universal Cartier cover
$(X^{\dag},\pi^{\dag},i_{X^{\dag}})\in \mbox{Ob }{\rm Cart}^{\dag}G(X;D)^p$ by Proposition~\ref{pr:fcc} which is also an initial object of ${\rm Cart}^{\dag}(X;D)^p$ and hence
isomorphic to $(X^{\dag\prime},\pi^{\dag\prime},i_{X^{\dag\prime}})$. Thus we conclude that $\hat\pi_1^{{\rm loc}}(\mbox{Reg }X^{\dag\prime})=\{1\}$.
\hfill\bsquare
\vskip 5mm
The following Lemma is an algebraic generalization of Brieskorn's fundamental lemma.

\begin{lm}[c.f., \cite{brieskorn}, Lemma 2.6]\label{lm:brieskorn's Lemma}
Let $f:(X,p)\rightarrow (Y,q)$ be a finite morphism between germs of normal complex analytic spaces $(X,p)$ and $(Y,q)$. Then for any $D\in \mbox{\rm Div}_{\mbox{\boldmath $Q$}}Y$,
there exists a canonical homomorphism $f_{\ast}:\hat\pi_{1,X,p}^{{\rm loc}}[f^{\ast}D]\rightarrow \hat\pi_{1,Y,q}^{{\rm loc}}[D]$ which satisfies 
$|\hat\pi_{1,Y,q}^{{\rm loc}}[D]:\mbox{\rm Im }f_{\ast}|\leq \deg f$. In particular, if $\hat\pi_{1,X,p}^{{\rm loc}}[f^{\ast}D]$ is finite, so is 
$\hat\pi_{1,Y,q}^{{\rm loc}}[D]$.
\end{lm}

{\it Proof. } For a given pointing $i_Y:{\cal M}_Y\rightarrow \overline{\cal M}_Y$, we choose a pointing $i_X:{\cal M}_X\rightarrow \overline{\cal M}_X=\overline{\cal M}_Y$ such that
$i_X\circ f^{\ast}=i_Y$. Take any $(Z,\alpha,i_Z)\in \mbox{Ob }{\rm Cart}^{\dag}G(Y;D)^p$. Let $\nu:W\rightarrow Y$ be the normalization of $Y$ in 
$i_X({\cal M}_X)\lor i_Z({\cal M}_Z)$. We note that there exist morphisms $\beta:W\rightarrow X$ and $\gamma:W\rightarrow Z$ such that $\alpha\circ\gamma=f\circ\beta=\nu$. Since
$\beta^{\ast}f^{\ast}D=\gamma^{\ast}\alpha^{\ast}D\in \mbox{Div }W$ and $\beta:W\rightarrow X$ is Galois (see, for example, \cite{nagata}, Theorem 3.6.3), we have $(W, \beta, i_W)\in
\mbox{Ob }{\rm Cart}^{m}G(X;f^{\ast}D)^p$ for a suitable pointing
$i_W$.  Let $(\tilde X,\pi_X, i_{\tilde X})$ ( resp. $(\tilde Y,\pi_Y, i_{\tilde Y})$ ) be a final object of ${\rm Cart}^{m}G(X;f^{\ast}D)^p$ 
( resp. ${\rm Cart}^{m}G(Y;D)^p$ ). Since $(\tilde X, f\circ \pi_X, i_{\tilde X})\in \mbox{Ob }{\rm Cart}^{m}(Y;D)^p$, there exists a morphism $\omega_{\tilde X}(i_{\tilde
X},i_{\tilde Y}):(\tilde X, f\circ \pi_X, i_{\tilde X})\rightarrow (\tilde Y,\pi_Y, i_{\tilde Y})$ in ${\rm Cart}^{m}(Y;D)^p$. Let $\nu^{\natural}:Y^{\natural}\rightarrow Y$ be the
normalization of $Y$ in $i_Z({\cal M}_Z)\cap i_{\tilde X}({\cal M}_{\tilde X})$. Since the induced finite morphism $\delta:Z\rightarrow Y^{\natural}$ is Galois, $i_Z({\cal M}_Z)$
and
$i_{\tilde X}({\cal M}_{\tilde X})$ are linearly disjoint over $i_{Y^{\natural}}({\cal M}_{Y^{\natural}})$, that is, $i_Z({\cal M}_Z)\otimes _{i_{Y^{\natural}}({\cal
M}_{Y^{\natural}})}i_{\tilde X}({\cal M}_{\tilde X})\simeq i_W({\cal M}_W)$ (see, for example, \cite{nagata}, Exercise 4.2.3). Let
$\eta\in W$ be the generic point of a prime divisor on $W$ and $\xi\in Z$ (resp. $\tilde\xi\in \tilde X$, resp. $\xi^{\natural}\in Y^{\natural}$) be its image on $Z$ (resp. $\tilde
X$, resp. $Y^{\natural}$). Consider the canonical morphism 
$\kappa:i_Z({\cal O}_{Z,\xi})\otimes _{i_{Y^{\natural}}({\cal O}_{Y^{\natural},\xi^{\natural}})} i_{\tilde X}({\cal O}_{\tilde X, \tilde \xi})\rightarrow i_W({\cal O}_{W,\eta})$ and put
$S:=i_{Y^{\natural}}({\cal O}_{Y^{\natural}})\setminus \{0\}$. We note that since $i_Z({\cal O}_{Z,\xi})$ is flat over $i_{Y^{\natural}}({\cal
O}_{Y^{\natural},\xi^{\natural}})$ by our construction,  $i_Z({\cal O}_{Z,\xi})\otimes _{i_{Y^{\natural}}({\cal
O}_{Y^{\natural},\xi^{\natural}})} i_{\tilde X}({\cal O}_{\tilde X, \tilde \xi})$ is a free $i_{\tilde X}({\cal O}_{\tilde X, \tilde \xi})$-module, in particular, a torsion free
$i_{Y^{\natural}}({\cal O}_{Y^{\natural},\xi^{\natural}})$-module. Since 
$S^{-1}\kappa:S^{-1}(i_Z({\cal O}_{Z,\xi})\otimes _{i_{Y^{\natural}}({\cal O}_{Y^{\natural},\xi^{\natural}})} i_{\tilde X}({\cal O}_{\tilde X, \tilde \xi}))\simeq S^{-1}(i_Z({\cal
O}_{Z,\xi}))\otimes _{i_{Y^{\natural}}({\cal M}_{Y^{\natural}})} S^{-1}(i_{\tilde X}({\cal O}_{\tilde X, \tilde
\xi}))\rightarrow i_W({\cal M}_W)$ is injective by the previous argument, so is $\kappa$, hence, in particular, $\mbox{Im }\kappa$ is a normal subring of $i_W({\cal O}_{W,\eta})$
whose total quotient ring coincides $i_W({\cal M}_W)$ which implies that $\mbox{Im }\kappa=i_W({\cal M}_W)$. Thus we conclude that $\kappa$ is an isomorphism and that $i_W({\cal
O}_{W,\eta})$ is flat over $i_{\tilde X}({\cal O}_{\tilde X,\tilde \xi})$, which implies that $(W,\beta,i_W)\in \mbox{Ob }{\rm Cart}^{\dag}G(X;f^{\ast}D)^p$. The canonical inclusion
$\mbox{Gal }(W/X)\simeq \mbox{Gal }(i_Z({\cal M}_Z)/i_Z({\cal M}_Z)\cap i_X({\cal M}_X))\rightarrow \mbox{Gal }(Z/Y)$ induces a homomorphism $f_{\ast}:\hat\pi_{1,X,p}^{{\rm
loc}}[f^{\ast}D]\rightarrow \hat\pi_{1,Y,q}^{{\rm loc}}[D]$. Since we have $[\mbox{Gal }(Z/Y):\mbox{Gal }(W/X)]=[i_Z({\cal M}_Z)\cap i_X({\cal M}_X):i_Y({\cal M}_Y)]\leq \deg f$, we
get the assertion (see also \cite{bourbaki}, \S 7.1, Corollaire 3).
\hfill\bsquare

\subsection{Lefshetz type theorem for $D$-local algebraic fundamental groups}
The aim of this section is to state and prove the Lefshetz type theorem for $D$-local algebraic fundamental groups.  

\begin{lm}[c.f., \cite{kawamata crep.}, Corollary 10.8]\label{lm:inv} Take any $D\in \mbox{\em Div}_{\mbox{\boldmath $Q$}}X\cap \mbox{\em Weil }X$ and let $\pi:\tilde X\rightarrow X$ be
the index one cover with respect to D. Assume that there exists a normal prime divisor $\Gamma$ passing through
$p\in X$ such that the following three conditions hold. 
\begin{description}
\item[$(1)$ ] \ $\tilde \Gamma:=\pi^{-1}\Gamma$ is normal, 
\item[$(2)$ ] \ $\Gamma$ does not contained in $\mbox{\em Supp }D$, 
\item[$(3)$ ] \ there exists an analytic closed subset
$\Sigma \subset X$ with $\mbox{\em codim}_X\Sigma\geq 2$ and $\mbox{\em
codim}_{\Gamma}(\Sigma\cap \Gamma)\geq 2$ such that $D|_U$ is Cartier and
$D_{\Gamma}:=j_{\Gamma}^{\ast}D\in \mbox{\em Div }\Gamma$, where $U:=X\setminus \Sigma$
and $j_{\Gamma}:\Gamma_0:=\Gamma\setminus \Sigma\rightarrow \Gamma$ is the natural embedding.  
\end{description} Then $D\in \mbox{\em Div }X$.
\end{lm}

{\it Proof.} Consider the exact sequence: 
$$
0\longrightarrow {\cal O}_{\tilde X}(\pi^{\ast}D-\tilde\Gamma)\longrightarrow {\cal O}_{\tilde X}(\pi^{\ast}D)\longrightarrow {\cal
O}_{\tilde\Gamma}(\pi^{\ast}D)\longrightarrow 0. 
$$
Note that ${\cal O}_{\tilde \Gamma}(\pi^{\ast}D)$ and $\pi^{\ast}{\cal O}_{\Gamma}(D_{\Gamma})$ are both invertible and coincide on $\tilde\Gamma\setminus \pi^{-1}(\Sigma)$, hence we
have
${\cal O}_{\tilde
\Gamma}(\pi^{\ast}D)=\pi^{\ast}{\cal O}_{\Gamma}(D_{\Gamma})$ by the normality of $\tilde\Gamma$. Since
$\pi_{\ast}^G:=\Gamma_X^G\circ\pi_{\ast}$ is an exact functor, where $G:=\mbox{Gal }(\tilde X/X)$, the above exact sequence induces a surjective map 
$$
\alpha:{\cal O}_X(D)=\pi^G_{\ast}{\cal O}_{\tilde X}(\pi^{\ast}D)\rightarrow \pi_{\ast}^G{\cal
O}_{\tilde\Gamma}(\pi^{\ast}D)={\cal
O}_{\Gamma}(D_{\Gamma})\otimes(\pi|_{\tilde\Gamma})_{\ast}^G{\cal O}_{\tilde\Gamma}={\cal O}_{\Gamma}(D_{\Gamma}).
$$ Take $\varphi_{\Gamma}\in {\cal
M}_{\Gamma}$ such that $\mbox{div }\varphi_{\Gamma}=-D|_{\Gamma}$. By the above argument, we have
$\varphi\in{\cal M}_X$ such that $\alpha(\varphi)=\varphi_{\Gamma}$. Since $(D+\mbox{div
}\varphi)|_{\Gamma}=D_{\Gamma}+\mbox{div }\varphi_{\Gamma}=0$ and $D+\mbox{div }\varphi$
is {\bf Q}-Cartier, we deduce that $\Gamma\cap\mbox{Supp }(D+\mbox{div
}\varphi)=\emptyset$, that is, $D+\mbox{div }\varphi=0$, hence $D\in \mbox{Div }X$.
\hfill\bsquare

\begin{lm}[c.f., \cite{reid c}, Lemma 1.12, \cite{seidenberg}]\label{lm:seidenberg} 
Let $X$ be a normal complex analytic space embedded in some domain in $\mbox{\boldmath
$C$}^n$. Consider the hypersurfaces
$H_{\tau}$ on $\mbox{\boldmath $C$}^n$ parametrized by $\tau\in \mbox{\boldmath
$P$}^n$ which is defined by a linear equation
$\tau_0+\sum_{i=1}^n\tau_iz_i=0$, where $z_1$, $\dots$, $z_n$ is a complete coordinate system of $\mbox{\boldmath
$C$}^n$. Then there exist a non-empty open subset $U\subset \mbox{\boldmath
$P$}^n$ and a countable union $Z$ of closed analytic subsets of $U$ such that for any $\tau\in U\setminus Z$, $H_{\tau}\cap X$ is a normal hypersurface on $X$. 
\end{lm}

{\it Proof.} Take an analytic open subset $U\subset \mbox{\boldmath
$P$}^n$ such that for any $\tau\in U$, $\bar H_{\tau}:=H_{\tau}\cap X$ is non empty and $H_{\tau}$
does not contain $X$. Since the base point free linear system $\{\bar H_{\tau}\}_{\tau\in U}$ on $X$
induces a base point free linear system on
$\mbox{Reg }X$, we have $\mbox{Sing }\bar H_{\tau}\subset \mbox{Sing } X$ and 
$\mbox{codim }_{\bar
 H_{\tau}}\mbox{Sing }\bar H_{\tau}\geq 2$ for any
$\tau\in U\setminus Z$, where $Z$ is a countable union of closed analytic subsets of $U$ by Bertini's theorem. Moreover, we may assume that for $k=1,\dots, d-2$, $\bar H_{\tau}$ does
not contain any maximal dimensional components of $(\mbox{Sing }X)\cap S_{k+1}({\cal O}_X)$, where
$d:=\dim X$ and $S_{k}({\cal O}_X)$ is a closed analytic set consisting of points at which the 
profoundity of ${\cal O}_X$ does not
exceed $k$. Since we have
$(\mbox{Sing }\bar H_{\tau})\cap S_{k}({\cal O}_{\bar H_{\tau}})\subset (\mbox{Sing }X)\cap
S_{k+1}({\cal O}_X)$, we see that $\dim (\mbox{Sing }\bar H_{\tau})\cap S_{k}({\cal O}_{\bar
H_{\tau}})\leq k-2$ for any
$k$, hence $\bar H_{\tau}$ is normal for any $\tau\in U\setminus Z$ by \cite{fischer},  2.27, Theorem.
\hfill\bsquare

\begin{rk}{\em Let $X$ be a normal Stein space. For
{\bf Q}-divisor $\Delta$ on $X$, let
$\mbox{Mult}_X(\Delta)\subset
\mbox{\boldmath $Q$}$ denote the subset consisting of all the multiplicities of $\Delta$ at prime
divisors on $X$. We note that for general normal hyperplanes $\bar H_{\tau}$, we have 
$\mbox{Mult}_{\bar H_{\tau}}(\mbox{Diff}_{\bar H_{\tau}}(\Delta))\subset \mbox{Mult}_X(\Delta)$.
}\end{rk}

To state the Lefshetz type theorem, we need to fix some sort of general conditions. We shall
consider the following conditions assuming $\dim X\geq 2$. 
\begin{description}
\item[ $(M1)$ ] $\Delta$ is a standard {\bf Q}-boundary. 
\item[ $(M2)$ ] $(X,\Delta)$ is divisorially log terminal.
\item[ $(M2)^{\ast}$ ] $(M2)^{\alpha}$ $(X,\Delta)$ is divisorially log terminal and $\{\Delta\}=0$ or $(M2)^{\beta}$ $(X,\Delta)$ is purely log terminal.
\item[ $(M3)$ ] There exists an irreducible component
$\Gamma$ of $\lfloor\Delta\rfloor$ passing through $p\in X$ such that $K_X+\Gamma$ is {\bf Q}-Cartier.
\end{description}

\begin{rk} {\em $(M2)^{\ast}$ is a slightly stronger condition than $(M2)$.}

\end{rk}

\begin{pr}\label{pr:ind} Assume the conditions $(M1)$, $(M2)$ and $(M3)$. Then 
$$
\mbox{\em ind}_p(K_X+\Delta)=\mbox{\em ind}_p(K_{\Gamma}+\mbox{\em Diff}_{\Gamma}(\Delta-\Gamma)).
$$
\end{pr}

{\it Proof.} Put $r_{\Gamma}:=\mbox{ind}_p(K_{\Gamma}+\mbox{Diff}_{\Gamma}(\Delta-\Gamma))$. 
Firstly, we note that $(X,\Gamma)$ is purely log terminal and that $\Gamma\cap \mbox{Supp }(\Delta-\Gamma)$ is purely one codimensional in $\Gamma$ since $\Delta-\Gamma$ is 
{\bf
Q}-Cartier by the conditions $(M2)$ and $(M3)$. We show that
$r_{\Gamma}(K_X+\Delta)$ is an integral divisor on
$X$ and is Cartier at general points of any prime divisors on $\Gamma$. By taking general
hyperplane sections, we only have to check that if $\dim X=2$, then $r_{\Gamma}(K_X+\Delta)$ is
Cartier. This can be checked by the classification of log canonical singularities with a standard {\bf Q}-boundary due to S. Nakamura
(see, \S 3.1 or 
\cite{kobayashi}, Theorem 3.1), but we can also argue in this way as follows. We note that $p\in X$ is a cyclic quotient singular point with the order, say, 
$n$ by the condition $(M2)$. If $\Gamma\cap(\lfloor\Delta\rfloor-\Gamma)\neq \emptyset$, then $X$ is smooth, hence this case is
trivial. Assume that $\Gamma\cap(\lfloor\Delta\rfloor-\Gamma)=\emptyset$. Since we can see that $(X,\Delta)$ is purely log terminal in this
case from the condition $(M2)$, we can write $\Delta=\Gamma+d\Xi$ for a prime divisor $\Xi$ such that
$(\Gamma,\Xi)_p=1$ and for some
$d=(l-1)/l$, where $l$ is a natural number and we have $\mbox{mult}_p\mbox{Diff}_{\Gamma}(\Delta-\Gamma)=(nl-1)/(nl)$, as in \cite{shokurov-cp}, Lemma
2.25, which implies $r_{\Gamma}(K_X+\Delta)\in \mbox{Div }X$. Going back to the general case, we see that
$D:=K_X+\Delta\in \mbox{Div}_{\mbox{\boldmath $Q$}}X$ and $\Gamma$ satisfies the conditions in Lemma~\ref{lm:inv} using
\cite{shokurov}, Corollary 2.2 and Lemma 3.6, hence we conclude that $r_{\Gamma}(K_X+\Delta)\in \mbox{Div }X$. 
\hfill\bsquare

\begin{rk} {\em We note that $\mbox{Diff}_{\Gamma}(\Delta-\Gamma)$ is also a standard {\bf Q}-boundary, 
since $\mbox{Diff}_{\tilde\Gamma}((\Delta-\Gamma)_{\tilde X})$ is a {\bf Q}-boundary (see \cite{shokurov}, (2.4.1)).

}\end{rk}
\begin{ex} {\em Let $X$ be the germ of $\mbox{\boldmath $C$}^2$ at the
origin and put $\Gamma:=\mbox{div }z$ and  
$\Delta:=\mbox{div }z+(1/n)\mbox{div }w+(1/n)\mbox{div }(z+w)$, where $(z,w)$ is a system of coordinates and $n\in \mbox{\boldmath $N$}$. Then we have $\mbox{ind}_0(K_X+\Delta)=n$ while
$\mbox{ind}_0(K_{\Gamma}+\mbox{Diff}_{\Gamma}(\Delta-\Gamma))=n/2$ (resp. $n$) if
$n$ is even (resp. if $n$ is odd), which explains why we need the assumptions in Proposition~\ref{pr:ind}. }
\end{ex}

A directed set $(\Lambda, \geq)$ naturally forms a cofilterd projective system assuming that for $\lambda$, $\mu\in \Lambda$, 
$\mbox{Card Hom}_{\Lambda}(\lambda,\mu)=1$ if and only if $\lambda\geq \mu$. We call this projective system $\Lambda$ a {\it cofilterd index
projective system}. Let us recall the following basic result (see for example, \cite{ribes}).

\begin{lm}\label{lm:topgrp} Let $\phi:\Lambda^{\prime}\rightarrow \Lambda$ be a covariant functor between cofilterd index
projective systems and $G:\Lambda\rightarrow \mbox{\em (Top. groups)}$, $H:\Lambda^{\prime}\rightarrow \mbox{\em (Top. groups)}$ be
two covariant functors to the category of topological groups. Assume that the following three conditions $(a)$, $(b)$ and $(c)$
hold.
\begin{description}
\item[$(a)$ ] $G_{\lambda}:=G(\lambda)$ and
$H_{\lambda^{\prime}}:=H(\lambda^{\prime})$ are compact for any $\lambda\in \mbox{\em Ob }\Lambda$ and
$\lambda^{\prime}\in\mbox{\em Ob }\Lambda^{\prime}$.
\item[$(b)$ ] $G(\lambda\rightarrow \mu)$ and $H(\lambda^{\prime}\rightarrow \mu^{\prime})$ are all surjective.
\item[$(c)$ ] There exists a natural transformation $\Psi:G\circ\phi\rightarrow H$ such that
$\Psi(\lambda^{\prime}):G_{\phi(\lambda^{\prime})}\rightarrow H_{\lambda^{\prime}}$ are surjective for any 
$\lambda^{\prime}\in\mbox{\em Ob
}\Lambda^{\prime}$.
\end{description}
Then there exists a canonical surjective morphism in \mbox{\em (Top. groups)} $:$
$$
\psi:\mbox{\em projlim }_{\lambda\in {\rm Ob }\Lambda}G_{\lambda}\rightarrow \mbox{\em projlim }_{\lambda^{\prime}\in {\rm Ob
}\Lambda^{\prime}}H_{\lambda^{\prime}}.
$$
\end{lm}

\vskip 5mm
Let $\Delta$ be a {\bf Q}-divisor on $X$ such that $K_X+\Delta$ is {\bf Q}-Cartier. In what follows, we put 
$$
{\cal I}^{\dag(m)}_1(G)(X,\Delta)^{(p)}:={\cal I}^{\dag(m)}_1(G)(X;K_X+\Delta)^{(p)}.
$$

\begin{th}\label{th:pi1surj} Assume the conditions $(M1)$, $(M2)^{\ast}$ and $(M3)$. Then there exists a canonical continuous 
surjective homomorphism $:$
$$
\psi_{\Gamma}:\hat\pi^{{\rm loc}}_{1,\Gamma,p}[\mbox{\em
Diff}_{\Gamma}(\Delta-\Gamma)]\rightarrow \hat\pi^{{\rm loc}}_{1,X,p}[\Delta]. 
$$\end{th}

{\it Proof. }For 
any
$(Y,f)\in \mbox{Ob }{\rm Cart}^m(X,\Delta)$, $\Gamma$ and $\Gamma_Y:=f^{-1}\Gamma$ are normal by
\cite{shokurov}, Lemma 3.6 and Corollary 2.2, hence they are irreducible since $f^{-1}(p)$ consists
of just one point. A canonical inclusion 
${\cal O}_{\Gamma}\rightarrow  \mbox{injlim}_{(Y,f,i_Y)\in {\rm Ob }{\rm Cart}^m(X,\Delta)^p}{\cal O}_{\Gamma_Y}$ extends to an inclusion 
$i_{\Gamma}:{\cal O}_{\Gamma}\rightarrow \overline{\cal M}_{\Gamma}$ and we fix this $i_{\Gamma}$. Then we have a canonical functor
$
\phi^{(p)}_{\Gamma}:{\rm Cart}^m(X,\Delta)^{(p)}\rightarrow {\rm Cart}^m(\Gamma,\mbox{Diff}_{\Gamma}(\Delta-\Gamma))^{(p)}
$
such that $\phi((Y,f))=(\Gamma_Y,f_{\Gamma})$, where $f_{\Gamma}:=f|_{\Gamma_Y}$. Take any $(Y,f)\in \mbox{Ob }{\rm Cart}^{\dag}(X,\Delta)$. 
We shall show that $(\Gamma_Y,f_{\Gamma})\in \mbox{Ob }{\rm Cart}^{\dag}(\Gamma,\mbox{Diff}_{\Gamma}(\Delta-\Gamma))$.
Let $\tilde X\rightarrow X$ be the log canonical cover with respect to $K_X+\Delta$. Note that $\Gamma_{\tilde X}$ is also normal and 
$\pi_{\Gamma}:=\pi|{\Gamma_Y}:\tilde\Gamma:=\Gamma_{\tilde X}\rightarrow \Gamma$ is the log canonical with respect to $K_{\Gamma}+\mbox{Diff}_{\Gamma}(\Delta-\Gamma)$ by
Proposition~\ref{pr:ind}. Take any pointings $i_{\tilde \Gamma}$ and $i_{\Gamma_Y}$ for 
$(\tilde \Gamma,\pi_{\Gamma})$ and $(\Gamma_Y,f_{\Gamma})\in \mbox{Ob
}{\rm Cart}^m(\Gamma,\mbox{Diff}_{\Gamma}(\Delta-\Gamma))$. We note also that $\varpi_{\Gamma_Y}(i_{\Gamma_Y},i_{\tilde \Gamma})=\varpi_Y(i_Y,i_{\tilde X})|_{\Gamma_Y}$ for some
pointings $i_Y$ and $i_{\tilde X}$ of $(Y,f)$ and $(\tilde X,\pi)\in \mbox{Ob }{\rm Cart}^{\dag}(X,\Delta)$. $(\Gamma_Y,f_{\Gamma})\in \mbox{Ob }{\rm Cart}^{\dag}(\Gamma,\mbox{Diff}_{\Gamma}(\Delta-\Gamma))$. By the covering theorem
in \cite{endre}, $(\tilde X,\tilde \Delta)$ is divisorially log terminal of index one, which
implies that $\tilde X$ is smooth in codimension two, hence, in particular, we have $\mbox{codim }_{\tilde\Gamma}(\mbox{Sing }\tilde X\cap\tilde\Gamma)\geq 2$. 
Since $\varpi_Y(i_Y,i_{\tilde X})$ is \'etale over
$\mbox{Reg }\tilde X$, we conclude that $\varpi_{\Gamma_Y}(i_{\Gamma_Y},i_{\tilde \Gamma})$ is \'etale in codimension one and $(\Gamma_Y,f_{\Gamma})\in \mbox{Ob }{\rm Cart}^{\dag}(\Gamma,\mbox{Diff}_{\Gamma}(\Delta-\Gamma))$. In other words, $\phi^{(p)}$ induces a functor
$
\phi^{(p)}_{\Gamma}:{\rm Cart}^{\dag}(G)(X,\Delta)^{(p)}\rightarrow {\rm Cart}^{\dag}(G)(\Gamma,\mbox{Diff}_{\Gamma}(\Delta-\Gamma))^{(p)},
$
where we used the same notation $\phi^{(p)}$. Consider the two functors
$$
G_X:{\rm Cart}^{\dag}G(X,\Delta)^p\rightarrow \mbox{(Top. groups)}\mbox{ and } G_{\Gamma}:{\rm Cart}^{\dag}G(\Gamma,\mbox{Diff}_{\Gamma}(\Delta-\Gamma))^p\rightarrow \mbox{(Top.
groups)}
$$
such that $G_X((Y,f,i_Y))=\mbox{Gal }(Y/X)$ and $G_{\Gamma}((\Gamma^{\prime},g,i_{\Gamma^{\prime}}))=\mbox{Gal }(\Gamma^{\prime}/\Gamma)$. Since $f$ is \'etale over a general points of
$\Gamma$ for any $(Y,f)\in \mbox{Ob }{\rm Cart}^{\dag}G(X,\Delta)$, there exists a natural equivalence $\Psi_{\Gamma}:G_{\Gamma}\circ\phi^p\rightarrow G_X$, which induces the desired
surjection 
$\psi_{\Gamma}:\hat\pi^{{\rm loc}}_{1,\Gamma,p}[\mbox{Diff}_{\Gamma}(\Delta-\Gamma)]\rightarrow \hat\pi^{{\rm loc}}_{1,X,p}[\Delta]$ by Lemma~\ref{lm:topgrp}.
\hfill\bsquare

\begin{rk}{\em Assume the conditions $(M1)$, $(M2)^{\alpha}$ and $(M3)$. Then, combined with Remark\ref{rk:integral case}, Theorem~\ref{th:pi1surj} says that there
exists a surjection 
$
\psi_{\Gamma}:\hat\pi^{{\rm loc}}_{1,\Gamma,p}[\mbox{Diff}_{\Gamma}(\Delta-\Gamma)]\rightarrow \hat\pi_1^{{\rm loc}}(\mbox{Reg }X).
$
For example, if $\dim X=4$, $\hat\pi^{{\rm loc}}_{1,\Gamma,p}[\mbox{Diff}_{\Gamma}(\Delta-\Gamma)]$ is finite under the assumptions as explained in Remark~\ref{rk:s-w}, hence so
is 
$\hat\pi_1^{{\rm loc}}(\mbox{Reg }X)$. 
}\end{rk}

Letting notation and assumptions be as in Theorem~\ref{th:pi1surj}, we obtain the following corollary. 

\begin{co}\label{co:induction} Assume that the universal Cartier cover of $\Gamma$ with respect to
$K_{\Gamma}+\mbox{\em Diff}_{\Gamma}(\Delta-\Gamma)$ exists. Then there exists the universal Cartier cover of $X$ with respect to $K_X+\Delta$. Moreover, there exists the following
exact sequence $:$
\begin{eqnarray}\label{eq:canexseq}
\{1\}\longrightarrow \hat\pi_1^{{\rm loc}}(\mbox{\em Reg
}\Gamma_{X^{\dag}})\longrightarrow \hat\pi^{{\rm loc}}_{1,\Gamma,p}[\mbox{\em Diff}_{\Gamma}(\Delta-\Gamma)]\longrightarrow \hat\pi^{{\rm loc}}_{1,X,p}[\Delta]\longrightarrow
\{1\},
\end{eqnarray}
where $\pi^{\dag}:X^{\dag}\rightarrow X$ is the universal Cartier cover of $X$ with respect to $K_X+\Delta$.
\end{co}

{\it Proof.} The first assertion follows from Proposition~\ref{pr:fcc} and Proposition~\ref{pr:fund.ext seq}. As for the last statement, let
$\pi^{\dag}_{\Gamma}:\Gamma^{\dag}\rightarrow \Gamma$ be  the universal Cartier cover of $\Gamma$ with respect to
$K_{\Gamma}+\mbox{Diff}_{\Gamma}(\Delta-\Gamma)$. Then the induced morphism $\tau^{\dag}_{\Gamma}:\Gamma^{\dag}\rightarrow \Gamma_{X^{\dag}}$ is the universal Cartier cover of $\Gamma_{X^{\dag}}$ since $\tau^{\dag}_{\Gamma}$ is \'etale in codimension one, which implies that $\mbox{Gal }(\Gamma^{\dag}/\Gamma_{X^{\dag}})\simeq
\hat\pi_1^{{\rm loc}}(\mbox{Reg }\Gamma_{X^{\dag}})$, hence we obtain the desired exact sequence. 
\hfill\bsquare

\begin{rk}\label{rk:3-plt}{\em Let notation be as above. Assume that $(X,p)$ is a three dimensional {\bf Q}-Gorenstein singularity and that $(X,\Gamma)$ is purely log terminal with
$\mbox{Sing X}\subset
\Gamma$. Then we see that
$(\tilde X,\tilde p)$ has only terminal singularities and that  
$(X^{\dag},p^{\dag})$ is an isolated compound Du Val singularity (see, \cite{milnor}, Theorem 5.2). We also note that $\Gamma^{\dag}$ is smooth and that $\Gamma_{X^{\dag}}\in
|-K_{X^{\dag}}|$ is a Du Val element.  Moreover, the above
exact sequence (\ref{eq:canexseq}) reduces to the following exact sequence:
\begin{eqnarray}\label{eq:extseq3}
\{1\}\longrightarrow \pi_1^{{\rm loc}}(\mbox{Reg
}\Gamma_{X^{\dag}})\longrightarrow \pi^{{\rm loc}}_{1,\Gamma,p}[\mbox{Diff}_{\Gamma}(0)]\longrightarrow \pi^{{\rm loc}}_1(\mbox{Reg }X)\longrightarrow \{1\},
\end{eqnarray}
which enables us to calculate the local fundamental group of the germ $(X,p)$, since $\pi_1^{{\rm loc}}(\mbox{Reg }\Gamma_{X^{\dag}})$ and $\pi^{{\rm
loc}}_{1,\Gamma,p}[\mbox{Diff}_{\Gamma}(0)]$ have faithful representations to the special unitary group $SU(2,\mbox{\boldmath $C$})$ and the unitary group
$U(2,\mbox{\boldmath
$C$})$ respectively, both of which are classified. It is important to determine the pair $(X^{\dag}, \pi^{{\rm loc}}_1(\mbox{Reg }X))$ which will lead us to the classification
$3$-dimensional purely log terminal singularities. }\end{rk}

\section{Types of degenerations of algebraic surfaces with Kodaira dimension zero}

\begin{df}\begin{em} (Minimal Semistable Degeneration)  
A minimal model $X\rightarrow {\cal D}$ obtained from a projective semistable degeneration of surfaces with non-negative Kodaira dimension 
$g:Y\rightarrow {\cal D}$ by applying the
Minimal Model Program is called a projective minimal semistable degeneration of surfaces. 
\end{em}\end{df}
A projective log minimal degeneration of Kodaira dimension zero is related to a minimal semistable degeneration as in the following way.

\begin{lm}\label{lm:semistable reduction} 
Let $f:X\rightarrow {\cal D}$ be a projective log minimal degeneration of surfaces with non-negative Kodaira dimension. Then there exists a finite covering $\tau:{\cal D}^{\sigma}\rightarrow {\cal D}$, a projective minimal semistable degeneration 
$f^{\sigma}:X^{\sigma}\rightarrow {\cal D}^{\sigma}$ which is bimeromorphically equivalent to $X\times_{{\cal D}}{\cal D}^{\sigma}$ over 
${\cal D}^{\sigma}$ and a generically finite morphism $\pi:X^{\sigma}\rightarrow X$ such that 
$f\circ \pi=\tau\circ f^{\sigma}$ and $K_{X^{\sigma}}+\Theta^{\sigma}=\pi^{\ast}(K_X+\Theta)$, where $\Theta^{\sigma}:=f^{\sigma \ast}(0)$.

\[\begin{array}{ccc}
X^{\sigma} &\mapright{\pi}&X\\
\mapdown{f^{\sigma}}&&\mapdown{f} \\
{\cal D}^{\sigma}&\mapright{\tau} & {\cal D}
\end{array}\]
\end{lm}

{\it Proof.}\ We use the idea explained in \cite{shokurov}, \S $2$. Let $\mu:Y\rightarrow X$ be a projective resolution of $X$ 
such that the support of the singular fiber of the induced morphism 
$g:Y\rightarrow {\cal D}$ has only simple normal crossings as its singularities. 
By the semistable reduction theorem (\cite{mumford}), there exists a finite covering 
$\tau:{\cal D}^{\sigma}\rightarrow {\cal D}$ and a projective resolution 
$Y^{\sigma}\rightarrow Y\times_{{\cal D}}{\cal D}^{\sigma}$ such that the induced degeneration 
$g^{\sigma}:Y^{\sigma}\rightarrow {\cal D}^{\sigma}$ is semistable. Let $\pi^{\prime}:X^{\prime}\rightarrow X$ be the normalization of 
$X$ in the meromorphic function field of $Y^{\sigma}$ and let $\varphi:X^{\sigma}\rightarrow X^{\prime}$ be a minimal model over 
$X^{\prime}$ obtained by applying the Minimal Model Program to the induced morphism $Y^{\sigma}\rightarrow X^{\prime}$. Then, 
by the ramification formula, we have $K_{X^{\prime}}+\Theta^{\prime}=\pi^{\prime \ast}(K_X+\Theta)$, where 
$\Theta^{\prime}:=\pi^{\prime -1}\Theta$. Since $(X^{\prime},\Theta^{\prime})$ is log canonical, we infer that 
$K_{X^{\sigma}}+\Theta^{\sigma}=\varphi^{\ast}(K_{X^{\prime}}+\Theta^{\prime})$. The induced morphism $f^{\sigma}:X^{\sigma}\rightarrow {\cal D}^{\sigma}$ as the 
Stein factorization of the morphism $X^{\sigma}\rightarrow {\cal D}$ gives the desired minimal 
semistable degeneration. \hfill \bsquare

\begin{df}[cf. Definition~\ref{df:123}]\label{df:type.deg}\begin{em}
A log minimal degeneration $f:X\rightarrow {\cal D}$ of surfaces of
Kodaira dimension zero is said to be {\it of type {\rm I}} (resp. {\it of type {\rm II}}, resp. {\it of type {\rm III}}), if there exists an irreducible component $\Theta_i$ of
$\Theta$ such that $(\Theta_i,
\mbox{Diff}_{\Theta_i}(\Theta-\Theta_i))$ is of type I (resp. of type II, resp. of type III).
\end{em}\end{df}

For a projective log minimal degenerations of
surfaces of Kodaira dimension zero $f:X\rightarrow {\cal D}$, take a projective minimal semistable
degeneration $f^{\sigma}:X^{\sigma}\rightarrow {\cal D}^{\sigma}$ obtained from $f$ as in Lemma~\ref{lm:semistable
reduction}. Then the following holds.

\begin{pr}\label{pr:deg.type}
$f$ is of  type
{\em I} $($ resp. of type {\em II}, resp. of type {\em III} $)$ if
and only if
$f^{\sigma}$ is of type {\em I} $($ resp. of type {\em II}, resp. of
type
{\em III}
$)$. Moreover two projective log minimal
degenerations
$f_j:X_j\rightarrow {\cal D}$
$(j=1,2)$ which are bimeromorphically equivalent to each other over
${\cal D}$ have exactly the same types as each other, i.e., types
{\em I}, {\em II} and {\em III} are bimeromorphic notion which are independent from the choice
of log minimal models.
\end{pr}
{\it Proof.} Let $\Theta_i$ be an irreducible component of $\Theta$ and 
$\Theta_i^{\sigma}$ be an irreducible component of $\pi^{-1}(\Theta_i)$
dominating $\Theta_i$. Since we have
$K_{X^{\sigma}}+\Theta^{\sigma}=\pi^{\ast}(K_{X}+\Theta)$, we have 
$K_{\Theta_i^{\sigma}}+\Delta_i^{\sigma}=\pi^{\ast}(K_{\Theta_i}+\Delta_i)$,
where
$\Delta_i^{\sigma}:=\mbox{Diff}_{\Theta_i^{\sigma}}(\Theta^{\sigma}-\Theta_i^{\sigma})$
and  
$\Delta_i:=\mbox{Diff}_{\Theta_i}(\Theta-\Theta_i)$,
hence $\lfloor\Delta_i^{\sigma}\rfloor=\pi^{-1}(\lfloor\Delta_i\rfloor)$ and
$$
\lfloor\mbox{Diff}_{\lfloor\Delta_i^{\sigma}\rfloor^{\nu}}(\Delta_i^{\sigma}-\lfloor\Delta_i^{\sigma}\rfloor)\rfloor 
=\pi^{\nu -1}
\lfloor\mbox{Diff}_{\lfloor\Delta_i\rfloor^{\nu}}
(\Delta_i-\lfloor\Delta_i\rfloor)\rfloor,
$$
where $\pi^{\nu}:\lfloor\Delta_i^{\sigma}\rfloor^{\nu}\rightarrow
\lfloor\Delta_i\rfloor^{\nu}$ is the induced morphism by $\pi$ between
the normalization of
$\lfloor\Delta_i^{\sigma}\rfloor$ and $\lfloor\Delta_i\rfloor$.  Let
$\pi^{\dag}:X^{\dag}\rightarrow X^{\sigma}$ be the index one cover of
$X^{\sigma}$ with respect to
$K_{X^{\sigma}}$ and $f^{\dag}:X^{\dag}\rightarrow {\cal D}^{\dag}$ be the 
degeneration obtained by Stein factorization. Then putting
$\Theta^{\dag}:=\pi^{\dag -1}(\Theta^{\sigma})$, $f^{\dag}$ is also a
projective minimal semistable degeneration with
$K_{X^{\dag}}+\Theta^{\dag}$ being Cartier. As in the same way as above, 
letting $\Theta_i^{\dag}$ be an irreducible component of
$\pi^{\dag -1}(\Theta_i^{\sigma})$ dominating $\Theta_i^{\sigma}$,
we have
$\lfloor\Delta_i^{\dag}\rfloor=\pi^{\dag
-1}(\lfloor\Delta_i^{\sigma}\rfloor)$ and
$
\lfloor\mbox{Diff}_{\lfloor\Delta_i^{\dag}\rfloor^{\nu}}(
\Delta_i^{\dag}-\lfloor\Delta_i^{\dag}\rfloor)\rfloor 
=\pi^{\dag \nu -1}
\lfloor\mbox{Diff}_{\lfloor\Delta_i^{\sigma}\rfloor^{\nu}}
(\Delta_i^{\sigma}-\lfloor\Delta_i^{\sigma}\rfloor)\rfloor,
$
where
$\Delta_i^{\dag}:=\mbox{Diff}_{\Theta_i^{\dag}}(\Theta^{\dag}-\Theta_i^{\dag})$
and $\pi^{\dag \nu}:\lfloor\Delta_i^{\dag}\rfloor^{\nu}\rightarrow
\lfloor\Delta_i^{\sigma}\rfloor^{\nu}$ is the induced morphism by
$\pi^{\dag}$ between the normalization of
$\lfloor\Delta_i^{\dag}\rfloor$ and $\lfloor\Delta_i^{\sigma}\rfloor$. By
Lemma~\ref{lm:ind 1 ls}, if $f^{\dag}$ is of type I (resp. of type II,
resp. of type III), then for any irreducible component
$\Theta_i^{\dag}$ of $\Theta^{\dag}$,
$(\Theta_i^{\dag},
\mbox{Diff}_{\Theta_i^{\dag}}(\Theta^{\dag}-\Theta_i^{\dag}))$ is of type I (resp. of type II,
resp. of type III). Thus we infer the first assertion. As for the last
assertion, construct two projective minimal semistable degeneration,
$f^{\sigma}_j:X^{\sigma}_j\rightarrow {\cal D}^{\sigma}$ $(j=1,2)$ as in
Lemma~\ref{lm:semistable reduction} from $f_j:X_j\rightarrow {\cal D}$
$(j=1,2)$ such that $f^{\sigma}_1$ and $f^{\sigma}_2$ are
bimeromorphically equivalent over ${\cal D}^{\sigma}$. Since there
exists a sequence of flops between $f^{\sigma}_1$ and
$f^{\sigma}_2$ ( see \cite{kollar}, Theorem 4.9 ), it is easily seen that $f^{\sigma}_1$ and
$f^{\sigma}_2$ have the same type, so we get the last assertion. 
\hfill\bsquare

\begin{rk}\begin{em}
From Proposition~\ref{pr:deg.type}, for a projective log minimal
degenerations of surfaces of Kodaira dimension zero $f:X\rightarrow {\cal
D}$, we can see that if $f$ is of type I (resp. of type II, resp. of type
III), then for any irreducible component
$\Theta_i$ of $\Theta$,
$(\Theta_i,
\mbox{Diff}_{\Theta_i}(\Theta-\Theta_i))$ is of type I (resp. of type II,
resp. of type III).
\end{em}\end{rk}

\section{Non-semistable degenerations of abelian or hyperelliptic surfaces}
In this section we prove Theorem~\ref{th:main th.ab}. Let $f:X\rightarrow {\cal D}$ be a projective log minimal degeneration of surfaces with Kodaira dimension zero. Take the log
canonical cover
$\pi:\tilde X\rightarrow X$ with respect to
$K_X+\Theta$, where $\Theta:=f^{\ast}(0)_{\rm red}$ and let $\tilde f:\tilde X\rightarrow \tilde{\cal D}$ be the induced degeneration via Stein factorization. For any irreducible
component
$\Theta_i$ of $\Theta$, $(X, \Theta_i)$ is purely log terminal, hence so is $(\tilde X,
\pi^{-1}\Theta_i)$. Thus the irreducible decomposition $\pi^{-1}\Theta_i=\sum_{j}\tilde \Theta_j^i$ is disjoint, which implies that each component of $\tilde\Theta:=\tilde
f^{\ast}(0)_{\rm red}$ is {\bf Q}-Cartier hence $\tilde X$ is {\bf Q}-Gorenstein. By \cite{endre},
Covering Theorem, we see that $(\tilde X,\tilde
\Theta)$ is divisorially log terminal. Thus we conclude that
$\tilde f:\tilde X\rightarrow
\tilde{\cal D}$ is a projective log minimal degeneration of surfaces with Kodaira dimension zero with $K_{\tilde X}+\tilde \Theta$ being Cartier. Let $\pi_i:\tilde\Theta_i\rightarrow
\Theta_i$ be the log canonical cover with respect to $K_{\Theta_i}+\Delta_i$, where
$\Delta_i:=\mbox{Diff}_{\Theta_i}(\Theta-\Theta_i)$.

\begin{lm}\label{lm:et} There exists an \'etale morphism $\tilde \Theta_j^i\rightarrow \tilde\Theta_i$.  
\end{lm}

{\it Proof.} Put $r_i:=\mbox{Min}\{n\in \mbox{\boldmath $N$}|n(K_{\Theta_i}+\Delta_i)\sim 0\}$.  By Proposition~\ref{pr:ind}, we infer that there exists an open neighbourhood $U$ of
$X$ containing $\Theta_i$ such that ${\cal O}_U(r_i(K_U+\Theta|_U))\in {\mbox{Tor Pic}}^{\circ}U$. Let $m_i$ be the order of ${\cal O}_U(r_i(K_U+\Theta|_U))$. For any connected
component
$V$ of $\pi^{-1}U$, $\pi|_V$ factors into $\pi|_V=\beta\circ\alpha$, where $\alpha:V\rightarrow W$ is \'etale of degree $m_i$ and $\beta:W\rightarrow U$ is finite of degree $r_i$.
We note $\pi|_V$ also factors into $\pi_i\circ\omega$, where $\omega:\pi^{-1}\Theta_i\rightarrow \tilde\Theta_i$ is finite. Since $\pi|_V$ is cyclic, we see that
$\beta^{-1}\Theta_i\simeq
\tilde\Theta_i$ and that $\omega$ is \'etale. Thus we get the assertion.
\hfill\bsquare 
\vskip 5mm
Assume that $e_{\mbox{top}}(X_t)=0$ for 
$t\in {\cal D}^{\ast}$. From the first part of Corollary~\ref{co:chi0} and Lemma~\ref{lm:et}, 
we obtain $e_{\mbox{orb}}(\tilde\Theta_i\setminus \tilde\Delta_i)=0$, where $\tilde\Delta_i:=\pi^{-1}_i\lfloor\Delta_i\rfloor$. Let $\pi_p:(\tilde X,\tilde p)\rightarrow (X,p)$ be the
log canonical cover of the germ of $X$ at $p\in \Theta_i\setminus \lfloor\Delta_i\rfloor$ with respect to $K_X+\Theta_i$. The last part of Corollary~\ref{co:chi0} says that
$(\tilde X,\pi_p^{-1}\Theta_i)$ has only singularity of type
$V_1(r;a,-a,1)$ at $\tilde p\in \tilde X$, where $(r,a)=1$. Thus by the exact sequence~(\ref{eq:extseq3}) in the previous section, we have $\pi^{{\rm
loc}}_1(\mbox{Reg }(X,p))\simeq\pi^{{\rm loc}}_{1,\Theta_i,p}[\mbox{Diff}_{\Theta_i}(0)]$. So when we want to calculate the local fundamental group $\pi^{{\rm loc}}_1(\mbox{Reg
}(X,p))$, we only have to calculate $\pi^{{\rm loc}}_{1,\Theta_i,p}[\mbox{Diff}_{\Theta_i}(0)]$. We note that the proof of the first assertion of Theorem~\ref{th:main th.ab} is
straightforward  since
$(X,p)$ has the universal Cartier cover $\pi^{\dag}:(X^{\dag},p^{\dag})\rightarrow (X,p)$ with respect to $K_X+\Theta_i$ such that $(X^{\dag},p^{\dag})$ is smooth
( see Proposition~\ref{pr:fcc} ).

\subsection{Case Type II}
In this section, we prove Theorem~\ref{th:main th.ab} in the case of type II. From proposition~\ref{pr:cltypeII}, For $p\in \Theta_i\setminus\Delta_i$, possible $\pi^{{\rm
loc}}_{1,\Theta_i,p}[\mbox{Diff}_{\Theta_i}(0)]$ is calculated to be
$\mbox{\boldmath
$Z$}/n\mbox{\boldmath $Z$}$ or
$\mbox{\boldmath $Z$}/2\mbox{\boldmath $Z$}\oplus \mbox{\boldmath
$Z$}/n\mbox{\boldmath $Z$}$, where $n=2,3,4$ or $6$. In applying the log minimal model program on $f:X\rightarrow {\cal D}$ with respect to $K_X$, we see that each extremal
contraction contracts a prime divisor to a curve and reducing to the surface case, that contracting generic curve does not intersect
$\mbox{Supp }\{\mbox{Diff}_{\Theta_i}(\Theta-\Theta_i)\}$. So the last assertion in the case of type II follows from the following Lemma.
\begin{lm}\label{lm:contII} Let $X$ be a normal $\mbox{\boldmath $Q$}$-Gorenstein $3$-fold and $S_1$, $S_2$, $E$ be mutually distinct {\bf Q}-Cartier prime divisors on
$X$ such that
$S_1\cap S_2=\emptyset$. Putting $D:=S_1+S_2+E$, assume that $(X,D)$ is divisorially log terminal and that 
there exists an extremal contraction $\varphi:X\rightarrow X^{\natural}$ to a normal $3$-fold $X^{\natural}$ such
that $S_i$ is $\varphi$-ample for $i=1$, $2$ and the following
$(1)$,
$(2)$ and $(3)$ hold.

\begin{description}
\item[$(1)$] $-K_X$ is $\varphi$-ample, 
\item[$(2)$] $K_X+D$ is numerically trivial over $X^{\natural}$, 
\item[$(3)$] $E$ is contracted to a curve on $X^{\natural}$ and general fibres of the induced morphism
$\varphi:E\rightarrow \varphi(E)$ does not intersects $\mbox{\em Supp }\{\mbox{\em Diff}_E(D-E)\}$. 
\end{description}
Then
$(X^{\natural}, D^{\natural})$ has only divisorially log terminal singularities, where $D^{\natural}:=\varphi_{\ast}D$.
\end{lm}

{\it Proof.} Let $F$ be any exceptional divisor in the function field of $X^{\natural}$ centered on $X^{\natural}$ 
whose log discrepancy with respect to $K_{X^{\natural}}+ D^{\natural}$
is non-positive. If $F=E$, the conditions $(1)$ and $(3)$ imply that $X^{\natural}$ is smooth and 
$D^{\natural}$ has simple normal crossings at the generic point of the center of $F$.  If $F\neq E$, log discrepancy at $F$ with respect to $K_X+ D$ is non-positive from the
condition $(2)$,  hence the support of $F$ at $X$ is a curve $C$ contained in $S_1\cap E$. If $C$ is contracted to a point by $\varphi$, then $(S_2, C)>0$ which contradicts the
assumption $S_1\cap S_2=\emptyset$. So $X^{\natural}$ is smooth and 
$D^{\natural}$ has simple normal crossings at the generic point of the center of $F$ also in this case. Thus we get the assertion by \cite{endre}, Divisorial Log Terminal Theorem.
\hfill\bsquare

\subsection{Case Type III}
The the results of Theorem~\ref{th:main th.ab} in the case of type III follows from the Propositions~\ref{pr:cltypeIII-1},
\ref{pr:cltypeIII-2.1}, \ref{pr:cltypeIII-2.2}, \ref{pr:cltypeIII-2.3}, \ref{pr:cltypeIII-2.4}, \ref{pr:cltypeIII-2.5} and the following lemma.

\begin{lm}\label{lm:no cusp} Let $(X,p)$ be a three dimensional {\bf Q}-Gorenstein singularity and $\Gamma$ be a prime divisor on $X$ passing through $p\in X$ such that $(X,\Gamma)$ is
purely log terminal. Let $\pi:(\tilde X,\tilde p)\rightarrow (X,p)$ be the log canonical cover with respect to $K_X+\Gamma$ and put $\tilde\Gamma:=\pi^{-1}\Gamma$. Assume
that
$(\Gamma,p)\simeq (\mbox{\boldmath $C$}^2,0)$ and
$\mbox{\em Diff}_{\Gamma}(0)=(1/2)\mbox{\em div}(z^2+w^n)$ $($ $n\geq 2$ $)$, where $(z,w)$ is a system of coordinate of $\Gamma$ at $p\in \Gamma$ and that $(\tilde X,\tilde \Gamma)$
has singularity of type
$V_1(r;a,-a,1)$ at $\tilde p\in \tilde X$, where $(r,a)=1$. Then we have $n=2$.
\end{lm}
{\it Proof.} It can be easily checked that $(X^{\dag},p^{\dag})$
and
$(\Gamma^{\dag},p^{\dag})$ are both smooth and that $\hat\pi^{{\rm loc}}_{1,\Gamma,p}[\mbox{Diff}_{\Gamma}(0)]\simeq \hat\pi^{{\rm loc}}_1(\mbox{Reg }X)\simeq G$, where $G$ is the
dihedral group of the order $2n$ (see also \cite{nakano}). Let $G=<a,b;a^n=1,b^2=1,b^{-1}ab=a^{-1}>$ be a presentation of $G$ and $\rho_{\Gamma}:G\rightarrow U(2, \mbox{\boldmath $C$})$
be a corresponding representation with respect to $(\Gamma,\mbox{Diff}_{\Gamma}(0))$ defined as follows. 

$$
\rho_{\Gamma}(a)=\left(
\begin{array}{cc}
e^{2\pi i/n} & 0\\
0 & e^{-2\pi i/n}
\end{array}
\right),
\quad
\rho_{\Gamma}(b)=
\left(
\begin{array}{cc}
0 & 1\\
1 & 0
\end{array}
\right).
$$

Let $\rho_X:G\rightarrow U(3,\mbox{\boldmath $C$})$ be a corresponding faithful representation with respect to $(X,p)$. Since $\Gamma^{\dag}\subset X^{\dag}$ is invariant under the
action of
$G$ through $\rho_X$, $\rho_X$ is equivalent to $\rho_{\Gamma}\oplus \chi$ for some character $\chi:G\rightarrow \mbox{\boldmath $C$}^{\times}$. Let $K$ be the kernel of the character 
$G\rightarrow \mbox{Aut }{\cal O}_{X^{\dag}}(K_{X^{\dag}}+\Gamma^{\dag})/m_{p^{\dag}}{\cal O}_{X^{\dag}}(K_{X^{\dag}}+\Gamma^{\dag})$ induced by $\rho_X$. Then we see that $K=\mbox{Ker
 det }\rho_{\Gamma}$. Since we have $\mbox{det }\rho_{\Gamma}(a)=1$ and $\mbox{det }\rho_{\Gamma}(b)=-1$, we have $<a>\subset K$ and $b\notin K$. Noting that we have  
$2=[G:<a>]=[G:K][K:<a>]$, we obtain $K=<a>$. We also note that we have $\mbox{ord }\chi(a)=n$ since $\tilde p\in\tilde X\simeq \mbox{\boldmath $C$}^3/K$ is isolated. On the other hand,
since we have
$\chi(b^{-1}ab)=\chi(a^{-1})$, we get $\chi(a)^2=1$. Thus we conclude that $n=2$.
\hfill\bsquare

The last assertion of the Theorem~\ref{th:main th.ab} in the case of type III follows from \cite{reid c}, Theorem 3.1.

\section{Canonical bundle formulae}
\subsection{Review on Fujino-Mori's canonical bundle formula} 
Let $f:X\rightarrow B$ be a morphism from a normal projective variety $X$ with $\dim X=m+1$ onto a smooth irreducible curve $B$ defined over the complex number field. 
Assume that $X$
has only canonical singularities and that a general fibre $F$ of $f$ is an irreducible variety with Kodaira dimension zero. Let $b$ be the smallest positive
integer such that the $b$-th plurigenera of
$F$ is not zero. According to \cite{fujino-mori}, there exists a {\bf Q}-divisor $L_{X/B}$ on $B$ with non-negative degree such that there exists a ${\cal O}_B$-algebra
isomorphism 
$$\oplus _{i\geq 0}{\cal O}_B(\lfloor iL_{X/B}\rfloor)\simeq \oplus _{i\geq 0}f_{\ast}{\cal O}_X(ibK_X).$$ $L_{X/B}$ is well defined in the sense
that $L_{X/B}$ is unique modulo linear equivalence\footnote{Any two {\bf Q}-divisors $D_1$ and $D_2$ on a variety
are said to be linearly equivalent to each other if $D_1-D_2$ is a principal divisor.}. Mori also defined a {\bf Q}-divisor $L^{ss}_{X/B}$ on $B$ as a
\lq\lq moduli contribution'' to a canonical bundle formula as follows. 

\begin{pr}[\cite{fujino-mori}, Corollary 2.5]\label{pr:moduli factor} There exists
a {\bf Q}-divisor $L^{ss}_{X/B}(\leq L_{X/B})$ with $\deg L^{ss}_{X/B}\geq 0$ such that 
\begin{description}
\item[(i)] $\tau^{\ast}L^{ss}_{X/B}\leq L_{X^{\prime}/B^{\prime}}$ for any finite surjective morphism 
$\tau:B^{\prime}\rightarrow B$ from
an irreducible smooth curve $B^{\prime}$, and that 

\item[(ii)] $\tau^{\ast}L^{ss}_{X/B}=L_{X^{\prime}/B^{\prime}}$ at $p^{\prime}\in B^{\prime}$ if $X\times _B B^{\prime}$ has a semistable resolution over a
neighbourhood of $p^{\prime}\in B^{\prime}$ or $f^{\prime \ast}(p^{\prime})$ has only canonical singularities, 

\end{description}

where $f^{\prime}:X^{\prime}\rightarrow B^{\prime}$ is a fibration by taking a
non-singular model of the second projection $X\times _B B^{\prime}\rightarrow
B^{\prime}$. 
\end{pr} 

\begin{rk}\begin{em} Since $L_{X/B}$  and 
$L^{ss}_{X/B}(\leq L_{X/B})$ depend only on the birational equivalence class of $X$ over $B$, we can define these {\bf Q}-divisors even if the singularity of $X$ is worse than
canonical by passing to a non-singular model.
\end{em}\end{rk}
Let $\pi:\tilde F\rightarrow F$ be a proper surjective morphism from a smooth variety $\tilde F$ obtained by taking $b$-th root of the unique element of
$|bK_F|$ and desingularization. Put $N(x):=\mbox{L.C.M.}\{n\in \mbox{\boldmath $N$}|\varphi(n)\leq x\}$, where $\varphi$ denotes the Euler's function. Let
$B_m$ be the m-th Betti number of $\tilde F$. The following theorem says, coefficients appearing in canonical bundle formulae can be
well controlled.

\begin{th}[\cite{fujino-mori}, Proposition 2.8 and Theorem 3.1]\label{th:mori's can.bundle formula}
\begin{description}
\item[(1)] $N(B_m) L^{ss}_{X/B}$ is a Weil divisor.
\item[(2)] Assume that $N L^{ss}_{X/B}$ is a Weil divisor. Then we have 
$bK_X=f^{\ast}(bK_B+L^{ss}_{X/B}+\sum _{p\in B} s_p p)+E,$
where $s_p\in \mbox{\boldmath $Q$}$ and $E$ is an effective {\bf Q}-divisor such that 
\begin{description}
\item[(i)] for each point $p\in B$, there exist positive integers $u_p$, $v_p$ such that
$0<v_p\leq bN$ and 
$$s_p=\frac{bNu_p-v_p}{Nu_p},$$
\item[(ii)] $s_p=0$ if $f^{\ast}(p)$ has only canonical singularity or if $f:X\rightarrow B$ has a semistable resolution in a neighbourhood of $p$, and 
\item[(iii)] $f_{\ast}{\cal O}_X(\lfloor nE\rfloor)={\cal O}_B$ for any $n\in \mbox{\boldmath $N$}$.
\end{description}
\end{description}
\end{th}

\subsection{Canonical Bundle Formulae and Log Minimal Models}

\begin{lm}\label{lm:lcbf}
Let $f:X\rightarrow B$ be a proper morphism from a normal variety $X$ onto a smooth irreducible projective curve $B$ defined over the complex number field whose general fibre $F$ of $f$
is a normal variety with only canonical singularity whose Kodaira dimension is zero. Let $b$ be the smallest positive integer such that the $b$-th plurigenera of $F$ is not
zero as in the previous section. Let $\Sigma\subset B$ be a finite set of points which consists of all the points $p\in B$, such that $f^{\ast}(p)$ is not a normal variety with only
canonical singularity. Assume that
$(X,\Theta)$ is log canonical where 
$\Theta:=(f^{\ast}\Sigma)_{\rm red}$. Then there exists $d\in \mbox{\boldmath $N$}$, such that 
$$f_{\ast}{\cal O}_X(n(K_X+\Theta))={\cal O}_B(n(K_B+(1/b)L^{ss}_{X/B}+\Sigma))$$ 
for any $n\in d\mbox{\boldmath $N$}$.
\end{lm}

{\it Proof.} Take a finite Galois cover $\tau:B^{\sigma}\rightarrow B$ from a smooth projective curve $B^{\sigma}$ such that $X\times _B B^{\sigma}\rightarrow B^{\sigma}$ has a
semistable resolution $g^{\sigma}:Y^{\sigma}\rightarrow B^{\sigma}$. Let $\pi:X^{\sigma}\rightarrow X$ be the normalization in the function field of $Y^{\sigma}$ and let
$f^{\sigma}:X^{\sigma}\rightarrow B^{\sigma}$ be the induced fibration. Put $\Theta^{\sigma}:=\pi^{-1}\Theta$. Looking at
the formula to be proved, we may assume that $\tau$ is \'etale over $B\setminus \mbox{Supp }\Sigma$. Taking sufficiently
divisible $n\in \mbox{\boldmath $N$}$, we have 

\begin{eqnarray*}\tau_{\ast}f^{\sigma}_{\ast}{\cal
O}_{X^{\sigma}}(n(K_{X^{\sigma}}+\Theta^{\sigma}))
&=&\tau_{\ast}f^{\sigma}_{\ast}\pi^{\ast}{\cal O}_X(n(K_X+\Theta))\\
&=&f_{\ast}\pi_{\ast}\pi^{\ast}{\cal O}_X(n(K_X+\Theta))\\
&=&f_{\ast}({\cal O}_X(n(K_X+\Theta)) \otimes \pi_{\ast}{\cal O}_{X^{\sigma}}).
\end{eqnarray*}
On the other hand, since $X^{\sigma}\setminus \mbox{Supp }\Theta^{\sigma}$ has only
canonical singularity, $\Theta^{\sigma}$ is Cartier and $(X^{\sigma},\Theta^{\sigma})$ is log
canonical, we can duduce that $X^{\sigma}$ has only canonical singularity. Therefore, we have 
\begin{eqnarray*}\tau_{\ast}f^{\sigma}_{\ast}{\cal
O}_{X^{\sigma}}(n(K_{X^{\sigma}}+\Theta^{\sigma}))
&=&
\tau_{\ast}{\cal O}_{B^{\sigma}}(n(K_{B^{\sigma}}+(1/b)L^{ss}_{X^{\sigma}/B^{\sigma}}+(\tau^{\ast}\Sigma)_{\rm red}))\\
&=&
\tau_{\ast}\tau^{\ast}{\cal O}_B(n(K_B+(1/b)L^{ss}_{X/B}+\Sigma))\\
&=&
{\cal O}_B(n(K_B+(1/b)L^{ss}_{X/B}+\Sigma))\otimes \tau_{\ast}{\cal O}_{B^{\sigma}}
\end{eqnarray*}
and hence 
$$f_{\ast}({\cal O}_X(n(K_X+\Theta)) \otimes \pi_{\ast}{\cal O}_{X^{\sigma}})={\cal O}_B(n(K_B+(1/b)L^{ss}_{X/B}+\Sigma))\otimes
\tau_{\ast}{\cal O}_{B^{\sigma}}.$$ Taking the invariant part under the action of 
$\mbox{\rm Gal }(B^{\sigma}/B)$, we obtain 
$$f_{\ast}{\cal O}_X(n(K_X+\Theta))={\cal O}_B(n(K_B+(1/b)L^{ss}_{X/B}+\Sigma)).$$ 
\hfill\bsquare

\begin{rk}\begin{em}\label{rk:cbf}Let $f^s:X^s\rightarrow B$ be a strictly log minimal fibration ( or degeneration ) of surfaces with Kodaira dimension zero projective over $B$ as in
Definition~\ref{df:lmf}. Since
$K_{X^s}$ is numerically trivial over
$B$, there exists a positive integer $\ell_p\in { \mbox{\boldmath $N$} }$ such that $f^{s\ast}(p)=\ell_p\Theta^s_p$ for any $p\in B$, where $\Theta^s_p:=f^{s\ast}(p)_{\rm
red}$. Let $\mu:Y\rightarrow X^s$ be a minimal model over $X^s$, that is, $\mu$ is a projective birational morphism from a normal
{\bf Q}-factorial $Y$ with only terminal singularity to $X^s$ such that $K_Y$ is $\mu$-nef. By running the minimal model program
over $B$ starting from the induced morphism 
$g:=f^s\circ\mu:Y\rightarrow B$, we obtain a minimal fibration
$h:Z\rightarrow B$ and a dominating rational map $\lambda:Y-\rightarrow Z$ over $B$. Since $X^s$  has only log terminal singularity, there exists an effective {\bf Q}-divisor $\Delta$
with
$\lfloor
\Delta\rfloor=0$ on $Y$ such that $$K_Y+\Delta=\mu^{\ast}K_{X^s}.$$ Since $K_Y+\Delta$, $K_Z+\lambda_{\ast}\Delta$ and $K_Z$ are all numerically trivial over $B$, there exists a
non-negative rational number $\mu_p\in { \mbox{\boldmath $Q$} }$ such that 
\begin{equation}
\lambda_{\ast}\Delta_p=\mu_ph^{\ast}(p),
\end{equation}
where $\Delta_p$ denotes the restriction of $\Delta$ in a neighbourhood of the fibre over $p\in B$. When $B$ is complete, a canonical bundle formula can be calculated by using
Lemma~\ref{lm:lcbf} as follows

\begin{equation}
K_Z=h^{\ast}(K_B+\frac{1}{b}L^{ss}_{Z/B}+\sum_{p\in B}(\frac{\ell_p-1}{\ell_p}-\mu_p)p).
\end{equation}

Define $s_p\in { \mbox{\boldmath $Q$} }$ by $$s_p:=b(\frac{\ell_p-1}{\ell_p}-\mu_p).$$
We can check Mori's estimate of $s_p$ in Theorem~\ref{th:mori's can.bundle formula} as in the following way. Firstly, we note that we may assume that there exists a prime
$\mu$-exceptional divisor over $p\in B$ such that $E\subset \mbox{Supp }\Delta$ and 
$\lambda_{\ast}E\neq 0$, since otherwise, $\mu_p=0$ and hence $s_p=(\ell_p-1)/\ell_p$. 
Put 
$$I_p:=\mbox{\rm Min}\{n\in \mbox{\boldmath $N$}|n(K_{X^s}+\Theta^s_p) \ \mbox{\rm is Cartier in a neighbourhood of the fibre } f^{s\ast}(p)\}.
$$
Since $(X^s,\Theta^s)$ is log canonical, we have 
$$
0\leq I_pa_l(E;X^s,\Theta^s_p)=I_p(a(E;X^s)-\mbox{mult}_E\mu^{\ast}\Theta^s_p+1)\in \mbox{\boldmath $Z$}.
$$
Thus if we put $v^{\prime}_p:=I_p(-a(E;X^s)+\mbox{mult}_E\mu^{\ast}\Theta^s_p)$, we have $v^{\prime}_p\in \mbox{\boldmath $N$}$ and $0<v^{\prime}_p\leq I_p$. Put 
$$
u^{\prime}_p:=\ell_p\mbox{mult}_E\mu^{\ast}\Theta^s_p=\mbox{mult}_Eg^{\ast}(p)\in \mbox{\boldmath $N$}.
$$
Then we have $\mu_p=-a(E;X^s)/u^{\prime}_p$ and $s_p=b(I_pu^{\prime}_p-v^{\prime}_p)/I_pu^{\prime}_p$. So the estimation of $s_p$ reduces to the estimation of $I_p$ in our case. By
passing to a divisorially log terminal model of
$(X^s,\Theta^s)$ and using Proposition~\ref{pr:ind} and \cite{tsunoda1}, Proof of Theorem 2.1, we get $I_p|N(21)$.
\end{em}
\end{rk}

\begin{rk}\begin{em}We should note that $\ell_p\in \mbox{\boldmath $N$}$ and $\mu_p\in \mbox{\boldmath $Q$}$ depends on the choice of strictly log minimal model. For example, consider a
degeneration of elliptic curve whose singular fibre is of type $_m\mbox{\rm I}_1$. Obviously, the minimal model is a strictly log minimal model in this case and we obtain $\ell_p=m$
and $\mu_p=0$. But blowing up the node of the singular fibre and blowing down the exceptional divisor we obtain another strictly log minimal model. When we use this model, we get
$\ell_p=2m$ and $\mu_p=1/(2m)$. 
\end{em}\end{rk}
To avoid this indeterminacy, we shall introduce the notion of moderate log canonical singularity.

\begin{df}\begin{em}\label{df:moderate} Let $(X,\Delta)$ be a normal log variety and let $\Delta^{\prime}$ be a boundary on $X$ such that $\Delta^{\prime}\leq \Delta$. Assume that
$(X,\Delta)$ is log canonical and that $K_X+\Delta^{\prime}$ is {\bf Q}-Cartier. $(X,\Delta)$ is said to be {\it moderately log canonical with respect $K_X+\Delta^{\prime}$} if for any
exceptional prime divisor
$E$ of the function field of $X$ with $a_l(E;X,\Delta)=0$, the inequality $a_l(E;X,\Delta^{\prime})>1$ holds.
\end{em}\end{df}

\begin{rk}\begin{em} If $X$ is {\bf Q}-Gorenstein and $(X,\Delta)$ is divisorially log terminal, then $(X,\Delta)$ is moderately log canonical with respect to $K_X$ by \cite{endre},
Divisorial Log Terminal Theorem. Thus it is easy to see that for a strictly log minimal fibration (or degeneration) $f^s:X^s\rightarrow B$ constructed in such a way as explained in
Remark~\ref{rk:construction of slmm}, $(X^s,\Theta^s)$ is moderate with respect to $K_{X^s}$, where $\Theta^s:=\sum _{p\in \Sigma}\Theta^s_p$. 
\end{em}\end{rk}

\begin{lm}\label{lm:l_p} Let $f^s_i:X^s_i\rightarrow B$ $(i=1,2)$ be two strictly log minimal fibration $($or degeneration$)$ of surfaces with Kodaira dimension zero which are
birationally equivalent to each other over $B$ and let
$l^{(i)}_p$ be a positive integer such that $f^{s\ast}_i(p)=\ell^{(i)}_p\Theta^s_{i,p}$, where $\Theta^s_{i,p}:=f^{s\ast}_i(p)_{\mbox{\em red }}$ for $i=1,2$. Assume that
$(X^s_1,\Theta^s_{1,p})$ is moderately log canonical with respect to $K_{X^s_1}$. Then $\ell^{(1)}_p\leq \ell^{(2)}_p$. Moreover, if $(X^s_2,\Theta^s_{2,p})$ is also
moderately log canonical with respect to $K_{X^s_2}$, then $f^s_1:X^s_1\rightarrow B$ and $f^s_2:X^s_2\rightarrow B$ are isomorphic in
codimension one to each other over a neighbourhood of $p\in B$ and $\ell^{(1)}_p=\ell^{(2)}_p$.
\end{lm}

{\it Proof.} Take a desingularization $\alpha_i:W\rightarrow X^s_i$ of $X^s_i$ and let $\omega:W\rightarrow B$ be the induced morphism. 
Let $G^{(i)}$ be a $\alpha_i$-exceptional effective {\bf
Q}-divisor defined by
$G^{(i)}:=K_W+\Theta^s_{i,W}-\alpha_i^{\ast}(K_{X^s_i}+\Theta^s_{i,p})$. We note that 
$G^{(1)}-G^{(2)}\in \omega^{\ast}(\mbox{\rm Div}(B)\otimes \mbox{\boldmath $Q$})$
by their definitions. Since $G^{(i)}$ is effective and $\mbox{\rm Supp }G^{(i)}$ does not contain the support of the fibre $\omega^{\ast}(p)$ entirely for $i=1,2$, 
we have $G^{(1)}=G^{(2)}$. If we assume that any $\alpha_1$-exceptional divisor contained in a fibre $\omega^{-1}(p)$ is $\alpha_2$-exceptional, then we have
$\ell^{(1)}_p=\ell^{(2)}_p$ obviously,  so we may assume that there exists a $\alpha_1$-exceptional prime divisor $E\subset \omega^{-1}(p)$ which is not $\alpha_2$-exceptional.
If we assume that $a_l(E;X^s_1, \Theta^s_{i,p})>0$, then $E\subset \mbox{\rm Supp }G^{(1)}=\mbox{\rm Supp }G^{(2)}$, which is a contradiction. Thus we have 
$a_l(E;X^s_1, \Theta^s_{i,p})=0$ and hence $a(E;X^s_1)>0$ by the assumption. Therefore we deduce that
$$
\ell^{(2)}_p=\mbox{\rm mult}_E\omega^{\ast}(p)=\ell^{(1)}_p\mbox{\rm mult}_E\alpha_1^{\ast}\Theta^s_{1,p}=\ell^{(1)}_p(a(E;X^s_1)+1)>\ell^{(1)}_p.
$$
The last assertion also follows from the above argument.
\hfill\bsquare

\begin{pr}[c.f. \cite{kollar}, Theorem 4.9]\label{pr:uniqueness of slmm} Let $f^s_i:X^s_i\rightarrow B$ $(i=1,2)$ be two strictly log minimal fibration $($or degeneration$)$ of surfaces
with Kodaira dimension zero projective over $B$ which are birationally equivalent to each other over $B$. Assume that
$(X^s_i,\Theta^s_{i,p})$ is moderately log canonical with respect to $K_{X^s_i}$ for $i=1,2$. Then $f^s_1:X^s_1\rightarrow B$ and $f^s_2:X^s_2\rightarrow B$ are connected by a sequence
of log flops over a neighbourhood of $p\in B$, that is, there exist birational morphisms between normal threefolds over a neighbourhood of $p\in B$ which are isomorphic in codimension
one$:$  
$$
X^s_1:=X^{(0)}\rightarrow Z^{(0)}\leftarrow X^{(1)}\rightarrow Z^{(1)}\cdots \leftarrow X^{(n)}=:X^s_2,
$$
where $X^{(k)}$ is {\bf Q}-factorial for $k=0,1,\dots n$. 
\end{pr}
{\it Proof. }Take a relatively ample effective divisor $H$ on $X^s_2$ over $B$ and let $H^{\prime}$ be the strict transform of $H$ on $X^s_1$. Applying the log minimal model program on
$X^s_1$ over $B$ with respect to $K_{X^s_1}+\varepsilon H^{\prime}$, where $\varepsilon$ is sufficiently small positive rational number, we may assume at first that $H^{\prime}$ is
$f^s_1$-nef since contraction morphisms appearing in the log minimal model program do not contract divisors. By the Base Point Free Theorem (\cite{nakayama}), some multiple of
$H^{\prime}$ defines a birational morphism $\gamma:X^s_1\rightarrow X^s_2$ over $B$ which is isomorphic in codimension one. Since $X^s_1$ and
$X^s_1$ are both {\bf Q}-factorial, $\gamma$ is an isomorphism and thus we get the assertion.
\hfill\bsquare

\begin{df}\begin{em}
Let $f^s:X^s\rightarrow B$ be a strictly log minimal fibration $($or degeneration$)$ of surfaces
with Kodaira dimension zero projective over $B$ such that
$(X^s,\Theta^s_{p})$ is moderately log canonical with respect to $K_{X^s}$ and let $\ell_p\in \mbox{\boldmath $N$}$ and $\mu_p\in \mbox{\boldmath $Q$}$ be as defined in
Remark~\ref{rk:cbf}. We define $\ell^{\ast}_p\in \mbox{\boldmath $N$}$, $\mu^{\ast}_p\in \mbox{\boldmath $Q$}$ and $s^{\ast}_p\in \mbox{\boldmath $Q$}$ as 
$\ell^{\ast}_p:=\ell_p$, $\mu^{\ast}_p:=\mu_p$ and 
$$s^{\ast}_p:=b(\frac{\ell^{\ast}_p-1}{\ell^{\ast}_p}-\mu^{\ast}_p).$$
\end{em}\end{df}

Proposition~\ref{pr:uniqueness of slmm} give the following:
\begin{co}\label{co:birational invariants} $\ell^{\ast}_p\in \mbox{\boldmath $N$}$ and $\mu^{\ast}_p\in \mbox{\boldmath $Q$}$ are birational $($or bimeromorphic$)$ invariants of
germs of singular fibres over $p\in B$ and hence so is $s^{\ast}_p$.
\end{co}

\begin{ex}\label{ex:invariants in the elliptic case}\begin{em}
For degenerations of elliptic curves, one can define invariants $\ell^{\ast}_p$, $\mu^{\ast}_p$ and $s^{\ast}_p$in the same way and it can be checked that $\ell^{\ast}_p$ coincides with
the multiplicty if the singular fibre is of type $_m\mbox{I}_b$ or otherwise, with the order of the semisimple part of the monodromy group around the singular fibre. We can also obtain
the following well known table:

\begin{center}
Table V\\
\begin{tabular}{c|c|c|c|c|c|c|c|c}\hline
                & $_m\mbox{I}_b$  &  $\mbox{I}_b^{\ast}$ &  $\mbox{II}$  &  $\mbox{II}^{\ast}$ & $\mbox{III}$ & $\mbox{III}^{\ast}$  & $\mbox{IV}$ &
$\mbox{IV}^{\ast}$
\\
\hline\hline
$\ell^{\ast}_p$ & m        &    2          &  6     &   6          & 4     & 4             & 3    & 3  \\ \hline
$\mu^{\ast}_p$         & 0        &    0          &  2/3   &   0          & 1/2   & 0             & 1/3  & 0  \\ \hline
$s^{\ast}_p$    & (m-1)/m  &    1/2        &  1/6   &   5/6        & 1/4   & 3/4           & 1/3  & 2/3  \\ \hline
\hline
\end{tabular}
\end{center}
\end{em}\end{ex}
Here we are using the Kodaira's notation (\cite{kodaira}). See also \cite{fujita2}.

\begin{lm}\label{lm:s=0} $s^{\ast}_p=0$ if and only if $\ell^{\ast}_p=1$.
\end{lm}

{\it Proof. } Firstly, assume that $\ell^{\ast}_p=1$, then we have obviously $s^{\ast}_p=0$ since $s^{\ast}_p\geq 0$ and $\mu^{\ast}_p\geq 0$. Secondarily, assume that $s^{\ast}_p=0$. 
Let $f^s:X^s\rightarrow B$ be a strictly log minimal fibration $($or degeneration$)$ of surfaces
with Kodaira dimension zero projective over $B$ such that
$(X^s,\Theta^s_{p})$ is moderately log canonical with respect to $K_{X^s}$. We may assume that the singularity of $X^s$ is worse than canonical and let $\mu:Y\rightarrow X^s$ and $E$ be
as in Remark~\ref{rk:cbf}. Then from the assumption, we have
$$
s^{\ast}_p=b(\frac{\ell^{\ast}_p-1}{\ell^{\ast}_p}+\frac{a(E;X^s)}{l^{\ast}_p\mbox{\rm mult}_E\mu^{\ast}\Theta^s_p})=0, 
$$
hence $a(E;X^s,\Theta^{\ast}_p)=1-l^{\ast}_p\mbox{\rm mult}_E\mu^{\ast}\Theta^s_p$. Since $a(E;X^s,\Theta^{\ast}_p)\geq 0$ and $l^{\ast}_p\mbox{\rm mult}_E\mu^{\ast}\Theta^s_p\in
\mbox{\boldmath $N$}$, we have $a(E;X^s,\Theta^{\ast}_p)=0$. From the definition of moderately log canonical singularity, we have $a(E;X^s)>0$ and hence $\mu^{\ast}_p<0$, which is
a contradiction.
\hfill\bsquare

\begin{lm}\label{lm:mu=0} Let $f^s:X^s\rightarrow B$ be a strictly log minimal fibration $($or degeneration$)$ of surfaces
with Kodaira dimension zero projective over $B$. Then $\mu_p=0$ if and only if $X^s$ has only canonical singularity over a neighbourhood of $p\in B$.
\end{lm}
{\it Proof. }We shall use notation in Remark~\ref{rk:cbf}. Assume that $\mu_p=0$ and that singularity of $X^s$ is worse than canonical over a neighbourhood of $p\in B$. Let
$W$ be a resolution of the graph of $\lambda$ and let $\alpha:W\rightarrow Y$ and $\beta:W\rightarrow Z$ be projections. Since $\beta_{\ast}\Delta^W=\lambda_{\ast}\Delta=0$ over $p\in
B$, we have $K_W+\Delta^W=\alpha^{\ast}(K_Y+\Delta)=\beta^{\ast}K_Z$ over $p\in B$. Since $Z$ has only terminal singularity by its construction, we obtain $-\Delta^W\geq 0$,
which is absurd. 
\hfill\bsquare

\begin{lm}\label{lm:ell-covering }Let $f^s:X^s\rightarrow {\cal D}$ be a strictly log minimal degeneration of surfaces
with Kodaira dimension zero projective over a unit complex disk ${\cal D}$ with the origin $p:=0\in {\cal D}$ and let $\pi:\tilde {\cal D}\rightarrow {\cal D}$ be a cyclic covering
from another unit disk $\tilde {\cal D}$ with the order $\ell_p$ which is \'etale over ${\cal D}^{\ast}:={\cal D}\setminus \{p\}$. Let $\tilde f:\tilde X\rightarrow \tilde {\cal D}$ be
a relatively minimal degeneration which is bimeromorphic to 
$X^s\times_{{\cal D}}\tilde {\cal D}\rightarrow \tilde {\cal D}$ over $\tilde {\cal D}$. Then $\tilde f^{\ast}(\tilde p)$ is reduced and $(\tilde X, \tilde f^{\ast}(\tilde p))$ is log
canonical and in particular
$s^{\ast}_{\tilde p}=0$, where $\tilde p:=0\in \tilde {\cal D}$.
\end{lm}
{\it Proof. }Let $\tilde X^{\prime}$ be a {\bf Q}-factorization of $(X^s\times_{{\cal D}}\tilde {\cal D})^{\nu}$ (see \cite{kawamata crep.}, \S 6, page 120). Then, as in the
same way as in the proof of Lemma~\ref{lm:semistable reduction}, we see that the induced degeneration
$\tilde f^{\prime}:\tilde X^{\prime}\rightarrow \tilde {\cal D}$ is strictly log minimal, projective over $\tilde {\cal D}$ and that $f^{\prime}(\tilde p)$ is reduced, which imply that
$\tilde X^{\prime}$ has only canonical singularity. Let $\tilde X$ be a minimal model over $\tilde X^{\prime}$, then the induced degeneration $\tilde f:\tilde X\rightarrow \tilde {\cal
D}$ turns out to be a minimal degeneration projective over $\tilde {\cal D}$ and it is easily seen that $\tilde f^{\ast}(\tilde p)$ is reduced and $(\tilde X, \tilde f^{\ast}(\tilde
p))$ is log canonical . Since minimal models are unique up to flops, we get the first assertion. As for the last assertion, we only have to check that $(\tilde X, \tilde f^{\ast}(\tilde
p))$ is moderately log canonical with respect to $K_{\tilde X}$, but which is trivial.
\hfill\bsquare

\vskip 5mm
The degree of the moduli contribution to a  canonical bundle formula can be calculated in a certain condition. 

\begin{pr}\label{pr:s=0 fibration} Let $f:X\rightarrow B$ be a proper surjective morphism from a normal algebraic threefolds $X$ with only canonical singularity onto a smooth projective
curve
$B$ whose general fibre is a surface with Kodaira dimension zero. Assume that $s^{\ast}_p=0$ for any $p\in B$. Then $\deg L^{ss}_{X/B}=\deg f_{\ast}{\cal O}_X(bK_{X/B})$.
\end{pr}
{\it Proof. }Let $f^s:X^s\rightarrow B$ be a strictly log minimal model of $f$ such that
$(X^s,\Theta^s_{p})$ is moderately log canonical with respect to $K_{X^s}$ for any $p\in B$. By Lemma~\ref{lm:s=0}, we have $\ell^{\ast}_p=1$ and hence $\mu^{\ast}_p=0$ for
any $p\in B$, from which we infer that $f^{s\ast}(p)$ is reduced and that $X^s$ has only canonical singularity by Lemma~\ref{lm:mu=0}. By the argument in \cite{mori.cl}, proof of
Definition-Theorem (1.11), there exists a Cartier divisor $\delta\in \mbox{\rm Div} B$ such that $bK_{X^s/B}\sim f^{s\ast}\delta$. Thus we have 
$$f_{\ast}{\cal O}_X(bK_{X/B})=f^s_{\ast}{\cal O}_{X^s}(bK_{X^s/B})\simeq {\cal O}_B(\delta).$$ On the other hand, since $bK_{X^s/B}\sim_{\mbox{\boldmath $Q$}}f^{s\ast}L^{ss}_{X^s/B}$, 
we have $\deg \delta=\deg L^{ss}_{X^s/B}=\deg L^{ss}_{X/B}$. Thus we get the assertion.
\hfill\bsquare

\begin{lm}\label{lm:ell-global covering} Let $f:X\rightarrow B$ be a proper surjective morphism from a normal algebraic threefolds $X$ onto a smooth projective curve $B$
whose general fibre is a surface with Kodaira dimension zero. Assume that $B\simeq \mbox{\boldmath $P$}^1$. Then there exists a Kummer covering $\pi:B^{\prime}\rightarrow B$ from a
smooth projective curve $B^{\prime}$ with $\mbox{\em Gal }(B^{\prime}/B)\simeq \oplus_{p\in B}\mbox{\boldmath $Z$}/\ell^{\ast}_p\mbox{\boldmath $Z$}$ such that
$s^{\ast}_{p^{\prime}}=0$ for any $p^{\prime}\in B^{\prime}$ for the second
projection $p_2:(X\times_{B}B^{\prime})^{\nu}\rightarrow B^{\prime}$
\end{lm}
{\it Proof. }The assertion follows from Lemma~\ref{lm:ell-covering } using the argument in \cite{oguiso}, \S4, page 112.
\hfill\bsquare

\begin{df}\begin{em} Let $f:X\rightarrow B$ be a proper surjective morphism from a normal algebraic threefolds $X$ onto $B\simeq \mbox{\boldmath $P$}^1$
whose general fibre is a surface with Kodaira dimension zero. Let ${\cal C}_f$ be the set of all the pair $(B^{\prime},\pi)$ 
which consists of a smooth projective curve $B^{\prime}$ and a finite surjective morphism $\pi:B^{\prime}\rightarrow B$ such that $s^{\ast}_{p^{\prime}}=0$ for any $p^{\prime}\in
B^{\prime}$ for the second
projection $p_2:(X\times_{B}B^{\prime})^{\nu}\rightarrow B^{\prime}$. For $f$, we define a positive integer $d(f)\in
\mbox{\boldmath $N$}$ as 
$d(f):=\mbox{\rm Min}\{\deg \pi| (B^{\prime},\pi)\in {\cal C}_f\}$.
\end{em}\end{df}

\begin{df}\begin{em}Let ${\cal CY}_B^3$ be the set of all the triple $(X,f,B)$ where $X$ is a normal projective threefold $X$ with only canonical singularity
whose canonical divisor $K_X$ is numerically trivial and $f:X\rightarrow B$ is a projective connected morphism onto $B\simeq \mbox{\boldmath $P$}^1$.
\end{em}\end{df}

The following conjecture is important for the bounding problem of Calabi-Yau threefolds.

\begin{conj}\begin{em}\label{conj:d-conj} There exists $d\in \mbox{\boldmath $N$}$ such that $d(f)\leq d$ for any $(X,f,B)\in {\cal CY}_B^3$.
\end{em}\end{conj}

By Lemma~\ref{lm:ell-global covering}, Conjecture~\ref{conj:d-conj} reduces to the following:

\begin{conj}\begin{em}\label{conj:l-conj} There exists $\ell\in \mbox{\boldmath $N$}$ such that $\prod_{p\in B}\ell^{\ast}_p\leq \ell$ for all 
$(X,f,B)\in {\cal CY}^3_B$.
\end{em}\end{conj}

The following proposition is an important step toward Conjecture~\ref{conj:l-conj}, which can be deduced from Theorem~\ref{th:mori's can.bundle formula}.

\begin{pr}\label{pr:bounding sp} There exists a finite subset ${\cal S}\subset \mbox{\boldmath $Q$}$ and a positive integer $\nu$ such that for all 
$(X,f,B)\in {\cal CY}^3_B$, $\{s^{\ast}_p|p\in B\}\subset {\cal S}$ and $\mbox{\em Card }\{p\in B|s^{\ast}_p>0\}\leq \nu$. 
\end{pr}

By Proposition~\ref{pr:bounding sp}, Conjecture~\ref{conj:l-conj} reduces to the following conjecture on degenerations:

\begin{conj}\begin{em}\label{conj:c-conj} Put $c^{\ast}_p:=\mu^{\ast}_p\ell^{\ast}_p$. There exists a finite subset ${\cal C}\subset \mbox{\boldmath $Q$}$ such that for any
degeneration of algebraic surfaces with Kodaira dimension zero over a one-dimensional complex disk, $c^{\ast}_p\in {\cal C}$.
\end{em}\end{conj}

\begin{rk}\begin{em} As we have seen in Example~\ref{ex:invariants in the elliptic case}, for any degenerations of elliptic curves, we have $c^{\ast}_p\in \{0,1,2,4\}$.
One can state the analogue of conjecture~\ref{conj:c-conj} in higher dimensional cases using the conjectural higher dimensional Log Minimal Model Program. In this case,
the problem involves another problem such as the boundedness of varieties with Kodaira dimension zero.  

\end{em}\end{rk}

\subsection{Abelian Fibred Case}
In this section, we prove the following theorem and apply this to Conjecture~\ref{conj:d-conj} in abelian fibred cases.

\begin{th}\label{th:mu.ab. case} For any degeneration of abelian surfaces over a one-dimensional complex disk, all the possible values of the invariant $\mu^{\ast}_p$ for 
the singular fibres can be listed in Table \mbox{\rm VI} and \mbox{\rm VII} except the case $\mu^{\ast}_p=0$. In particular,we have
$$c^{\ast}_p\in \{0,1/5,1/4,1/3,2/5,1/2,2/3,1,3/2,2,3,4,5,6\}.$$
\end{th}

\begin{df}\begin{em}Let ${\cal CY}_{B,\mbox{\rm ab}}^3$ be the set of all the triple $(X,f,B)$ where $X$ 
is a normal projective threefold $X$ with only canonical singularity
whose canonical divisor $K_X$ is numerically trivial and $f:X\rightarrow B$ is a projective connected morphism 
onto $B\simeq \mbox{\boldmath $P$}^1$ whose geometric generic fibre is an abelian
surface.
\end{em}\end{df}

By Theorem~\ref{th:mu.ab. case}, we can give a positive answer to Conjecture~\ref{conj:d-conj} using the argument in the previous section.

\begin{co}\label{co:co of main-th2}There exists $d\in \mbox{\boldmath $N$}$ such that $d(f)\leq d$ for any $(X,f,B)\in {\cal CY}_{B,\mbox{\rm ab}}^3$.
\end{co}

For the proof of Theorem~\ref{th:mu.ab. case}, we need the following:

\begin{lm}\label{lm:ab.typeII} Let $f:X\rightarrow {\cal D}$ be a log minimal Type {\rm II} degeneration of abelian surfaces. Then local fundamental groups of $X$ at any point in $X$
is cyclic. The same holds also for a strictly log minimal model $f^s:X^s\rightarrow {\cal D}$ obtained by applying the log minimal model program on $f$ with respect to $K_X$.
\end{lm}

{\it Proof. }We use notation in Lemma~\ref{lm:semistable reduction}. Since $X^{\sigma}$ is smooth and the support of the singular fibre $f^{\sigma}$ has only normal crossing
singularity with all of the components being relatively elliptic ruled surfaces, Galois group $G:=\mbox{\rm Gal }({\cal D}^{\sigma}/{\cal D})$ acts biregulary on $X^{\sigma}$ and
$\pi:X^{\sigma}\rightarrow X$ factors into $\pi=\pi_1\circ\pi_2$ where $\pi_1:X^{\sigma}/G\rightarrow X$ is a bimeromorphic morphism and $\pi_2:X^{\sigma}\rightarrow X^{\sigma}/G$ is
the quotient map. Let $\{E_j|j\in J\}$ be the set of all the $\pi_1$-exceptional divisors. Since both of $X^{\sigma}/G$ and $X$ are {\bf Q}-factorial, 
we have $\mbox{\rm Exc }\pi_1=\cup_{j\in J}E_j$, hence $\pi_1$ induces an isomorphism 
$X^{\sigma}/G\setminus \cup_{j\in J}E_j\rightarrow  X\setminus \cup_{j\in J}\mbox{\rm Center}_{X}(E_j)$. So we only have to check that $\mbox{\rm Center}_{X}(E_j)$ is contained in a
double curve of $\Theta$ for any $j\in J$ but which is immediate by \cite{endre}, Divisorial Log Terminal Theorem, since we have $a_l(E_j;X,\Theta)=0$ obviously.
Thus we get the assertion.
\hfill\bsquare

\vskip 5mm
{\it Proof of Theorem~\ref{th:mu.ab. case}. }If there exists a strictly log minimal model $f^s:X^s\rightarrow {\cal D}$ of the degeneration to be considered such that
$(X^s,\Theta^s)$ is moderately log canonical with respect to $K_{X^s}$ and that $X^s$ has only canonical singularity, then we have $\mu^{\ast}_p=0$ by its definition and there is
nothing we have to   prove. In particular,  we do not have to care about the Type III degeneration by Theorem~\ref{th:main th.ab}. Consider a strictly log minimal model
$f^s:X^s\rightarrow {\cal D}$ obtained by applying the log minimal model program on a type I or II log minimal degeneration $f:X\rightarrow {\cal D}$ with respect to $K_X$. Assume that
there exists a point
$x\in X^s$ such that
$(X^s,x)$ is not canonical. By Theorem~\ref{th:main th.ab}, Lemma~\ref{lm:ab.typeII} and \cite{oguiso}, Proposition 3.6 and 3.8,types of singularity of $(X^s,\Theta^s)$ at $x\in
X^s$ are of type $V_1(r;a_0,a_1,a_2)$ with
$r=3, 4, 5, 6, 8, 10, 12$ or
$V_2(r;a_0,a_1,a_2)$ with $r=3,4,6$.  In what follows, we use notation in Remark~\ref{rk:cbf}. By \cite{reid2}, \S 4, $E$ corresponds to some primitive vector
$(1/r)(\overline{ka_0},\overline{ka_1},\overline{ka_2})$ and we have 
$$
a(E;X^s)=(1/r)(\overline{ka_0}+\overline{ka_1}+\overline{ka_2})-1,\quad \mbox{\rm mult}_E\mu^{\ast}\Theta^s=(1/r)\overline{ka_2},
$$ hence 
$$
\mu^{\ast}_p=\frac{r-(\overline{ka_0}+\overline{ka_1}+\overline{ka_2})}{\ell^{\ast}_p\overline{ka_2}}
$$
in the case $V_1(r;a_0,a_1,a_2)$. In the same way, we obtain 
$$
\mu^{\ast}_p=\frac{r-(\overline{ka_0}+\overline{ka_1}+\overline{ka_2})}{\ell^{\ast}_p(\overline{ka_0}+\overline{ka_1})}
$$
in the case $V_2(r;a_0,a_1,a_2)$. Determination of possible primitive vectors in the case $V_1(r;a_0,a_1,a_2)$ with
$r=5, 8, 10, 12$ was essentially done in \cite{oguiso}, since degenerations in these cases are Type I in our terminology which coincides \lq\lq moderate" in Oguiso's terminology. When
we determine possible primitive vectors in other cases, we only have to note that
$\sum_{i=0}^2\overline{ka_i}<r$ which is obvious restriction from the assumption and that 
$(r,\overline{ka_2})=1$ since
$\mbox{\rm Sing }X^s\subset \Theta^s$ in the case $V_1(r;a_0,a_1,a_2)$ and $X^s$ is smooth at the generic point of double curves in the case $V_2(r;a_0,a_1,a_2)$. 
\hfill\bsquare

\begin{rk}\begin{em} Theorem~\ref{th:mu.ab. case} implies the inequality $s_p\geq 1/6$ holds for any degeneration of abelian surfaces which has been obtained in the
case of Type I degenerations by Oguiso (\cite{oguiso}, Main Theorem).
\end{em}\end{rk}

\begin{center}
Table VI, Case $V_1(r;a_0,a_1,a_2)$\\
\begin{tabular}{c|c|c|c|c}\hline
&\mbox{\rm primitive vectors}      &  $\mu^{\ast}_p$ &   $s^{\ast}_p$ & \mbox{\rm divisibility of }$\ell^{\ast}_p$\\ \hline\hline
(1)&(1/3)(1,0,1) & $1/\ell^{\ast}_p$  &   $(\ell^{\ast}_p-2)/\ell^{\ast}_p$ & $3|\ell^{\ast}_p$\\\hline
(2)&(1/4)(1,1,1) & $1/\ell^{\ast}_p$ & $(\ell^{\ast}_p-2)/\ell^{\ast}_p$ & $4|\ell^{\ast}_p$\\ \hline
(3)&(1/4)(0,1,1) & $2/\ell^{\ast}_p$ & $(\ell^{\ast}_p-3)/\ell^{\ast}_p$  & $4|\ell^{\ast}_p$\\ \hline
(4)&(1/5)(1,2,1) & $1/\ell^{\ast}_p$ & $(\ell^{\ast}_p-2)/\ell^{\ast}_p$  & $5|\ell^{\ast}_p$\\ \hline
(5)&(1/6)(3,1,1) & $1/\ell^{\ast}_p$ & $(\ell^{\ast}_p-2)/\ell^{\ast}_p$ & $6|\ell^{\ast}_p$\\ \hline
(6)&(1/6)(2,1,1) & $2/\ell^{\ast}_p$ & $(\ell^{\ast}_p-3)/\ell^{\ast}_p$ & $6|\ell^{\ast}_p$\\ \hline
(7)&(1/6)(1,1,1) & $3/\ell^{\ast}_p$ & $(\ell^{\ast}_p-4)/\ell^{\ast}_p$ & $6|\ell^{\ast}_p$\\ \hline
(8)&(1/6)(1,0,1) & $4/\ell^{\ast}_p$ & $(\ell^{\ast}_p-5)/\ell^{\ast}_p$ & $6|\ell^{\ast}_p$\\ \hline
(9)&(1/8)(5,1,1) & $1/\ell^{\ast}_p$ & $(\ell^{\ast}_p-2)/\ell^{\ast}_p$ & $8|\ell^{\ast}_p$\\ \hline
(10)&(1/8)(3,1,1) & $3/\ell^{\ast}_p$ & $(\ell^{\ast}_p-4)/\ell^{\ast}_p$ & $8|\ell^{\ast}_p$\\ \hline
(11)&(1/8)(3,1,3) & $1/(3\ell^{\ast}_p)$ & $(3\ell^{\ast}_p-4)/(3\ell^{\ast}_p)$ & $8|\ell^{\ast}_p$\\ \hline
(12)&(1/10)(7,1,1) & $1/\ell^{\ast}_p$ & $(\ell^{\ast}_p-2)/\ell^{\ast}_p$ & $10|\ell^{\ast}_p$\\ \hline
(13)&(1/10)(3,1,1) & $5/\ell^{\ast}_p$ & $(\ell^{\ast}_p-6)/\ell^{\ast}_p$ & $10|\ell^{\ast}_p$\\ \hline
(14)&(1/10)(3,1,3) & $1/\ell^{\ast}_p$ & $(\ell^{\ast}_p-2)/\ell^{\ast}_p$ & $10|\ell^{\ast}_p$\\ \hline
(15)&(1/12)(7,1,1) & $3/\ell^{\ast}_p$ & $(\ell^{\ast}_p-4)/\ell^{\ast}_p$ & $12|\ell^{\ast}_p$\\ \hline
(16)&(1/12)(4,3,1) & $4/\ell^{\ast}_p$ & $(\ell^{\ast}_p-5)/\ell^{\ast}_p$ & $12|\ell^{\ast}_p$\\ \hline
(17)&(1/12)(5,1,1) & $5/\ell^{\ast}_p$ & $(\ell^{\ast}_p-6)/\ell^{\ast}_p$ & $12|\ell^{\ast}_p$\\ \hline
(18)&(1/12)(3,2,1) & $6/\ell^{\ast}_p$ & $(\ell^{\ast}_p-7)/\ell^{\ast}_p$ & $12|\ell^{\ast}_p$\\ \hline
(19)&(1/12)(5,1,5) & $1/(5\ell^{\ast}_p)$ & $(5\ell^{\ast}_p-6)/(5\ell^{\ast}_p)$ & $12|\ell^{\ast}_p$\\ \hline
(20)&(1/12)(3,2,5) & $2/(5\ell^{\ast}_p)$ & $(5\ell^{\ast}_p-7)/(5\ell^{\ast}_p)$ & $12|\ell^{\ast}_p$\\ \hline

\hline
\end{tabular}
\end{center}
\vskip 0.3in

\begin{center}
Table VII, Case $V_2(r;a_0,a_1,a_2)$\\
\begin{tabular}{c|c|c|c|c}\hline
&\mbox{\rm primitive vectors}      &  $\mu^{\ast}_p$ &   $s^{\ast}_p$ & \mbox{\rm divisibility of }$\ell^{\ast}_p$\\ \hline\hline
(21)&(1/3)(1,0,1) & $1/\ell^{\ast}_p$  &   $(\ell^{\ast}_p-2)/\ell^{\ast}_p$ & $3|\ell^{\ast}_p$\\\hline
(22)&(1/4)(1,0,1) & $2/\ell^{\ast}_p$ & $(\ell^{\ast}_p-3)/\ell^{\ast}_p$ & $4|\ell^{\ast}_p$\\ \hline
(23)&(1/4)(1,1,1) & $1/(2\ell^{\ast}_p)$ & $(2\ell^{\ast}_p-3)/(2\ell^{\ast}_p)$ & $2|\ell^{\ast}_p$\\ \hline
(24)&(1/6)(1,0,1) & $4/\ell^{\ast}_p$ & $(\ell^{\ast}_p-5)/\ell^{\ast}_p$ & $6|\ell^{\ast}_p$\\ \hline
(25)&(1/6)(1,1,1) & $3/(2\ell^{\ast}_p)$ & $(2\ell^{\ast}_p-5)/(2\ell^{\ast}_p)$ & $3|\ell^{\ast}_p$\\ \hline
(26)&(1/6)(1,2,1) & $2/(3\ell^{\ast}_p)$ & $(3\ell^{\ast}_p-5)/(3\ell^{\ast}_p)$ & $2|\ell^{\ast}_p$\\ \hline
(27)&(1/6)(1,3,1) & $1/(4\ell^{\ast}_p)$ & $(4\ell^{\ast}_p-5)/(4\ell^{\ast}_p)$ & $3|\ell^{\ast}_p$\\ \hline

\hline
\end{tabular}
\end{center}
\section{Bounding the number of singular fibres of Abelian Fibered Calabi Yau threefolds}

Let $B$ be a smooth projective irreducible curve defined over the complex number field and let $K$ denote the rational function field of $B$. 
\begin{df}\label{df:Sigma-f}\begin{em} 
We define the loci of singular fibres $\Sigma_f$, $\Sigma_{\varphi}$ as follows.
\begin{description}
\item[ ($1$) ] Let $f:X\rightarrow B$ is a projective connected morphism from a normal {\bf Q}-factorial projective variety $X$ with only
canonical singularity onto a $B$. We define the subset of closed points of $B$, $\Sigma_f$ by
$$
\Sigma_f:=\{p\in B |\mbox{ $f$ is not smooth over a neighbourhood of $p$ }\}.
$$
\item[ ($2$) ] Let $A_K$ be an abelian variety over $K$,
and let $\varphi:A\rightarrow B$ be the N\'eron model of $A_K$. 
We define the subset of closed points of $B$, $\Sigma_{\varphi}$ by
$$
\Sigma_{\varphi}:=\{p\in B|\mbox{ $A_K$ does not have good reduction at $p$ }\}.
$$ 
\end{description}
\end{em}\end{df}

\begin{rk}\begin{em} Consider two morphisms $f_i:X_i\rightarrow B$ ($i=1,2$) as in Definition~\ref{df:Sigma-f}, (1), such that $\dim X_i=3$ and the that the geometric generic fiber of
$f_i$ is an abelian variety for $i=1,2$.
 Assume that $K_{X_i}$ is $f_i$-nef for $i=1,2$ and that $X_1$ is birationally equivalent to $X_2$ over $B$. Then we see that $\Sigma_{f_1}=\Sigma_{f_2}$ and hence $\mbox{Card
}\Sigma_1=\mbox{Card }\Sigma_2$, for, as is well known, birationally equivalent minimal models are connected by a sequence of flops while abelian varieties do not contain rational
curves. Moreover, both of the defintions of (1) and (2) are compatible in the case $\dim A_K=2$, that is, for a
birational model
$f:X\rightarrow B$ of
$A_K$ such as in Definition~\ref{df:Sigma-f}, (1), $\Sigma_f$ coincides with $\Sigma_{\varphi}$.
\end{em}\end{rk}

The aim of this section is to prove the following theorem which gives a positive answer to Oguiso's question communicated to the author; Are the numbers
of singular fibres of abelian fibered Calabi Yau threefolds are bounded?

\begin{th}\label{th:bouding the number of singular fibres}
There exists $s\in \mbox{\boldmath $N$}$, such that for any triple $(X,f,B)\in {\cal CY}_{B,\mbox{\rm ab}}^3$,
$$
s_f:=\mbox{\em Card }\Sigma_f\leq s.
$$
\end{th}

\vskip 6mm
Let $f:X\rightarrow B$ be a projective connected morphism from a normal {\bf Q}-factorial projective variety $X$ with only
canonical singularity onto $B$. Put $B_{\circ}:=B\setminus \Sigma_f$ and let $f_{\circ}:X_{\circ}\rightarrow B_{\circ}$ be 
the restriction of $f$ to $X_{\circ}:=f^{-1}(B_{\circ})$. According to the theory of relative Picard schemes and relative Albanese maps due to Grothendiek (see \cite{grothendieck-f} and
see also \cite{fujiki} working on the analytic category), there exists the relative Picard schemes $\mbox{Pic }(X_{\circ}/B_{\circ})\rightarrow B_{\circ}$ representing the Picard
functor whose connected component containing the unit section, denoted by $\mbox{Pic}^{\circ}(X_{\circ}/B_{\circ})\rightarrow B_{\circ}$, is a projective abelian scheme. Moreover,
there exists a morphism
$\alpha:X_{\circ}\rightarrow
\mbox{Alb}^1(X_{\circ}/B_{\circ})$ over $B_{\circ}$ where
$\mbox{Alb}^1(X_{\circ}/B_{\circ})$ is $B_{\circ}$-torsor under the dual projective abelian scheme of $\mbox{Pic}^{\circ}(X_{\circ}/B_{\circ})$ denoted by 
$\mbox{Alb}^0(X_{\circ}/B_{\circ})$ satisfying the
universal mapping property. $\alpha$ is called {\it the relative Albanese map associated with $f$}. If we assume that the geometric generic fibre of $f$ is an abelian variety, then 
$\alpha$ is an isomorphism over $B_{\circ}$ by the universal mapping property.

\begin{df}\begin{em} Let $X_{\eta}$ be the generic fibre of $f$. 
The N\'eron Model $\varphi:{\cal A}\rightarrow B$ of $\mbox{Alb}^0(X_{\eta})$ is an extension of the abelian scheme $\varphi_{\circ}:\mbox{Alb}^0(X_{\circ}/B_{\circ})\rightarrow
B_{\circ}$. We call $\varphi$, {\it the Albanese group scheme associated with $f$}.
\end{em}
\end{df}

\begin{co}There exists $s\in \mbox{\boldmath $N$}$, such that for any triple $(X,f,B)\in {\cal CY}_{B,\mbox{\rm ab}}^3$,
$$
s_{\varphi}:=\mbox{\em Card }\Sigma_{\varphi}\leq s,
$$
where $\varphi:{\cal A}\rightarrow B$ is the Albanese group scheme associated with $f$.
\end{co}

{\it Proof.} Since ${\cal A}\times_{B}B_{\circ}\simeq \mbox{Alb}^0(X_{\circ}/B_{\circ})$, we have $\Sigma_{\varphi}\subset \Sigma_f$ and hence $s_{\varphi}\leq s_f$.
Thus the assertion follows immediately from Theorem~\ref{th:bouding the number of singular fibres}.

\hfill\bsquare

The following Lemma can be deduced immediately from the property of torsors.

\begin{lm}\label{lm:duality}
Let 
$f:X\rightarrow B$ be a projective connected morphism from a normal {\bf Q}-factorial projective variety $X$ with only
canonical singularity onto $B$ whose geometric generic fiber is an abelian variety and let $\varphi:{\cal A}\rightarrow B$ be the Albanese group scheme associated to $f$. 
\begin{description}
\item[ (i) ] we have
$L^{ss}_{X/B}\sim_{\mbox{\boldmath $Q$}}L^{ss}_{{\cal A}/B}$ and 
\item[ (ii) ] moreover, if we assume that $\dim X=3$ and that $f:X\rightarrow B$ admits an analytic local section in a neighbourhood of any closed points $p\in B$,
then we have $\Sigma_f=\Sigma_{\varphi}$ and $s_p^{\ast}(f)=s_p^{\ast}(\varphi)$, where $s_p^{\ast}(f)$ and $s_p^{\ast}(\varphi)$ are the analytic local bimeromorphic invariants
$s_p^{\ast}$ of the fibres of $f$ and
$\varphi$ over
$p\in B$ defined in Lemma~$\ref{lm:s=0}$ respectively. In particular, the Kodaira dimensions of $X$ and ${\cal A}$ are the same.
\end{description}
\end{lm}

By Corollary~\ref{co:co of main-th2}, there exists $d\in \mbox{\boldmath $N$}$ such that $d(f)\leq d$ for any $(X,f,B)\in {\cal CY}_{B,\mbox{\rm ab}}^3$
and in what follows, we fix such $d\in \mbox{\boldmath $N$}$. Firstly, take any $(X,f,B)\in {\cal CY}_{B,\mbox{\rm ab}}^3$. Then there exists a finite
surjective morphism $\tau:B^{\prime}\rightarrow B$ from a projective smooth irreducible curve $B^{\prime}$ with $\deg \tau\leq d$ which is \'etale over any
closed points $p\in B$ with $s_p^{\ast}(f)=0$ and $p\neq p_{\circ}$, where $p_{\circ}$ is a certain closed point of $B$ which is not contained in
$\Sigma_f$, such that $s_{p^{\prime}}^{\ast}(f^{\prime})=0$ for any closed points $p^{\prime}\in B^{\prime}$, where $f^{\prime}:X^{\prime}\rightarrow
B^{\prime}$ is a minimal model of $p_2:(X\times_{B}B^{\prime})^{\nu}\rightarrow B^{\prime}$. Since the number of
closed points of $B$ with $s_p^{\ast}(f)>0$ is bounded by Proposition~\ref{pr:bounding sp}, we only have to bound $s_{f^{\prime}}:=\mbox{Card
}\Sigma_{f^{\prime}}$ to prove Theorem~\ref{th:bouding the number of singular fibres}. Let $X^{\prime}_{\eta^{\prime}}$ be the generic fibre of
$f^{\prime}$ over the generic point $\eta^{\prime}$ of $B^{\prime}$ and let $\varphi^{\prime}:{\cal A}^{\prime}\rightarrow B^{\prime}$ be the N\'eron model of
$\mbox{Alb}^0(X^{\prime}_{\eta^{\prime}})$. Since $\ell^{\ast}_{p^{\prime}}(f^{\prime})=1$ for any
closed points $p^{\prime}\in B^{\prime}$ by Lemma~\ref{lm:s=0}, $f^{\prime}:X^{\prime}\rightarrow B^{\prime}$ admits an analytic local section in a
neighbourhood of any closed points $p^{\prime}\in B^{\prime}$, so we only have to bound $s_{\varphi^{\prime}}:=\mbox{Card }\Sigma_{\varphi^{\prime}}$ by 
Lemma~\ref{lm:duality}. From Proposition~\ref{pr:s=0 fibration}, we see that 
$$
\deg e^{\prime\ast}\det\Omega^1_{{\cal A}^{\prime}/B^{\prime}}=\deg L^{ss}_{{\cal A^{\prime}}/B^{\prime}}=\deg L^{ss}_{X^{\prime}/B^{\prime}}=(\deg \tau)(\deg L^{ss}_{X/B})\leq 2d,
$$
where $e^{\prime}:B^{\prime}\rightarrow {\cal A^{\prime}}$ is the unit section of $\varphi^{\prime}$.

\begin{lm}\label{lm:stable reduction}
$\varphi^{\prime}:{\cal A}^{\prime}\rightarrow B^{\prime}$ has semi-abelian reduction at all closed points $p^{\prime}\in B^{\prime}$.
\end{lm}
{\it Proof. }
From Proposition~\ref{pr:s=0 fibration}, we see that, for $\mbox{Alb}^0(X^{\prime}_{\eta^{\prime}})$, the stable height equals to the differential height, hence by \cite{moret-baily},
Ch.X, proposition 2.5 (ii), we get the result.
\hfill\bsquare

\vskip 5mm

Let $K^{\prime}$ be the rational function field of $B^{\prime}$. Define an abelian variety $Z^{\prime}$ over $K^{\prime}$, by 
$$
Z^{\prime}:=(\mbox{Alb}^0(X^{\prime}_{\eta^{\prime}})\times_{K^{\prime}}\mbox{Alb}^0(X^{\prime}_{\eta^{\prime}})^t)^4,
$$ where $t$ denotes its dual. Let $\zeta^{\prime}:{\cal Z}^{\prime}\rightarrow B^{\prime}$ be the N\'eron model of $Z^{\prime}$.

\begin{lm}\label{lm:stable reduction zeta}
Put $B_{\circ}^{\prime}:=B^{\prime}\setminus
\Sigma_{\varphi^{\prime}}$ and let $\zeta^{\prime}_{\circ}:{\cal Z}^{\prime}_{\circ}\rightarrow B^{\prime}_{\circ}$ denote the restriction of
$\zeta^{\prime}$ to $\zeta^{\prime -1}(B^{\prime}_{\circ})$. Then,
\begin{description}
\item[ (i) ] $\Sigma_{\zeta^{\prime}}=\Sigma_{\varphi^{\prime}}$ and $\zeta^{\prime}_{\circ}:{\cal Z}^{\prime}_{\circ}\rightarrow B^{\prime}_{\circ}$ is a principally polarized abelian
scheme.
\item[ (ii) ] $\zeta^{\prime}:{\cal Z}^{\prime}\rightarrow B^{\prime}$ has semi-abelian reduction at all closed points $p^{\prime}\in B^{\prime}$ and 
\item[ (iii) ] $\deg \zeta^{\prime}_{\ast}\det\Omega^1_{{\cal Z}^{\prime}/B^{\prime}}=8\deg e^{\prime\ast}\det\Omega^1_{{\cal A}^{\prime}/B^{\prime}}\leq 16d$
\end{description}
\end{lm}

{\it Proof. }Since
$\varphi_{\circ}^{\prime}:\mbox{Alb}^0(X_{\circ}^{\prime}/B_{\circ}^{\prime})\rightarrow B_{\circ}^{\prime}$ is a projective abelian scheme,
there exists a polarization
$\lambda^{\prime}_{\circ}:\mbox{Alb}^0(X_{\circ}^{\prime}/B_{\circ}^{\prime})\rightarrow \mbox{Alb}^0(X_{\circ}^{\prime}/B_{\circ}^{\prime})^t$ which is an isogeny over
$B_{\circ}^{\prime}$. Let $\varphi^{\prime t}:{\cal A}^{\prime t}\rightarrow B^{\prime}$ be the N\'eron Model of $\mbox{Alb}^0(X^{\prime}_{\eta^{\prime}})^t$. Then
$\lambda^{\prime}_{\circ}$ extends to an isogeny $\lambda^{\prime}:{\cal A}^{\prime}\rightarrow {\cal A}^{\prime t}$ over $B^{\prime}$, hence we have 
$\Sigma_{\varphi^{\prime}}=\Sigma_{\varphi^{\prime t}}$ and ${\cal
A}^{\prime t}\rightarrow B^{\prime}$ also has semi-abelian reduction at all closed points $p^{\prime}\in B^{\prime}$ (see \cite{BLR}, \S 7.3, Proposition 6, Corollary 7 ). Therefore, 
using \cite{BLR}, \S 7.4 Proposition 3, we conclude that ${\cal Z}^{\prime \circ}=({\cal A}^{\prime \circ}\times_{B^{\prime}}{\cal A}^{\prime t\circ})^4\rightarrow B^{\prime}$ is a
semi-abelian sheme and $\Sigma_{\zeta^{\prime}}=\Sigma_{\varphi^{\prime}}$, where ${\cal Z}^{\prime \circ}$, ${\cal A}^{\prime \circ}$ and ${\cal A}^{\prime t\circ}$ 
denotes connected components of ${\cal Z}^{\prime}$, ${\cal A}^{\prime}$ and ${\cal A}^{\prime t}$ containing their unit respectively. It is well known that
$\zeta^{\prime}_{\circ}:{\cal Z}^{\prime}_{\circ}\rightarrow B^{\prime}_{\circ}$ has a principal polarization and that we have the equality $\deg
\zeta^{\prime}_{\ast}\det\Omega^1_{{\cal Z}^{\prime}/B^{\prime}}=8\deg e^{\prime\ast}\det\Omega^1_{{\cal A}^{\prime}/B^{\prime}}$ (see \cite{faltings-chai}, Ch I, \S 5, Lemma 5.5,
\cite{moret-baily}, Ch IX and \cite{zarkhin}).
\hfill\bsquare

\vskip 5mm

Let $H_g$, $\mbox{Sp}(2g,\mbox{\boldmath $Z$})$ and $\Gamma_{g,n}$ be the Siegel's upper half-plane, integral symplectic group of degree $2g$ and the congruence subgroup of level $n$
of 
$\mbox{Sp}(2g,\mbox{\boldmath $Z$})$ respectively as usual. For $n\geq 3$, ${\cal A}_{g,n}:=H_g/\Gamma_{g,n}$ is known to be smooth quasi projective scheme over the complex number field
$\mbox{\boldmath $C$}$. According to $\cite{faltings}$, there exists a normal projective variety ${\cal A}_{g,n}^{\ast}$ over $\mbox{\boldmath $C$}$ which contains ${\cal A}_{g,n}$ as a
open subsheme and an ample invertible $\omega_{{\cal A}_{g,n}^{\ast}}$ such that for any $g$-dimensional principally polarized abelian 
variety $A_K$ defined over $K$ with a level $n$ structure, the period map $\phi:B_{\circ}:=B\setminus \Sigma_{\varphi}\rightarrow {\cal A}_{g,n}$ associated with 
$\varphi_{\circ}:{\cal A}_{\circ}:=\varphi^{-1}(B_{\circ})\rightarrow B_{\circ}$ extends to a morphism $\phi:B\rightarrow {\cal A}_{g,n}^{\ast}$ which gives an isomorphism
$e^{\ast}\det\Omega^1_{{\cal A}/B}\simeq \phi^{\ast}\omega_{{\cal A}_{g,n}^{\ast}}$, where $\varphi:{\cal A}\rightarrow B$ is the N\'eron Model of $A_K$ and $e$ is its unit section.
Take a natural number
$N=N(g,n)\in
\mbox{\boldmath $N$}$ depending only on $g$ and $n$ such that $\omega_{{\cal A}_{g,n}^{\ast}}^{\otimes N}$ is very ample. We may assume that there exists an effective divisor $H$ on
${\cal A}_{g,n}^{\ast}$ such that ${\cal O}_{{\cal A}_{g,n}^{\ast}}(H)\simeq \omega_{{\cal A}_{g,n}^{\ast}}^{\otimes N}$ and that ${\cal A}_{g,n}^{\ast}\setminus
{\cal A}_{g,n}\subset\mbox{Supp H}$.  Take any $g$-dimensional principally polarized abelian 
variety $A_K$ defined over $K$ with a level $n$ structure. Under the assumption that $n\geq 3$, its N\'eron Model $\varphi:{\cal A}\rightarrow B$ is known to have semi-abelian reduction
at any closed points of $B$ so we have 
$$\Sigma_{\varphi}\subset \phi^{-1}({\cal A}_{g,n}^{\ast}\setminus
{\cal A}_{g,n})\subset  \phi^{-1}(\mbox{Supp }H)
$$
and hence,
\begin{equation}\label{equation:bounding s}
s_{\varphi}:=\mbox{Card }\Sigma_{\varphi}\leq \deg \phi^{\ast}H=N(g,n)\deg \phi^{\ast}\omega_{{\cal A}_{g,n}^{\ast}}=N(g,n)\deg e^{\ast}\det\Omega^1_{{\cal A}/B}.
\end{equation}
If $g$-dimensional principally polarized abelian 
variety $A_K$ does not have a level $n$ structure, there exists an finite extension $K^{\prime}$ of $K$ with 
$[K^{\prime}:K]\leq \mbox{Card }\mbox{Sp}(2g,\mbox{\boldmath $Z$}/n\mbox{\boldmath $Z$})$ such that $A_{K^{\prime}}:=A_{K}\times_{K}K^{\prime}$ has a level $n$-structure.
Let $\tau:B^{\prime}\rightarrow B$ be the finite surjective morphism from the smooth projective model $B^{\prime}$ of $K^{\prime}$ and let 
$\varphi^{\prime}:{\cal A}^{\prime}\rightarrow B^{\prime}$ be the N\'eron Model of $A_{K^{\prime}}$. If we assume that $\varphi:{\cal A}\rightarrow B$ has semi-abelian reduction
at any closed points of $B$, then we have ${\cal A}^{\prime}={\cal A}\times_{B}B^{\prime}$ and hence, $\Sigma_{\varphi^{\prime}}=\tau^{-1}(\Sigma_{\varphi})$. Thus,
by (\ref{equation:bounding s}), we have

\begin{eqnarray}
s_{\varphi}\leq s_{\varphi^{\prime}}&\leq & N(g,n)\deg e^{\prime\ast}\det\Omega^1_{{\cal A}^{\prime}/B^{\prime}}\\
&=&N(g,n)\deg \tau\deg e^{\ast}\det\Omega^1_{{\cal A}/B}\\
&=&N(g,n)\mbox{ Card }\mbox{Sp}(2g,\mbox{\boldmath $Z$}/n\mbox{\boldmath $Z$})\deg e^{\ast}\det\Omega^1_{{\cal A}/B}.
\end{eqnarray}

\vskip 5mm

{\it Proof of Theorem~\ref{th:bouding the number of singular fibres}.}
As we preiviously remarked, we only have to bound $s_{\varphi^{\prime}}$. From Lemma~\ref{lm:stable reduction zeta}, we have 
\begin{eqnarray}
s_{\varphi^{\prime}}&=&s_{\zeta^{\prime}}\\
&\leq & N(16, 3)\mbox{ Card }\mbox{Sp}(32,\mbox{\boldmath $Z$}/3\mbox{\boldmath $Z$})\deg \zeta^{\prime}_{\ast}\det\Omega^1_{{\cal Z}^{\prime}/B^{\prime}}\\
&\leq &16dN(16, 3)\mbox{ Card }\mbox{Sp}(32,\mbox{\boldmath $Z$}/3\mbox{\boldmath $Z$}).
\end{eqnarray}

Thus we get the assertion.
\hfill\bsquare

Department of Mathematics Graduate School of Science, \\Osaka University, \\Toyonaka, Osaka, 560-0043, Japan
\\{\it E-mail address\/}: koji@math.sci.osaka-u.ac.jp


\vskip 5mm













\begin{thebibliography}{1}

\bibitem{abhyankar} S. S. Abhyankar, {\em Local Analytic Geometry}, Academic Press, New York-London, 1964.

\bibitem{artin-mazur} M. Artin and B. Mazur, {\em Etale Homotopy}, Lecture Notes in Math. {\bf 100}, Springer Heidelberg, 1969.

\bibitem{an} V. A. Alekseev and V. V. Nikulin, {\em Classification of del Pezzo Surfaces with
Log-Terminal Singularities of Index $\leq 2$, and Involutions on $K3$ Surfaces}, Soviet Math. Dokl.
{\bf 39}, (1989), No.3, pp. 507-511.


\bibitem{BLR} S. Bosh, W. L\"utkebohmert and M. Raynaud, {\em N\'eron Models}, Ergeb. Math. Grenzgeb (3), 
{\bf 21}, Spriger-Verlag, Berlin-Heidelberg-New York-London-Paris-Tokyo-Hong
Kong, (1980).



\bibitem{bourbaki} N. Bourbaki, {\em Topologie G\'en\'erale}, Herman, Paris, 1961.


\bibitem{brieskorn} E. Brieskorn and H. Kn\"orrer, {\em Plane Algebraic Curves}, Birkh\"auser,
1986.

\bibitem{catanese} F. Catanese, {\em Automorphisms of Rational Double Points and Moduli Spaces of
Surfaces of General Type}, Compos. Math. {\bf 61}, (1987), pp. 81-102.


\bibitem{corti/utah} A. Corti, {\em Adjunction of log divisors,} \lq\lq Flips  and
abundunce for algebraic threefolds'',
  Ast\'erisque {\bf 211}, A summer seminar at the university of Utah Salt Lake City. 1992.
 

\bibitem{crauder} B. Crauder and D. Morrison, {\em Triple point free degenerations of surfaces with
Kodaira number zero}, 
\lq \lq The Birational Geometry of Degenerations" Progress in Math. {\bf 29}, Birkh\"auser, 1983, 
pp.353-386. 

\bibitem{crauder2} B. Crauder and D. Morrison, {\em Minimal Models and Degenerations of surfaces
with Kodaira number Zero}, Trans. of the Amer. Math. Soc. {\bf 343}, No. 2, (1994), pp. 525-558.

\bibitem{dimca} A. Dimca, {\em Singularities and topology of hypersurfaces}, Springer-Verlag, 1992.

\bibitem{faltings} G. Faltings, {\em Arakelov's theorem for Abelian Varieties}, Invent.math., {\bf 73}, (1983), pp. 337-347.



\bibitem{faltings-chai} G. Faltings and C-L. Chai, {\em Degeneration of Abelian Varieties}, Ergeb. Math. Grenzgeb (3), 
{\bf 22}, Spriger-Verlag, Berlin-Heidelberg-New
York-London-Paris-Tokyo-Hong Kong-Barcelona, (1990).



\bibitem{fong-mc/utah} L.-Y. Fong and J. Mckernan, {\em Log abundance for surfaces,} \lq\lq Flips  and
abundunce for algebraic threefolds'',
  Ast\'erisque {\bf 211}, A summer seminar at the university of Utah Salt Lake City. 1992.
 

\bibitem{fischer} G. Fisher, {\em Complex Analytic Geometry}, Lecture Notes in Math. {\bf 538},
Springer-Verlag, Berlin-Heidelberg-New York, 1976.


\bibitem{fujiki} A. Fujiki, {\em Relative Algebraic Reduction and Relative Albanese Map for a Fiber Space in ${\cal C}$},
Publ. RIMS, Kyoto Univ. {\bf 19}, (1983), pp. 207-236.


\bibitem{fujino-mori} O. Fujino and S. Mori, 
{\em A Canonical Bundle Formula,} J. Differential Geometry {\bf 56} (2000), pp.167-188.


\bibitem{fujita} T. Fujita, {\em Fractionally logarithmic canonical rings of algebraic surfaces},
J. Fac. Sci. Univ. Tokyo Sect. IA {\bf 30}, (1984), pp. 685-696.

\bibitem{fujita2} T. Fujita, {\em Zariski decomposition and canonical rings of elliptic threefolds},
J. Math. Soc. Japan {\bf 38}, No.1, (1986), pp. 19-37.

\bibitem{furushima} M. Furushima, {\em Singular del Pezzo surfaces and analytic compactifications
of $3$-dimensional complex  affine space ${ \mbox{\boldmath $C$} }^3$,} Nagoya Math. J. {\bf
104}, (1986), pp. 1-28. 

\bibitem{grauert-remmert} H. Grauert and R. Remmert, {\em Komplexe R\"aume}, Math. Ann. {\bf 136}, (1958), pp. 245-318.


\bibitem{grothendieck} A. Grothendieck, {\em Techniques de Construction en G\'eom\'etrie Analytique}, S\'eminaire Henri Cartan,13ieme ann\'ee (1960/1).

\bibitem{grothendieck-f} A. Grothendieck, {\em Foundations de la G\'eom\'etrie Analytique} [Extraits du S\'eminaire Bourbaki 1957-1962], Paris 1962.


\bibitem{sga1} A. Grothendieck et. al., {\em Rev\^etements Etales et Groupe Fondamental}, Lecture Notes in Math. {\bf 224},  Springer-Verlag,
Berlin-Heidelberg-New York, 1971. 


\bibitem{sga4} A. Grothendieck et. al., {\em T\'eorie des Topos et Cohomologie \'Etale des Sch\'emas},Lecture Notes in Math. {\bf 269},  Springer-Verlag, Berlin-Heidelberg-New York,
1972. 



\bibitem{hironaka} H. Hironaka, {\em Resolution of Singularities of an Algebraic Variety over a
Field of Characteristic Zero I, II}, Ann. of Math. {\bf 79}, (1964), pp. 109-326.

\bibitem{hunt} B. Hunt, {\em A bound on the Euler number for certain Calabi-Yau 3-folds}, J. reine angew. Math. {\bf 411}, (1990), pp.137-170.

\bibitem{kato} M. Kato, {\em On the uniformizations of orbifolds}, Adv. Stud. in Pure Math. {\bf 9}, (1986), pp.149-172.


\bibitem{kawamata crep.}Y. Kawamata, 
{\em Crepant blowing-up of $3$-dimensional canonical singularities 
and its application to degeneration of surfaces,} Jour. of Ann. of Math. {\bf 127},
(1988), pp.93-163. 


\bibitem{kawamata ab.} Y. Kawamata,
{\em Abundance theorem for minimal threefolds,} Jour. of Invent. math. {\bf 108}, (1992), pp. 229-246.

\bibitem{kawamata term.} Y. Kawamata,
{\em Termination of Log Flips for Algebraic 3-folds,} International Jour. Math. {\bf 30}, No.5 (1992)  pp. 653-659.

\bibitem{mumford} G. Kempf, F. Knudsen, D. Mumford and B. Saint-Donat, {\em Toroidal embeddings.
I},  Lecture Notes in Math. {\bf 339}, Springer-Verlag, Berlin, 1973.


\bibitem{kns} R. Kobayashi, S. Nakamura and F. Sakai, {\em A Numerical Characterization of Ball
Quotients for Normal Surfaces with Branch Loci}, Proc. Japan Acad., {\bf 65}, Ser. A (1989),
pp.238-241.

\bibitem{kobayashi} R. Kobayashi, {\em Uniformization of Complex Surfaces,} 
Adv. Stud. in Pure Math. {\bf 18-II}, (1990), pp. 313-394.


\bibitem{kodaira} K. Kodaira, {\em On compact analytic surfaces II}, Ann. of Math. {\bf 77},
 (1963), pp.563-626.


\bibitem{kollar} J. Koll\'ar, {\em Flops,} Nagoya Math. J. {\bf 113} (1989), pp.15-36.


\bibitem{kollar/utah} J. Koll\'ar, {\em Log flips and abundance $:$ An overview,} \lq\lq Flips  and
abundunce for algebraic threefolds'',
  Ast\'erisque {\bf 211}, A summer seminar at the university of Utah Salt Lake City. 1992.




\bibitem{kulikov} V. Kulikov, {\em Degeneration of K$3$ surfaces and Enriques surfaces,} Jour.
of Math.USSR Izv. {\bf 11}, (1977), pp.957-989.



\bibitem{lipman-sommese} J. Lipman and A. J. Sommese, {\em On blowing down projective spaces in
singular varieties}, J. Reine. Angew. Math. {\bf 362}, (1985), pp. 51-62. 

\bibitem{matsumura} H. Matsumura, {\em Commutative Algebra}, Second Edition, The Benjamin/Cummings Publishing Company, Inc., 1980.



\bibitem{megyesi} G. Megyesi, {\em Chern classes of Q-sheaves}, \lq\lq Flips  and
abundunce for algebraic threefolds'',
  Ast\'erisque {\bf 211}, A summer seminar at the university of Utah Salt Lake City. 1992.

\bibitem{milnor} J. Milnor, {\em Singular points of complex hypersurfaces}, Princeton Univ. Press, Princeton, New Jersey and Univ. of Tokyo Press, Tokyo, 1968.


\bibitem{miranda-persson} R. Miranda and U. Persson, {\em On Extremal Rational Elliptic Surfaces},
Math. Z. {\bf 193}, (1986), pp.537-558.

\bibitem{mitchell} B. Mitchell, {\em Theory of Categories}, Academic Press, New York and London, 1965.



\bibitem{miyanishi-zhang} M. Miyanishi and D. Q. Zhang, {\em Gorenstein log del pezzo surfaces of
rank one,} Jour. of Algebra {\bf 118}, No. 1, (1988), pp.63-84.



\bibitem{miyaoka} Y. Miyaoka,  {\em The maximal number of quotient singularities on surfaces 
with given numerical invariants,} Math. Ann. {\bf 268}, (1984), pp.159-171.



\bibitem{moret-baily} L. Moret-Bailly, {\em Pinceaux de vari\'et\'es ab\'eliannes,} Ast\'erisque {\bf 129} (1985).




\bibitem{morrison} D. Morrison,  {\em Semistable degenerations of Enriques' and Hyperelliptic surfaces,} Duke Math. Jour. {\bf 48}, No. 1 (1981), pp.197-249.




\bibitem{mori.cl} S. Mori, {\em Classification of higher-dimensional varieties,} Proc. of Symp. in Pure Math. {\bf 40}, (1987), pp. 269-331.



\bibitem{mori} S. Mori, {\em Flip theorem and the existence of minimal models,} Jour. of the
AMS. {\bf 1}, (1988), pp. 117-253.



\bibitem{murre} J. Murre, {\em Lectures on an Introduction to Grothendiek's Theory of the Fundamental Group},
 Lecture Notes, Tata Institute of Fundamental Research, Bombay, 1967.



\bibitem{nagata} M. Nagata, {\em Field Theory}, Pure and Applied Mathematics, A series of Monographs and Textbooks,
Marcel Dekker, Inc.,  New York and Basel, 1977.

\bibitem{nakano} T. Nakano and K. Tamai, {\em On Some Maximal Galois Coverings over Affine and Projective Planes}, Osaka J. Math. {\bf 33} (1996), 347-364.


\bibitem{nakayama}N. Nakayama, 
{\em The lower semi-continuity of the plurigenera of complex varieties} Algebraic Geometry, 
Sendai 1985, Adv. Stud. Pure Math. {\bf 10}, Kinokuniya, Tokyo, and North-Holland, Amsterdam, 1987, 
pp.551-590.

\bibitem{namba} M. Namba, {\em Branched coverings and algebraic functions}, Pitman Research Notes in Math. Series {\bf 161}, 
Longman Scientific and Technical, Harlow, John Wiley and Sons, New York, 1987.



\bibitem{oguiso-shioda} K. Oguiso and T. Shioda, {\em The Mordell-Weil Lattice of a
Rational Elliptic Surface},  Comment. Math. Univ. St. Pauli {\bf 40}, No.1, (1991), pp. 83-99.

\bibitem{oguiso} K. Oguiso, {\em A Note on Moderate Abelian Fibrations},  Contemporary Math. {\bf 207}, (1997), pp.101-118.


\bibitem{ohno} K. Ohno, {\em Toward determination of the singular fibers of minimal degeneration 
of surfaces with $\kappa =0$}, 
 Osaka J. Math. {\bf 33} No.1, (1996), pp. 235-305.


\bibitem{persson-pinkham} U. Persson and H. Pinkham, {\em Degeneration of surfaces with trivial canonical divisor}, Ann. of Math. {\bf 113} (1981), pp. 45-66.


\bibitem{persson} U. Persson, {\em Configurations of Kodaira fibers on rational elliptic
surfaces},  Math. Z. {\bf 205}, (1990), pp. 1-49.




\bibitem{prill} D. Prill, {\em Local Classification of Quotients of Complex Manifolds by
Discontinuous Groups}, Duke Math. J. {\bf 34}, (1967) pp. 375-386.


\bibitem{reid c} M. Reid, {\em Canonical threefold}, G\'eom\'etry Alg\'ebrique Angers, A. Beauville
ed., Sijthoff \& Noordhoff, 1980, pp. 273-310.

\bibitem{reid} M. Reid, {\em Minimal models of canonical $3$-folds,}  Algebraic Varieties and
Analytic Varieties, Adv. Stud. in Pure Math. {\bf 1}, 1983, pp.131-180.

\bibitem{reid2} M. Reid, {\em Young person's guide to canonical singularities,} 
Proc. Symp. in Pure Math. {\bf 46}, 1987, pp.345-414.

\bibitem{ribes} L. Ribes, {\em Introduction of Profinite Groups and Galois Cohomology}, Queen's Papers in Pure and Applied Mathemathecs-No. {\bf 24}, Queen's University, Kingston,
Ontario, 1970.



\bibitem{sakai} F. Sakai, {\em Weil divisors on Normal surfaces,} 
Duke Math. J. {\bf  51}, No. 4, (1984), pp. 877-887.

\bibitem{seidenberg} A. Seidenberg, {\em The hyperplane sections of normal varieties}, Trans. of the
Amer. Math. Soc., {\bf 69}, No. 2, (1950), pp. 357-386.

\bibitem{seidenberg2} A. Seidenberg, {\em Saturation of an analytic ring}, Amer. Jour. of Math., {\bf 94}, (1972), pp. 424-430.

\bibitem{serre} J.P. Serre, {\em Topics in Galois Theory}, Jones and Bartlett Publishers, Boston-London, 1992.


\bibitem{shepherd-wilson} N. I. Shepherd-Barron and P. M. H. Wilson, {\em Singular $3$-folds with
Numerically Trivial First and Second Chern Classes}, J. Algebraic Geometry {\bf 3}, (1994),
pp.265-281.


\bibitem{shioda} T. Shioda, {\em On the Mordell-Weil Lattices,} Comment. Math. Univ. St.
Pauli {\bf 39}, (1990), pp.211-240.



\bibitem{shokurov} VV. Shokurov, {\em $3$-fold Log Flips,} Russian Acad. Sci. Izv.
Math. {\bf 40}, (1993), pp.95-202.

\bibitem{shokurov2} VV. Shokurov, {\em Semistable $3$-fold Flips,} Russian Acad. Sci. Izv. Math. {\bf 42}, No.2, (1993), pp.371-425.

\bibitem{shokurov-cp} VV. Shokurov, {\em Complements on surfaces}, preprint, 1997.

\bibitem{steenrod} N. Steenrod, {\em The Topology of Fibre Bundles}, Princeton Math. Series {\bf 14}, Princeton UP, 1951.



\bibitem{endre} E. Szab\'o, {\em Divisorial log terminal singularities}, J. Math. Sci. Univ. Tokyo
{\bf 1}, (1994), pp. 631-639.


\bibitem{gabi} G. Teodosiu, {\em A Class of Analytic Coverings Ramified over $u^3=v^2$}, J. London Math. Soc. (2) {\bf 38}, (1988), pp. 231-242.


\bibitem{tsunoda1} S. Tsunoda, {\em Structure of open algebraic surfaces, I}, J. Math. Kyoto
Univ. {\bf 23}, No. 1, (1983), pp. 95-125.


\bibitem{ueno1} K. Ueno, {\em On fibre spaces of normally polarized abelian varieties of
dimension $2$, I}, Jour. Fac. Sci. Univ. Tokyo, Sect. IA {\bf 18}, (1971), pp.37-95. 

\bibitem{ueno2} K. Ueno, {\em On fibre spaces of normally polarized abelian varieties of
dimension $2$, II}, Jour. Fac. Sci. Univ. Tokyo, Sect. IA {\bf 19}, (1972)  pp.163-199.

\bibitem{ueno3} K. Ueno, {\em On algebraic fibre spaces abelian varieties}, Math. Ann. {\bf 237}, (1978)  pp.1-22.


\bibitem{zariski} O. Zariski, {\em On the problem of two variables possessing a given branch
curve}, Amer. J. Math., {\bf 51}, (1929) pp. 305-328.

\bibitem{zarkhin} Yu. G. Zarkhin, {\em Endomorphisms of abelian varieties and points of finite order in characteristic $p$}, Math. Notes of the Academy of Sciences of the USSR, A
translation of Mathematicheskie Zametki, {\bf 21}, (1977) pp. 415-419.


\end{thebibliography}
\end{document}